\documentclass[final,faculty=fw,department=wis,phddegree=mat]{adsphd}

\usepackage{nomencl}   

\title{The monodromy property for \texorpdfstring{$K3$}{K3} surfaces allowing a triple-point-free model}

\author{Annelies}{Jaspers}
\supervisor{Prof.~dr.~L.~H.~Halle}{University of Copenhagen}
\supervisor{Prof.~dr.~J.~Nicaise}{Imperial College London \& KU Leuven}

\president{Prof.~dr.~J.~Quaegebeur}
\jury{}
\externaljurymember{Prof.~dr.~N.~Budur}{KU Leuven}
\externaljurymember{Prof.~dr.~A.~H\"oring}{Universit\'e de Nice}
\externaljurymember{Prof.~dr.~D.~Petersen}{University of Copenhagen}
\externaljurymember{Prof.~dr.~W.~Veys }{KU Leuven}

\researchgroup{Algebraic Geometry and Number Theory Group}
\website{} 
\email{} 

\address{Celestijnenlaan 200B box 2400}
 \addresscity{B-3001 Leuven} 

\date{}
\copyyear{2017}

 \setcustomcoverpage{mycoverpage.tex} 

\usepackage{textcomp} 
\usepackage{pifont}   
\usepackage{booktabs}
\usepackage{color}
\usepackage{url}
\usepackage{caption}
\usepackage{enumerate}
\usepackage{amsmath,amsfonts,amsthm,amssymb}
\usepackage{graphicx,xcolor}
\usepackage{pdfpages}
\usepackage{array}
\usepackage{multirow}
\usepackage{subcaption}
\usepackage[all,cmtip]{xy}
\usepackage{longtable}
\usepackage{tikz}
\usepackage{mathrsfs}
\usepackage{amscd,graphics,amsbsy}
\usepackage{subfiles}
\usepackage{latexsym}
\usepackage{mathrsfs}
\usepackage{nth} 
\usepackage{parskip}
\usepackage{todonotes}
\usepackage{float}
\usepackage{stmaryrd}
\usepackage{listings}
\usepackage{verbatim}
\usepackage{rotating}
\usepackage{anyfontsize}
\usetikzlibrary{positioning}
\usetikzlibrary{shapes,arrows}

\renewcommand{\O}{\mathcal{O}}
\newcommand{\fiber}{\mathcal{F}}

\newcommand{\Z}{\mathbb{Z}}
\newcommand{\Q}{\mathbb{Q}}
\newcommand{\R}{\mathbb{R}}
\newcommand{\C}{\mathbb{C}}
\newcommand{\A}{\mathbb{A}}
\renewcommand{\L}{\mathbb{L}}
\renewcommand{\P}{\mathbb{P}}
\newcommand{\X}{\mathcal{X}}
\newcommand{\Y}{\mathcal{Y}}
\newcommand{\p}{\mathfrak{p}}
\newcommand{\q}{\mathfrak{q}}
\newcommand{\m}{\mathfrak{m}}
\newcommand{\Gal}{\mathrm{Gal}}
\newcommand{\Aut}{\mathrm{Aut}}
\newcommand{\Sch}{\mathrm{Sch}}

\DeclareMathOperator{\Spec}{Spec}
\DeclareMathOperator{\ord}{ord}
\DeclareMathOperator{\Pic}{Pic}
\DeclareMathOperator{\Frac}{Frac}
\DeclareMathOperator{\lcm}{lcm}

\newtheorem{theorem}{Theorem}[section]
\newtheorem*{theorem*}{Theorem}
\newtheorem{lemma}[theorem]{Lemma}
\newtheorem*{lemma*}{Lemma}
\newtheorem{proposition}[theorem]{Proposition}
\newtheorem*{proposition*}{Proposition}
\newtheorem{corollary}[theorem]{Corollary}

\newtheorem{definition}[theorem]{Definition}
\newtheorem{prop-def}[theorem]{Proposition -- Definition}
\newtheorem{conjecture}[theorem]{Conjecture}
\theoremstyle{definition}
\newtheorem{example}[theorem]{Example}
\newtheorem{countercandidate}[theorem]{Combinatorial Countercandidate}
\newtheorem{remark}[theorem]{Remark}

\makeatletter
\def\thm@space@setup{%
  \thm@preskip=\parskip \thm@postskip=0pt
}
\makeatother

\definecolor{gray}{gray}{0.5}
\newcommand{\note}[2][]{\ignorespaces} 

\begin{document}

\makefrontcoverXII

\maketitle

\frontmatter 


\chapter*{Dankwoord}                                  \label{ch:acknowledgement}
\addcontentsline{toc}{chapter}{Dankwoord}
\markboth{DANKWOORD}{DANKWOORD}

Na vier jaar werken is deze thesis eindelijk klaar. Dit heb ik zeker niet alleen gedaan: veel mensen hebben een steentje bijgedragen.

I would like to start by thanking my two excellent supervisors: Johannes and Lars. 
Johannes, bedankt voor alle tijd die je gespendeerd hebt in de laatste vier jaar met uitleggen, nalezen en coachen, voor alle wiskundige hulp en suggesties.
Bedankt om altijd in mij te geloven, ook als ik dat zelf niet deed, bedankt om mij de kans te geven om Londen te verkennen, bedankt om mij de vrijheid te geven om zo veel conferenties bij te wonen en bedankt dat ik af en toe in Champaign mocht werken. 
En ten slotte, bedankt om gewoon interesse te tonen in mij als persoon: mijn weekendplannen, mijn reizen, mijn ambities voor na het doctoraat.

Lars, every visit to Copenhagen was amazing.  Whenever I was there, my research leapt forward immensely. I am absolutely grateful for the opportunity to work together intensely. Working with you is such a pleasure and I enjoyed our many fruitful blackboard discussions a lot. 
I would also like to thank you for your optimism and enthusiasm. There were times that I lost the joy of doing mathematical research, but your motivating words and your faith in me encouraged me to continue and finish this PhD. 

I am honored that Nero Budur, Andreas H\"oring, Dan Petersen, Johan Quaegebeur, and Wim Veys accepted to be members of my examination committee. Thanks for reading this manuscript, because this thesis was certainly improved by your suggestions and comments.

One of the things I enjoyed the most during these past four years were the many stays abroad.
Thanks to the department of mathematics of the university of Copenhagen for hosting me three times. Every detail of my stay was always perfectly planned: an office, a place to stay, a bike, subscriptions to the mailing lists, $\ldots$. Thanks for making my visits so easy and pleasant. And of course, the colleagues there had an important contribution to my well-being too. Thanks to Dino, Fabien, Giacomo, Nadim, Olivier, Raymond, and Sho for organising interesting reading seminars, for the Thursday cake club, for the many dinner parties and nice lunch breaks, and for even occasionally obliging me to speak French.

My six month visit to Imperial College in London wouldn't have been the same without the presence of the many nice people I met there. I would especially like to mention Arne, Enrica, Jacub, Peter, and Thomas. Thanks, Enrica, for your sincere friendship and for inviting me to all the social activities of LSGNT. Thanks, Peter, for offering me an office \emph{with} a window, for your insights in American politics, and for listening to my worries and doubts. Thanks, Arne, for your continuous interest in my work and my well-being. Thanks, Jacub, for the weekly farmers' market lunches and for the interesting mathematical discussions. Thanks, Thomas, for making the hunt for an affordable place to stay in London unexpectedly easy, and for all the useful tips to discover the city.

Many of my dearest memories during this PhD happened during conferences and summer schools: I met many interesting mathematicians, I learned a lot during talks, and I felt part of a community. But I will probably remember most what happened in the evenings and during free afternoons: returning after midnight from the mountains in Nordfjordeid, tasting salami, cheese and wine in an osmiza in Trieste, drinking tequila in Guanajuato, visiting the not-that-exciting China Lights exhibition and afterwards laughing for hours in our hotel room in Utrecht, hiking the white rocks of the Calanques in Luminy, being impressed by the presence of many Fields medaillist in Salt Lake City, swimming in the Donau in Regensburg, and many more. Thanks to all the people who shared these experiences with me!

Although some would doubt it, I still spent most of my time in the department in Leuven. Thanks to my colleagues for being excellent colleagues! Thanks for the ice cream breaks, for the cookie breaks, for the cake breaks, for the board game nights, for the movie nights, for the risotto nights and for the cocktail nights. In het bijzonder zou ik graag de collega's van het gelijkvloers willen bedanken: Alexander, Christophe, Hans, Jasper, Lore en Marta. Wat waren wij toch een gezellige bende! Samen elke minuut van de regeringsvorming volgen, samen wiskunde-olympiadevraagjes oplossen, samen video's van schattige panda's kijken, samen klagen over de relevantie van ons werk, en samen vreugdedansjes doen bij doorbraken. Maar uiteraard steunden we elkaar ook op wiskundig vlak en konden we altijd bij elkaar terecht met vragen, hoe dom ook. Bedankt voor alles!

Twee collega's liggen me extra nauw aan het hart: Jasper en Marta. 
Jasper, samen begonnen aan het doctoraat, en samen eindigen, samen vier jaar lang een kantoor delen. Bedankt voor alle raad, voor alle gesprekken over politiek, over feminisme en recent over de zoektocht naar een job. 

Marta, I was always welcome in your office to chat, to gossip, to laugh, for a tea break, and for a hug when my research didn't progress as hoped for. I would also like to thank you for all the amazing experiences during the many conferences together. You were not only a colleague, but your were also one of my best friends. I missed you a lot this past year!

En dan zijn er nog de mensen die er voor gezorgd hebben dat ik geen complete vakidioot ben geworden en die er altijd voor mij waren om samen plezier te beleven.
Bedankt aan Janne, Lynn en Val\'erie voor de woensdagse stepaerobic-afleiding, lange Caverna-avonden, lunch-onderonsjes en klimsessies. Bedankt aan Sofie voor onze tijd samen in Kopenhagen. Het was geweldig om met jou de stad te verkennen, af en toe nog eens Nederlands te praten, en tips uit te wisselen over hoe al ons bezoek rond te leiden. Bedankt aan Eveline voor de maandelijkse eerste-van-de-maandmailtjes en de bezoekjes in zowel Kopenhagen als Londen. Bedankt aan Pieter om regelmatig binnen te springen op mijn kantoor voor wat broodnodige afleiding, en bedankt om zo een goede fondue-partner te zijn in de chocoladebar. Ik ben een bofkont met zo'n schitterende vrienden! 

Ik ben ook mijn ouders ontzettend dankbaar. Jullie staan altijd, maar dan ook altijd, voor mij klaar: met een luisterend oor, met een theetje of een chocolaatje, met raad en met begrip dat je niet hoeft te vragen waarover mijn doctoraat gaat. Het is fijn om zo te kunnen rekenen op twee schatten die zo trots op mij zijn. Ik weet dat jullie me veel gemist hebben en ik weet dat jullie misschien stiekem hoopten dat ik wat meer in Belgi\"e bleef, maar daar hebben jullie nooit iets van laten merken. Integendeel, jullie hebben me altijd nog extra aangemoedigd om mijn ding te doen, hoe avontuurlijker, hoe beter. Ik ben jullie daar eeuwig dankbaar voor.

En ten slotte, wil ik nog mijn lieve Daan bedanken. Gelukkig heb je, door onze dagelijkse Skype-gesprekken, geen enkel euforisch moment moeten missen: wiskundige doorbraken, een eerste geaccepteerde publicatie, complimentjes van een promotor, een thesis die naar de jury kon doorgestuurd worden. Maar jij hebt ook alle frustratie, alle tranen, alle twijfels van dichtbij meegemaakt en gedeeld. Jij was de enige die mij het gevoel gaf mij echt te begrijpen. Bedankt om mij kalm te maken als het nodig was, bedankt om ingenieuze back-upsystemen te verzinnen, bedankt om mij urenlang over mijn onzekerheden te laten praten, bedankt om mij aan te moedigen om door te zetten, bedankt om programmeerraad te geven, en bedankt om klaar te staan met een glas cava als ik het \'echt verdiend heb.
Je maakt mij intens gelukkig, maar dat weet je.

\cleardoublepage


\chapter*{Abstract}                                 \label{ch:abstract}
\addcontentsline{toc}{chapter}{Abstract}
\markboth{}{}

The aim of this thesis is to study under which conditions $K3$ surfaces allowing a triple-point-free model satisfy the monodromy property. This property is a quantitative relation between the geometry of the degeneration of a Calabi-Yau variety $X$ and the monodromy action on the cohomology of $X$: a Calabi-Yau variety $X$ satisfies the monodromy property if poles of the motivic zeta function~$Z_{X,\omega}(T)$ induce monodromy eigenvalues on the cohomology of $X$.

Let $k$ be an algebraically closed field of characteristic 0, and set $K=k(\!( t )\!)$.
In this thesis, we focus on $K3$ surfaces over $K$ allowing a triple-point-free model, i.e., $K3$ surfaces allowing a strict normal crossings model such that three irreducible components of the special fiber never meet simultaneously. Crauder and Morrison classified these models into two main classes: so-called  flowerpot degenerations and chain degenerations. This classification is very precise, which allows to use a combination of geometrical and combinatorial techniques to check the monodromy property in practice.

The first main result is an explicit computation of the poles of $Z_{X,\omega}(T)$ for a $K3$ surface $X$ allowing a triple-point-free model and a volume form $\omega$ on $X$. We show that $Z_{X,\omega}(T)$ can have more than one pole. This is in contrast with previous results: so far, all Calabi-Yau varieties known to satisfy the monodromy property have a \emph{unique} pole.

We prove that $K3$ surfaces allowing a flowerpot degeneration satisfy the monodromy property. We also show that the monodromy property holds for $K3$ surfaces with a certain chain degeneration.
We don't know whether all $K3$ surfaces with a chain degeneration satisfy the monodromy property, and we investigate what characteristics a $K3$ surface $X$ not satisfying the monodromy property should have. We prove that there are 63 possibilities for the special fiber of the Crauder-Morrison model of a $K3$ surface~$X$ allowing a triple-point-free model that does not satisfy the monodromy property.

\cleardoublepage


\chapter*{Samenvatting}
\addcontentsline{toc}{chapter}{Samenvatting}
\markboth{}{}

In deze thesis bestuderen we onder welke voorwaarden $K3$-oppervlakken voldoen aan de monodromie-eigenschap. Deze eigenschap is een kwantitatieve relatie tussen de degeneratie van een Calabi-Yau-vari\"eteit $X$ en de monodromieactie op de cohomologie van $X$: een Calabi-Yau-vari\"eteit voldoet aan de monodromie-eigenschap als elke pool van de motivische zetafunctie $Z_{X,\omega}(T)$ overeenkomt met een monodromie-eigenwaarde op de cohomologie van $X$.

Zij $k$ een algebra\"isch gesloten veld van karakteristiek 0 en neem ${K=k(\!( t)\!)}$. In deze thesis concentreren we ons op $K3$-oppervlakken over $K$ met een strikt-normale-kruisingenmodel zonder drievoudige punten. Dit zijn $K3$-oppervlakken met een strikt-normale-kruisingenmodel zodat drie irreduciebele componenten nooit een gemeenschappelijke intersectie hebben. Crauder en Morrison klassificeerden deze modellen in twee klassen: de zogenaamde bloempotdegeneraties en ketendegeneraties. Deze klassificatie is heel precies zodat we een combinatie van meetkundige en combinatorische technieken kunnen gebruiken om de monodromie-eigenschap na te gaan.

Het eerste grote resultaat in deze thesis is een expliciete berekening van de polen van $Z_{X,\omega}(T)$ voor een $K3$-oppervlak $X$ met een strikt-normale-kruisingenmodel zonder drievoudige punten en een volumevorm $\omega$ op $X$. We tonen aan dat $Z_{X,\omega}(T)$ meer dan \'e\'en pool kan hebben. Dit is in contrast met eerdere resultaten: alle Calabi-Yau-vari\"eteiten waarover we tot nu toe wisten dat ze aan de monodromie-eigenschap voldeden, hadden een unieke pool.

We bewijzen verder dat $K3$-oppervlakken met een bloempotdegeneratie aan de monodromie-eigenschap voldoen. We tonen ook aan dat de monodromie-eigenschap geldt voor sommige $K3$-oppervlakken met een ketendegeneratie. We weten echter niet of alle $K3$-oppervlakken met een ketendegeneratie deze eigenschap hebben, en we bestuderen $K3$-oppervlakken die er misschien niet aan voldoen. 
We geven een exhaustieve lijst van 63 mogelijkheden voor de speciale vezel van het Crauder-Morrisonmodel van een $K3$-oppervlak dat niet aan de mononodromie-eigenschap voldoet.
\cleardoublepage


\setcounter{tocdepth}{2}
\tableofcontents


\mainmatter 

\chapter*{Introduction}\label{ch:introduction}
\addcontentsline{toc}{chapter}{Introduction}
\markboth{INTRODUCTION}{INTRODUCTION}

In this thesis, we study the monodromy property for a specific type of Calabi-Yau variety: for $K3$ surfaces allowing a triple-point-free model. The monodromy property is a precise relation between the geometry of $snc$-models of a Calabi-Yau variety $X$ and the cohomology of $X$. Inspired by the $p$-adic and motivic monodromy conjecture, Halle and Nicaise formulated this property in \cite{HalleNicaiseAbelian}. 

There is a good reason why we study the specific class of $K3$ surfaces allowing a triple-point-free model: Crauder and Morrison classified triple-point-free $snc$-models of $K3$ surfaces, and such a classification is an extremely useful tool to verify the monodromy property in practice. To our knowledge, this is the only systematic classification of $snc$-models of Calabi-Yau varieties without a semi-stability condition.

Fix an algebraically closed field $k$ of characteristic 0, and set $K=k(\!(t)\!)$ and $R=k\llbracket t \rrbracket$.
Let us now explore the results obtained in this thesis in more detail.

\section*{The first ingredient: the monodromy property for Calabi-Yau varieties}

\subsubsection*{\texorpdfstring{$p$}{p}-adic zeta functions and the \texorpdfstring{$p$}{p}-adic monodromy conjecture}

For a very long time, mathematicians have been interested in solutions $x\in \Z^m$ of congruences of the form
\[f(x) \equiv 0 \mod n,\]
for $f$ a polynomial in $\Z[x_1, \ldots, x_m]$ and $n$ an integer. For example Euler, Legendre and Gauss made significant progress for quadratic polynomials~$f$ in one variable. In 1796, Gauss proved the law of quadratic reciprocity, which makes it possible to determine for any odd prime $p$, and for any quadratic equation of the form $x^2\equiv a \mod p$, whether it has a solution or not.

Another natural question mathematicians wonder about is: \emph{how many} solutions does the equation 
$f(x) \equiv 0 \mod n$
 have?
This problem can be simplified, thanks to the Chinese remainder theorem: we actually only need to study the case where $n$ is a power of a prime. Let's denote by $N_d$ the number of solutions of $f(x)\equiv 0 \mod p^{d+1}$.
One way to study the number $N_d$ for $d$ big, is to study the \emph{Poincar\'e series} $P_f(T)$, which is defined as the generating power series
\[P_f(T) = \sum_{d\in \Z_{\geq 0}} N_d T^d.\]

Another, closely related function is the \emph{$p$-adic zeta function} $Z_{f}(s)$, which also encodes the information of all $N_d$ for $d\in \Z_{\geq 0}$.

In \cite{BorevichShafarevich}, Borevich and Shafarevich conjectured that $P_f(T)$ is a rational function. In 1974, Igusa gave a proof of this fact in \cite{Igusa74} and \cite{Igusa75}, by proving that $Z_f(s)$ is a rational function in $p^{-s}$, which implies that $P_f(T)$ is rational too. The rationality of these functions implies that there is a meromorphic continuation to $\C$, and therefore we can talk about poles of $P_f(T)$ and $Z_f(s)$. The poles of $Z_f(s)$ give valuable information on the asymptotic behaviour of $N_d$ when $d$ tends to infinity, see for example \cite{Segers} for a nice explanation.

With the \emph{$p$-adic monodromy conjecture}, Igusa formulated a fascinating conjecture, linking the arithmetical concept of the $p$-adic zeta function to something of more differential-topological nature: the local monodromy eigenvalues of the analytic function $f\colon \C\to \C$. Monodromy eigenvalues give valuable information on the singularities of the hypersurface in $\C^m$ defined by $f=0$. The $p$-adic monodromy conjecture states that every pole of the $p$-adic zeta function corresponds to a local monodromy eigenvalues.

The $p$-adic monodromy conjecture has been proven in some significant cases, for example for polynomials in two variables \cite{Loeser88}, and for homogeneous polynomials in three variables \cite{Rodrigues-Veys01} \cite{Artal-Cassou-et-al}, but the general case remains wide open.

\subsubsection*{Motivic zeta functions and the motivic monodromy conjecture for hypersurface singularities}

Inspired by Kontsevich's theory of motivic integration, Denef and Loeser introduced in \cite{DenefLoeser98} the \emph{motivic zeta function} $Z_f^{mot}(T)$ associated with a polynomial $f\in k[x_1, \ldots, x_m]$. This function is some kind of `super $p$-adic zeta function', because it captures the $p$-adic zeta functions $Z_f(s)$, for almost all primes $p$. 
It is a more geometrical invariant of the polynomial $f$ than the $p$-adic zeta function. 

The motivic zeta function $Z_f^{mot}(T)$ is a formal power series with coefficients in a certain Grothendieck ring of varieties. The Denef-Loeser formula, proved in \cite{DenLoe01}, shows that the motivic zeta function is in fact a rational function. This formula makes it even possible to explicitly compute the motivic zeta function  from the data of an embedded resolution of singularities of the hypersurface in $\A_k^m$ defined by $f=0$.

The $p$-adic monodromy conjecture got a motivic upgrade. This conjecture is called Denef and Loeser's \emph{motivic monodromy conjecture for hypersurface singularities}, and states that every pole of the motivic zeta function corresponds to a local monodromy eigenvalue.

A proof of the motivic monodromy conjecture would immediately imply a proof of the $p$-adic monodromy conjecture. And so far, proofs of specific cases of the $p$-adic monodromy conjecture seem to generalize to the motivic monodromy conjecture without major adaptations.

For more information on the $p$-adic and motivic monodromy conjecture, we refer to Section~\ref{sect:p-adic-motivic-monodromy}.

\subsubsection*{Motivic zeta functions for Calabi-Yau varieties}

For every pair $(X,\omega)$, where $X$ is a Calabi-Yau variety over $K$ and $\omega$ a volume form on $X$, Halle and Nicaise formulated in \cite{HalleNicaiseAbelian} the invariant $Z_{X,\omega}(T)$, called the \emph{motivic zeta function} associated with $(X,\omega)$. It is defined as a formal power series
\[Z_{X,\omega}(T) = \sum_{d\in \Z_{>0}} N_d T^d,\]
where $N_d$ is an element in a certain localized Grothendieck ring of varieties. When $K=\C(\!(t)\!)$, the coefficient~$N_d$ measures the space of $\C(\!(\sqrt[d]{t})\!)$-rational points on $X$.
The motivic zeta function $Z_{X,\omega}(T)$ is a natural analog of Denef and Loeser's motivic zeta function for hypersurface singularities. 

The Denef-Loeser formula, which allows to compute $Z_f^{mot}$ in terms of the data of an embedded resolution of singularies of $f$, has an analog in the context of Calabi-Yau varieties. Nicaise and Sebag in \cite{NicaiseSebag_motivic_serre}, and Bultot and Nicaise in \cite{BultotNicaise}, proved a formula for $Z_{X,\omega}(T)$ in terms of the degeneration of $X$ at $t=0$, or more precisely, in terms of an $snc$-model of $X$. An $snc$-model of $X$ is a regular, proper $R$-scheme $\X$ such that $\X\times_R K \simeq X$ and such that the special fiber $\X_k = \X \times_R k$ is a strict normal crossings divisor on $\X$. This formula immediately implies that $Z_{X,\omega}(T)$ is a rational function. Moreover, from the data of an $snc$-model of $X$, it is possible to explicitly compute the motivic zeta function $Z_{X,\omega}(T)$ with this formula.

As in the context of hypersurface singularities, we are interested in the poles of $Z_{X,\omega}(T)$. So far, we have results for abelian varieties \cite{HalleNicaiseAbelian}, and for Calabi-Yau varieties with a so-called equivariant Kulikov model \cite{HalleNicaiseKulikov}. In both cases, the motivic zeta function $Z_{X,\omega}(T)$ has a unique pole. In this thesis, we will see that the motivic zeta function of a $K3$ surface with a triple-point-free model can have more than one pole.

\subsubsection*{The monodromy property for Calabi-Yau varieties}

Of course, we wonder whether there is an analog of the motivic monodromy conjecture for hypersurface singularities in the context of Calabi-Yau varieties as well. A Calabi-Yau variety is said to satisfy the \emph{monodromy property}, if every pole of the motivic zeta function $Z_{X,\omega}(T)$ corresponds to a monodromy eigenvalue. Halle and Nicaise proved that abelian varieties, and Calabi-Yau varieties with an equivariant Kulikov model, satisfy the monodromy property, in \cite{HalleNicaiseAbelian} and \cite{HalleNicaiseKulikov} respectively. In this thesis, we will see more examples of Calabi-Yau varieties that satisfy the monodromy property.

For more information on the monodromy property, we refer to Chapter~\ref{ch:GMP}.

\section*{The second ingredient: the Crauder-Morrison classification}

In \cite{CrauMor}, Crauder and Morrison study triple-point-free $snc$-models of proper, smooth surfaces with trivial pluricanonical bundle. These models are strict normal crossings models such that three irreducible components of the special fiber never meet simultaneously. They show that such a model can always be birationally modified to a so-called Crauder-Morrison model, which is classified in one of the following three classes: 
\begin{itemize}
\item flowerpot degenerations,
\item cycle degenerations, and
\item chain degenerations. 
\end{itemize}

The Crauder-Morrison classification specifies `building blocks' of which the special fiber of a Crauder-Morrison model is made: flowers, flowerpots, cycles and chains. They also classify the flowers into 21 classes. Below, we find a picture of the dual graph of a possible chain degeneration with three flowers.

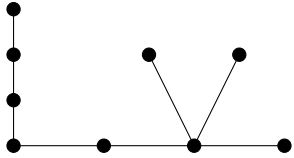
\begin{center}
\begin{tikzpicture}[scale=0.6]

\node [draw,shape=circle, fill, inner sep=0pt,minimum size=5pt] at (0,0) {};

\node [draw,shape=circle, fill, inner sep=0pt,minimum size=5pt] at (2,0) {};

\node [draw,shape=circle, fill, inner sep=0pt,minimum size=5pt] at (4,0) {};

\node [draw,shape=circle, fill, inner sep=0pt,minimum size=5pt] at (6,0) {};

\node [draw,shape=circle, fill, inner sep=0pt,minimum size=5pt] at (0,1) {};

\node [draw,shape=circle, fill, inner sep=0pt,minimum size=5pt] at (0,2) {};

\node [draw,shape=circle, fill, inner sep=0pt,minimum size=5pt] at (0,3) {};

\node [draw,shape=circle, fill, inner sep=0pt,minimum size=5pt] at (3,2) {};

\node [draw,shape=circle, fill, inner sep=0pt,minimum size=5pt] at (5,2) {};

\draw (0,0) -- (6,0)
(0,0) -- (0,3)
(4,0) -- (3, 2)
(4,0) -- (5, 2)
;
\end{tikzpicture}
\end{center}

For more information on the Crauder-Morrison classification, we refer to Chapter~\ref{ch:CM} of this thesis.

%
%

\section*{The monodromy property for $K3$ surfaces allowing a triple-point-free model}

The Denef-Loeser formula for Calabi-Yau varieties makes it possible to explicitly compute the motivic zeta function $Z_{X,\omega}(T)$ in terms of the data of an $snc$-model of the Calabi-Yau variety $X$. This makes a classification of $snc$-models like the one of Crauder and Morrison so useful: it gives us plenty of information that can be used to verify the monodromy property in practice.
Therefore, we are interested in Calabi-Yau varieties to which we can apply the Crauder-Morrison classification. These are exactly the abelian surfaces allowing a triple-point-free model, and $K3$ surfaces allowing a triple-point-free model. 

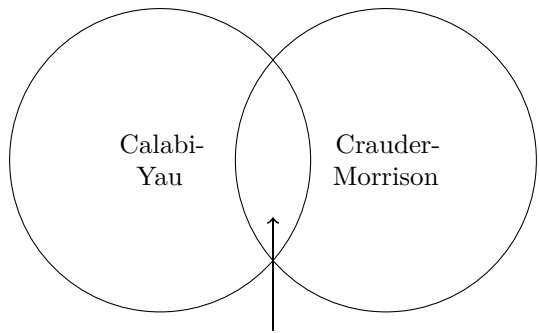
\def\firstcircle{(0,0) circle (2cm)}
\def\secondcircle{(3cm,0) circle (2cm)}

\begin{center}
\begin{tikzpicture}
\draw \firstcircle node [text width = 1.5 cm, align = center]{Calabi-Yau};
\draw \secondcircle node [text width = 1.5 cm, align = center] {Crauder-Morrison};

\draw [->, line width = 0.75pt](1.5cm,-2.25cm) to (1.5cm,-0.75cm); 

\node at (1.75cm,-3.25cm) {\parbox{5.75cm}{\begin{itemize} \item abelian surfaces with a\\ triple-point-free model \item $K3$ surfaces with a\\ triple-point-free model \end{itemize}}};

\end{tikzpicture}
\end{center}

We already know that abelian surfaces satisfy the monodromy property by \cite{HalleNicaiseAbelian}, so we will focus on $K3$ surfaces allowing a triple-point-free model.

Since the Crauder-Morrison classification is formulated for more general surfaces, we investigate whether the classification can be refined for $K3$ surfaces. Our most important result is the following theorem.

\textbf{Theorem~\ref{thm:CM-K3-no-cycle}.} \emph{Let $X$ be a $K3$ surface over $K$ and let $\X$ be a Crauder-Morrison model of $X$. Then $\X$ is either a flowerpot degeneration or a chain degeneration, but not a cycle degeneration.}

The Crauder-Morrison classification and the description of the building blocks are very precise, which allows us to use a combination of geometrical and combinatorial arguments to verify the monodromy property. The first major result is the computation of the poles of the motivic zeta function $Z_{X,\omega}(T)$ for a $K3$ surface $X$ allowing a triple-point-free model and a volume form $\omega$ on $X$.

\textbf{Theorem~\ref{thm:poles}.}
\emph{Let $X$ be a $K3$ surface over $K$.  Let $\X$ be a Crauder-Morrison model of $X$ over $R$ with special fiber $\X_k=\sum_{i\in I} N_i E_i$. Let $\omega$ be a volume form on $X$ and let $(N_i, \nu_i)$ be the numerical data of $E_i$, for every $i\in I$.}

\emph{The rational number $q\in \Q$ is a pole of $Z_{X,\omega}(T)$ if and only if there exists an element $i\in I$ with $q=-\nu_i/N_i$ and such that 
\begin{enumerate}[(i)]
\item either $\rho_i$ is minimal,
\item or $E_i$ is the top of a conic-flower.
\end{enumerate}
Moreover, in case (i), $q$ is a pole of order 1 if $\X$ is a flowerpot degeneration, and of order 2 if $\X$ is a chain degeneration. In case (ii), $q$ is a pole of order 1.}

For the definition of $\nu_i, \rho_i$, and a conic-flower, we refer to the main text of this thesis.

In particular, we see that if the Crauder-Morrison model $\X$ has a conic-flower, then $Z_{X,\omega}(T)$ has more than one pole. This is in contrast to abelian varieties and Calabi-Yau varieties with an equivariant Kulikov-model, because these kinds of Calabi-Yau varieties have a motivic zeta function with a unique pole.

To verify the monodromy property for $K3$ surfaces allowing a triple-point-free model, we need to check whether every pole of the motivic zeta function $Z_{X,\omega}(T)$ induces a monodromy eigenvalue. The Crauder-Morrison classification also serves as a tool to compute monodromy eigenvalues, and we find the following theorem.

\textbf{Theorem~\ref{thm:holds-pot}.} \emph{Let $X$ be a $K3$ surface over $K$ with Crauder-Morrison model~$\X$. If $\X$ is a flowerpot degeneration, then $X$ satisfies the monodromy property.} 

For chain degenerations, the situation is more complex. However, we have a partial result:

\textbf{Theorem~\ref{thm:holds-chain}.} \emph{Let $X$ be a $K3$ surface over $K$ with Crauder-Morrison model~$\X$. If $\X$ is a chain degeneration satisfying some extra conditions, the $K3$ surface~$X$ satisfies the monodromy property.}

Our results are summarized below.

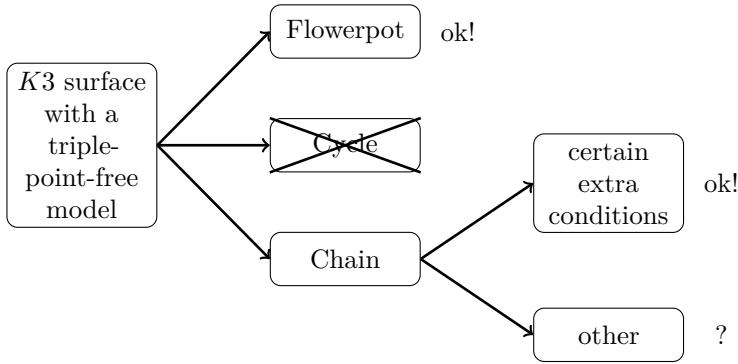
\tikzstyle{block} = [rectangle, draw, 
    text width=5em, text centered, rounded corners, minimum height=2em]
    \begin{center}
\begin{tikzpicture}[node distance = 1.5cm]
\node[block](K3){$K3$ surface with a triple-point-free model};

\node[block, right of=K3, node distance = 3.5cm](cycle){Cycle}; 
\node[block, above of =cycle](flowerpot){Flowerpot};  
\node[block, below of = cycle](chain){Chain}; 
\draw[->, line width = 1pt] (K3.east) to (cycle.west);
\draw[->, line width = 1pt] (K3.east) to (flowerpot.west);
\draw[->, line width = 1pt] (K3.east) to (chain.west);

\draw [line width = 1pt] (cycle.north east) to (cycle.south west);
\draw [line width = 1pt] (cycle.north west) to (cycle.south east);

\node[right of = flowerpot, node distance = 1.5cm]{{ok!}};

\node[right of=chain, node distance = 3.5cm](ch-mult){}; 
\node[block, above of= ch-mult, node distance = 1cm](ch-rational){certain extra conditions}; 
\node[block, below of=ch-mult, node distance = 1 cm](ch-rest){other};

\draw[->, line width = 1pt] (chain.east) to (ch-rest.west);
\draw[->, line width = 1pt] (chain.east) to (ch-rational.west);

\node[right of = ch-rational, node distance = 1.5cm]{{ok!}};
\node[right of = ch-rest, node distance = 1.5cm]{{?}};  
\end{tikzpicture}
\end{center}

We don't know whether all $K3$ surfaces with a chain degeneration satisfy the monodromy property. The difficulty is that we don't know which special fibers can occur as a chain degeneration of a $K3$ surface. In the last chapter, we will produce a list of 63 `combinatorial countercandidates'. These are combinatorial descriptions of special fibers such that if there exists a $K3$ surface $X$ over $K$ with an $snc$-model $\X$ such that the special fiber satisfies the properties of one of these 63 combinatorial countercandidates, then $X$ does not satisfy the monodromy property. Conversely, if there exists a $K3$ surface $X$ over $K$ allowing a triple-point-free model that does not satisfy the monodromy property, then the special fiber $\X_k$ of its Crauder-Morrison model $\X$ satisfies the properties of one of these 63 combinatorial countercandidates.

\section*{Outline of the thesis}

Chapter 1. \emph{Preliminaries}

We start by fixing notation and terminology, and we explain the minimal background needed for this thesis.

Chapter 2. \emph{The monodromy property for Calabi-Yau varieties}

In the second chapter, we give an introduction to the monodromy property for Calabi-Yau varieties. First, we will motivate the subject by giving some historical context: we will discuss the $p$-adic and motivic monodromy conjecture for hypersurface singularities. Then we will introduce the two main ingredients of the monodromy property, namely the motivic zeta function, and monodromy eigenvalues. We will formulate the monodromy property and give the results known so far.

\vspace{1cm}
Chapter 3. \emph{The Crauder-Morrison classification}

In the third chapter, we will give the Crauder-Morrison classification of triple-point-free $snc$-models of smooth, proper surfaces with trivial pluricanonical bundle. This classification is one of the major tools in this thesis. We will also refine the classification specifically for $K3$ surfaces.

Chapter 4. \emph{Poles of the motivic zeta function}

This chapter is devoted to the first main theorem in this thesis: for a $K3$ surface~$X$ allowing a triple-point-free model and a volume form $\omega$ on $X$, we find a set $S^\dagger\subset \Z\times\Z_{>0}$ such that the motivic zeta function $Z_{X,\omega}(T)$ is an element of $\mathcal{M}_k^{\hat{\mu}}\left[T, \frac{1}{1-\L^a T^b}\right]_{(a,b)\in S^\dagger}$. Moreover, we prove that if $(a,b)\in S^\dagger$, then $a/b$ is a pole of the motivic zeta function, and therefore no element of $S^\dagger$ can be omitted. The main tools in the proof are the Denef-Loeser formula and the Crauder-Morrison classification. We explicitly compute the contribution of a flower to the motivic zeta function by writing Python code, which can be found in Appendix~\ref{ch:appendix1}.

Chapter 5. \emph{$K3$ surfaces satisfying the monodromy property}

The main result in this chapter is `$K3$ surfaces with a Crauder-Morrison model that is a flowerpot degeneration satisfy the monodromy property'. For most flowerpots, this is not very difficult to prove, but flowerpots that are rational, non-minimal ruled surfaces require some more work. We also prove that $K3$ surfaces with a Crauder-Morrison model that is a chain degeneration satisfy the monodromy property, under certain extra conditions.

Chapter 6. \emph{Future research: a proof or a counterexample?}

The last chapter is a guideline for further research and we present a variety of possibly useful techniques. We introduce the concept of combinatorial countercandidate: a combinatorial description of a special fiber that could lead to a counterexample of the monodromy property. We generate a list of 63 combinatorial countercandidates by writing Python code, which can be found in Appendix~\ref{app:combinatorial-countercandidate}.

\cleardoublepage


\chapter{Preliminaries}\label{ch:preliminaries}

The aim of this chapter is twofold: we will fix notation used in this thesis, and we will briefly introduce the minimal background needed. For the reader who would like to get a more in-depth introduction, we will give some references.
\section{Notation and conventions}

For any field $k$, a \emph{$k$-variety} is a reduced, separated $k$-scheme of finite type. A \emph{surface} over $k$ is a projective, geometrically connected $k$-scheme of dimension 2.

Unless stated otherwise, $k$ will be an algebraically closed field of characteristic~0 and we set $K=k(\!( t )\!)$ and $R=k\llbracket t \rrbracket$. We fix an algebraic closure $K^{alg}$ of $K$.

For any integer $d\geq 1$, the notation $\mu_d$ is used for the \emph{group of $d$-th roots of unity} in $k$. The \emph{profinite group of roots of unity} in $k$ will be denoted by
\[\hat{\mu}=\varprojlim \mu_d.\]

For $D$ and $D'$ two divisors on an integral, regular, projective scheme, we write $D\sim D'$ if $D$ and $D'$ are \emph{linearly equivalent} and $D\equiv D'$ if $D$ and $D'$ are\emph{ numerically equivalent}. For a regular, projective scheme $\X$ over $R$, we write $\omega_{\X/R}$ for the \emph{relative canonical sheaf} of $\X$ over $R$ and $K_{\X/R}$ for a divisor associated with $\omega_{\X/R}$. Note that $K_{\X/R}$ is only defined up to linear equivalence.

A \emph{ruled surface} is a surface $X$, together with a surjective morphism $\pi\colon X\to C$ to a non-singular curve $C$, such that all fibers are connected and that for all but finitely many points $y\in C$, the fiber $X_y$ is isomorphic to $\P^1$.
We say $X$ is \emph{rational or elliptic ruled} if $C$ is a rational or elliptic curve respectively.
A ruled surface is said to be a \emph{minimal ruled surface}, if all fibers are isomorphic to $\P^1$. Notice that what is called a ruled surface in \cite[Section~V.2]{Hartshorne} and \cite[Section~V.4]{BHPV}, is called a minimal ruled surface here. It can be proven that a smooth, projective surface $X$ is a ruled surface if and only if there exists a composition of blow-ups $X\to X'$, where $X'$ is a minimal ruled surface.
As a consequence, if $X$ is a ruled surface with reducible fiber $\fiber$, then any irreducible component of $\fiber$ is rational with negative self-intersection.

\section{Calabi-Yau varieties and \texorpdfstring{$K3$}{K3} surfaces }

\begin{definition}
A \emph{Calabi-Yau variety} is a smooth, proper, geometrically connected variety with trivial canonical sheaf.
\end{definition}

Contrary to what is often done, we do not require that $h^{i,0}(X)$ vanishes for $0<i<\dim X$ in the definition of a Calabi-Yau variety $X$.

A well-known class of examples of Calabi-Yau varieties are the abelian varieties. In this thesis, we will focus on $K3$ surfaces.

\begin{definition}
A \emph{$K3$ surface} $X$ is a 2-dimensional Calabi-Yau variety with $H^1(X,\O_X)=0$.
\end{definition}

The theory on $K3$ surfaces is vast and we will just state some very basic properties.

\begin{proposition} \label{thm:pre-K3-cohomology}
Let $X$ be a $K3$ surface, the dimension of the $\ell$-adic cohomology is
\[\dim\, H^m(X, \Q_\ell)=
\begin{cases}
1 & \text{for } m=0, 4,\\
22 & \text{for } m=2,\\
0 & \text{otherwise.}\\
\end{cases}
\]

\end{proposition}
For a proof of this proposition, we refer to \cite[Section~1.3.3]{HuybrechtsK3}.

\begin{corollary}
Let $X$ be a $K3$ surface. Then $\chi(X)=24$, where $\chi$ denotes the topological Euler characteristic.
\end{corollary}

\begin{example}
Let $f\in k[x,y,z,w]$ be a homogeneous polynomial of degree~4. If the surface $X$ defined by $f=0$ in $\P^3_k$ is smooth, then $X$ is a $K3$ surface.
\end{example}

We refer to \cite{HuybrechtsK3} for more information on $K3$ surfaces.

\section{Models and the Kulikov classification} \label{sect:models}

For a smooth, proper $K$-scheme $X$, a \emph{model} of $X$ is a flat $R$-scheme~$\X$ endowed with an isomorphism $\X\times_R K \simeq X$.
The base change $\X_k=\X\times_R k$ is called the \emph{special fiber} of the model $\X$.

 A model $\X$ is called an \emph{$snc$-model}, if it is regular and proper and if $\X_k=\sum_{i\in I} N_iE_i$ is a strict normal crossings divisor. Because we assumed that $k$ has characteristic~0, such a model always exists by Nagata's compactification theorem \cite{Nagata1}\cite{Nagata2} and Hironaka's resolution of singularities \cite{Hironaka}.
 
An $snc$-model $\X$ with a reduced special fiber is said to be \emph{semi-stable}. By Mumford's semi-stable reduction theorem \cite[Section~4.3]{SemistableReductionTheorem}, we have that for every proper $K$-scheme $X$, there exists a finite extension $L$ of $K$ such that $X\times_K L$ has a semi-stable model.

An $snc$-model $\X$ of $X$ is said to be \emph{triple-point-free} if no three distinct irreducible components in the special fiber intersect. For such a model, we can define the \emph{dual graph} $\Gamma$ in the following way. Write $\X_k=\sum_{i\in I} N_i E_i$. The vertices of the graph $\Gamma$ are $\{v_i\}_{i\in I}$, and for every irreducible component of $E_i\cap E_j$, where $i\neq j$, there is an edge between $v_i$ and $v_j$. Note that for $snc$-models that are not triple-point-free, it is possible to define the dual complex of $\X$, see for example \cite[Definition~12]{DualComplex}. When $\X$ is triple-point-free, the dual complex and the dual graph coincide.
 


\subsubsection*{The Kulikov classification}

Semi-stable $snc$-models of $K3$ surfaces are completely classified by Kulikov in \cite{Kulikov1} and Persson and Pinkham \cite{Persson-Pinkham}.

\begin{theorem} \label{thm:pre-Kulikov}
Let $X$ be a $K3$ surface over $K$ and let $\X$ be a semi-stable $snc$-model of $X$ over $R$. There exists a birational modification $\widetilde{\X}$ of $\X$ such that $\widetilde{\X}$ is a semi-stable $snc$-model of $X$ with $K_{\widetilde{\X}/R}$ trivial. 
Moreover, one of the following properties hold.
\begin{enumerate}[(i)]
\item The special fiber $\widetilde{\X}_k$ is smooth.
\item The special fiber $\widetilde{\X}_k$ is a chain of surfaces $V_0 + V_1 + \cdots + V_k + V_{k+1}$. The surfaces $V_1, \ldots, V_k$ are elliptic ruled, and $V_0$ and $V_{k+1}$ are rational surfaces. The intersection curves are elliptic.
\item The special fiber $\widetilde{\X}_k$ is a union of rational surfaces, and its dual complex is a triangulated 2-sphere.

\end{enumerate}
\end{theorem}

The semi-stable $snc$-model $\widetilde{\X}$ of $X$ from the previous theorem, is called the \emph{Kulikov model} of $X$

To construct $\widetilde{X}$ from $\X$, one has to run the Minimal Model Program on $\X$ and then resolve the singularities by means of small resolutions. These resolutions force us to allow $\widetilde{\X}$ to be an algebraic space.

Note that for every $K3$ surface $X$, there exists a finite extension $L$ of $K$ such that $X \times_K L$ has a Kulikov model, because of the semi-stable reduction theorem.

\section{The Grothendieck ring of varieties}

The Grothendieck ring of varieties is one of the most important concepts in this thesis. The motivic zeta function, one of the major ingredients in the monodromy property, is defined over a localization of this ring. But it is also an interesting topic in itself and a lot of research is still going on to investigate this mysterious object the Grothendieck ring still is. For a nice survey on this subject and some open questions, we refer to \cite{NicaiseSebag_Grothendieck_ring} and \cite[Chapter~1]{MotivicIntegration}

\subsection{The Grothendieck ring of varieties}

Let $k$ be any field.

\begin{definition}
The \emph{Grothendieck ring of $k$-varieties} $K_0(Var_k)$ is the abelian group generated by the isomorphism classes $[X]$ of separated $k$-schemes $X$ of finite type, modulo the scissor relations
\[[X]=[Y]+[X\setminus Y],\]
for any closed subscheme $Y$ of $X$. The group $K_0(Var_k)$ is endowed with a ring structure by considering the unique multiplication such that
\[[X]\cdot [X'] = [X\times_k X']\]
for every pair $(X, X')$ of separated $k$-schemes of finite type.
\end{definition}

We immediately see that $[\emptyset]$ is the zero element of $K_0(Var_k)$ and $[\Spec (k)]$ is the unit element.
We introduce the symbol $\L$, which stands for the class of $\A^1_k$. We denote by $\mathcal{M}_k$ the localized Grothendieck ring $K_0(Var_k)[\L^{-1}]$.

\begin{remark}
The scissor relations tell us that $[X]=[X_{red}]$ for every separated $k$-scheme $X$ of finite type, where $X_{red}$ is the maximal reduced closed subscheme of $X$. This means that if we replace `separated $k$-scheme of finite type' by `$k$-variety' in the definition of Grothendieck ring, we actually define the same ring. But in that case, one should be careful with the definition of the ring multiplication: when $k$ is imperfect and $X$ and $X'$ are $k$-varieties, the product $X\times_k X'$ is not necessarily reduced.
\end{remark}

The relation $[X]=[Y]+[X\setminus Y]$ allows to cut and paste varieties. If we cut a variety $X$ into two pieces $Y$ and $X\setminus Y$, adding (pasting) the classes of these two pieces in the Grothendieck ring, gives the class of the original variety $X$ again. This is why we call this relation the scissor relation.

The scissor relation can be very helpful to compute the class of a variety in practice. For instance, the projective line $\P^1_k$ has a subvariety isomorphic to $\A^1_k$ and the complement of this subvariety is a point. This means that
\[[\P^1_k]=\L+1.\]
Another example: the projective plane $\P^2_k$ has a subvariety isomorphic to the affine plane $\A_k^2$ and its complement is a projective line. Since $\A^2_k\simeq \A^1_k\times_k \A^1_k$, we have
\[[\P^2_k]=\L^2+\L+1.\]
By induction, one finds that
\[[\P^n_k]=\L^n+\cdots + \L+1.\]

\begin{remark} A subset $U$ of a $k$-variety $X$ is called \emph{locally closed} if it is open in its closure in $X$. A locally closed subset $U$ has a unique structure of subvariety of $X$. We call a subset $C$ of $X$ a \emph{constructible} subset, when it is the finite union of locally closed subsets of $X$. Although a constructible subset doesn't necessarily have a unique structure of subvariety of $X$, its class in the Grothendieck ring of varieties is well defined. Indeed, if we write $C$ as a finite disjoint union of locally closed subsets $U_1, \ldots, U_n$, then the class
\[[C]=\sum_{i=1}^n [U_i]\]
does not depend on the choice of partition into locally closed subsets.  
\end{remark} 

\note{A constructible subset is not necessarily a subscheme. As it is the disjoint union of locally closed subschemes, it is a scheme, but it is not a subscheme. For example in $\A^2$, when we take the constructible subset that is the plane minus a coordinate axis plus the origin, then this is not a subscheme of $\A^2$.}

Because of the scissor relations, taking the class $[X]$ of a variety $X$ in the Grothendieck ring satisfies all the requirements a reasonable \emph{measure} should have: if we cut a variety into two subvarieties, their sizes should add up to the original size. Moreover, the size of a product variety should be the product of the sizes of its factors. 
So every reasonable way to measure the size of a variety~$X$ factors through its class $[X]$ in the Grothendieck ring of varieties. Or said even differently, the measure $X\to [X]$ specializes to any other (more intuitive) measure.  That's why the Grothendieck ring of $k$-varieties should be thought of as the \emph{universal measure}.

Let's look at some examples of such specialization morphisms.

\begin{example}\label{ex:specialization-Grothendieck-ring}
\hfill
\begin{enumerate}[(i)]
\item Fix a finite field $k=\mathbb{F}_q$ and consider the invariant $\sharp$ that counts the number of $\mathbb{F}_q$-rational points on a variety over $\mathbb{F}_q$. Since counting points is clearly additive and multiplicative, we have a well-defined morphism
\[\sharp\colon  K_0(Var_{\mathbb{F}_q})\to \Z\colon  [X]\mapsto \sharp\, X(\mathbb{F}_q).\]
This means in particular that two varieties having the same class in $K_0(Var_{\mathbb{F}_q})$ have the same number of rational points.
This morphism localizes to a ring morphism
$\sharp\colon \mathcal{M}_{\mathbb{F}_q}\to \Z[q^{-1}],$
because 
\[\sharp\,\L= \sharp \A^1_{\mathbb{F}_q}(\mathbb{F}_q) = q.\]

\item For $k=\C$, an important additive, multiplicative invariant of a complex variety $X$ is the compactly supported Euler characteristic $\chi_c$ of $X(\C)$ with respect to the complex topology. This means that there exists a well-defined morphism
\[\chi_c\colon K_0(Var_k)\to \Z\colon  [X]\mapsto \chi_c( X).\]

For complex varieties, we have that $\chi_c(X)=\chi(X)$, where $\chi(X)$ is the topological Euler characteristic of $X(\C)$ with respect to the complex topology. This was proven by Laumon in \cite{Laumon}. Therefore, there also exists a well-defined morphism
\[\chi\colon K_0(Var_k)\to \Z\colon  [X]\mapsto \chi( X).\]

This last morphism, for example, allows to easily compute the Euler characteristic of projective space:
\[\chi(\P^n_\C) = \chi(\L^n+\cdots + \L+1) = \chi(\L)^n+\cdots+\chi(\L)+1 = n+1.\]
The morphism $\chi$ localizes to a ring morphism $\chi\colon \mathcal{M}_\C\to \Z$, since $\chi(\L)=1$.

For an arbitrary field $k$, one can consider the $\ell$-adic Euler characteristic $\chi$, with $\ell$ an invertible prime in $k$. One can show that $\chi$ does not depend on the choice of $\ell$, see for example \cite[point~(2a) p.~28]{Katz-ell-adic-cohomology}. 

\item For $k$ a field of characteristic 0, there exists a unique ring morphism $P\colon K_0(Var_k)\to \Z[v]$ such that, when $X$ is a smooth and proper $k$-variety, the class $[X]$ gets mapped to the Poincar\'e polynomial
\[P(X;v) = \sum_{i\geq 0} (-1)^i b_i(X)\, v^i,\]
where $b_i(X)= \dim H^i(X,\Q_\ell)$ is the $i$-th $\ell$-adic Betti number of $X$. 
\note{There is unique such ring morphism because $K_0(Var_k)$ is generated by all smooth and proper varieties by Hironaka's resolution of singularies.}

In particular, if a class in the Grothendieck ring has a non-zero Poincar\'e polynomial, then the class is non-zero itself. But we also have some converse result. 
In \cite[Proposition~8.7]{NicaiseTrace}, it is proven that for any $k$-variety $X$, the degree of $P(X;v)$ is $2\dim(X)$ and its leading coefficient equals the number of irreducible components of $X\times k^{alg}$ of dimension $\dim(X)$. Therefore, a (non-empty) variety has a non-zero Poincar\'e polynomial. 

Notice that $P(X;1)=\chi(X)$.
Furthermore, we have $P(\L)=v^2$ and hence, there is a well-defined localized ring morphism
\[P\colon \mathcal{M}_k\to \Z[v,v^{-1}].\]

The Poincar\'e polynomial can also be defined over a field $k$ of positive characteristic $p$. For any smooth, proper $k$-variety $X$, it still holds that
\[P(X;v) = \sum_{i\geq 0} (-1)^i b_i(X)\, v^i,\]
where $b_i(X)= \dim H^i(X,\Q_\ell)$ is the $i$-th $\ell$-adic Betti number of $X$ with $\ell$ a prime invertible in $k$. However, we don't know whether $P\colon K_0(Var_k)\to \Z$ is uniquely defined by this property, as is the case for characteristic 0.
For more information on the Poincar\'e polynomial and how to define it for fields of positive characteristic, we refer to \cite[Section 8]{NicaiseTrace}. \note{We can define the Poincar\'e polynomial also for fields of positive characteristic. When $X$ is smooth and proper, we still have \[P(X;v) = \sum_{i\geq 0} (-1)^i b_i(X) v^i.\] But we don't know whether this uniquely characterizes $P$ on $K_0(Var_k)$, unless we assume resolution of singularities in positive characteristic.}

\item For $k$ a field of characteristic zero, there exists a unique ring morphism $HD\colon K_0(Var_k)\to \Z[u,v]$ such that for a smooth and proper $k$-variety $X$, the class $[X]$ gets mapped to the Hodge-Deligne polynomial
\[HD(X;u,v) = \sum_{p,q\geq 0} (-1)^{p+q} h^{p,q}(X)\, u^p v^q,\]
where $h^{p,q}(X) = \dim_k H^p(X, \Omega^q_X)$ is the $(p,q)$-th Hodge number of $X$. Notice that $HD(X;v,v) = P(X;v)$. Since $HD(\L; u,v) = uv$, we can localize this morphism to
$HD\colon \mathcal{M}_k \to \Z[u,u^{-1},v,v^{-1}]$.

The Hodge-Deligne polynomial of a complex variety $X$ equals
\[HD(X;u,v) = \sum_{p,q\geq 0} \sum_{i\geq 0} (-1)^i h^{p,q}(H^i_c(X(\C),\C)) u^p v^q,\]
where $h^{p,q}(H^i_c(X(\C),\C))$ is the dimension of the $(p,q)$-component of Deligne's mixed Hodge structure on $H^i_c(X(\C),\C)$.


\end{enumerate}
\end{example}

The Grothendieck ring $K_0(Var_k)$ is still a mysterious object and a lot of researchers are investigating the properties of this ring. 
In 2002, Poonen proved in \cite{Poonen02} that $K_0(Var_k)$ is not a domain, if $k$ has characteristic 0. We even know that the class of $\L$ is a zero-divisor in $K_0(Var_\C)$, as shown by Borisov in \cite{Borisov15}. Moreover, when $k$ has characteristic zero, $\mathcal{M}_k$ is not a domain either. This can be proven with the same argument as Ekedahl uses in \cite{Ekedahl}. Cauwbergs writes this down in more detail in \cite{Cauwbergs}.

\subsection{The {equivariant} Grothendieck ring of varieties}

The equivariant Grothendieck ring of varieties $K_0^G(Var_k)$ is a variant of the Grothendieck ring of varieties, where we take into account a group action of a fixed finite group $G$ on the varieties. This subsection is based on \cite[Sections~3~and~4]{Hartmann_quotient}, and we refer to this paper for more details.

\subsubsection{Notation}
Let $k$ be a field and fix a finite group $G$. Let $X$ be a separated scheme of finite type over $k$ endowed with a group action of $G$. Unless explicitely stated otherwise, we assume that groups act on schemes from the left. We say the group action is \emph{good}, if every orbit of this action is contained in an affine open subscheme of $X$. We denote by $(\Sch_{k,G})$ the category whose objects are separated schemes of finite type over $k$ with a good $G$-action, and whose morphisms are $G$-equivariant morphisms of $k$-schemes. One can check that the fiber product exists in this category, see for example \cite[Section~2]{Hartmann_quotient}. 

\subsubsection{Affine bundles}
\begin{definition}
Let $S$ be a $k$-scheme. An \emph{affine bundle over $S$ of rank $d$} is an $S$-scheme $A$ with a vector bundle $V\to S$ of rank $d$ and a morphism of $S$-schemes $\phi\colon V\times_S A\to A$ such that $\phi\times p_A\colon  V\times_S A \to A\times_S A$ is an isomorphism of $S$-schemes, where $p_A$ is the projection to $A$.
\end{definition}

\begin{remark} \label{rmk:affine-bundle}
The inspiration behind the definition of affine bundle is that it should be a generalization of the following fact in linear algebra: assume $A$ is an affine space of dimension $n$ over $k$ and $V=k^n$. We think of $A$ as $k^n$ where we `forget' the origin. Then we can look at the map
\[V\times A\to A\colon (v,x)\to x+v,\]
which is just adding a vector to a point as in high-school geometry. This map induces a bijection
\[V\times A\to A\times A\colon (v,x)\to (x+v,x),\]
which represents the easy fact that for every two points $x$ and $x'$ in $A$, the vector~$x-x'$ exists in $V$. We see that $A$, together with $V$, satisfies the definition of an affine bundle.
\end{remark}

\begin{remark}
It can be shown that every affine bundle of rank $d$ over $S$ is a locally trivial fibration over $S$ with fiber $\A_k^d$. 

\note{$A$ is a torsor under $V$. In Milne's book on \'etale cohomology, p.76, it is proven that the torsors (=principal homogenious spaces) under $V$ are parametrized by the first (flat) cohomology of $V$. Since Zariski and flat cohomology coincide for $V$, I think the $H^1$ is trivial, so there is just one torsor, so it is a locally trivial fibration.}
\end{remark}

If there is an action on an affine bundle $A$, ideally, the action on $A$ should be compatible with the action on $V$. 
\begin{definition}
The affine bundle $A$ over $S$ is said to be \emph{$G$-equivariant}, if $A$ and $S$ are in $(\Sch_{k,G})$ and $A\to S$ is $G$-equivariant. We call the $G$-action on $A\to S$ \emph{affine} if there is a $G$-action on $V$, linear over the action on $S$, such that $\phi\colon V\times_S A\to A$ is $G$-equivariant.
\end{definition}

\begin{remark}
In the special case we discussed in Remark \ref{rmk:affine-bundle}, a $G$-action on $A$ is affine if $g\cdot (x+v) = g\cdot x + g \cdot v$ for every $g\in G$, $x\in A$ and $v\in V$.
\end{remark}

\subsubsection{The equivariant Grothendieck ring \texorpdfstring{$K_0^G(Var_k)$}{}}

\begin{definition}
The \emph{ equivariant Grothendieck ring of $k$-varieties} $K_0^G(Var_k)$ is the abelian group generated by the isomorphism classes $[X]$ of objects $X\in (Sch_{k,G})$, modulo the following relations:
\begin{enumerate}
\item $[X]=[Y]+[X\setminus Y]$, whenever $Y$ is a closed $G$-equivariant subscheme of $X$ (scissor relations), and
\item $[A]=[S \times_k \A^d_k]$, whenever $A\to S$ is a $G$-equivariant affine bundle of rank~$d$ over $S\in (\Sch_{k,G})$ with affine $G$-action and with trivial $G$-action on $\A^d_k$. 
\end{enumerate}
We endow $K_0^G(Var_k)$ with a ring structure by putting the unique multiplication such that
\[[X]\cdot [X'] = [X\times_k X']\]
for all $X, X'\in (\Sch_{k,G})$ and where the fiber product is taken in the category~$(\Sch_{k,G})$.
\end{definition}

Let $\L\in K_0^G(Var_k)$ be the class of the affine line with trivial action.
In particular, we see that, if $A\to S$ is a $G$-equivariant affine bundle of rank $d$ over $S$, then $[A]=[S]\L^d$.
 We denote by $\mathcal{M}^G_k$ the localization $K_0^G(Var_k)[\L^{-1}]$.

\begin{remark}
When $G=\{e\}$ is the trivial group, we have that the equivariant Grothendieck ring $K_0^G(Var_k)$ equals the usual Grothendieck ring~$K_0(Var_k)$. \note{Relation 2. becomes indeed trivial. (think about this)} Moreover, $\mathcal{M}_k^G$ equals $\mathcal{M}_k$ in this case.
\end{remark}

\begin{remark}
A morphism of finite groups $G'\to G$ induces ring morphisms
\[K_0^G(Var_k)\to K_0^{G'}(Var_k) \text{ and } \mathcal{M}_{k}^G\to \mathcal{M}_{k}^{G'}.\]
\end{remark}

\begin{remark}
In the literature, other definitions of the equivariant Grothendieck ring are often used. We refer to \cite[Remark~4.5, Remark 4.6]{Hartmann_quotient} for a discussion.
\end{remark}

In this thesis, we will mostly be interested in actions of the profinite group $\hat{\mu}=\varprojlim \mu_d$ of roots of unity of the field $k$. The ring $K_0^{\hat{\mu}}(Var_k)$ can be defined as follows.

\begin{definition}
Let $\widehat{G}=\varprojlim_{i\in I} G_i$ be a profinite group, where all $G_i$ are finite groups. Then we define
\[K_0^{\widehat{G}}(Var_k) = \varinjlim_{i\in I} K_0^{G_i}(Var_k) \quad\text{ and }\quad \mathcal{M}_k^{\widehat{G}} = \varinjlim_{i\in I} \mathcal{M}_k^{G_i}.\]
\end{definition}

\section{Motivic integration on Calabi-Yau varieties}
\label{sect:motivic_integration}


Let $X$ be a Calabi-Yau variety over $K$, and let $\omega$ be a volume form on $X$, i.e., a nowhere-vanishing differential form of maximal rank.
For every integer~$d\geq 1$, set $R(d)=R[\pi]/(\pi^d-t)$ and $K(d)$ as the fraction field of $R(d)$. The field~$K(d)$ is the unique totally ramified extension of $K$ in $K^{alg}$ of degree $d$. We define
\[X(d)=X\times_{K} K(d),\]
and $\omega(d)$ as the pullback of $\omega$ to $X(d)$. There is a left group-action of $\mu_d$ on $X(d)$. 
 Let $\Y$ be an equivariant, weak N\'eron model of $X(d)$. This means that $\Y$ is a separated, smooth $R(d)$-scheme, endowed with a good $\mu_d$-action and a $\mu_d$-equivariant isomorphism $\Y\times_{R(d)} K(d) \simeq X(d)$, such that the natural map $\Y(R(d))\to X(K(d))$ is a bijection. It can be proven that such an equivariant, weak N\'eron model always exists. 
 
 For every connected component $C$ of $\Y_k$, we denote by $\ord_C\omega(d)$ the unique integer $a$ such that $t^{-a/d}\omega(d)$ is a generator of $\omega_{\Y/R(d)}$, locally at the generic point of $C$.
\begin{definition}
Let $X$ be a Calabi-Yau variety over $K$, and let $\omega$ be a volume form on $X$. Fix an integer $d\geq 1$ and let $\Y$ be an equivariant, weak N\'eron model of $X(d)$.
The \emph{motivic integral} of $\omega(d)$ on $X(d)$ is defined as
\[\int_{X(d)} |\omega(d)| = \sum_{C\in \pi_0(\Y_k)} [C]\, \L^{-\ord_C \omega(d)} \in \mathcal{M}_k^{\hat{\mu}},\]
where $\pi_0(\Y_k)$ is the set of connected components of $\Y_k$.
\end{definition}

It is non-trivial to prove that this definition does not depend on the choice of equivariant, weak N\'eron model $\Y$. Without the $\hat{\mu}$-action, this is done by Loeser and Sebag in \cite[Proposition~4.3.1]{LoeserSebag}. The generalization to the equivariant version follows from the change of variables formula for equivariant motivic integrals, proven by Hartmann in \cite{Hartmann_motivic}.

\note{
We omit the factor $\L^{-m}$ from the definition given in \cite{HalleNicaiseAbelian}. This is done because we also omit the factor $\L^m$ in the definition of motivic zeta function and hence they still define the same motivic zeta function.  The omission of $\L^m$ and $\L^{-m}$ was also done in \cite{BultotNicaise} and \cite{HalleNicaiseKulikov}.
}

\note{
The integral $\int_{X(d)}|\omega(d)|$ is a motivic integral measuring the space of $K(d)$-rational points on $X$. 
It can be computed using a weak N\'eron model of $X(d)$ as explained in Section~\ref{sect:motivic_integration}.
}

\section{Algebraic spaces}

In Chapter~\ref{ch:CM}, we allow models $\X$ of smooth, proper $K$-varieties to be \emph{algebraic spaces}. An algebraic space is a generalization of the notion of a scheme. For the reader who is not familiar with the concept of algebraic space, we refer to \cite{Knutson} and \cite[Tag 025R, Tag 03BO and Tag 03H8]{stacks-project}.

The reason why we need to work with algebraic spaces is that in the category of algebraic spaces, more contraction results hold. In particular, we refer to the work of Artin in \cite[Theorem~6.2]{Artin}. 

Let $X$ be a smooth, proper surface over $K$ with $\omega_{X/K}^{\otimes m}\simeq \O_K$ for some integer $m\geq 1$. In Chapter~\ref{ch:CM}, we define a Crauder-Morrison model $\X$ of $X$ as an algebraic space over $R$ such that $\X$ is a triple-point-free $snc$-model of $X$, where no component in the special fiber $\X_k$ can be contracted such that the result is again a triple-point-free $snc$-model. It is essential that these contractions are considered in the category of \emph{algebraic spaces}. If we consider only contractions in the category of schemes, then the classification of Crauder-Morrison models in Theorem~\ref{thm:CM-classification} does not necessarily hold.

An algebraic space of finite type over a field $k$ has a well-defined class in the Grothendieck ring of $k$-varieties $K_0(Var_k)$, see for example \cite[(7.1.3)]{HalleNicaiseKulikov}. The theory of motivic integration can also be extended to algebraic spaces, as done by Halle and Nicaise in \cite[Section~7]{HalleNicaiseKulikov}.

\cleardoublepage


\chapter{The monodromy property for Calabi-Yau varieties}\label{ch:GMP}

In this introductory chapter, we will explain the monodromy property for Calabi-Yau varieties. In the first section, we will motivate the subject by giving the historical context: we will discuss the $p$-adic and motivic monodromy conjecture for hypersurface singularities. In the second and third section, we will define the two main ingredients in the monodromy property, namely the motivic zeta function and monodromy eigenvalues. In the last section, we will be ready to define the monodromy property. We will also discuss what is known about varieties satisfying the monodromy property.
\section{History: the \texorpdfstring{$\lowercase{p}$}{p}-adic and motivic monodromy conjecture} \label{sect:p-adic-motivic-monodromy}

Before we start defining the main concepts of this thesis, let us first look at some history. The monodromy property for Calabi-Yau varieties is a variant of the motivic monodromy conjecture for hypersurface singularities, which, in its turn, is a generalization of the $p$-adic monodromy conjecture. To motivate and better understand the monodromy property for Calabi-Yau varieties, it is useful to first explain the $p$-adic and motivic monodromy conjecture for hypersurface singularities. For a nice introduction and more details regarding these conjectures, we refer to \cite{NicaisePadicMotivicMonConj} and \cite[Chapter~0]{Bart-thesis}. It is important to keep in mind that there are no direct implications between the motivic monodromy conjecture for hypersurface singularities and the monodromy property for Calabi-Yau varieties: a proof for the motivic monodromy conjecture for hypersurface singularities would not imply that the monodromy property holds for certain types of Calabi-Yau varieties and vice versa.

\subsection{The \texorpdfstring{$p$}{p}-adic monodromy conjecture} \label{sect:p-adic_monodromy}

\subsubsection{The \texorpdfstring{$p$}{p}-adic zeta function and the Poincar\'e series}

Let $f$ be a non-constant polynomial in $\Z[x_1, \ldots, x_m]$ for some integer $m$, and fix a prime number $p$. It is a very natural question to ask how many solutions the congruence $f\equiv 0 \mod p^{d+1}$ has for $d\geq 0$. Let's denote by $N_d$ the number of such solutions, so 
\[N_d = \sharp \left\{a\in (\Z/p^{d+1})^m \mid f(a)\equiv 0 \mod p^{d+1}\right\}.\]

If we use these numbers as coefficients in a generating series,
we get the Poincar\'e series of $f$. 

\begin{definition}
Let $p$ be a prime number and let $f\in \Z[x_1,\ldots, x_m]\setminus \Z$ be a non-constant polynomial, for some $m\in \Z_{>0}$. Let $N_d$ be the number of solutions~$x\in(\Z/p^{d+1})^m$  of $f(x)\equiv 0\mod p^{d+1}$, for any integer $d\geq 0$. The \emph{Poincar\'e series of $f$} is the generating series
\[P_f(T)=\sum_{d\geq 0} N_d T^d \in \Z\,\llbracket T\rrbracket .\]
\end{definition}

\note{Note the difference in definition of the Poincar\'e polynomial in \cite{NicaisePadicMotivicMonConj} and \cite{HalleNicaiseMotivicCY}. The relation is $Q(f;T) = T\cdot P(T)$. We use the definition of \cite{HalleNicaiseMotivicCY}. 
}

\begin{remark}
This power series has a radius of convergence $R\geq 1/p^m$, since $N_d\leq p^{(d+1)m}$. \end{remark}
\note{radius of convergence $= 1/\limsup_{d\to +\infty} \sqrt[d]{|c_d|}$ and $\limsup_{d\to +\infty} \sqrt[d]{|c_d|} \leq \limsup_{d\to +\infty} p^{(d+1)m/d}=p^m$.}

The Poincar\'e series of a polynomial $f$ is closely related to the $p$-adic zeta function of $f$. 

\begin{definition}
Let $| \cdot |_p$ denote the $p$-adic absolute value on $\Q_p$. The \emph{$p$-adic zeta function of $f$} is defined by the $p$-adic integral
\[Z_f(s) = \int_{\Z_p^n}| f(x)|_p^s\, d\mu,\]
for every $s\in \C$ with $\textrm{Re}(s)>0$ and where $\mu$ is the usual Haar measure on $\Z_p^n$.
\end{definition}

By definition of $p$-adic integration, the function $Z_f(s)$ can be written as a power series in $p^{-s}$:
\begin{eqnarray} \label{eq:padic-zeta-function}
Z_f(s) &= \sum_{d\geq 0} \mu\left(\left\{a\in \Z_p^m\mid v_p(f(a)) = d\right\}\right)\, p^{-ds},
\end{eqnarray}

with $v_p$ the $p$-adic valuation. Readers not familiar with $p$-adic integration, can take \eqref{eq:padic-zeta-function} as the definition of the $p$-adic zeta function.

\begin{remark}
The Poincar\'e series $P_f(T)$ and the $p$-adic zeta function $Z_f(s)$ are closely related and their relation
\begin{equation}
P_f(p^{-n}T) = \frac{p^n\left(1-Z_f(s)\right)}{1-T}, \label{eq:relation-poinc-series-zetafunction}
\end{equation}
can be shown by direct computation, where $T=p^{-s}$. From this, we see that $Z_f$ encodes the number of solutions of $f$ modulo powers of $p$ as well.
\end{remark}

In problem 9 of Section~1.5 of \cite{BorevichShafarevich}, Borevich and Shafarevich asked whether the Poincar\'e series is a rational function. This question got an affirmative answer in 1974, when Igusa proved that the $p$-adic zeta function of $f$ is a rational function in $p^{-s}$ by using Hironaka's resolution of singularities. The original proof can be found in the papers \cite{Igusa74} and \cite{Igusa75}, and a simplified and more readable proof can be found in the appendix of \cite{Igusa77}. The rationality of $Z_f$ immediately implies the rationality of $P_f$ by formula \eqref{eq:relation-poinc-series-zetafunction}. In \cite[Theorem 3.2]{Denef84}, Denef gave a completely different proof of the rationality of the $p$-adic zeta function using cell decomposition.

The rationality of $Z_f(s)$ and $P_f(T)$ implies that they admit meromorphic continuations to $\C$. This means that we can talk about poles of $Z_f(s)$ and~$P_f(T)$.
From formula \eqref{eq:relation-poinc-series-zetafunction}, we immediately see that the poles of $P_f$ can be computed from $Z_f$ and vice versa. \note{Note that in formula \eqref{eq:relation-poinc-series-zetafunction} the variables are different ($T$ versus $s$) and that there is a factor $T-1$, so a pole $1$ disappears.} This means that knowing the poles of $Z_f$ gives information on the asymptotic behaviour of $N_d$ when $d$ tends to infinity. We refer to \cite{Segers} for a nice explanation of this fact.

\note{From Manu's thesis: "Rationality of the series $\sum_{l=0}^\infty a_l T^l$ is equivalent to the sequence $(a_l)_l$ being a linear recurrence sequence."}

Igusa's proof of the rationality of $Z_f(s)$ even gives a finite set containing the real parts of all the poles of $Z_f(s)$. 
Let $g\colon Y\to \A_{\Q_p}^m$ be an embedded resolution of singularities of the zero locus $f^{-1}(0)\subset \A_{\Q_p}^m$. If $\alpha$ is a pole of $Z_f(s)$, then $Re(\alpha)$ is of the form $-\nu_i/N_i$, where $N_i$ and $\nu_i-1$ are the multiplicities
of $f\circ g$ and $g^*dx$ along an irreducible component~$E_i$ of $g^{-1}(f^{-1}(0))$. It is remarkable that a lot of irreducible components of $g^{-1}(f^{-1}(0))$ do not contribute to the poles of $Z_f(s)$, a phenomenon that is intimately connected with the $p$-adic monodromy conjecture.

\subsubsection{Monodromy eigenvalues}

Let $f\colon \C^m\to \C$ be a non-constant analytic map. We denote by $X_0$ the analytic space defined by $f=0$ and let $x$ be any point of $X_0$. We would like to study the topology of $X_0$ in a neighbourhood of $x$, especially when $X_0$ is not smooth at $x$. 
The idea is the following: locally around $x$, we will study the analytic space $X_{\delta}$ defined by $f=\delta$ instead of $f=0$, for complex numbers $\delta$ close enough to 0. For such `small' $\delta$, the analytic space $X_\delta$ is smooth and all these spaces~$X_\delta$ are diffeomorphic for every choice of $\delta$. When we let the space $X_\delta$ turn once around the space $X_0$, we get an automorphism of $X_\delta$, which is called the monodromy transformation. This automorphism gives valuable information about the singularity $x$ of $X_0$. Let us now make this idea more concrete.

Let $B=B(x,\varepsilon)$ be an open ball around $x$ in $\C^m$ with radius $\varepsilon$ and let $D=D(0,\eta)$ be an open disc around 0 in $\C$ with radius $\eta$. By $D^*$, we denote the punctured disc $D\setminus\{0\}$. For a suitable choice of $0<\eta \ll \varepsilon \ll 1$, we have that the map
\[f_x\colon f^{-1}(D^*)\cap B \to D^*\]
is a locally trivial fibration. 
\note{
This means that all fibers are smooth and diffeomorphic to each other and that for every $t\in D^*$, there is an open neighborhood $U_t$ aroung $t$ such that $f_x^{-1}(U_t)$ is isomorphic to $U_t\times X$, where $X$ is a fiber. According to Bart Bories' thesis, this is a theorem of Hamm, L\^{e} and Milnor. The fact that all fibers are smooth for a suitable choice for $\delta$ can be understood as follows: we prove that $f$ has only a finite number of singular fibers. Suppose $f$ has a singularity in $x$, then the gradient of $f$ at $x$ is zero. Let $A$ be a connected component of the singular locus of $f$. Then $\nabla f=0$ on $A$ and hence $f$ is constant on $A$. This means that $A$ is a part of a fiber of $f$. Since the singular locus has only a finite number of connected components, there are only a finite number of singular fibers of $f$.
}
This fibration is called the \emph{Milnor fibration} of $f$ at $x$. The fact that the Milnor fibration is locally trivial implies in particular that all fibers of $f_x$ are diffeomorphic to one another. The fiber of $f_x$ over any point of $D^*$ is called the \emph{Milnor fiber} $F_x$ of $f$ at $x$, which is well defined up to diffeomorphism. 

To study the topological structure of $X_0$ around $x$, we let the Milnor fiber $F_x$ turn once around $X_0$ and we look at the transformation $F_x$ underwent. More formally, we choose a loop $\gamma$ in $D^*$ around the origin. Over each point of this loop, we have a fiber that is diffeomorphic to the Milnor fiber. Since $f_x$ is locally trivial, we can move the Milnor fiber in a continuous way through open neighbourhoods where $f_x$ is trivial. This results in an automorphism  of $F_x$, which we call the \emph{monodromy transformation} $M_x$. It can be shown that, up to homotopy, the monodromy transformation is well defined, i.e., it does not depend on the chosen loop $\gamma$ and the chosen trivialization of $f_x$.

We can describe the monodromy transformation even more formally, which will turn out to be useful when we will give a more algebraic definition of monodromy in Section~\ref{sect:MonodromyEigenvalues}.
First, notice that the Milnor fiber is homotopy-equivalent to
\[\left(f^{-1}(D^*)\cap B\right) \times_{D^*} \widetilde{D^*},\]
where $\widetilde{D^*}$ is the universal covering space of $D^*$. This is true because $\widetilde{D^*}$ is contractible. The group $\pi_1(D^*)\simeq \Z$ of covering transformations of $\widetilde{D^*}$ over $D^*$, acts naturally on~$\widetilde{D^*}$ and this action can be lifted canonically to $(f^{-1}(D^*)\cap B)\times_{D^*} \widetilde{D^*}$. The action of the canonical generator of $\pi_1(D^*)$, namely a positively oriented loop turning once around the origin, is homotopic to the monodromy transformation $M_x$. 

The monodromy transformation $M_x$ induces linear maps
\[M_{x,i}\colon  H^i(F_x,\C)\to H^i(F_x,\C)\]
on the singular cohomology spaces $H^i(F_x,\C)$ for all $i\geq 0$. These are uniquely determined, since $M_x$ is well defined up to homotopy. 

\begin{definition} \label{def:monodromy-eigenvalue}
Let $f\colon \C^m\to \C$ be a non-constant analytic map. We denote by $X_0$ the analytic space defined by $f=0$ and let $x$ be any point of $X_0$. Let $F_x=\left(f^{-1}(D^*)\cap B\right) \times_{D^*} \widetilde{D^*}$ be the Milnor fiber of $f$ at $x$ and let $M_{x,i}$ be the action of the canonical generator of $\pi_1(D^*,\eta)$ on $H^i(F_x,\C)$.
We say that a complex number $\alpha\in \C$ is a \emph{monodromy eigenvalue of $f$ at $x$}, if $\alpha$ is an eigenvalue of $M_{x,i}$ for some $i\geq 0$.
\end{definition}

All monodromy eigenvalues are roots of unity \cite[Th\'eor\`eme de mono\-dromie~2.1]{SGA7}.
An important tool to compute monodromy eigenvalues in practice is the monodromy zeta function of $f$ at $x$.

\begin{definition} \label{def:monodromy-zeta-p-adic}
The \emph{monodromy zeta function} of $f$ at $x$ is the alternating product of the characteristic polynomials of the monodromy transformations $M_{x,i}$ on the cohomology spaces $H^i(F_x, \C)$:
\[\zeta_{f,x}(T) = \prod_{i\geq 0} \det \left(Id - T\cdot M_{x,i}\mid H^i(F_x,\C)\right)^{(-1)^{i+1}}.\]
\end{definition}

\begin{remark}
Zeroes and poles of $\zeta_{f,x}$ are monodromy eigenvalues of $f$ at $x$. A priori, not every monodromy eigenvalue of $f$ at $x$ is necessarily a zero or pole of $\zeta_{f,x}$, since cancellations may occur in the definition of $\zeta_{f,x}$. Denef proved in \cite[Lemma~4.6]{Denef93} that every monodromy eigenvalue of $f$ at $x$ is a zero or pole of $\zeta_{f,y}$ for some, possibly different, point $y\in X_0$.
\end{remark}

\note{The monodromy zeta function is defined as $\zeta_{f,x}(T) = \prod_{i\geq 0} \det (Id - T\cdot M_{x,i})^{(-1)^{i+1}}$, but classically, the characteristic polynomial of $M_{x,i}$ should be $\det(T\cdot Id- M_{x,i})$. So to get the alternating product of the characteristic polynomials, we actually need the substitution $T\to 1/T$ and the zeros and poles of $\zeta_{f,x}$ are actually inverses of the monodromy eigenvalues. But $\zeta_{f,x}$ is a rational function with coefficients in $\Q$ (A'Campo), so zeros and poles occur in conjugate pairs. Since monodromy eigenvalues are roots of unity, the inverse is exactly its complex conjugate and we conclude that zeros and poles of $\zeta_{f,x}$ are monodromy eigenvalues.}

As already mentioned before, the monodromy zeta function can be useful to compute monodromy eigenvalues in practice. A'Campo gave a particularly nice formula for the monodromy zeta function in terms of an embedded resolution of singularities $g$ of $f$.  

\begin{theorem}[{\cite[Th\'eor\`eme~3]{Acampo}}]
Let $f\colon \C^m\to \C$ be a non-constant analytic map. We denote by $X_0$ the analytic space defined by $f=0$ and let $x$ be any point of $X_0$.
Let $g\colon Y\to \A^m_\C$ be an embedded resolution of singularities of~$f$ and denote by $\{E_i\}_{i\in I}$ the irreducible components of $g^{-1}(X_0)$. Set $N_i$ as the multiplicity of $E_i$ and $E_i^\circ = E_i\setminus \bigl(\bigcup_{j\in I\setminus\{i\}} E_j\bigr)$.
The monodromy zeta function can be expressed as
\[\zeta_{f,x}(T) = \prod_{i\in I} \left(1-T^{N_i}\right)^{-\chi\left(E_i^\circ \cap g^{-1}(x)\right)},\]
where $\chi(\cdot)$ is the topological Euler characteristic.
\end{theorem}

\note{
The monodromy transformation and the Milnor fiber give information on the singularity of $X_0$ in $x$. For example, the analytic Milnor fiber gives $\widehat{\O}_{X_0,x}$ and hence the exact singularity of $X_0$ in $x$. In dimension 1, the topological Milnor fiber also gives the exact singularity of $X_0$ in $x$, up to topological equivalence.
}

\subsubsection{The $p$-adic monodromy conjecture}

The monodromy conjecture is a fascinating conjecture relating the poles of the $p$-adic zeta function of a polynomial $f$ - an arithmetical concept - to the monodromy eigenvalues of $f$ considered as an analytic function - a concept that is differential-topological in nature.

\begin{conjecture}[The $p$-adic monodromy conjecture]
Let $f$ be a non-constant polynomial in $\Z[x_1, \ldots, x_m]$ for some integer $m$. For almost all primes $p$, we have the following: if $\alpha$ is a pole of the $p$-adic zeta function $Z_f(s)$, then $\exp(2\pi i \textrm{Re}(\alpha))$ is a monodromy eigenvalue of $f$ at some point $x\in X_0$.
\end{conjecture}

The $p$-adic monodromy conjecture has been proven, for example, for polynomials in two variables in \cite{Loeser88} and homogeneous polynomials in three variables in \cite{Rodrigues-Veys01} and \cite{Artal-Cassou-et-al}. However, the general case remains wide open. 

\subsection{The motivic monodromy conjecture} \label{sect:motivic_monodromy}

In 1995, Kontsevich launched the idea of motivic integration, which can be thought of as a geometrized version of $p$-adic integration. The key idea is to replace the finite field $\Z/p\Z$ by $\C$ and the $p$-adic integers $\Z_p$ by the field of formal complex power series $\C\,\llbracket t\rrbracket $. Instead of using the Haar measure as in $p$-adic integration, for motivic integration, the size of a variety is measured by its class in the Grothendieck ring. Kontsevich's goal was to prove that two birationally equivalent complex Calabi-Yau varieties have the same Hodge numbers. 

Kontsevich never published his ideas, but the theory of motivic integration was further developed by other mathematicians, such as Cluckers, Denef, Nicaise, Loeser and Sebag. In \cite{DenefLoeser98}, Denef and Loeser introduce the motivic zeta function associated with a polynomial $f\in k[x_1,\ldots, x_m]$ and the motivic monodromy conjecture. A more modern version of the motivic zeta function appears in \cite{DenLoe01}, and it is this version that we will use in this thesis. The motivic zeta function can be thought of as the \emph{universal} zeta function, since it specializes to the $p$-adic zeta function for almost all primes $p$. This means that if a statement is proven for the motivic zeta function, we immediately have a proof for the $p$-adic equivalent as well.

\subsubsection{The motivic zeta function}

Let $k$ by any field. We will now formulate a motivic counterpart of the $p$-adic zeta function ${Z_f(s) = \int_{\Z_p^n}| f(x)|_p^s\, d\mu}$.  

\begin{definition} \label{def:motivic-zeta-function}
Let $f\in k[x_1,\ldots, x_m]$ be a polynomial. The \emph{motivic zeta function} of $f$ is
\begin{align*}
Z_f^{mot}(T)
&= \sum_{d\geq 0} [\X_{d,1}]\, \L^{-md}\,T^d  \in \mathcal{M}_{k}\,\llbracket T\rrbracket,
\end{align*}
where $\X_{d,1} = \left\{\varphi \in \left(k[t]/(t^{d+1})\right)^m \mid f(\varphi)\equiv t^d \mod (t^{d+1})\right\}$ for every~${d\in \Z_{\geq 0}}$.
\end{definition}
\note{
In \cite{DenLoe01}, the motivic zeta function with angular component is actually only defined equivariantly (with the $\hat{\mu}$-action). And without the extra factor $\L^{-m}$.
}

Let us now have a look at how the $p$-adic and motivic zeta functions are related. The $p$-adic zeta function compares best with the \emph{naive motivic zeta function}, which is a slightly different version than the definition we gave before.
Remember that using the definition of $p$-adic integration, we could write the $p$-adic zeta function as
\[Z_f(s)=\sum_{d\geq 0} \mu\left(\left\{a\in \Z_p^m\mid v_p(f(a)) = d\right\}\right)\, p^{-ds}.\]

The key-idea of motivic integration is to replace $\Z_p$ by $k\,\llbracket t\rrbracket $, so the set~${\left\{a\in \Z_p^m\mid v_p(f(a)) = d\right\}}$ changes to
\[\left\{\varphi \in k\,\llbracket t\rrbracket ^m \mid ord_t(f(\varphi))=d\right\},\]
where $ord_t(\sum_{i\geq 0} a_it^i)=\min \{i\geq 0 \mid a_i\neq 0\}$. Compare this to the $p$-adic valuation, which is defined as $v_p(\sum_{i\geq 0}a_ip^i)=\min \{i\geq 0 \mid a_i\neq 0\}$. Note that $ord_t(f(\varphi))=d$ if and only if $f(\varphi) \equiv a t^d \mod (t^{d+1})$  for some $a\in k \setminus \{0\}$.
 
In the motivic setting, the role of $p$ will be played by the class of the affine line~$\L$ and $p^{-s}$ is replaced by $T=\L^{-s}$. Putting all these ideas together, the {naive motivic zeta function} associated with $f$ should be
\[Z_f^{mot, naive}(T)=\sum_{d\geq 0} \mu\left(\left\{\varphi \in k\,\llbracket t\rrbracket ^m \mid ord_t(f(\varphi))=d\right\}\right) T^d \in \mathcal{M}_{k}\,\llbracket T\rrbracket ,\] 
where $\mu$ is the \emph{motivic measure}. To define this motivic measure, we need, for every integer $d\geq 0$, the \emph{truncation map} $\pi_d$:
\[\pi_d\colon \left(k\,\llbracket t\rrbracket \right)^m\to \left(k[t]/(t^{d+1})\right)^m \colon  \Big(\sum_{j\geq 0} a_{ij}t^j\Big)_i\mapsto \Big(\sum_{j=0}^{d} a_{ij}t^j\Big)_i.\]
A subset $A\subset (k\,\llbracket t\rrbracket )^m$ is called \emph{cylindric}, if there exists a constructible subset~${C\subset \left(k[t]/(t^{d+1})\right)^m}$, when $\left(k[t]/(t^{d+1})\right)^m$ is viewed as $k^{(d+1)m}$, such that $A=\pi_d^{-1}(C)$, for some integer $d$. The motivic measure $\mu$ assigns a `size' to such cylindric subsets:
\[\mu(A) := [C]\,\L^{-m(d+1)} \in \mathcal{M}_k.\]
It can be easily shown that $\mu(A)$ is independent of the choice of $d$.

This is why the naive motivic zeta function can be rewritten to
\[Z^{mot, naive}_f(T) = \L^{-m}\sum_{d\geq 0} [\X_d]\, \L^{-md}\,T^d,\]
where 
\[\X_d = \left\{\varphi \in \left(k[t]/(t^{d+1})\right)^m \mid f(\varphi) \equiv a t^d \mod (t^{d+1}) \text{ for some } a\in k\setminus\{0\} \right\}.\]

The factor $\L^{-m}$ is often omitted, like we also did in Definition~\ref{def:motivic-zeta-function}. Unfortunately, by discarding $\L^{-m}$, the psychological link with the $p$-adic zeta function gets weakened.

\begin{remark}
The reader who is familiar with motivic integration, will recognize that the naive motivic zeta function can also be defined as the motivic integral
\[Z_f^{mot, naive}(s)=\int_{\mathcal{L}(\A^m_k)}\L^{-ord_t(f) \,  s}.\]
\end{remark}

It can often be difficult to compute the motivic zeta function by using just the definition. Inspired by the $p$-adic case, Denef and Loeser \cite{DenLoe01} proved an amazing formula for $Z_f^{mot}$ in terms of an embedded resolution of singularities for $f$, for fields $k$ with $\text{char}(k)=0$. This formula turns out to be extremely useful in practice.

\begin{theorem} \label{thm:Denef-Loeser}
Assume $k$ is a field of characteristic 0 and let $f\in k[x_1,\ldots, x_m]$ be a polynomial.
Let $g\colon Y\to \A_k^m$ be an embedded resolution of singularities for $f$ and let $\{E_i \mid i\in I\}$ be the irreducible components of the divisor $g^{-1}(f^{-1}(0))$. Define $N_i$ and $\nu_i-1$ to be the multiplicities of $f\circ g$ and~$g^*dx$ along the component $E_i$. \note{where $dx$ is a local non-vanishing volume form at any point of $g(E_i)$, i.e. a local generator of the sheaf of differential forms of maximal degree.} 
\note{
Manu defines $\nu_i$ by the equality $K_{Y/X}=\sum_{i\in I}(\nu_i-1)E_i$.
}
For any non-empty subset $J\subseteq I$, set $E_J=\cap_{j\in J} E_j$ and $E_J^\circ=E_J\setminus \left(\cup_{i\in I\setminus J} E_i\right)$. The motivic zeta function can be written as
\[Z_f^{mot}(T) = \sum_{J\subset I} (\L-1)^{|J|-1}[\widetilde{E_J^\circ}]\prod_{j\in J}\frac{\L^{-\nu_j}T^{N_j}}{1-\L^{-\nu_j}T^{N_j}},\]
where $[\widetilde{E_J^\circ}]$ is a certain Galois cover of $E_J^{\circ}$ of degree $\gcd_{j\in J}N_j$, and is defined in \cite[Section~3.3]{DenLoe01}.
\end{theorem}
\note{Let $m_J = gcd(N_j)_{j\in J}$. We introduce an unramified Galois cover $\widetilde{E_J^\circ}$ of $E^{\circ}_J$, with Galois group $\mu_{m_J}$, as follows. Let $U$ be an affine Zariski open subset of $Y$, such that, on $U$, $f \circ g = uv^{m_J}$, with $u$ a unit on $U$ and $v$ a morphism from $U$ to $\A_k^1$. Then the restriction of $\widetilde{E_J^\circ}$ above $E_J^\circ \cap U$, denoted by $\widetilde{E_J^\circ} \cap U$, is defined as
\[\{(z, y) \in  \A^1_k \times (E^\circ_J \cap U) \mid z^{m_J} = u ^{-1}\}.\]
Note that $E^\circ_J$ can be covered by such affine open subsets $U$ of $Y$. Gluing together the covers $\widetilde{E_J^\circ} \cap U$, in the obvious way, we obtain the cover $\widetilde{E_J^\circ}$ of $E^\circ_J$ which has a natural $\mu_{m_J}$-action (obtained by multiplying the $z$-coordinate with the elements
of $\mu_{m_J}$). This $\mu_{m_J}$-action on $\widetilde{E_J^\circ}$ induces a $\hat{\mu}$-action on $\widetilde{E_J^\circ}$ in the obvious way.
}

In particular, this formula implies that $Z_f^{mot}$ is rational over $\mathcal{M}_k$. More precisely, 
\[Z_f^{mot}(T)\in \mathcal{M}_k\left[T,\frac{1}{1-\L^{-\nu_i}T^{N_i}}\right]_{i\in I}\subset \mathcal{M}_k\,\llbracket T\rrbracket .\]

\subsubsection{The motivic monodromy conjecture}

Since the motivic zeta function specializes to the $p$-adic zeta function for almost all $p$, one hopes that a motivic upgrade of the $p$-adic monodromy conjecture exists as well. One has to be careful though when formulating the motivic monodromy conjecture since $\mathcal{M}_k$ is not a domain and hence the concept of pole is a priori not defined. We will discuss this issue more in-depth in Section~\ref{sect:def-poles}. 
Denef and Loeser formulated the motivic monodromy conjecture as follows:

\begin{conjecture}[The motivic monodromy conjecture for hypersurface singularities]
Let $k$ be a subfield of $\C$ and
let $f\in k[x_1, \ldots, x_m]$ be a non-constant polynomial. Set $X_0$ to be the complex analytic space defined by the equation $f=0$ in $\C^m$. There exists a finite set ${S}\subset \Z_{<0}\times \Z_{>0}$ such that
\[Z_f^{mot}(T) \in \mathcal{M}_k\left[T, \frac{1}{1-\L^{a}T^{b}}\right]_{(a,b)\in {S}},\]
and such that for each $(a,b)\in {S}$, the complex number $\exp(2\pi i a/b)$ is a monodromy eigenvalue of $f$ at some point of $X_0$.
\end{conjecture}

\note{One can drop the condition that $k$ is a subfield of $\C$ by using $\ell$-adic nearby cycles to define the notion of local monodromy eigenvalue. This does not yield a more general conjecture, since by the Lefschetz principle, one can always reduce to the case where $k$ is a subfield of $\C$.}

 So far, it seems that almost all known partial proofs of the $p$-adic monodromy conjecture could be naturally adapted to the motivic setting. Since the motivic zeta function specializes to the $p$-adic zeta function for almost all primes $p$, researchers try to solve this more general conjecture instead of the $p$-adic one.

\section{The motivic zeta function for Calabi-Yau varieties} \label{sect:motivic-zeta-CY}

Fix an algebraically closed field $k$ of characteristic zero. Define $R=k\llbracket t\rrbracket$ and $K=k(\!( t) \!)$. Fix an algebraic closure $K^{alg}$ of $K$.

\note{
We could define $R$ to be any complete discrete valuation ring, $K$ its fraction field and $k$ its residue field. We assume $k$ to be of characteristic zero. But this definition is not more general than the one we gave since there is an isomorphism $R\simeq k\llbracket  t \rrbracket$, since $R$ has equal characteristic zero. However, this morphism is not canonical and depends on the choice of uniformizer and the choice of $k$-algebra structure $k\to R$.
}

\note{
In principle, we don't need $R$ to be complete, Henselian should be enough. We want $R$ to be Henselian since we want to know all the finite extensions of $K$. But asking $R$ to be complete doesn't really matter since we can always take the completion and the motivic zeta function will not change.
}

Let $X$ be a Calabi-Yau variety over $K$ of dimension $m$. 
Let $\omega$ be a volume form on $X$, i.e., a nowhere-vanishing differential form of maximal rank.
For every integer $d\geq 1$, define $K(d)$ as the unique totally ramified extension of $K$ in $K^{alg}$ of degree $d$. \note{When $k$ is not algebraically closed, such a totally ramified extension of $K$ is not necessarily unique.} \note{If we don't fix an algebraic closure $K^{alg}$ of $K$, the extension $K(d)$ is uniquely defined modulo isomorphism. But this isomorphism is not canonical. If we fix the algebraic closure, then everything is uniquely defined.} For example, when $K=\C(\!(t)\!)$, we have $K(d) = \C(\!(\sqrt[d]{t})\!)$. We define
\[X(d)=X\times_{K} K(d),\]
and $\omega(d)$ as the pullback of $\omega$ to $X(d)$.

We start by defining the motivic zeta function associated with $X$ and $\omega$. In \cite[Definition~6.1.4]{HalleNicaiseMotivicCY}, the motivic zeta function is defined for a specific choice of $\omega$, but this condition on $\omega$ can be omitted.

\begin{definition}[{\cite[Definition~6.1.4]{HalleNicaiseMotivicCY}}] \label{def:motivic-zeta-function-CY}
The \emph{motivic zeta function} $Z_{X,\omega}(T)$ of $X$ is defined as
\[Z_{X,\omega}(T) = \sum_{d\in \Z_{> 0}}\left( \int_{X(d)}|\omega(d)|\right)T^d \in \mathcal{M}_k^{\hat{\mu}}\,\llbracket T\rrbracket .\]
The integral $\int_{X(d)}|\omega(d)|$ is a motivic integral measuring the space of $K(d)$-rational points on $X$. 
\end{definition}

\note{
In some references (for example \cite{HalleNicaiseMotivicCY}) a factor $\L^m$ is added. But this gives the same object because then in the definition of $\int_{X(d)}|\omega(d)|$ a factor $\L^{-m}$ is added as well. They cancel each other out. We use the formulation as used in \cite{BultotNicaise} and \cite{HalleNicaiseKulikov}.
}

\note{
In principal, we could also define the motivic zeta function for varieties other than Calabi-Yau's, but in that case, we don't expect a link with monodromy anymore.
}

\note{
We need a volume form to gauge the measure. The volume form makes the integral $\int_{X}|\omega|=\int_{Gr(\X)} \L^{-\ord \omega}$ independent of the chosen model. $Gr(\X)$ is the Greenberg scheme of $\X$ and the schematic structure of $Gr(\X)$ changes when $\X$ changes. 
}

\note{
This is because $\int_{X}|\omega|=\int_{Gr(\X)} \L^{-\ord \omega}$ and $Gr(\X(d)) = \X(d) \big( K(d) \big) = X(d) \big( K(d) \big)$ because $\X(d)$ is a proper model.
}

To understand the inspiration for the definition of the motivic zeta function for Calabi-Yau varieties, we should have a look at an alternative interpretation of~$Z_{f,x}^{mot}(T)$ that Nicaise and Sebag give in \cite{NicaiseSebag_motivic_serre} and \cite[Corollory~9.6]{Nicaise_trace-formula-rigid-varieties} in terms of non-archimedean geometry. The function $Z_{f,x}^{mot}(T)$ is the \emph{local} motivic zeta function and is defined as follows:
let $f\in k[x_1,\ldots, x_m]$ be a polynomial and let $x\in \A_k^m$ be a point such that $f$ is smooth at $x$. We define
\[Z_{f,x}^{mot}(T)=\sum_{d\geq 0} [\X_{d,1,x}] \L^{-md}T^d \in \mathcal{M}_{k}\,\llbracket T\rrbracket ,\]
where
\[\X_{d,1,x}=\{\varphi \in \X_{d,1} \mid \varphi(0)=x \}.\]
Nicaise and Sebag gave an alternative formula for $Z_{f,x}^{mot}(T)$ in terms of the analytic Milnor fiber $\mathscr{F}_x$ of $f$ at $x$ and the Gelfand-Leray form $\omega/df$ on~$\mathscr{F}_x$ associated with a volume form $\omega$ on $\A^n_k$:
\[Z_{f,x}^{mot}(\L T)=\sum_{d\in \Z_{>0}}\left(\int_{\mathscr{F}_x(d)} \left|\frac{\omega}{df}(d)\right|  \right) T^d.\]

This interpretation gave the inspiration to formulate the version of the motivic zeta function discussed for Calabi-Yau varieties. 

\subsubsection{Denef-Loeser formula}

The formula in Theorem~\ref{thm:Denef-Loeser} can be extremely useful to compute the motivic zeta function $Z_f^{mot}(T)$ in practice, when an embedded resolution of $f$ is known. Luckily, we have an analogue for zeta functions of Calabi-Yau varieties, proven by Bultot and Nicaise in \cite{BultotNicaise}, that requires the data of an $snc$-model $\X$ of $X$. We also refer to \cite[Corollary~7.7]{NicaiseSebag_motivic_serre}, where Nicaise and Sebag prove a similar result for formal schemes.

Let $X$ be a Calabi-Yau variety over~$K$ with a volume form $\omega$ and $snc$-model $\X$.
Denote by $\X_k=\sum_{i\in I} N_i E_i$ the special fiber of the model $\X$. 
For every $i\in I$, we define $\nu_i$ as the order of $\mathrm{div}(\omega)$ along $E_i$, when $\omega$ is viewed as a rational section of the relative log-canonical bundle $\omega_{\X/R}(\X_{k,red}-\X_k)$. The couple $(N_i,\nu_i)$ is called the \emph{numerical data} of $E_i$.

\begin{theorem}[{\cite[Corollary~4.3.2]{BultotNicaise}}]
\label{thm:Denef-Loeser-CY}
Let $X$ be a Calabi-Yau variety over~$K$ with a volume form $\omega$ and $snc$-model $\X$.
Denote by $\X_k=\sum_{i\in I} N_i E_i$ the special fiber of the model $\X$ and let $(N_i, \nu_i)$ be the numerical data of $E_i$. 

For every non-empty subset $J\subseteq I$, we define
$E_J = \cap_{j\in J} E_j$ and \[{E_J^\circ = E_J\setminus \left(\cup_{i\in I\setminus J} E_i\right)}.\]
Let $\widetilde{E^\circ_J}$ be the finite \'etale cover $E_J^\circ \times_\X \mathcal{Y}$ of $E_J^\circ$, where $\mathcal{Y}$ is the normalization of 
$\X\times_{R} R[x]/(x^{N_J}- t)$ with $N_J = \gcd\{N_j\mid j\in J\}$.
The motivic zeta function of $X$ can be expressed as
\[Z_{X, \omega}(T)=\sum_{\emptyset \neq J \subseteq I} (\L-1)^{|J|-1}[\widetilde{E^\circ_J}]\prod_{j\in J}\frac{\L^{-\nu_j}T^{N_j}}{1-\L^{-\nu_j}T^{N_j}}\]
in $\mathcal{M}_k^{\hat{\mu}}\llbracket T\rrbracket$.
\end{theorem}

\begin{remark}
In Chapter \ref{ch:poles}, we need a stronger version of Theorem~\ref{thm:Denef-Loeser-CY}, where we weaken the assumption of $\X$ being a scheme to $\X$ being an algebraic space. Theorem~\ref{thm:Denef-Loeser-CY} still holds for algebraic spaces because of \cite[Proposition~7.2.2]{HalleNicaiseKulikov} and the fact that normalized base change commutes with formal completion.

\end{remark}
\note{We need to say why the DL-formula holds for algebraic spaces. The idea is that, when $\X$ is an algebraic space, there exists a disjoint union of schemes $\X_k=\coprod E_{\alpha}$. The formal completion $\widehat{\X}_{E_\alpha}$ of $\X$ is a formal scheme. (formal algebraic spaces and formal schemes are the same.) In \cite{BultotNicaise}, $Z_{E_\alpha}(T)$ is defined, using $\widehat{\X}_{E_\alpha}$. We will then have $Z_\X(T) = \sum Z_{E_\alpha}(T)$ and this does not depend on the chosen decomposition of $\X_k$ in disjoint schemes. For each $Z_{E_\alpha}(T)$, the DL-formula holds. We still need to argue why this implies that the DL-formula also holds for the sum.}

\begin{remark}
Theorem \ref{thm:Denef-Loeser-CY} immediately implies that $Z_{X,\omega}(T)$ is a rational function over $\mathcal{M}_k^{\hat{\mu}}$.
Although the definition of a pole will not be given until Section~\ref{sect:def-poles}, we can intuitively see that all poles of $Z_{X,\omega}(T)$ are of the form $-\nu_i/N_i$ for some $i\in I$. Since a normal crossings model is not unique, one cannot expect all `candidate poles' $-\nu_i/N_i$ to be actual poles of 
$Z_{X,\omega}(T)$. But even candidate poles that appear in \emph{every} model, will not always be actual poles. This phenomenon is intimately related with the monodromy property for Calabi-Yau varieties.
\end{remark}

\section{Monodromy eigenvalues}
 \label{sect:MonodromyEigenvalues}

In this subsection, we will give a more algebraic definition of monodromy eigenvalues. In Definition~\ref{def:monodromy-eigenvalue}, we defined the monodromy transformation~$M_{x,i}$ of a hypersurface singularity to be the action of the canonical generator of~$\pi_1(D^*)$ on $H^i(F_x,\C)$, the singular cohomology of the Milnor fiber $F_x$. Let us have a look at how we can mimic this construction for a smooth and proper $K$-variety~$X$. 


In Definition~\ref{def:monodromy-eigenvalue}, we considered a small open disc $D$ around the origin and the punctured disc $D^*$. One can think of $\Spec(R)$ as a small open disc around the origin, therefore, it plays the role of $D$. Furthermore, $\Spec(K)$ can be thought of as the punctured disc $D^*$. 

The algebraic analogue of $\widetilde{D^*}$ is $\Spec(K^{alg})$. This can be understood as follows: the universal covering space $\widetilde{D^*}$ is the inverse limit of the projective system of finite covers $\{ D^*\to D^*\colon  t\mapsto t^d\}_{d\in \Z_{\geq 0}}$ ordered by divisibility. The right algebraic analogue of a cover is a finite \'etale morphism. So we would like to look at the inverse limit of the projective system of finite \'etale morphisms $\{\Spec (K(d))\to \Spec (K)\colon  { t}'\to  t\}_{d\in \Z_{\geq 0}}$, where $K(d)$ is the unique totally ramified extension of $K$ of degree $d$ in $K^{alg}$ with uniformizer $ t'$. This inverse limit is exactly $\Spec(K^{alg})$, which means that it is indeed the appropriate analogue of $\widetilde{D^*}$. 

The topological monodromy group $\pi_1(D^*)$ is equal to the group of deck transformations $\Aut(\widetilde{D^*}/D^*)$. 
So the algebraic counterpart of $\pi_1(D^*)$ is $\Aut(\Spec(K^{alg})/\Spec(K))$, which is exactly $\Gal(K^{alg}/K)$.
The Galois group $\Gal(K^{alg}/K)$ is isomorphic to the profinite group of roots of unity $\hat{\mu}=\varprojlim \mu_d$. 
\note{$\Gal(K^{alg}/K)\simeq \hat{\mu}$ is true because $k$ has characteristic zero and is algebraically closed.}
\note{$\hat{\mu}$ is not necessarily the group of \emph{complex} roots of unity. It is the group of roots of unity in $k$. This is why there is no canonical topological generator when $k\neq \C$.}

In the topological setting, we noticed that the Milnor fiber $F_x$ is homotopy equivalent to $(g^{-1}(D^*)\cap B) \times_{D^*} \widetilde{D^*}$. So the algebraic match of $F_x$ is the base change $X\times_K K^{alg}$.

We put all the algebraic counterparts in the following table:

\begin{center}
\begin{tabular}{|l|l|}
\hline 
Complex topological & Algebraic counterpart
\\ 
\hline 

$D$ & $\Spec (R)$ \\ 

$D^*$ & $\Spec(K)$ \\ 

$\widetilde{D^*}$ & $\Spec(K^{alg})$\\

$\pi_1(D^*)$ \tiny{ $(=\Aut(\tilde{D}^*/D^*))$} & $\Gal(K^{alg}/K)$\tiny{ $(=\Aut(\Spec(K^{alg})/\Spec(K)))$}\\
 
$F_x$ \tiny{ $(=(g^{-1}(D^*)\cap B) \times_{D^*} \tilde{D}^*)$} & $X\times_K K^{alg}$ \\ 
\hline 
\end{tabular} 
\end{center}

When we interpret $X$ as a family of varieties over the punctured disc, the monodromy action can be thought of as the automorphism on a variety in the family, when travelling once around the origin. Let's define this more formally.

\begin{definition}
Let $X$ be a smooth, proper variety over $K$ and let $\sigma$ be a topological generator of $\Gal(K^{alg}/ K)$.
For every $i\geq 0$, the \emph{monodromy transformation} $M_{X,i}$ is the action of $\sigma$ on the $\ell$-adic cohomology space $H^i(X\times_K K^{alg}, \Q_\ell)$. A \emph{monodromy eigenvalue} is an eigenvalue of $M_{X,i}$ for some $i\geq 0$.
\end{definition}

\note{When $k\neq \C$, then there is no canonical topological generator of $\Gal(K^{alg}/ K)$. This means that the monodromy transformation change when we choose a different generator. It should be possible to show that the set of monodromy eigenvalues is independent of the choice of generator, but it has never been written down explicitly. In any case, there is no problem for $K3$ surfaces because of A'Campo.}

\note{If $k$ is not of characteristic zero, then we should take $\Gal(K^{sep}/ K)$ instead of $\Gal(K^{alg}/ K)$.}

\subsubsection{Monodromy zeta function}

\begin{definition} \label{def:GMP-monodromy-zeta}
The \emph{monodromy zeta function} is the alternating product of the characteristic polynomials of the monodromy transformations $M_{X,i}$ on the $\ell$-adic cohomology spaces $H^m(X\times_K K^{alg}, \Q_\ell)$:
\[\zeta_X(T) = \prod_{m\geq 0} \Bigl(\det \left( T\cdot Id - \sigma\mid H^m(X\times_K K^{alg}, \Q_\ell)\right)\Bigr)^{(-1)^{m+1}} \in \Q_\ell(T).\]
\end{definition}

Remember that in Section~\ref{sect:p-adic-motivic-monodromy}, we found that both the motivic zeta function and the monodromy zeta function of a hypersurface singularity can be computed from the data of an embedded resolution of singularities. In Section~\ref{sect:motivic-zeta-CY}, we found that the motivic zeta function of a Calabi-Yau variety can be computed when the data of an $snc$-model of the Calabi-Yau variety is given. It turns out that the monodromy zeta function from Definition~\ref{def:GMP-monodromy-zeta} can be computed from the data of an $snc$-model as well. 

Nicaise proved in \cite[Theorem~2.6.2]{NicaiseTameRam} the following simple expression for the monodromy zeta function.

\begin{proposition} \label{thm:GMP-ACampo-dvr}
Let $X$ be a smooth, proper $K$-variety with
$snc$-model $\X$. Denote by $\X_k=\sum_{i\in I} N_i E_i$ the special fiber of the model $\X$. 

The monodromy zeta function of $X$ can be written as
\[\zeta_X(T) = \prod_{i\in I} \left(T^{N_i}-1\right)^{-\chi(E_i^\circ)},\]
where $\chi(E_i^\circ)$ is the topological Euler characteristic of $E_i^\circ= E_i\setminus \left( \cup_{j\in I\setminus \{ i\}} E_j\right)$.
\end{proposition}
This is a variant of A'Campo's formula for the local monodromy zeta function of a hypersurface singularity. 

\note{From this proposition, it is plausible that monodromy eigenvalues give information about which models are possible. For example, if there is an eigenvalue $\neq 1$, then it is not possible to have a semistable model. Another example: if there exists a Jordan block of rang 3, then all models must have triple points.}

By the formula of Proposition~\ref{thm:GMP-ACampo-dvr}, we can find some of the monodromy eigenvalues from the data of an $snc$-model of $X$, but in general we won't find all of them, since cancellations may occur: a monodromy eigenvalue is not necessarily a zero or pole of $\zeta_X(T)$. 
Note that for a $K3$ surface, all monodromy eigenvalues appear as poles of the monodromy zeta function, since the cohomology spaces of a $K3$ surface are trivial in odd degree.

\section{The monodromy property for Calabi-Yau varieties} \label{sect:global_monodromy_property}

Informally, the monodromy property for Calabi-Yau varieties expresses that poles of the motivic zeta function correspond to monodromy eigenvalues. By Theorem~\ref{thm:Denef-Loeser-CY}, one can interpret the monodromy property for a Calabi-Yau variety $X$ as a precise relation between its cohomology and the geometry of its $snc$-models. The following formulation of the monodromy property appeared in \cite [Definition~6.4.1]{HalleNicaiseMotivicCY} and a formulation of the equivariant monodromy property can be found in \cite[Definition~2.3.5]{HalleNicaiseKulikov}

\begin{definition}\label{def:GMP}
Let $X$ be a Calabi-Yau variety over~$K$ with volume form $\omega$, and let $\sigma$ be a topological generator of the monodromy group $\textnormal{Gal}(K^{alg}/K)$. 
We say that $X$ satisfies the \emph{monodromy property} if there exists a finite subset $S$ of $\Z\times \Z_{>0}$ such that
\[Z_{X, \omega}(T)\in \mathcal{M}_k^{\hat{\mu}}\left[T, \frac{1}{1-\L^{a}T^b}\right]_{(a,b)\in S},\] 
and such that for each $(a,b)\in S$, we have that $\exp(2\pi i a/b)$ is an eigenvalue of $\sigma$ on $H^m(X\times_K K^{alg}, \Q_\ell)$, for some $m\geq 0$, and for every embedding of $\Q_\ell$ into~$\C$. 

\end{definition}

\note{$H^m(X\times_K K^{alg}, \Q_\ell)$ is a $\Q_\ell$-vector space, so its eigenvalues are $\ell$-adic numbers. But when we embed $\Q_\ell$ in $\C$, this eigenvalue gets mapped to $\exp(2\pi i \nu/N)$. The embedding of $\Q_\ell$ into $\C$ goes as follows: first we take the metric completion $\C_\ell$ of an algebraic closure of $\Q_\ell$. By the axiom of choice, it is possible to prove that $\C_\ell$ is isomorphic to $\C$. This embedding cannot be constructed explicitly because we use the axiom of choice.}

%

\begin{remark}
Whether a Calabi-Yau variety $X$ satisfies the monodromy property does not depend on the choice of volume form $\omega$.
Indeed, let $\omega'$ be another volume form on $X$. Since $X$ has trivial canonical bundle, there exists a unit $u\in K^\times$ such that $\omega'=u\cdot \omega$.
It follows immediately from the definition of the motivic zeta function that
\begin{align}\label{eq:zeta-function-different-volumeform}
Z_{X,\omega'}(T) = Z_{X,\omega}(\mathbb{L}^{-ord_t(u)}\,T).
\end{align}
\note{
Let $\mathcal{Y}$ be a weak Neron model of $X(d)$. For every irreducibel component $C$ of $\mathcal{Y}_k$, the order $\ord_C(u\omega)$ is defined as the unique integer $a$ such that $t^{-a/d} u \omega$ is a generator of $\omega_{\mathcal{Y}/R(d)}$ locally at the generic point of $C$. We have $u = t^{\ord_t(u)}\cdot u'$, with $\ord_t(u')=1$, i.e.\ $u'\in R^\times$. Therefore $t^{-a/d} u \omega = t^{\frac{-a + d\cdot \ord_t (u)}{d}} u ' \omega$.
This implies that \[\ord_C(\omega) = \ord_C(\omega') - d\cdot \ord_t(u).\]
Therefore
\begin{align*}
\int_{X(d)} |\omega'(d)| &= \sum_{C\in \pi_0(\mathcal{Y})} [C]\L^{-\ord_C \omega'}\\
&=\sum_{C\in \pi_0(\mathcal{Y})} [C]\L^{-\ord_C \omega -d\cdot \ord_t(u)}\\
&= \L^{-d\cdot \ord_t(u)}\int_{X(d)} |\omega(d)|
\end{align*}

We conclude
\begin{align} 
Z_{X,\omega'}(T) &= \sum_{d>0} \left( \int_{X(d)} |\omega'(d)| \right) T^d\\
&= \sum_{d>0} \left( \L^{-d\cdot \ord_t(u)}\int_{X(d)} |\omega(d)|\right) T^d\\
&= Z_{X,\omega}(\L^{-\ord_t(u)} T).
\end{align}
}
Let $S$ and $S'$ be as in Definition~\ref{def:GMP} for $Z_{X,\omega}(T)$ and $Z_{X,\omega'}(T)$ respectively. Equation \eqref{eq:zeta-function-different-volumeform} implies that if $(a,b)\in S$, then $(a-\ord_t(u)\cdot b,b) \in S'$. But if $\exp(2\pi i a/b)$ is a monodromy eigenvalue, then so is $\exp(2\pi i (a/b+\ord_t u))$. Therefore, the monodromy property does not depend on the chosen volume form.

\end{remark}

In \cite[Theorem~8.5]{HalleNicaiseAbelian}, Halle and Nicaise prove that abelian varieties satisfy the monodromy property, when the $\hat{\mu}$-action is ignored. The proof uses in an essential way properties of abelian varieties and their N\'eron models. In \cite[Theorem~4.2.2]{HalleNicaiseKulikov}, they upgrade the result so that the $\hat{\mu}$-action is taken into account.

\begin{theorem} \label{thm:GMP-abelian}
Let $X$ be an abelian variety over $K$ and let $\omega$ be a volume form on $X$.
The motivic zeta function $Z_{X,\omega}(T)$ has a unique pole and the monodromy property holds for $X$.
\end{theorem}
\begin{remark}
For a suitable choice of $\omega$, the unique pole of $Z_{X,\omega}(T)$ coincides with Chai's base change conductor.
\end{remark}


\note{It is not proven that the equivariant monodromy property holds for abelian varieties. The first part of the proof can be just copied word for word, but in the second part, there are still details to be checked.}

Halle and Nicaise also proved the following result in \cite[Corollary~5.3.3]{HalleNicaiseKulikov}:

\begin{theorem} \label{thm:GMP-kulikov}
Let $X$ be a Calabi-Yau variety over $K=\C(\! ( t )\! )$ with volume form $\omega$, and assume that $X$ has an equivariant Kulikov
model over $R(d)=\C\llbracket \sqrt[d]{t} \rrbracket$ for some $d > 0$. The motivic zeta function $Z_{X,\omega}(T)$ has a unique pole and the monodromy property holds for $X$.
\end{theorem}


So far, all Calabi-Yau varieties known to satisfy the monodromy property, have a motivic zeta function with a unique pole. We will give the first class of Calabi-Yau varieties satisfying the monodromy property with a motivic zeta function with more than one pole.

In this thesis, we will focus on Calabi-Yau varieties of dimension 2. To investigate the monodromy property in dimension 2, the only remaining case is that of $K3$ surfaces, i.e., 2-dimensional Calabi–Yau varieties $X$ with $H^1(X, \O_X) =0$. Indeed, if $X$ is a 2-dimensional Calabi–Yau variety with $H^1(X, \O_X) \neq 0$, then $X$ is an abelian surface and hence the monodromy property holds for $X$ by Theorem~\ref{thm:GMP-abelian}.

Using Kulikov's classication, Stewart and Vologodsky give in \cite{Stewart-Vologodsky} an explicit formula for $Z_{X, \omega}(T)$ when $X$ is a $K3$ surface allowing a semi-stable model. From their formula, it can be deduced that the motivic zeta function $Z_{X,\omega}(T)$ has only one pole. From Proposition~\ref{thm:Denef-Loeser-CY}, we immediately see that this pole is an integer and therefore, $X$ trivially satisfies the monodromy property. 

In this thesis, we will therefore focus on $K3$ surfaces without a semi-stable model. In Theorem~\ref{thm:Denef-Loeser-CY} and Theorem~\ref{thm:GMP-ACampo-dvr}, we explained why information about an $snc$-model $\X$ is extremely useful to verify the monodromy property for $X$ in practice. To our knowledge, the only classification of $snc$-models of $K3$ surfaces without the condition of semi-stability, is the classification of Crauder and Morrison \cite{CrauMor} where they classify triple-point-free $snc$-models of surfaces with a numerically trivial canonical bundle. We will study and refine this classification in Chapter~\ref{ch:CM}.

\cleardoublepage


\chapter{The Crauder-Morrison classification}\label{ch:CM}

In this chapter, we will discuss and extend the results of \cite{CrauMor}. In this paper, Crauder and Morrison classify triple-point-free $snc$-models of smooth, proper surfaces with trivial pluricanonical bundle. The tools to prove this classification where developed by Crauder in \cite{Crauder}, and errata were published in \cite{CrauMor2}.

In Chapter~\ref{ch:GMP}, in particular in Theorem~\ref{thm:Denef-Loeser-CY} and Theorem~\ref{thm:GMP-ACampo-dvr}, we explained why information about an $snc$-model $\X$ of a Calabi-Yau variety $X$ is extremely useful to verify the monodromy property for $X$ in practice. This is the reason why the Crauder-Morrison classification is one of the major tools in this thesis: it gives \emph{concrete} information about $snc$-models of $K3$ surfaces allowing a triple-point-free model.

In Section~\ref{sect:CM-notation}, we fix notation for this chapter and discuss some basic facts.
The main theorem of \cite{CrauMor} is given in Section~\ref{sect:CM-classification}. The statement there is a rather rough classification and a lot more information can be given. The details that are relevant for this thesis are formulated in sections~\ref{sect:flowers} to \ref{sect:Euler-characteristics}. Most of these statements appeared in \cite{CrauMor}, but some are original. In particular, Theorem~\ref{thm:CM-K3-no-cycle} has been announced in \cite{Jaspers}.

\section{Notation and basic facts} \label{sect:CM-notation}

\subsubsection{Notation}

We fix an algebraically closed field $k$ of characteristic zero. Define $R=k\llbracket t\rrbracket$ and $K=k(\!(t) \!)$, the fraction field of $R$.
Let $X$ be a smooth, proper surface over $K$ with $\omega_{X/K}^{\otimes m}\simeq \O_X$ for some $m\geq 1$. \note{When we write $\omega_X$, this means $\omega_{X/K}$.} Let $\X$ be an $snc$-model of $X$ with special fiber $\X_k=\sum_{i\in I}N_i E_i$, where we allow $\X$ to be an algebraic space. Suppose moreover that $\X$ is triple-point-free, i.e., 
$E_i \cap E_j \cap E_k = \emptyset$,
whenever $i$, $j$ and $k$ are pairwise distinct. Assume furthermore that $\X$ is relatively minimal for this property, i.e., it is not possible to contract components in the special fiber (in the category of algebraic spaces) such that the result is still a triple-point-free $snc$-model. A relatively minimal, triple-point-free $snc$-model is called a \emph{Crauder-Morrison model}.

Let $\omega$ be a nowhere-vanishing section of $\omega_{X/K}^{\otimes m}$. 
We can view $\omega$ as a rational section of $\omega_{\X/R}(\X_{k, red}-\X_k)^{\otimes m}$, which is a multiple of the relative log-canonical bundle. Define $\nu_i$ such that $m\nu_i$ is the multiplicity of $\mathrm{div}(\omega)$ along $E_i$. The couple $(N_i,\nu_i)$ is called the \emph{numerical data} of $E_i$.

We define \emph{the weight} $\rho_i$ of the component $E_i$ to be
\begin{align}
\rho_i=\frac{\nu_i}{N_i}+1. \label{def:CM-weight}
\end{align}


Notice that the definition of $\nu_i$ and $\rho_i$ depends on the choice of $\omega$. If $\omega'$ is another nowhere-vanishing section of $\omega_{X/K}^{\otimes m}$, inducing $\nu_i'$ and $\rho_i'$, then $\nu_i = \nu_i' + c N_i$ and $\rho_i = \rho_i' + c$ for a fixed constant $c\in \Z$.
\note{Because $X$ is proper, reduced and geometrically connected, the ring of regular functions is $k$. (If not geometrically connected, but just connected, then the ring of regular functions is a finite extension of $k$. (Why?)) $\omega$ and $\omega'$ differ by a regular function (a constant), and $c$ is the valuation of this constant.}

We define $\Gamma$ to be the \emph{dual graph} of $ \mathcal{X}_k $. Notice that the dual complex of $\X_k$ is indeed a graph, because $\X_k$ is triple-point-free. Denote by $\Gamma_{min}$ the subgraph of $\Gamma$ spanned by the vertices corresponding to components $E_i$ with minimal weight~$\rho_i$. 

\subsubsection{Basic facts}

Let $K_{\X/R}$ be the relative canonical divisor of $\X$ over $R$.
By definition of $\nu_i$, we have
\begin{equation}
K_{\X/R} \equiv \sum_{i\in I} (\nu_i+N_i-1) E_i. \label{eq:canonical-bundle-component-nu}
\end{equation}

Since the special fiber $ \X_k $ is a numerically trivial divisor on $\X$, we have the relation
\[ E_j \equiv - \frac{1}{N_j} \sum_{i \neq j} N_i E_i. \]
Together with the adjunction formula \cite[Proposition~5.73]{KollarMori}, this relation yields
\begin{align} 
 K_{E_j/k} &\equiv \left((K_{\X/R} + E_j) \cdot E_j\right)_{E_j} \nonumber\\
 &\equiv \sum_{i \neq j} \left(\nu_i - \frac{N_i}{N_j}\nu_j - 1\right) (E_i\cdot E_j)_{E_j}\nonumber \\
 &\equiv \sum_{i \neq j} \left(N_i (\rho_i - \rho_j) - 1\right) (E_i \cdot E_j)_{E_j}.\label{eq:canonical-bundle-component-rho}
\end{align}
\note{
$E_i\cdot E_j$ is an element in the Chow group and $E_i\cap E_j$ is a subscheme.
}
\note{
Adjunction formula: $K_{E_j/k} = (K_{\X/R}+E_j)\cdot E_j$. So $K_{E_j/k} = (\sum_{i\in I}r_i E_i +E_j)\cdot E_j = \left(\sum_{i\neq j} r_i E_i - \frac{r_j+1}{N_j} \sum_{i \neq j} N_i E_i\right)\cdot E_j=\sum_{i \neq j} (N_i (\rho_i - \rho_j) - 1) E_i \cdot E_j = \sum_{i \neq j} (N_i ((r_i+1)/N_i - (r_j+1)/N_j) - 1) E_i \cdot E_j = \sum_{i\neq j} (r_i+1 - \frac{r_j+1}{N_j}N_i -1)E_i\cdot E_j$
}
\note{
Note that this expression is independent of our choice of $\omega$. Indeed, $\rho_i-\rho_j = (\rho_i' + c) - (\rho_j'+c).$
}
\note{
To prove the adjunction formula for algebraic spaces: An algebraic space is locally in the \'etale topology a scheme, so then it is not that difficult to prove the adjunction formula for algebraic spaces.
}

For every $j\neq i$, we define
\[a_{ji} = 
\begin{cases}
N_i(\rho_i-\rho_j)-1 & \text{if } E_i\cap E_j\neq \emptyset,\\
0 & \text{otherwise}.
\end{cases}
\]
Then
\begin{equation}
K_{E_j/k}\equiv \sum_{i\neq j} a_{ji} (E_i\cdot E_j)_{E_j}, \label{eq:canonical-bundle-component-general}
\end{equation}
and we have
\begin{align}
a_{ji}<-1 & \qquad\text{iff}\qquad \rho_j>\rho_i, \nonumber\\
a_{ji}=-1 & \qquad\text{iff}\qquad \rho_j=\rho_i, \label{eq:CM-equal-weight} \\
a_{ji}>-1 & \qquad\text{iff}\qquad \rho_j<\rho_i. \nonumber
\end{align}

\section{The Crauder-Morrison classification} \label{sect:CM-classification}
Crauder and Morrison classified the special fiber~$\X_k$ of relatively minimal, triple-point-free $snc$-models $\X$ of a smooth, proper surface $X$ over $K$ with trivial pluricanonical bundle. 


\begin{theorem}[Crauder - Morrison Classification \cite{CrauMor}] \label{thm:CM-classification}
Let $X$ be a smooth, proper surface over $K$ with $\omega_{X/K}^{\otimes m}\simeq \O_X$ for some $m\geq 1$, and let $\X$ be a Crauder-Morrison model of $X$. Then $\X$ has the following properties:

\begin{enumerate}[(i)]
\item $\Gamma_{min}$ is a connected subgraph of $\Gamma$. It is either a single vertex, a cycle or a chain. We call $\X$ a flowerpot degeneration, a cycle degeneration or a chain degeneration respectively.

\item Each connected component of $\Gamma\setminus \Gamma_{min}$ is a chain (called a \emph{flower}) $F_0\text{---} F_1 \text{---} \cdots \text{---}F_{\ell}$, where only $F_\ell$ meets $\Gamma_{min}$, and $F_\ell$ meets a unique vertex of $\Gamma_{min}$. 
The weights strictly decrease: $\rho_0>\rho_1>\cdots>\rho_\ell$.
The surface $F_0$ is either minimal ruled or it is isomorphic to $\P^2$. If $F_0\simeq \P^2$, then $F_0\cap F_1$ is either a line or a conic on $F_0$. The surface $F_i$ is minimal ruled with sections $F_{i-1}\cap F_i$ and $F_i \cap F_{i+1}$ for $1\leq i \leq \ell$.

\item Suppose $\Gamma_{min}$ is a single vertex $P$ (called a \emph{flowerpot}). The surface $P$ is isomorphic to $\P^2$, or it is a ruled surface, or $K_{P/k}\equiv 0$. 

\item Suppose $\Gamma_{min}$ is a cycle $V_1\text{---}V_2\text{---}\cdots \text{---} V_k$. Then there are no flowers and all components have the same multiplicity, i.e., there exists an integer $N \geq 1$ such that
$\X_k = N \left(\sum_{i=1}^k V_i\right)$.
Furthermore, for every $i=1,\ldots, k$, the component $V_i$ is an elliptic, minimal ruled surface with sections $V_{i-1}\cap V_i$ and $V_i\cap V_{i+1}$, where we identify $V_0=V_k$ and $V_{k+1}=V_1$. 

\item Suppose $\Gamma_{min}$ is a chain $V_0\text{---}V_1\text{---}\cdots \text{---} V_k\text{---} V_{k+1}$. If $i=1,\ldots, k$, the surface $V_i$ is an elliptic, ruled surface with sections $V_{i-1}\cap V_i$ and $V_i\cap V_{i+1}$.
If $i=0$ or $k+1$, the surface $V_i$ is either isomorphic to $\P^2$, or it is a rational or elliptic, ruled surface.
\end{enumerate}
\end{theorem}

\begin{remark}
Note that in Theorem~\ref{thm:CM-classification}, we do not assume the special fiber to be reduced. So in general,
the model $\X$ does not need to be semi-stable. A classification of semi-stable models of $K3$ surfaces has been given by Kulikov in \cite{Kulikov1} with corrections by Persson and Pinkham in \cite{Persson-Pinkham}. Crauder and Morrison weaken the hypothesis of semi-stability, but impose the extra condition of the model being triple-point-free.
\end{remark}

In the pictures below, you can find illustrations of possible dual graphs of models of surfaces considered in Theorem~\ref{thm:CM-classification}. We label the vertices of the graphs with $(N_i,\rho_i)$.

\begin{minipage}{0.5\textwidth}
\begin{figure}[H]
\centering
\begin{tikzpicture}[scale = 0.5]

\node [draw,shape=circle, fill, inner sep=0pt,minimum size=5pt] at (0,0) {};
\node [below] at (0,-0.2) {{$(3, {4}/{3})$}};
\node [draw,shape=circle, fill, inner sep=0pt,minimum size=5pt] at (-2,1) {};
\node [left] at (-2,1) {{$(2,{3}/{2})$\ \ }};
\node [draw,shape=circle, fill, inner sep=0pt,minimum size=5pt] at (-2,2) {};
\node [left] at (-2,2) {{$(1,2)$\ \ }};
\node [draw,shape=circle, fill, inner sep=0pt,minimum size=5pt] at (2,1) {};
\node [right] at (2,1) {{$\ (2,{5}/{2})$}};
\node [draw,shape=circle, fill, inner sep=0pt,minimum size=5pt] at (2,2) {};
\node [right] at (2,2) {{\ $(1,3)$}};
\node [draw,shape=circle, fill, inner sep=0pt,minimum size=5pt] at (2,3) {};
\node [right] at (2,3) {{\ $(1,4)$}};
\node [draw,shape=circle, fill, inner sep=0pt,minimum size=5pt] at (2,4) {};
\node [right] at (2,4) {{\ $(2,{9}/{2})$}};

\draw
(0,0) -- (-2,1)
(-2,1) -- (-2,2)
(0,0) -- (2,1)
(2,1) -- (2,2)
(2,2)--(2,3)
(2,3)--(2,4);
\end{tikzpicture}
\caption{A flowerpot degeneration}
\label{fig:flowerpot}
\end{figure}
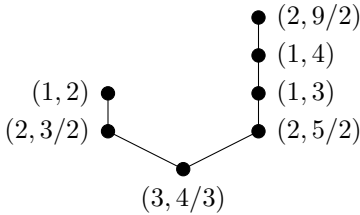
\end{minipage} \hfill
\begin{minipage}{0.45\textwidth}
In Figure \ref{fig:flowerpot}, you can see an example of a flowerpot degeneration with two flowers. 
Notice that although the dual graph is in fact a chain, there is only one surface with minimal weight $4/3$, and therefore it is a flowerpot degeneration. The weights of the components of the flowers are strictly decreasing towards the flowerpot. 
\end{minipage}\hfill

\begin{minipage}{0.5\textwidth}
In Figure~\ref{fig:cycle-degeneration}, you can see a cycle degeneration. All components have weight 2, and hence there are no flowers. In fact, a cycle degeneration never has flowers. Moreover, all components have the same multiplicity.
\end{minipage} \hfill
\begin{minipage}{0.4\textwidth}
\begin{center}
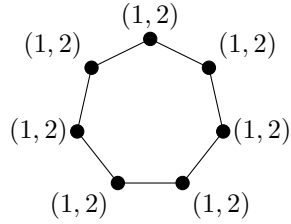
\begin{figure}[H]
\begin{tikzpicture}
\node [draw,shape=circle, fill, inner sep=0pt,minimum size=5pt] at (0,1) {};
\node [above] at (0,1) {$(1,2)$};
\node [draw,shape=circle, fill, inner sep=0pt,minimum size=5pt] at (0.78,0.62) {};
\node [above right] at (0.78,0.62) {$(1,2)$};
\node [draw,shape=circle, fill, inner sep=0pt,minimum size=5pt] at (0.97,-0.22) {};
\node [right] at (0.97,-0.22) {$(1,2)$};
\node [draw,shape=circle, fill, inner sep=0pt,minimum size=5pt] at (0.43,-0.9) {};
\node [below right] at (0.43,-0.9) {$(1,2)$};
\node [draw,shape=circle, fill, inner sep=0pt,minimum size=5pt] at (-0.43,-0.9) {};
\node [below left] at (-0.43,-0.9) {$(1,2)$};
\node [draw,shape=circle, fill, inner sep=0pt,minimum size=5pt] at (-0.97,-0.22) {};
\node [left] at (-0.97,-0.22) {$(1,2)$};
\node [draw,shape=circle, fill, inner sep=0pt,minimum size=5pt] at (-0.78,0.62) {};
\node [above left] at (-0.78,0.62) {$(1,2)$};

\draw (0,1) -- (0.78,0.62);
\draw (0.97,-0.22) -- (0.78,0.62);
\draw (0.97,-0.22) -- (0.43,-0.9);
\draw (0.43,-0.9) -- (-0.43,-0.9);
\draw (-0.43,-0.9) -- (-0.97,-0.22);
\draw (-0.97,-0.22) -- (-0.78,0.62);
\draw (-0.78,0.62) -- (0,1);
\end{tikzpicture}
\caption{A cycle degeneration}
\label{fig:cycle-degeneration}
\end{figure}
\end{center}
\end{minipage}\hfill

\begin{minipage}{0.5\textwidth}
\begin{figure}[H]
\centering
\begin{tikzpicture}[scale=0.6]

\node [draw,shape=circle, fill, inner sep=0pt,minimum size=5pt] at (0,0) {};
\node [below] at (0,0) {$(4,3/2)$};
\node [draw,shape=circle, fill, inner sep=0pt,minimum size=5pt] at (4,0) {};
\node [below] at (4,0) {$(2,3/2)$};

\node [draw,shape=circle, fill, inner sep=0pt,minimum size=5pt] at (0,2) {};
\node [above] at (0,2) {$(2,2)$};

\node [draw,shape=circle, fill, inner sep=0pt,minimum size=5pt] at (3,2) {};
\node [above] at (3,2) {$(4,7/4)$};

\node [draw,shape=circle, fill, inner sep=0pt,minimum size=5pt] at (5,2) {};
\node [above] at (5,2) {$(1,5/2)$};

\draw (0,0) -- (4,0)
(0,0) -- (0,2)
(4,0) -- (3, 2)
(4,0) -- (5, 2)
;
\end{tikzpicture}
\caption{A chain degeneration}\label{fig:chain-degeneration}
\end{figure}
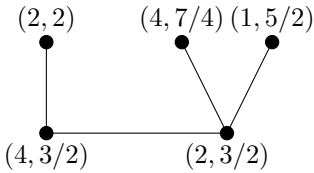
\end{minipage}\hfill
\begin{minipage}{0.45\textwidth}
In Figure~\ref{fig:chain-degeneration}, you can see an example of a possible chain degeneration. There are two components in the chain, both with weight $3/2$. Moreover, there are three flowers, all of length $1$. 
\end{minipage}

We also refer to Example~\ref{ex:poles-more-than-one} for an example of a $K3$ surface allowing a flowerpot degeneration.


\section{Flowers} \label{sect:flowers}


\subsection{Terminology}

\begin{definition}
Let $X$ be a smooth, proper surface over $K$ with $\omega_{X/K}^{\otimes m}\simeq \O_X$ for some $m\geq 1$, and let $\X$ be a Crauder-Morrison model of $X$. Suppose that the special fiber~$\X_k$ has a flower $N_0 F_0 + N_1 F_1 + \cdots + N_{\ell} F_\ell$, where $F_i\cap F_j=\emptyset$ if and only if $j\not\in \{i-1,i, i+1\}$. The top of the flower is $F_0$, and $F_\ell$ meets $\Gamma_{min}$ in $F_{\ell+1}$.

The double curve $F_{\ell}\cap F_{\ell+1}$ is called a \emph{flowercurve}.
The \emph{genus} of the flower is the genus of the flowercurve $F_\ell\cap F_{\ell+1}$ .
We call the flower \emph{rational} if $F_0$ is a rational surface.
\end{definition}

\begin{remark}\label{rmk:CM-genus-flower}
For $1\leq i\leq \ell$, the double curves $F_{i-1}\cap F_i$ and $F_{i}\cap F_{i+1}$ are sections of the minimal ruled surface~$F_i$ and therefore, they have the same genus.
So if a flower has genus $g$, all double curves $F_i\cap F_{i+1}$ are curves of genus $g$ for $0\leq i\leq \ell$.
\end{remark}
\begin{remark}\label{rmk:CM-rational-flower}
 A flower is rational if and only if it has genus zero. Indeed, if the flower is rational, then either $F_0$ is a rational, ruled surface, or $F_0\simeq \P^2$. If $F_0$ is a rational, ruled surface, then $F_0\cap F_1$ is a rational curve, since it is a section on $F_0$. If $F_0\simeq \P^2$, then the curve $F_0\cap F_1$ is either a line or a conic, and therefore, it is a rational curve. Remark~\ref{rmk:CM-genus-flower} explains why the flower has genus 0. Conversely, if the flower is not rational, then $F_0$ is ruled over a curve of genus $g\geq 1$. Since $F_0\cap F_1$ is a section of the ruling on $F_0$, it has genus $g$, and hence the flower has genus $g\geq 1$, by Remark~\ref{rmk:CM-genus-flower}.
\end{remark}

\begin{definition}
If $F_0\simeq \P^2$, and if $F_0\cap F_1$ is a conic, then the flower is called a \emph{conic-flower}.
\end{definition}

\begin{prop-def}[{\cite[Theorem~3.2 and  Theorem~3.5]{CrauMor}}] 
\label{thm:CM-type-flower}
Let $X$ be a smooth, proper surface over $K$ with $\omega_{X/K}^{\otimes m}\simeq \O_X$ for some $m\geq 1$, and let $\X$ be a Crauder-Morrison model of $X$. Suppose that the special fiber~$\X_k$ has a flower $N_0 F_0 + N_1 F_1 + \cdots + N_{\ell} F_\ell$, where $F_i\cap F_j=\emptyset$ if and only if $j\not\in \{i-1,i, i+1\}$. The top of the flower is $F_0$, and $F_\ell$ meets $\Gamma_{min}$ in $F_{\ell+1}$.

There is an integer $M\in \{2,3,4,6,8, 12\}$, called the \emph{type} of the flower, such that
\[K_{F_{\ell+1}/k}\equiv (-1+\frac{2}{M}) (F_\ell\cdot F_{\ell+1})_{F_{\ell+1}}-D,\]
where $D$ is an effective divisor on $F_{\ell+1}$ disjoint from the flowercurve $F_{\ell}\cap F_{\ell+1}$. 
We also call $M$ the type of the flowercurve $F_{\ell}\cap F_{\ell+1}$.
\end{prop-def}

\begin{remark}\label{rmk:CM-canonical-bundle-min-weight} 
In the notation of equation \eqref{eq:canonical-bundle-component-general} with $E_i=F_\ell$ and $E_j=F_{\ell+1}$, we have $a_{ji} = -1+\frac{2}{M}$. So in particular $-1< a_{ji}\leq 0$ and $a_{ji}=0$ if and only if $M=2$.
\end{remark}

\subsection{Classification of the flowers}

Crauder and Morrison classified the flowers into 21 combinatorial classes.

\begin{theorem}[{\cite[Theorem~3.9]{CrauMor} and \cite[Appendix~2]{CrauMor2}}] \label{thm:CM-classification-flowers}
Let $X$ be a smooth, proper surface over $K$ with $\omega_{X/K}^{\otimes m}\simeq \O_X$ for some $m\geq 1$, and let $\X$ be a Crauder-Morrison model of $X$. Suppose that the special fiber~$\X_k$ has a flower $N_0 F_0 + N_1 F_1 + \cdots + N_{\ell} F_\ell$, where $F_i\cap F_j=\emptyset$ if and only if $j\not\in \{i-1,i, i+1\}$. The top of the flower is $F_0$, and $F_\ell$ meets $\Gamma_{min}$ in $F_{\ell+1}$. The multiplicity of $F_{\ell+1}$ is $N_{\ell+1}$.

Then there exists an integer $N\geq 1$, such that the flower satisfies one of the 21~possibilities described in Tables \ref{table:flowers-P2-line}, \ref{table:flowers-P2-conic} and~\ref{table:flowers-ruled}.

\begin{table}[H]
\begin{center}
\begin{tabular}{|c|c|c|c|c|}
\hline
Name & Type &Composition&$\ell$&$N_{F_{\ell+1}}$\\
\hline
2A & 2 & $NF_0$ & 0&$2N$\\
3A &3 & $NF_0+2NF_1$&1 & $3N$\\
3B &3 & $NF_0$ &0& $3N$\\ 
4A &4 & $NF_0+2NF_1+3NF_2$&2 & $4N$\\
4B&4 & $NF_0$ & 0&$4N$\\
6A&6 & $NF_0+2NF_1+3NF_2+4NF_3+5NF_4$&4 & $6N$\\
6B&6 & $NF_0$ & 0&$6N$\\
\hline
\end{tabular}
\caption{Flowers with $F_0\simeq \P^2$ and $F_0\cap F_1$ is a line}\label{table:flowers-P2-line}
\end{center}
\end{table}

\begin{table}[H]
\begin{center}
\begin{tabular}{|c|c|c|c|c|}
\hline
Name & Type &Composition&$\ell$&$N_{F_{\ell+1}}$\\
\hline
2B&2 & $2NF_0$  & 0& $N$ \\
2C&2 & $2NF_0+NF_1+\cdots + NF_\ell$ &$\ell$&  $N$ \\
4C&4 & $2NF_0+NF_1+\cdots + NF_\ell$ &$\ell$& $2N$ \\ 
4D&4 & $NF_0+NF_1+\cdots + NF_\ell$ &$\ell$&  $N$ \\
6C&6 & $2NF_0+NF_1+\cdots + NF_{\ell-1}+2NF_\ell$ &$\ell$& $3N$ \\
6D&6 & $2NF_0+NF_1+\cdots + NF_\ell$  &$\ell$& $3N$ \\
6E&6 & $2NF_0$& 0&$3N$\\
\hline
\end{tabular}
\caption{Flowers with $F_0\simeq \P^2$ and $F_0\cap F_1$ is a conic}\label{table:flowers-P2-conic}
\end{center}
\end{table}

\begin{table}[H]
\begin{center}
\begin{tabular}{|c|c|c|c|c|}
\hline
Name & Type &Composition&$\ell$&$N_{F_{\ell+1}}$\\
\hline
4$\alpha$&4 & $NF_0$ &0&$2N$ \\
6$\alpha$ &6& $NF_0+2NF_1$&1 & $3N$  \\
6$\beta$ &6& $NF_0$ &0& $3N$ \\ 
8$\alpha$ &8& $NF_0+2NF_1+3NF_2$&2 & $4N$\\
8$\beta$ &8& $NF_0$&0 & $4N$ \\
12$\alpha$ &12& $NF_0+2NF_1+3NF_2+4NF_3+5NF_4$ &4& $6N$\\
12$\beta$ &12& $NF_0$ &0& $6N$ \\
\hline
\end{tabular}
\caption{Flowers with $F_0$ minimal ruled}\label{table:flowers-ruled}
\end{center}
\end{table}
\end{theorem}

\subsection{Self-intersection of a double curve in the flower}

We will now prove some results concerning the self-intersection of the double curves $F_i\cap F_{i+1}$ for $i=0, \ldots, \ell$.


\begin{lemma} \label{thm:CM-self-intersection-flower-curve}
Let $X$ be a smooth, proper surface over $K$ with $\omega_{X/K}^{\otimes m}\simeq \O_X$ for some $m\geq 1$, and let $\X$ be a Crauder-Morrison model of $X$. Suppose that the special fiber~$\X_k$ has a flower $N_0 F_0 + N_1 F_1 + \cdots + N_{\ell} F_\ell$, where $F_i\cap F_j=\emptyset$ if and only if $j\not\in \{i-1,i, i+1\}$. The top of the flower is $F_0$, and $F_\ell$ meets $\Gamma_{min}$ in $F_{\ell+1}$.

Denote by $C$ the flowercurve $F_\ell\cap F_{\ell+1}$. Suppose the flower has genus $g$ and type $M$, then 
\[(C^2)_{F_{\ell+1}}=M(g-1).\]
So in particular, if the flower is rational, then $(C^2)_{F_{\ell+1}}=-M$.
\end{lemma}
\begin{proof}
By Proposition-Definition~\ref{thm:CM-type-flower}, we can write
\[K_{F_{\ell+1}/k}\equiv (-1+\frac{2}{M}) C-D,\]
\note{
It is not a problem to interchange $F_\ell\cap F_{\ell+1}$ and $F_\ell\cdot F_{\ell+1}$. The intersection cycle 
$F_\ell \cdot F_{\ell+1}$
is represented by the curve 
$F_\ell \cap F_{\ell+1}$
because the surfaces intersect transversally.
}
where $D$ is an effective divisor on $F_{\ell+1}$ disjoint from the flowercurve $C$.
The adjunction formula gives
\[2g-2=K_{F_{\ell+1}/k}\cdot C+ C^2 = \frac{2}{M}C^2. \]
It follows that $C^2=M(g-1)$.
\end{proof}

We also state and prove the following well-known relation:

\begin{lemma} \label{thm:CM-self-intersection-different-surface}
Let $X$ be a smooth, proper surface over $K$ with $\omega_{X/K}^{\otimes m}\simeq \O_X$ for some $m\geq 1$, and let $\X$ be a Crauder-Morrison model of $X$. Suppose that the special fiber~$\X_k$ has a flower $N_0 F_0 + N_1 F_1 + \cdots + N_{\ell} F_\ell$, where $F_i\cap F_j=\emptyset$ if and only if $j\not\in \{i-1,i, i+1\}$.

Let $C_i=F_{i-1}\cap F_i$. For $i=1, \ldots, \ell$, we have
\[ (C_{i+1}^2)_{F_i} = - (C_i^2)_{F_i}.\]
\end{lemma}

\begin{proof}
Since $F_i$ is minimal ruled, $\mathrm{Num}(F_i)$ is generated by a section $C$ and a fiber $\fiber$, and moreover, $\fiber^2=0$ and $C\cdot \fiber=1$ by \cite[Proposition~V.2.3]{Hartshorne}. As~$C_i$ and $C_{i+1}$ are sections, we can write
\[ C_i \equiv C + a \fiber\qquad \text{ and }\qquad C_{i+1} \equiv  C + b \fiber, \]
for some integers $a,b\in \Z$.
On the one hand,
\begin{align*}
C_{i+1}^2 + C_i^2
	&= (C+a\fiber)^2 + (C+b\fiber)^2\\
	&= C^2+2a(C\cdot \fiber) +\fiber^2 + C^2 + 2b(C\cdot \fiber) + \fiber^2\\
	&= 2(C^2+a+b).
\end{align*}

On the other hand,
\[0=C_i\cdot C_{i+1} = (C+a\fiber)\cdot (C+b\fiber) = C^2+a+b,\]
since there are no triple points. This implies $(C_{i+1}^2)_{F_i} = - (C_i^2)_{F_i}$.
\end{proof}

\subsection{Combinatorial relations on the numerical data}

For $i=0, \ldots, \ell+1$, let $(N_i, \nu_i)$ be the numerical data of the component $F_i$. We will derive combinatorial relations between the numerical data of consecutive components in a flower. Corollary~\ref{thm:CM-K3-no-4D} is an important consequence. 

\begin{lemma}\label{thm:CM-flowers-numerical-relations} Let $X$ be a smooth, proper surface over $K$ with $\omega_{X/K}^{\otimes m}\simeq \O_X$ for some $m\geq 1$, and let $\X$ be a Crauder-Morrison model of $X$. Suppose that the special fiber~$\X_k$ has a flower $N_0 F_0 + N_1 F_1 + \cdots + N_{\ell} F_\ell$, where $F_i\cap F_j=\emptyset$ if and only if $j\not\in \{i-1,i, i+1\}$. The top of the flower is $F_0$, and $F_\ell$ meets $\Gamma_{min}$ in $F_{\ell+1}$.
\begin{enumerate}[(i)]
\item Suppose $F_0\simeq \P^2$ and $F_0\cap F_1$ is a line on $F_0$. Then we have
\[\nu_1 = \frac{N_1}{N_0}\nu_0-2.\]
\item Suppose $F_0\simeq \P^2$ and $F_0\cap F_1$ is a conic on $F_0$. Then we have
\[\nu_1 = \frac{N_1}{N_0}\nu_0-1/2.\]
\item Suppose $F_0$ is a minimal ruled surface. Then we have
\[\nu_1 = \frac{N_1}{N_0}\nu_0-1.\]
\item For $j=1, \ldots, \ell$, we have
\[\nu_{j+1} = \frac{N_{j-1}+N_{j+1}}{N_j} \nu_j- \nu_{j-1}.\]
\end{enumerate}
\end{lemma}

\begin{proof}
\begin{enumerate}[(i)]

\item Suppose $C_1=F_0\cap F_1$ is a line on $F_0\simeq \P^2$. From \eqref{eq:canonical-bundle-component-rho}, it follows that
\[K_{F_0/k} \equiv \left(\nu_1 - \frac{N_1}{N_0}\nu_0 - 1\right) C_1.\]
Adjunction on $F_0$ gives
\[ -2 = (K_{F_0/k} + C_1) \cdot C_1 = \left(\nu_1 - \frac{N_1}{N_0}\nu_0 - 1\right)C_1^2 + C_1^2=\left(\nu_1 - \frac{N_1}{N_0}\nu_0\right)C_1^2. \]
Since $C_1$ is a line, $C_1^2=1$ and therefore $\nu_1 - \frac{N_1}{N_0}\nu_0=-2$.

\item Suppose $C_1=F_0\cap F_1$ is a conic on $F_0\simeq \P^2$. From \eqref{eq:canonical-bundle-component-rho}, it follows that
\[K_{F_0/k} \equiv \left(\nu_1 - \frac{N_1}{N_0}\nu_0-1\right)C_1.\]
Adjunction on $F_0$ gives
\[ -2 = (K_{F_0/k} + C_1) \cdot C_1 = \left(\nu_1 - \frac{N_1}{N_0}\nu_0-1\right)C_1^2 + C_1^2 = \left(\nu_1 - \frac{N_1}{N_0}\nu_0\right) C_1^2. \]
Since $C_1$ is a conic, $C_1^2=4$, and therefore $ \nu_1 - \frac{N_1}{N_0}\nu_0 = - 1/2 $.

\item Suppose $C_1=F_0\cap F_1$ is a section on the minimal ruled surface $F_0$. From \eqref{eq:canonical-bundle-component-rho}, it follows that
\[ K_{F_0/k} \equiv \left(\nu_1 - \frac{N_1}{N_0}\nu_0- 1\right)C_1. \]
Let $\fiber$ be a fiber of the ruling on $F_0$. The adjunction formula gives
\[-2=2g(\fiber)-2=K_{F_0/k}\cdot \fiber+\fiber^2=K_{F_0/k}\cdot \fiber.\]
 Hence, we get that
\[ - 2 = \left(\nu_1 - \frac{N_1}{N_0}\nu_0 - 1\right) C_1 \cdot \mathcal{F} = \nu_1 - \frac{N_1}{N_0}\nu_0- 1. \]
Therefore, $\nu_1 - \frac{N_1}{N_0}\nu_0=-1.$

\item Let $C_j=F_{j-1}\cap F_j$ and $C_{j+1}=F_j\cap F_{j+1}$. Both $C_j$ and $C_{j+1}$ are sections on the minimal ruled surface $F_j$. From \eqref{eq:canonical-bundle-component-rho}, it follows that
\begin{align*}
 K_{F_j/k} \equiv \left(\nu_{j-1} - \frac{N_{j-1}}{N_j}\nu_j - 1\right)C_j + \left(\nu_{j+1} - \frac{N_{j+1}}{N_j}\nu_j - 1\right)C_{j+1}. 
\end{align*}
Let $\fiber$ be a fiber of the minimal ruled surface $F_j$.
On the one hand, the adjunction formula applied to $\fiber$ gives
\[ K_{F_j/k} \cdot \fiber = (K_{F_j/k} + \mathcal{F}) \cdot \mathcal{F}= -2.\]
On the other hand, the expression for $K_{F_j/k}$ gives that
\begin{align*}
 K_{F_j/k}\cdot \fiber &= \left(\nu_{j-1} - \frac{N_{j-1}}{N_j}\nu_j - 1\right)C_j \cdot \fiber\\
&\qquad + \left(\nu_{j+1} - \frac{N_{j+1} }{N_j}\nu_j- 1\right)C_{j+1}\cdot \fiber\\
&= \left(\nu_{j-1} - \frac{N_{j-1}}{N_j}\nu_j\right)+ \left(\nu_{j+1} - \frac{N_{j+1}}{N_j}\nu_j\right)-2
\end{align*}
Hence, $\left(\nu_{j-1}-\frac{N_{j-1}}{N_j}\nu_j\right)+\left(\nu_{j+1}-\frac{N_{j+1}}{N_j}\nu_j\right) = 0$.
\end{enumerate}
\end{proof}

\begin{corollary} \label{thm:CM-K3-no-4D}
Suppose $X$ is a $K3$ surface or an abelian variety over $K$, with Crauder-Morrison model~$\X$. The special fiber $\X_k$ does not have flowers of type~$4D$.
\end{corollary}
\begin{proof}
Write $\X_k=\sum_{i\in I}N_i E_i$. Let $(N_i, \nu_i)$ be the numerical data of $E_i$.
Since $X$ is a $K3$ surface or an abelian variety, we have that $\omega_{X/k}\simeq \O_X$, which means that $m=1$. It follows that $\nu_i\in \Z_{>0}$, for all $i\in I$. For a flower of type $4D$, relation (ii) of Lemma~\ref{thm:CM-flowers-numerical-relations} gives
\[\nu_1=\nu_0-\frac{1}{2},\]
which is a contradiction.
\end{proof}

\begin{corollary}\label{thm:CM-flowers-relation-mu}
For every class in Theorem~\ref{thm:CM-classification-flowers}, there is a fixed relation between $\nu_{\ell+1}$ and $\nu_0$, which can be found in Table \ref{table:flowers-relation-mu}.

\begin{table}[H]
\begin{center}
\begin{tabular}{|c|c||c|c||c|c|}
\hline
Type &$\nu_{{\ell+1}}$ &Type &$\nu_{{\ell+1}}$&Type &$\nu_{{\ell+1}}$\\
\hline
2A & $2\nu_0-2$&2B & $(\nu_0-1)/2$&4$\alpha$ & $2\nu_0-1$\\
3A & $3\nu_0-4$ &2C & $(\nu_0-2\ell-1)/2$&6$\alpha$ & $3\nu_0-2$\\
3B & $3\nu_0-2$& 4C & $\nu_0-2\ell$&6$\beta$ & $3\nu_0-1$\\ 
4A & $4\nu_0-6$&4D & $(2\nu_0 - \ell-1)/{2}$&8$\alpha$ & $ 4\nu_0-3$\\
4B & $4\nu_0-2$&6C & $(3\nu_0-6\ell+5)/2$&8$\beta$ & $4\nu_0-1 $\\
6A & $6\nu_0-10$&6D & $(3\nu_0-6\ell+1)/2$&12$\alpha$ & $6\nu_0-5$\\
6B & $6\nu_0-2$&6E & $(3\nu_0-1)/2$&12$\beta$ & $6\nu_0-1$\\
\hline 
\end{tabular}
\caption{Relation between $\nu_{{\ell+1}}$ and $\nu_0$}\label{table:flowers-relation-mu}
\end{center}
\end{table}
\end{corollary}

\begin{proof}
By induction and Lemma~\ref{thm:CM-flowers-numerical-relations}, we can write $\nu_{\ell+1}$ in terms of $\nu_0$. The computations are straightforward and we will illustrate this fact by computing $\nu_{\ell+1}$ for two classes of flowers.

A flower of type $4A$ is given combinatorially by $NF_0+2NF_1+3NF_2$ for some integer~$N\geq 1$, and ${N_3=4N}$. 
By numerical relation (iv), we get
$\nu_2=2\nu_1-\nu_0$. 
Using the numerical relation (i), we find that
$\nu_1 = 2\nu_0-2$, and hence
$\nu_2 = 3\nu_0-4$.
By relation (iv), we have $\nu_3 = 2\nu_2-\nu_1$.
Substituting the expressions for $\nu_1$ and~$\nu_2$ gives
\[\nu_3 = 4\nu_0 - 6.\]

A flower of type $2C$ is given combinatorially by
$ 2N F_0 + N F_1 + \ldots + N F_\ell$, for some integer~$N\geq 1$, 
and $ N_{\ell+1}= N $. Relation (ii) gives
$ \nu_1 = (\nu_0-1)/2$.
Relation (iv) gives $\nu_2 = 3 \nu_1 - \nu_0$,
and~$ {\nu_{j+1} = 2 \nu_j - \nu_{j-1}}$,
for $ 2 \leq j \leq \ell $. An easy induction argument shows that~${ \nu_j = \frac{\nu_0 - 2j + 1}{2}, }$
for $1\leq j\leq \ell+1$. So in particular, we get
\[ \nu_{\ell+1} = \frac{\nu_0 - 2\ell- 1}{2}. \]
\end{proof}

\section{Flowerpots} \label{sect:flowerpots}

In this section, we will give more details about the structure of flowerpots and the possible combinations of flowers that a flowerpot degeneration can contain.

\subsection{Geometry of the pot}

\begin{proposition}\label{thm:CM-geometry-pot}
Assume that $X$ is a smooth, proper surface over $K$ with $\omega_{X/K}^{\otimes m}\simeq \O_X$ for some $m\geq 1$, such that it has a Crauder-Morrison model $\X$ where the subgraph $\Gamma_{min}$ of the dual graph $\Gamma$ has a unique vertex $P$.

The geometry of the flowerpot $P$ is one of the following:
\begin{enumerate}[(i)]
\item Either $ K_{P/k} \equiv 0 $,
\note{We actually have a stronger result: $N_P K_P \sim 0$.}
\item or  $ P \cong \mathbb{P}^2 $,
\item or $P$ is minimal ruled,
\item or $P$ is rational, non-minimal ruled.
\end{enumerate} 
\end{proposition} 
\begin{proof}
Combine the remark directly before \cite[Proposition 2.17]{Crauder} with \cite[Theorem 3.2]{CrauMor}.
In \cite[Lemma 4.1]{CrauMor}, it is stated that a non-minimal ruled pot is rational.
\end{proof}

\subsection{Combinations of flowers in a flowerpot degeneration}

\begin{lemma}\label{thm:CM-pot-hodge-index}
Assume that $X$ is a smooth, proper surface over $K$ with $\omega_{X/K}^{\otimes m}\simeq \O_X$ for some $m\geq 1$, such that it has a Crauder-Morrison model $\X$ where the subgraph $\Gamma_{min}$ of the dual graph $\Gamma$ has a unique vertex $P$.

Assume there are $\lambda$ flowers in $\X_k$, for some integer $\lambda\geq 1$. Denote by $C_1, \ldots, C_\lambda$ the flowercurves on $P$. Then one of the following holds:
\begin{enumerate}[(i)]
\item
There exists a $ j \in \{ 1, \ldots, \lambda \} $, such that the flowercurve $C_j$ has genus ${g(C_j) > 1} $, and for all $ i \in \{ 1, \ldots, \lambda \} \setminus \{ j\} $, the flowercurve $C_i$ is rational.
\item For all $ i \in \{ 1, \ldots, \lambda \} $, the flowercurve $C_i$ has genus $ g(C_i) \leq 1 $ .
\end{enumerate}
\end{lemma}
\begin{proof}
Define $ a_i = 1 - \frac{2}{M_i} $ for all $i= 1, \ldots, \lambda$, where $M_i$ is the type of the flowercurve $C_i$. We have
\[ K_{P/k} \equiv -\sum_{i=1}^\lambda a_i C_i,\]
by Remark~\ref{rmk:CM-canonical-bundle-min-weight}.
By adjunction and the fact that there are no triple points, we have for any $i=1, \ldots, \lambda$:
\[ 2 g(C_i) - 2 = (K_{P/k} + C_i)\cdot C_i = (1-a_i) C_i^2 .\]
Suppose (ii) does not hold. Hence there exists a $j\in \{1,\ldots, \lambda\}$ such that the flowercurve $C_j$ has genus $g(C_j) > 1 $. We will show that for $i\in \{1,\ldots, \lambda\}\setminus\{j\}$, the flowercurve $C_i$ is rational. By Theorem~\ref{thm:CM-classification-flowers}, the type of $C_j$ is ${M_j \geq 4}$. 
This implies that $C_j^2 > 0$, by Lemma~\ref{thm:CM-self-intersection-flower-curve}. Moreover, for any $ {i \neq j}$, we know that $ C_i \cdot C_j = 0 $, since there are no triple points. Hence the Hodge Index Theorem~\cite[Example~19.3.1]{Fulton-intersection} implies that $C_i^2 < 0$. Therefore,  $g(C_i) = 0$, by Lemma~\ref{thm:CM-self-intersection-flower-curve}.
\end{proof}

We will now discuss some restrictions on the combination of flowers in a flowerpot degeneration, depending on the geometry of the flowerpot.

\begin{proposition}\label{thm:CM-pot-NK=0}
Assume that $X$ is a smooth, proper surface over $K$ with $\omega_{X/K}^{\otimes m}\simeq \O_X$ for some $m\geq 1$, such that it has a Crauder-Morrison model $\X$ where the subgraph $\Gamma_{min}$ of the dual graph $\Gamma$ has a unique vertex $P$. 

If $K_{P/k} \equiv 0 $, then any flower in $\X_k$ has either type~2 or genus 1.
\end{proposition}
\begin{proof}
Denote by $C_1, \ldots, C_\lambda$ the flowercurves on $P$ with types $M_1, \ldots, M_\lambda$ respectively. We will show that for every $i=1, \ldots, \lambda$, we have either $g(C_i)=1$ or $M_i=2$.
By Remark~\ref{rmk:CM-canonical-bundle-min-weight}, we can write
\[ K_{P/k} \equiv \sum_{i = 1 }^\lambda \left(-1 + \frac{2}{M_i}\right) C_i. \]
Since $ K_{P/k} \equiv 0 $, it follows that $ K_{P/k} \cdot C_i = 0 $, and therefore, $C_i^2=0$ or $M_i=2$. Adjunction yields that if $C_i^2=0$, then $g(C_i) = 1 $.
\end{proof}

\begin{proposition}[{\cite[Lemma~3.6]{CrauMor}}]\label{thm:CM-pot-P2}
Assume that $X$ is a smooth, proper surface over $K$ with $\omega_{X/K}^{\otimes m}\simeq \O_X$ for some $m\geq 1$, such that it has a Crauder-Morrison model $\X$ where the subgraph $\Gamma_{min}$ of the dual graph $\Gamma$ has a unique vertex $P$.

If $P\simeq \P^2$, then there exists exactly one flower, and it has type 4 and genus 10, or type 8 and genus 3.
\end{proposition}

\begin{proposition}[{\cite[Lemma~3.7]{CrauMor}}]\label{thm:CM-pot-minimal-ruled}
Assume that $X$ is a smooth, proper surface over $K$ with $\omega_{X/K}^{\otimes m}\simeq \O_X$ for some $m\geq 1$, such that it has a Crauder-Morrison model $\X$ where the subgraph $\Gamma_{min}$ of the dual graph $\Gamma$ has a unique vertex $P$.

If $P$ is a minimal ruled surface, and if there is a rational flower, then there is exactly one other flower, and it has genus $g\geq 2$.
\end{proposition}

\section{Cycles} \label{sect:cycles}

In this section, we prove that cycle degenerations do not occur for $K3$ surfaces. This result has been announced in \cite{Jaspers}. 

\begin{lemma} \label{thm:CM-cycle-euler-characteristic}
Let $X$ be a smooth, proper surface over $K$ with $\omega_{X/K}^{\otimes m}\simeq \O_X$ for some $m\geq 1$, and let $\X$ be a Crauder-Morrison model of $X$. Assume that the subgraph $\Gamma_{min}$ of the dual graph $\Gamma$ is a cycle. Then $X$ has trivial $\ell$-adic Euler characteristic, i.e.\ $\chi(X)=0$.
\end{lemma}

\begin{proof}
By Theorem~\ref{thm:CM-classification} (iv), we can find an integer $N\geq 1$ such that ${\X_k = N \sum_{i=1}^k V_i}$, with $V_i$ an elliptic, minimal ruled surface with sections $V_{i-1}\cap V_i$ and $V_i\cap V_{i+1}$, for all $i=1, \ldots, k$. Here we identify $V_0=V_k$ and $V_{k+1} = V_1$. So for every $i=1, \ldots, k$, we have
\[\chi (V_i) = 0 \text{ and } \chi(V_{i-1}\cap V_i) = \chi (V_i\cap V_{i+1})=0.\]
Define $V_i^\circ = V_i \setminus \left(\bigcup_{j\neq i} V_j \right)$, then
\[\chi(V_i^\circ) = \chi (V_i) - \chi(V_{i-1}\cap V_i) -\chi (V_i\cap V_{i+1})=0.\]
By taking the valuation at infinity of both sides of the A'Campo formula in Proposition~\ref{thm:GMP-ACampo-dvr}, we get
\[\chi(X) = N \sum_{i=1}^k\, \chi(V_i^\circ)=0.\]
\end{proof}

\begin{theorem} \label{thm:CM-K3-no-cycle}
Let $X$ be a $K3$ surface over $K$ and let $\X$ be a Crauder-Morrison model of $X$. Then $\X$ is either a flowerpot degeneration or a chain degeneration, but not a cycle degeneration.
\end{theorem}
\begin{proof}
This follows from Lemma~\ref{thm:CM-cycle-euler-characteristic} and the fact that the $\ell$-adic Euler characteristic of a $K3$ surface is 24.
\end{proof}

\section{Chains} \label{sect:chains}

\subsection{Combinations of flowers in a chain degeneration}
\begin{proposition}[{\cite[Corollory~6.2]{CrauMor}}]
\label{thm:CM-chain-flowers}
Let $X$ be a smooth, proper surface over $K$ with $\omega_{X/K}^{\otimes m}\simeq \O_X$ for some $m\geq 1$, and let $\X$ be a Crauder-Morrison model of $X$. Suppose that the subgraph $\Gamma_{min}$ of the dual graph $\Gamma$ is a chain $V_0 - V_1 - \cdots - V_k - V_{k+1}$. The following properties hold.
\begin{enumerate}[(i)]
\item If $j\in \{1, \ldots, k\}$, then all flowers meeting $V_j$ are rational of type 2.

\item If $V_0\simeq \P^2$, then $V_0$ does not meet any flowers.

\item If $V_0$ is a rational, ruled surface, then all flowers meeting $V_0$ are rational of type 2. 

\item If $V_0$ is an elliptic, ruled surface, then flowers meeting $V_0$ are either rational of type 2 or non-rational of type $4\alpha$. There is either a unique non-rational flower of type $4\alpha$, or there are two non-rational flowers of type $4\alpha$. In the latter case, both these flowers are elliptic.
\end{enumerate}
The analogous result of (ii), (iii), and (iv) holds for $V_{k+1}$ as well.
\end{proposition}

\subsection{Structure of the chain}

\begin{proposition}[{\cite[Theorem~6.7 and Proposition~6.9]{CrauMor}}] \label{thm:CM-chain-good}
Let $X$ be a smooth, proper surface over $K$ with $\omega_{X/K}^{\otimes m}\simeq \O_X$ for some $m\geq 1$, and let $\X$ be a Crauder-Morrison model of $X$. Suppose that the subgraph $\Gamma_{min}$ of the dual graph~$\Gamma$ is a chain.
There exists a birational modification $\widetilde{\X}$ of $\X$ with the following properties:
\begin{enumerate}[(i)]
\item $\widetilde{\X}$ is a Crauder-Morrison model of $X$.
\item There exist integers $\alpha, \beta$ and $k$ with $0\leq \alpha\leq \beta\leq k+1$, such that $\widetilde{\Gamma}_{min}$ is $N_0V_0+N_1V_1+\cdots + N_kV_k+N_{k+1}V_{k+1}$
and
\[N_0>N_1>\cdots >N_\alpha = \cdots = N_{\beta}<\cdots<N_k<N_{k+1}.\]
\item The surfaces $V_1, \ldots, V_{\alpha-1}, V_{\beta+1}, \ldots, V_k$ are elliptic, minimal ruled surfaces, and don't meet any flowers.
\item If $\alpha\geq 1$, then there exists an integer $N\geq 1$ such that the chain 
\[V_0 - \cdots - V_{\alpha-1}\]
 satisfies one of the possibilities described in Table \ref{table:CM-chain-classification}. A similar result holds for the chain $V_{\beta+1} - \cdots - V_{k+1}$, if $\beta\leq k$.
\begin{table}[H]
\begin{center}
\begin{tabular}{|c|c|c|c|c| }
\hline
$V_0$&$N_0V_0+ \cdots + N_{\alpha-1} V_{\alpha-1}$ &$\alpha$& $N_{\alpha}$& Remark\\
\hline \hline
\multirow{3}{*}{$V_0\simeq \P^2$} & $3NV_0$ & 1& $N$ & \multirow{3}{*}{}\\\cline{2-4}
& $3NV_0+2NV_1$&2&$N$ & No flowers on $V_0$\\
\cline{2-4}
& $3NV_0$& 1&$2N$&\\ 
\hline
\multirow{4}{*}{$V_0\simeq \Sigma_2$}& $4NV_0$&1&$N$&\multirow{4}{*}{}\\
\cline{2-4}
& $4NV_0+3NV_1+2NV_2$&3&$N$& Unique flower on $V_0$\\
\cline{2-4}
&$4NV_0+3NV_1$&2&$2N$& of type $2A$\\
\cline{2-4}
& $4NV_0$&1&$3N$&\\ 
\hline
$V_0$ is ruled& $2NV_0$&1&$N$&\\ \hline
\end{tabular}
\caption{Classification of the chains $V_0-\cdots - V_{\alpha-1}$}\label{table:CM-chain-classification}
\end{center}
\end{table}
\end{enumerate}
\end{proposition}

\note{The birational modifications we need are flops and contractions. Take $i$ as small as possible with $V_i$ not minimal ruled and $N_i>N_{i+1}$ or $N_{i}<N_{i+1}$. It can be proved that $V_i$ meets only flowers of type $2A$ and every such flowercurve meets an exceptional curve. For every flower we do the following operation: we first blow up the exceptional curve on $V_i$, the exceptional divisor is $\P^1\times \P^1$, and then we contract along the other ruling of $\P^1\times  \P^1$. (This is a flop). Afterwards we contract the strict transform of this flower to its flowercurve. Performing this operation for every flower, we make $V_i$ minimal ruled without flowers. Note that this operation also altered $V_{i+1}$. We keep performing this operation for every non-minimally ruled $V_i$ with $N_i>N_{i+1}$ or $N_{i}<N_{i+1}$.}

From now on, we assume that a Crauder-Morrison model with a chain degeneration always satisfies the properties in Proposition~\ref{thm:CM-chain-good}.

\subsection{Geometry of the components in the chain}

\begin{proposition}\label{thm:CM-chain-structure-fibers}
Let $X$ be a smooth, proper surface over $K$ with $\omega_{X/K}^{\otimes m}\simeq \O_X$ for some $m\geq 1$, and let $\X$ be a Crauder-Morrison model of $X$. Suppose that the subgraph $\Gamma_{min}$ of the dual graph $\Gamma$ is a chain with components 
$V_0, V_1, \cdots , V_k, V_{k+1}$,
where $V_i\cap V_j\neq \emptyset$ if and only if $j\in \{i-1,i,i+1\}$.
\begin{enumerate}[(i)]
\item For $i= 1, \ldots, k$, the double curves $C_{i}=V_{i-1}\cap V_i$ and $C_{i+1}=V_i\cap V_{i+1}$ are sections on the elliptic, ruled surface $V_i$, and any reducible fiber of the ruling on $V_i$ is a chain
\[ \fiber_1 + \ldots + \fiber_s, \]
where $\fiber_1^2 = \fiber_s^2 = - 1 $ and $ \fiber_j^2 = - 2 $ for $j=2, \ldots, s-1$. Moreover, $\fiber_1$ meets $C_i$ transversally in a unique point, and  $\fiber_s$ meets $C_{i+1}$ transversally in a unique point as well. For $j = 2, \ldots, s-1$, the curve $\fiber_j$ does not intersect $C_i$, nor $C_{i+1}$.

\item Suppose $V_0$ is ruled over an elliptic curve $D$, and let $C_1=V_0\cap V_1$. One of the following holds:
\begin{enumerate}
\item The morphism $C_1\to D$ is an isomorphism. Any reducible fiber of the ruling on $V_0$ is a chain
\[  \fiber_1 + \ldots + \fiber_s, \]
where $\fiber_1^2 =\fiber_s^2 = - 1 $ and $ \fiber_j^2 = - 2 $ for $j=2, \ldots, s-1$. Moreover, $\fiber_1$ meets $C_1$ transversally in a unique point. For $j = 2, \ldots, s$, the curve $\fiber_j$ does not intersect $C_1$.
\item The morphism $C_1\to D$ is an \'etale morphism of degree 2. Any reducible fiber of the ruling on $V_0$ is a chain
\[ \fiber_1 + \ldots + \fiber_s, \]
where $\fiber_1^2 =\fiber_s^2 = - 1 $ and $ \fiber_j^2 = - 2 $ for $j=2, \ldots, s-1$. Moreover, $\fiber_1$ and $\fiber_s$ meet $C_1$ transversally in a unique point each. For $j = 2, \ldots, s-1$, the curve $\fiber_j$ does not intersect $C_1$.
\end{enumerate}
A similar result holds for $V_{k+1}$.

\item Suppose $V_0$ is ruled over a rational curve $D$. The double curve $C_1=V_0\cap V_1$ is an elliptic curve and the morphism $C_1\to D$ is a morphism of degree 2 which is branched above four points $\{d_1, \ldots, d_4\}\subset D$. 
 
Moreover, any reducible fiber of the ruling on $V_0$ above a point $P\not\in \{d_1, \ldots, d_4\}$, is a chain
\[ \fiber_1 + \ldots + \fiber_s ,\]
where $\fiber_1^2 =\fiber_s^2 = - 1 $ and $ \fiber_j^2 = - 2 $ for $j=2, \ldots, s-1$. Moreover, $\fiber_1$ and $\fiber_s$ meet $C_1$ transversally in a unique point each. For $j = 2, \ldots, s-1$, the curve $\fiber_j$ does not intersect $C_1$.

A reducible fiber of the ruling on $V_0$ above a point $P\in \{d_1, \ldots, d_4\}$, has one of the following forms:
\begin{enumerate}
\item 
$ 2 \fiber_1 + \ldots +2 \fiber_{s-2} + \fiber_{s-1} + \fiber_s$, where $\fiber_1, . . . , \fiber_{s-1}$ is a chain. The component $\fiber_s$ intersects $\fiber_{s-2}$ transversally in a unique point and is disjoint from $\fiber_j$ for $j=1,\ldots, s-3, s-1$. Moreover, $\fiber_1^2= -1$ and $\fiber^2_j = -2$ for $j=2, \ldots, s$. The component $\fiber_1$ meets $C_1$ in a unique point and $\fiber_j$ is disjoint from $C_1$ for $j=2,\ldots, s$.
\item $\fiber_1+\fiber_2$, where $\fiber_1^2=\fiber_2^2=-1$. The curve $C_1$ intersects $\fiber_1$ and $\fiber_2$ both transversally in $\fiber_1\cap \fiber_2$.
\end{enumerate}

A similar result holds for $V_{k+1}$.
\end{enumerate}
\end{proposition}

To prove this proposition, we first need several lemmas.

\begin{lemma}\label{thm:CM-chain-exceptional-curve}
Let $V$ be a surface over $k$ such that $ K_{V/k} \equiv \sum_{j=1}^n a_j C_j $ for rational numbers $ -1 \leq a_j < 0 $, with $a_1 = -1 $ and where all $C_j$ are smooth, pairwise disjoint curves. Then the curve $C_1$ has genus $ g(C_1) = 1 $.

Suppose moreover that $E \subset V$ is a rational curve with $E^2=-1$ meeting $C_1$ non-trivially. Then $E$ does not meet $ C_j $ for $ 2\leq  j \leq n $, and $E$ intersects $C_1$ transversally in a unique point.
\end{lemma}
\begin{proof}
We have that
\[ K_{V/k} \equiv -C_1 + \sum_{j=2}^n a_j C_j, \]
so the adjunction formula implies $ g(C_1) = 1 $.

Under the assumptions in the statement, we have that $ E \neq C_j $ for all $ 1 \leq j \leq n $, as all $C_j$ are pairwise disjoint. So $ C_1 \cdot E > 0 $, and $ C_j \cdot E \geq 0 $ for $j=2,\ldots, n$. The adjunction formula gives that $ - K_{V/k} \cdot E = 1 $. Therefore we find
\[ 1 = - K_{V/k} \cdot E =  C_1 \cdot E - \sum_{j=2}^n a_j C_j \cdot E, \]
which implies that $ C_1 \cdot E = 1 $ and $ \sum_{j=2}^n a_j C_j \cdot E = 0 $. We can conclude that $ C_j \cdot E = 0 $ for all $ 2 \leq j \leq n $, and that $C_1$ meets $E$ transversally in a unique point.
\end{proof}

\begin{lemma}\label{thm:CM-chain-canonical-divisor}
Let $V$ be a surface over $k$ such that $ K_{V/k} \equiv \sum_{j=1}^n a_j C_j $ for rational numbers $ -1 \leq a_j < 0 $, with $a_1 = -1 $ and where all $C_j$ are smooth, pairwise disjoint curves. Assume that there is a rational curve $ E \subset V $, with $E^2=-1$, meeting $C_1$ transversally in a unique point, and that $E$ meets no other $C_j$ for $2 \leq j \leq n $. Let $ f \colon  V \to \overline{V} $ be the contraction of $E$.

The curves $ \overline{C}_j = f(C_j) $ are smooth and pairwise disjoint, and
\[ K_{\overline{V}/k} \equiv \sum_{j=1}^n a_j \overline{C}_j.\]
\end{lemma}
\begin{proof}
Since $f$ is the contraction of $E$, it is an isomorphism on $V\setminus E$. Because $E$ is disjoint from $C_j$ for $2\leq j\leq n$, the curves $ \overline{C}_j $ are also smooth, pairwise disjoint curves. Moreover, since $C_1$ meets $E$ transversally in a unique point, it follows from \cite[Proposition V.3.6 and Proposition V.3.2]{Hartshorne} that $ \overline{C}_1 $ is smooth. \note{$f^*\overline{C} = C+rE$, with $r$ the multiplicity of $\overline{C}$ in $P=f(E)$ (V.3.6). By V.3.2, $0=f^*\overline{C}\cdot E = C\cdot E +rE^2 = 1-r$, so $r=1$. Remark V.3.5.2 explains why this implies $\overline{C}$ is smooth in $P$.
}

To prove the last claim, we have that $f^* \overline{C}_j = C_j$, for $ 2\leq j \leq n $, and that $f^* \overline{C}_1 = C_1 + E$, by \cite[Proposition~V.3.6]{Hartshorne}.
Therefore,
\[ f^* \Big(\sum_{j=1}^n a_j \overline{C}_j\Big) = \sum_{j=1}^n a_j C_j + a_1 E = K_{V/k} - E, \]
since $ a_1 = - 1 $.

By \cite[Proposition~V.3.3]{Hartshorne}, it holds that
$K_{V/k} - E = f^* K_{\overline{V}/k}$,
so
\[f^* K_{\overline{V}/k} = f^* \Big(\sum_{j=1}^n a_j \overline{C}_j\Big). \]
 The map $ f^*$ is injective, since 
$ f^*\colon  \Pic(\overline{V})\to \Pic(V)\simeq \Pic(\overline{V})\oplus \Z$
is an inclusion \cite[Proposition V.3.2]{Hartshorne}. Hence
\[ K_{\overline{V}/k} \equiv \sum_{j=1}^n a_j \overline{C}_j.\] 
\end{proof}

\begin{lemma}\label{thm:CM-chain-minimal-ruled}
Let $V$ be a ruled surface over $k$ such that $ K_{V/k} \equiv \sum_{j=1}^n a_j C_j $ for rational numbers $ -1 \leq a_j < 0 $, with $a_1 = -1 $, and where all $C_j$ are smooth, pairwise disjoint curves. Assume that $C_1$ does not meet any rational curve~$E$ with $E^2=-1$. 
Then the curve $C_1$ does not meet any smooth, rational curve with negative self-intersection. Moreover, $V$ is a minimal ruled surface.
\end{lemma}
\begin{proof}
Assume there exists a smooth, rational curve $E$ on $V$ with $ E^2 \leq -2 $ and $C_1\cdot E>0$. Adjunction implies that 
\[ C_1 \cdot E - \sum_{j=2}^n  a_j C_j \cdot E = - K_{V/k} \cdot E \leq 0. \]
Hence $ \sum_{j=2}^n a_j C_j \cdot E > 0 $.
However, the assumption $ C_1 \cdot E > 0 $ implies that $E\neq C_j$, for $2\leq j\leq n$, since the curve~$C_1$ is disjoint from $C_j$, for every $2\leq j \leq n$. So $C_j\cdot E \geq 0$ for every $2\leq j \leq n$.  But this is a contradiction, because $a_j<0$, for every $2\leq j \leq n$.

To prove the last claim, notice that if $V$ is a ruled surface that is not minimal ruled, then there exists a reducible fiber $b_1\fiber_1+ b_2 \fiber_2 + \cdots + b_s \fiber_s$ where each $\fiber_j$ is a smooth, rational curve with negative self-intersection. Since $g(C_1) = 1$, by Lemma~\ref{thm:CM-chain-exceptional-curve}, the curve~$C_1$ must be a horizontal component, and hence meets every fiber non-trivially. But we have shown that $C_1$ does not meet any smooth, rational curve with negative self-intersection, hence $V$ must be minimal ruled.
\end{proof}

\emph{Proof of Proposition~\ref{thm:CM-chain-structure-fibers}.}
Let $V=V_i$ be one of the surfaces in the chain for $0\leq j \leq k+1$. Then $ K_{V/k} \equiv \sum_{j=1}^n a_j C_j $ for rational numbers $ -1 \leq a_j < 0 $ by equations \eqref{eq:canonical-bundle-component-general} and \eqref{eq:CM-equal-weight}. After renumbering, we can assume $C_1$ is the double curve on $V$, where $V=V_i$ meets $V_{i-1}$ or $V_{i+1}$. Because $V_i$ and $V_{i\pm 1}$ have the same weight, $a_1 = -1$. Moreover, all $C_j$ are smooth, pairwise disjoint curves. 

If $V$ is not minimal ruled, then by Lemma~\ref{thm:CM-chain-minimal-ruled}, there is a rational curve $E$ with $E^2=-1$ meeting $C_1$. By Lemma~\ref{thm:CM-chain-exceptional-curve}, we know that $E\cdot C_j=0$, for $2\leq j \leq n$. Let $V\to V^{(1)}$ be the contraction of $E$. By Lemma~\ref{thm:CM-chain-canonical-divisor}
\[K_{V^{(1)}/k} \equiv \sum_{j=1}^n a_j C^{(1)}_j,\]
where $C^{(1)}_j$ is the image of $C_j$ under the contraction. Furthermore, all $C_j^{(1)}$ are smooth and pairwise disjoint, by Lemma~\ref{thm:CM-chain-canonical-divisor}. If $V^{(1)}$ is not minimal ruled, then there is a rational curve $E$ with $E^2=-1$ meeting $C_1^{(1)}$. We can apply the same reasoning as before, to get a contraction $V^{(1)}\to V^{(2)}$. 

We can repeat this process and obtain a series of morphisms
\[ V = V^{(0)} \to V^{(1)} \to \cdots \to V^{(t)} = \overline{V} \]
by successively contracting rational curves with self-intersection $-1$ meeting $C_1^{(l)}$. At any stage of the process, Lemma~\ref{thm:CM-chain-canonical-divisor} guarantees that $ K_{V^{(l)}/k} \equiv \sum_{j=1}^n a_j C_j^{(l)} $, where $C_j^{(l)}$ denotes the image of $C_j$ in $V^{(l)}$. The curves $C_j^{(l)}$ are smooth and pairwise disjoint and $C_1^{(l)}$ is an elliptic curve on the ruled surface $V^{(l)}$. \note{elliptic because of Lemma~\ref{thm:CM-chain-exceptional-curve}}

By Lemma~\ref{thm:CM-chain-minimal-ruled}, this process stops as soon as there are no rational curves with self-intersection $-1$ meeting the image of $C_1$ and the surface $\overline{V} = V^{(t)} $ is minimal ruled.

Going backwards in this process, we find that
\[ V^{(l-1)} \to V^{(l)} \]
is precisely blowing up a smooth point on $C_1^{(l)}$ for every $l=1, \ldots t$. 

\begin{enumerate}[(i)]
\item Assume $i=1, \ldots, k$, then $V=V_i$ is elliptic ruled, $ n \geq 2 $, and after renumbering, $a_1=a_2=-1$. The adjunction formula applied to a general fiber of the ruling on $V$, gives that $C_1$ and $C_2$ are both sections. Repeatedly blowing up smooth points on $C_1^{(l)}$, gives the result. We illustrate this in Figure \ref{fig:CM-structure-fiber-middle-chain}. You see a series of blow-ups where the centers are the marked points. The horizontal, dashed lines are $C_1^{(l)}$ and $C_2^{(l)}$ and the vertical lines are the components of a fiber of the ruling. The labels $-2, -1$ and $0$ indicate the self-intersections of the components of the fibers.

\begin{figure}[H]
\minipage{0.24\textwidth}
\centering
\begin{tikzpicture}[scale = 0.4]
\draw[dashed]	(-0,0.25) -- (2.5,0.25)
		(-0,4.75) -- (2.5,4.75)
;

\node [right] at (2.5,0.25) {\scriptsize{$C_2$}};
\node [right] at (2.5,4.75) {\scriptsize{$C_1$}};

\draw	(1.25,5) -- (2,3.75)
		(1.25,3) -- (2,4.25)
		(1.25,2) -- (2,0.75)
		(1.25,0) -- (2,1.25)
;
\node at (1.25,2.65) {\tiny{$\vdots$}};

\node [left] at (1.25,5) {\tiny{$-1$}};
\node [left] at (1.25,3) {\tiny{$-2$}};
\node [left] at (1.25,2) {\tiny{$-2$}};
\node [left] at (1.25,0) {\tiny{$-1$}};






\node [right] at (3,2.25) {$\longrightarrow \cdots$};

\end{tikzpicture}
\endminipage \hfill
\minipage{0.23\textwidth}
\begin{tikzpicture}[scale = 0.4]
\draw[dashed]	(0,0.25) -- (2.5,0.25)
		(0,4.75) -- (2.5,4.75)
;

\node [right] at (2.5,2.25) {$\longrightarrow$};
\node [left] at (0,2.25) {$\longrightarrow$};


\draw	(1.25,5) -- (2,3.5)
		(1.25,2.25) -- (2,4)
		(1.25,2.75) -- (2,1)
		(1.25,0) -- (2,1.5)
;

\node [left] at (1.25,5) {\tiny{$-1$}};
\node [left] at (1.25,3) {\tiny{$-2$}};
\node [left] at (1.25,2) {\tiny{$-2$}};
\node [left] at (1.25,0) {\tiny{$-1$}};

\node [shape=circle, fill, inner sep=0pt,minimum size=4pt] at (1.375,4.75) {};

\end{tikzpicture}
\endminipage\hfill
\minipage{0.17\textwidth}
\begin{tikzpicture}[scale = 0.4]
\draw[dashed]	(-0,0.25) -- (2.5,0.25)
		(-0,4.75) -- (2.5,4.75)
;

\node [right] at (2.5,2.25) {$\longrightarrow$};


\draw	(1.25,5) -- (2,3.08)
		(1.25,1.41) -- (2,3.58)
		(1.25,1.91) -- (2,0)
;

\node [left] at (1.25,5) {\tiny{$-1$}};
\node [left] at (1.25,1.91) {\tiny{$-2$}};
\node [left] at (2,0) {\tiny{$-1$}};

\node [shape=circle, fill, inner sep=0pt,minimum size=4pt] at (1.348,4.75) {};

\end{tikzpicture}
\endminipage\hfill
\minipage{0.17\textwidth}
\begin{tikzpicture}[scale = 0.4]
\draw[dashed]	(-0,0.25) -- (2.5,0.25)
		(-0,4.75) -- (2.5,4.75)
;

\node [right] at (2.5,2.25) {$\longrightarrow$};


\draw	(1.25,5) -- (2,2.25)
		(1.25,0) -- (2,2.75)
		
;

\node [left] at (1.25,5) {\tiny{$-1$}};
\node [left] at (1.25,0) {\tiny{$-1$}};

\node [shape=circle, fill, inner sep=0pt,minimum size=4pt] at (1.318,4.75) {};

\end{tikzpicture}
\endminipage\hfill
\minipage{0.17\textwidth}
\begin{tikzpicture}[scale = 0.4]
\draw[dashed]	(-0,0.25) -- (2.5,0.25)
		(-0,4.75) -- (2.5,4.75)
;

\node [right] at (2.5,0.25) {\scriptsize{$C_2^{(t)}$}};
\node [right] at (2.5,4.75) {\scriptsize{$C_1^{(t)}$}};

\draw	(1.25,5) -- (1.25,0)
;

\node [left] at (1.25,5) {\tiny{$0$}};

\node [shape=circle, fill, inner sep=0pt,minimum size=4pt] at (1.25,4.75) {};


\end{tikzpicture}
\endminipage\hfill
\caption{}
\label{fig:CM-structure-fiber-middle-chain}
\end{figure}

\item Assume $i=0$ or $k+1$, and that $V=V_i$ is a ruled surface over an elliptic curve $D$. First, suppose $n>1$, so there is at least one flower meeting $V$. The adjunction formula applied to a general fiber of the ruling on $V$, gives that $C_1$ is a section. Repeatedly blowing up smooth points on $C_1^{(l)}$, gives the result. We illustrate this in Figure~\ref{fig:CM-structure-fiber-end-section}. 
\begin{figure}[H]
\minipage{0.24\textwidth}
\centering
\begin{tikzpicture}[scale = 0.4]
\draw[dashed]	
		(-0,4.75) -- (2.5,4.75)
;

\node [right] at (2.5,4.75) {\scriptsize{$C_1$}};

\draw	(1.25,5) -- (2,3.75)
		(1.25,3) -- (2,4.25)
		(1.25,2) -- (2,0.75)
		(1.25,0) -- (2,1.25)
;

\node at (1.25,2.65) {\tiny{$\vdots$}};

\node [left] at (1.25,5) {\tiny{$-1$}};
\node [left] at (1.25,3) {\tiny{$-2$}};
\node [left] at (1.25,2) {\tiny{$-2$}};
\node [left] at (1.25,0) {\tiny{$-1$}};






\node [right] at (3,2.25) {$\longrightarrow \cdots$};

\end{tikzpicture}
\endminipage \hfill
\minipage{0.23\textwidth}
\begin{tikzpicture}[scale = 0.4]
\draw[dashed]	
		(0,4.75) -- (2.5,4.75)
;

\node [right] at (2.5,2.25) {$\longrightarrow$};
\node [left] at (0,2.25) {$\longrightarrow$};


\draw	(1.25,5) -- (2,3.5)
		(1.25,2.25) -- (2,4)
		(1.25,2.75) -- (2,1)
		(1.25,0) -- (2,1.5)
;

\node [left] at (1.25,5) {\tiny{$-1$}};
\node [left] at (1.25,3) {\tiny{$-2$}};
\node [left] at (1.25,2) {\tiny{$-2$}};
\node [left] at (1.25,0) {\tiny{$-1$}};

\node [shape=circle, fill, inner sep=0pt,minimum size=4pt] at (1.375,4.75) {};

\end{tikzpicture}
\endminipage\hfill
\minipage{0.17\textwidth}
\begin{tikzpicture}[scale = 0.4]
\draw[dashed]	
		(-0,4.75) -- (2.5,4.75)
;

\node [right] at (2.5,2.25) {$\longrightarrow$};


\draw	(1.25,5) -- (2,3.08)
		(1.25,1.41) -- (2,3.58)
		(1.25,1.91) -- (2,0)
;

\node [left] at (1.25,5) {\tiny{$-1$}};
\node [left] at (1.25,1.91) {\tiny{$-2$}};
\node [left] at (2,0) {\tiny{$-1$}};

\node [shape=circle, fill, inner sep=0pt,minimum size=4pt] at (1.348,4.75) {};

\end{tikzpicture}
\endminipage\hfill
\minipage{0.17\textwidth}
\begin{tikzpicture}[scale = 0.4]
\draw[dashed]	
		(-0,4.75) -- (2.5,4.75)
;

\node [right] at (2.5,2.25) {$\longrightarrow$};


\draw	(1.25,5) -- (2,2.25)
		(1.25,0) -- (2,2.75)
		
;

\node [left] at (1.25,5) {\tiny{$-1$}};
\node [left] at (1.25,0) {\tiny{$-1$}};

\node [shape=circle, fill, inner sep=0pt,minimum size=4pt] at (1.318,4.75) {};

\end{tikzpicture}
\endminipage\hfill
\minipage{0.17\textwidth}
\begin{tikzpicture}[scale = 0.4]
\draw[dashed]	
		(-0,4.75) -- (2.5,4.75)
;

\node [right] at (2.5,4.75) {\scriptsize{$C_1^{(t)}$}};

\draw	(1.25,5) -- (1.25,0)
;

\node [left] at (1.25,5) {\tiny{$0$}};

\node [shape=circle, fill, inner sep=0pt,minimum size=4pt] at (1.25,4.75) {};


\end{tikzpicture}
\endminipage\hfill
\caption{}
\label{fig:CM-structure-fiber-end-section}
\end{figure}

Suppose now that $n=1$, then the adjunction formula applied to a general fiber, gives that $C_1$ is a bisection, i.e., $C_1\cdot \fiber = 2$ for any fiber $\fiber$ of the ruling on $V$. By Lemma~\ref{thm:CM-chain-exceptional-curve}, we also know that $C_1$ is elliptic. Hurwitz' formula \cite[Corollary IV.2.4]{Hartshorne} asserts that $C_1\to D$ is \'etale. \note{$(0=0+\sum(e_P-1))$} Repeatedly blowing up smooth points on $C_1^{(l)}$, gives the result. We illustrate this in Figure~\ref{fig:CM-structure-fiber-end-bisection-etale}. 
\note{For Hurwitz formula, we like that characteristic is zero (or otherwise tame ramification, i.e.\ characteristic does not divide $e_P$).} 

\begin{figure}[H]
\minipage{0.24\textwidth}
\centering
\begin{tikzpicture}[scale = 0.4]
\draw[dashed]	(-0,0.25) -- (2.5,0.25)
		(-0,4.75) -- (2.5,4.75)
;

\node [right] at (2.5,0.25) {\scriptsize{$C_1$}};
\node [right] at (2.5,4.75) {\scriptsize{$C_1$}};

\draw	(1.25,5) -- (2,3.75)
		(1.25,3) -- (2,4.25)
		(1.25,2) -- (2,0.75)
		(1.25,0) -- (2,1.25)
;

\node at (1.25,2.65) {\tiny{$\vdots$}};

\node [left] at (1.25,5) {\tiny{$-1$}};
\node [left] at (1.25,3) {\tiny{$-2$}};
\node [left] at (1.25,2) {\tiny{$-2$}};
\node [left] at (1.25,0) {\tiny{$-1$}};






\node [right] at (3,2.25) {$\longrightarrow \cdots$};

\end{tikzpicture}
\endminipage \hfill
\minipage{0.23\textwidth}
\begin{tikzpicture}[scale = 0.4]
\draw[dashed]	(0,0.25) -- (2.5,0.25)
		(0,4.75) -- (2.5,4.75)
;

\node [right] at (2.5,2.25) {$\longrightarrow$};
\node [left] at (0,2.25) {$\longrightarrow$};


\draw	(1.25,5) -- (2,3.5)
		(1.25,2.25) -- (2,4)
		(1.25,2.75) -- (2,1)
		(1.25,0) -- (2,1.5)
;

\node [left] at (1.25,5) {\tiny{$-1$}};
\node [left] at (1.25,3) {\tiny{$-2$}};
\node [left] at (1.25,2) {\tiny{$-2$}};
\node [left] at (1.25,0) {\tiny{$-1$}};

\node [shape=circle, fill, inner sep=0pt,minimum size=4pt] at (1.375,4.75) {};

\end{tikzpicture}
\endminipage\hfill
\minipage{0.17\textwidth}
\begin{tikzpicture}[scale = 0.4]
\draw[dashed]	(-0,0.25) -- (2.5,0.25)
		(-0,4.75) -- (2.5,4.75)
;

\node [right] at (2.5,2.25) {$\longrightarrow$};


\draw	(1.25,5) -- (2,3.08)
		(1.25,1.41) -- (2,3.58)
		(1.25,1.91) -- (2,0)
;

\node [left] at (1.25,5) {\tiny{$-1$}};
\node [left] at (1.25,1.91) {\tiny{$-2$}};
\node [left] at (2,0) {\tiny{$-1$}};

\node [shape=circle, fill, inner sep=0pt,minimum size=4pt] at (1.348,4.75) {};

\end{tikzpicture}
\endminipage\hfill
\minipage{0.17\textwidth}
\begin{tikzpicture}[scale = 0.4]
\draw[dashed]	(-0,0.25) -- (2.5,0.25)
		(-0,4.75) -- (2.5,4.75)
;

\node [right] at (2.5,2.25) {$\longrightarrow$};


\draw	(1.25,5) -- (2,2.25)
		(1.25,0) -- (2,2.75)
		
;

\node [left] at (1.25,5) {\tiny{$-1$}};
\node [left] at (1.25,0) {\tiny{$-1$}};

\node [shape=circle, fill, inner sep=0pt,minimum size=4pt] at (1.318,4.75) {};

\end{tikzpicture}
\endminipage\hfill
\minipage{0.17\textwidth}
\begin{tikzpicture}[scale = 0.4]
\draw[dashed]	(-0,0.25) -- (2.5,0.25)
		(-0,4.75) -- (2.5,4.75)
;

\node [right] at (2.5,0.25) {\scriptsize{$C_1^{(t)}$}};
\node [right] at (2.5,4.75) {\scriptsize{$C_1^{(t)}$}};

\draw	(1.25,5) -- (1.25,0)
;

\node [left] at (1.25,5) {\tiny{$0$}};

\node [shape=circle, fill, inner sep=0pt,minimum size=4pt] at (1.25,4.75) {};


\end{tikzpicture}
\endminipage\hfill
\caption{}
\label{fig:CM-structure-fiber-end-bisection-etale}
\end{figure}

\item Assume $i=0$ or $k+1$ and that $V=V_i$ is a ruled surface over a rational curve~$D$. By Lemma~\ref{thm:CM-chain-exceptional-curve}, the curve $C_1$ is elliptic. The adjunction formula applied to a general fiber $\fiber$ of the ruling on $V$, gives that $C_1\cdot \fiber$ is either $1$ or $2$. But since $C_1$ is elliptic and $D$ is rational, $C_1\cdot \fiber=2$. The Hurwitz formula gives
\[2g(C_1)-2 = 2\cdot(2g(D)-2)+\sum_{P\in C_1}(e_P-1)=-4+\sum_{P\in C_1}(e_P-1),\]
where $e_P$ denotes the ramification index of a point ${P\in C_1}$. Since $C_1$ is elliptic, we find
\[\sum_{P\in C_1}(e_P-1) = 4.\]
The result of \cite[Formula~(4.8)~p.~290]{Liu} asserts that $e_P\leq 2$, for all ${P\in C_1}$, and hence $C_1$ is ramified above exactly four points $\{d_1, d_2, d_3, d_4\}\subset D$. 

Repeatedly blowing up smooth points on $C_1^{(l)}$ gives the result. For an illustration of the process above a point that is not a branch point, we refer to Figure~\ref{fig:CM-structure-fiber-end-bisection-etale}. The evolution of the fiber above a branch point is illustrated in Figure~\ref{fig:CM-structure-fiber-end-bisection-branch}.

\begin{figure}[H]
\minipage{0.24\textwidth}
\centering
\begin{tikzpicture}[scale = 0.4]
\draw[dashed]	
		(-0,4.75) -- (2.5,4.75)
;

\node [right] at (2.5,4.75) {\scriptsize{$C_1$}};

\draw	(1.25,5) -- (2,3.75)
		(1.25,3) -- (2,4.25)
		(1.25,2) -- (2,0.75)
		(1.25,0) -- (2,1.25)
		(1.1,0.25) -- (1.85,1.5)
;

\node at (1.25,2.65) {\tiny{$\vdots$}};

\node [left] at (1.25,5) {\tiny{$-1$}};
\node [left] at (1.25,3) {\tiny{$-2$}};
\node [left] at (1.25,2) {\tiny{$-2$}};
\node [left] at (1.25,-0.25) {\tiny{$-2$}};
\node [left] at (1.25,0.25) {\tiny{$-2$}};






\node [right] at (3,2.25) {$\longrightarrow \cdots$};

\end{tikzpicture}
\endminipage \hfill
\minipage{0.23\textwidth}
\begin{tikzpicture}[scale = 0.4]
\draw[dashed]	
		(0,4.75) -- (2.5,4.75)
;

\node [right] at (2.5,2.25) {$\longrightarrow$};
\node [left] at (0,2.25) {$\longrightarrow$};


\draw	(1.75,5) -- (1,3.08)
		(1.75,1.41) -- (1,3.58)
		(1.75,1.91) -- (1,0)
		(1.55,2.48866) -- (0.8,0.578666)
;

\node [left] at (1.75,5) {\tiny{$-1$}};
\node [left] at (1,3.58) {\tiny{$-2$}};
\node [left] at (1,0) {\tiny{$-2$}};
\node [left] at (0.8,0.578666) {\tiny{$-2$}};

\node [shape=circle, fill, inner sep=0pt,minimum size=4pt] at (1.6484,4.75) {};

\end{tikzpicture}
\endminipage\hfill
\minipage{0.17\textwidth}
\begin{tikzpicture}[scale = 0.4]
\draw[dashed]	
		(-0,4.75) -- (2.5,4.75)
;

\node [right] at (2.5,2.25) {$\longrightarrow$};

\draw	(1.25,5) -- (2.15,1.69118)
		(1.25,0) -- (2.15,3.30882)
		(1,0.83333) -- (1.9,4.13)
;

\node [left] at (1.25,5) {\tiny{$-1$}};
\node [left] at (1.25,0) {\tiny{$-2$}};
\node [left] at (1,0.83333) {\tiny{$-2$}};

\node [shape=circle, fill, inner sep=0pt,minimum size=4pt] at (1.348,4.75) {};

\end{tikzpicture}
\endminipage\hfill
\minipage{0.17\textwidth}
\begin{tikzpicture}[scale = 0.4]
\draw[dashed]	(0,2.5) -- (2.5,2.5);

\node [right] at (2.5,2.25) {$\longrightarrow$};


\draw	(1.25,5) -- (2.15,1.69118)
		(1.25,0) -- (2.15,3.30882)
		
;

\node [left] at (1.25,5) {\tiny{$-1$}};
\node [left] at (1.25,0) {\tiny{$-1$}};

\node [shape=circle, fill, inner sep=0pt,minimum size=4pt] at (1.93,2.5) {};

\end{tikzpicture}
\endminipage\hfill
\minipage{0.17\textwidth}
\begin{tikzpicture}[scale = 0.4]
\draw[dashed] (2.5,4.5) .. controls (0.833,2.5) .. (2.5,0.5);

\node [right] at (2.5,0.25) {\scriptsize{$C_1^{(t)}$}};

\draw	(1.25,5) -- (1.25,0)
;

\node [left] at (1.25,5) {\tiny{$0$}};

\node [shape=circle, fill, inner sep=0pt,minimum size=4pt] at (1.25,2.5) {};


\end{tikzpicture}
\endminipage\hfill
\caption{}
\label{fig:CM-structure-fiber-end-bisection-branch}
\end{figure}

\end{enumerate}

\begin{corollary} \label{thm:CM-relation-L-flower-elliptic}
Let $X$ be a smooth, proper surface over $K$ with $\omega_{X/K}^{\otimes m}\simeq \O_X$ for some $m\geq 1$, and let $\X$ be a Crauder-Morrison model of $X$. Suppose that the subgraph $\Gamma_{min}$ of the dual graph $\Gamma$ is a chain. Let $V_i$ be a component in the chain that is an elliptic, ruled surface.
 Let $L$ be the number of blow-ups in the contraction $V_i\to \overline{V}_i$ to the minimal ruled surface $\overline{V}_i$ and let $\lambda$ be the number of flowercurves on $V_i$ of type $2$.
Then 
\[L\geq 2\lambda.\]
\end{corollary}
\begin{proof}
Let $C_1, \ldots, C_\lambda$ be the flowercurves on $V_i$ of type $2$. For $j=1, \ldots, \lambda$, the curve $C_j$ is a rational curve, and therefore it is contained in a reducible fiber on $V_i$. 
Moreover, all $C_j$ are disjoint, because there are no triple points.
Lemma~\ref{thm:CM-self-intersection-flower-curve} implies that $C_j^2=-2$ for $j=1, \ldots, \lambda$. On the other hand, Proposition~\ref{thm:CM-chain-structure-fibers} describes the structure of the reducible fibers. It follows that $L\geq 2\lambda$.
\end{proof}

\section{Euler characteristics}\label{sect:Euler-characteristics}


In this section, we will compute $\chi(E_i^\circ)$, for certain components $E_i$ of the special fiber of a Crauder-Morrison model. From Proposition~\ref{thm:GMP-ACampo-dvr}, it is clear that this is useful for computing monodromy eigenvalues. 

\begin{lemma} \label{thm:CM-euler-characteristics}
Let $X$ be a smooth, proper surface over $K$ with $\omega_{X/K}^{\otimes m}\simeq \O_X$ for some $m\geq 1$, and let $\X$ be a Crauder-Morrison model of $X$. If $\X$ is a chain degeneration, we may assume that $\X$ satisfies the properties of Proposition~\ref{thm:CM-chain-good}.
 Write the special fiber as $\X_k=\sum_{i\in I}N_i E_i$. Denote by $E_i^\circ=E_i\setminus \left(\bigcup_{j\neq i} E_j\right)$. The following statements hold.
\begin{enumerate}[(i)]
\item If $E_i\simeq \P^2$ is a top of a flower, then $\chi(E_i^\circ)=1$.
\item If $E_i$ is a top of a flower of genus $g\geq 0$, and if $E_i$ is a ruled surface, then $\chi(E_i^\circ) = 2-2g$.
\item If $E_i$ is a component of a flower, but not the top component, then ${\chi(E_i^\circ)=0}$.  
\item If $\X$ is a chain degeneration, and if $E_i$ is a component of the chain, but not an end component of the chain, then $\chi(E_i^\circ)\geq 0$.
\item If $\X$ is a chain degeneration, and if $E_i$ is a component of the chain, but not an end component of the chain, and if the multiplicity $N_i$ is not minimal among the components in the chain, then $\chi(E_i^\circ)=0$.
\item If $\X$ is a chain degeneration, and if $E_i\simeq \P^2$ is an end component of the chain, then $\chi(E_i^\circ) = 3$.
\item If $\X$ is a chain degeneration, and if $E_i$ is an elliptic, ruled surface and an end component of the chain, then $\chi(E_i^\circ) \geq 0$.
\item If $\X$ is a chain degeneration, and if $E_i$ is a rational, minimal ruled surface and an end component of the chain, then $\chi(E_i^\circ)\geq 0$.
%
%
\end{enumerate} 
\end{lemma} 
\begin{proof}
\begin{enumerate}[(i)]
\item Let $E_i\simeq \P^2$ be a top of a flower. We have $\chi(E_i)=\chi(\P^2)=3$. Let $C$ be the double curve on $E_i$, which is either a line or a conic by Theorem~\ref{thm:CM-classification}. In either case, $\chi(C) = 2$.  So 
\[\chi(E_i^\circ) = \chi(E_i)-\chi(C) = 1.\]

\item Let $E_i$ be a top of a flower of genus $g\geq 0$, such that $E_i$ is a ruled surface. Denote by $C$ the double curve on $E_i$. By Remark~\ref{rmk:CM-genus-flower}, the curve~$C$ has genus~$g$, so $\chi(C) = 2-2g$. Moreover, $C$ is a section of $E_i$, and therefore $\chi(E_i) = \chi(\P^1)\cdot\chi(C) = 2\cdot(2-2g)$. We conclude
\[\chi(E_i^\circ) = \chi(E_i)-\chi(C) = 2-2g.\]

\item Let $E_i$ be a component of a flower, but not the top component. Let $C_i$ and $C_{i+1}$ be the two double curves on $E_i$. By Theorem~\ref{thm:CM-classification}, we know that $E_i$ is a minimal ruled surface with sections $C_i$ and $C_{i+1}$. So $\chi(C_i)=\chi(C_{i+1})$ and $\chi(E_i)=\chi(\P^1)\cdot \chi(C_i)= 2\cdot\chi(C_i)$. We conclude
\[\chi(E_i^\circ) = \chi(E_i) - \chi(C_i)-\chi(C_{i+1}) = 0.\]

\item Suppose $\X$ is a chain degeneration, and let $E_i$ be a component of the chain, but not an end component.
Let $L$ be the number of blow-ups in the contraction $E_i\to \overline{E}_i$ to the minimal ruled surface $\overline{E}_i$. We have $\chi(\overline{E}_i)=0$, because $\overline{E}_i$ is an elliptic, minimal ruled surface. Since blowing up increases the topological Euler characteristic by 1, we have $\chi(E_i) = L$. 
\note{
blowing-up is replacing a point by a $\P^1$. This has the effect on Euler characteristic: $-\chi(pt)+\chi(\P^1)=-1+2=1$.
}

Let $C_1$ and $C_2$ be the double curves on $E_i$, where $E_i$ meets another component in the chain. By Theorem~\ref{thm:CM-classification}, the double curves $C_1$ and $C_2$ are sections of $E_i$, and hence elliptic. So $\chi(C_1)=\chi(C_2) = 0$.

Let $C_3, \ldots, C_{\lambda+2}$ be the flowercurves on $E_i$, where $\lambda$ is the number of flowers meeting $E_i$. By Proposition~\ref{thm:CM-chain-flowers}, all flowers meeting $E_i$ are rational of type 2. This means, in particular, that $C_j$ is rational and hence, $\chi(C_j)=2$ for $j=3,\ldots, \lambda+2$. 
Therefore
\[\chi(E_i^\circ) = \chi(E_i) - \sum_{j=1}^{\lambda+2} \chi(C_j) = L - 2 \lambda\geq 0,\]
by Corollary~\ref{thm:CM-relation-L-flower-elliptic}.

\note{For every fiber, let $n$ be the number of components, $m$ be the number of components with self-intersection $-2$ and $p$ be the number of flowercurves in the fiber. So $p\leq m\leq n$ and $n=m+2$ and $p\leq (m+1)/2$ (no triple points). So $n\geq 2p+1$. Let $q$ be the number of reducible fibers on $E_i$. Then we have $L_i = \sum_{j=1}^q (n_j-1) = \sum n_j - q$ and hence $L_i\geq \sum_{j=1}^q  (2p_j + 1) - q =2 \sum_{j=1}^q p_j = 2\lambda_i$.}

\item Suppose $\X$ is a chain degeneration, and let $E_i$ be a component of the chain, but not an end component and  such that the multiplicity $N_i$ is not minimal among the components in the chain. By Proposition~\ref{thm:CM-chain-good}, the component~$E_i$ is an elliptic, minimal ruled surface, and doesn't meet any flowers. Therefore, there are exactly two double curves $C_1$ and $C_2$ and both $C_1$ and $C_2$ are sections on $E_i$. Since $E_i$ is an elliptic, minimal ruled surface, we have that $\chi(E_i)=0$. Moreover, $C_1$ and $C_2$ are elliptic curves, and hence $\chi(C_1)=\chi(C_2)=0$. We conclude
\[\chi(E_i^\circ) = \chi(E_i)-\chi(C_1)-\chi(C_2) = 0.\]

\item Suppose $\X$ is a chain degeneration, and let $E_i\simeq \P^2$ be an end component of the chain. So $\chi(E_i)=3$. By Proposition~\ref{thm:CM-chain-flowers}, the component $E_i$ does not meet any flowers and hence, there is exactly one double curve $C$. By Theorem~\ref{thm:CM-classification}, the curve $C$ is elliptic. So $\chi(C)=0$. Therefore
\[\chi(E_i^\circ) = \chi(E_i)-\chi(C) =3.\]

\item Suppose $\X$ is a chain degeneration, and let $E_i$ be an elliptic, ruled surface and an end component of the chain.
Let $L$ be the number of blow-ups in the contraction $E_i\to \overline{E}_i$ to the minimal ruled surface $\overline{E}_i$. We have $\chi(\overline{E}_i)=0$, because $\overline{E}_i$ is an elliptic, minimal ruled surface. Since blowing up increases the topological Euler characteristic by 1, we have $\chi(E_i) = L$. \note{blowing-up is replacing a point by a $\P^1$. This has the effect on Euler characteristic: $-\chi(pt)+\chi(\P^1)=-1+2=1$.}

 Let $C_1$ be the intersection curve of $E_i$ with a component in the chain. By Theorem~\ref{thm:CM-classification}, $C_1$ is a section of $E_i$ and hence elliptic, so $\chi(C_1) = 0$.
 
  By Proposition~\ref{thm:CM-chain-flowers}, all flowers meeting $E_i$ are rational of type 2 or non-rational of type 4. Denote by $\lambda$ the number of flowers of type 2 and by $\lambda'$ the number of flowers of type $4$. Let $C_2, \ldots, C_{\lambda+1}$ be rational flowercurves on $E_i$ of type 2. We have, in particular, that $\chi(C_j)= 2$ for $j=2,\ldots, \lambda+1$. Let $C_{\lambda+2}, \ldots, C_{\lambda+\lambda'+1}$ be non-rational flowercurves on $E_i$ of type 4. So $\chi(C_j)\leq 0$ for $j=\lambda+2,\ldots, \lambda+\lambda'+1$. 
Therefore
\[\chi(E_i^\circ) = \chi(E_i) - \sum_{j=1}^{\lambda+\lambda'+1} \chi(C_j) \geq L - 2 \lambda\geq 0,\]
by Corollary~\ref{thm:CM-relation-L-flower-elliptic}.


\item 
Suppose $\X$ is a chain degeneration, and let $E_i$ be a rational, minimal ruled surface and an end component of the chain.  We have $\chi({E}_i)=4$.
Let $C_1$ be the intersection of $E_i$ with a component in the chain. By Theorem~\ref{thm:CM-classification}, the curve $C_1$ is elliptic. So $\chi(C_1) = 0$. 

In \cite[Section~V.4]{BHPV}, it is explained that, on a rational, minimal ruled surface, there is at most one smooth curve  with strictly negative self-intersection.
By Proposition~\ref{thm:CM-chain-flowers}, all flowers meeting $E_i$ are rational of type 2 and therefore, Lemma~\ref{thm:CM-self-intersection-flower-curve} implies that they have strictly negative self-intersection.
Therefore, there is at most one flowercurve on $E_i$.

If there is no flowercurve on $E_i$, then
\[\chi(E_i^\circ)=\chi(E_i)-\chi(C_1) = 4.\]
If there is a flowercurve $C_2$ on $E_i$, then it is rational, and hence $\chi(C_2)=2$. Therefore,
\[\chi(E_i^\circ)=\chi(E_i)-\chi(C_1) - \chi(C_2) = 2.\]

\end{enumerate}
\end{proof}

\cleardoublepage


\chapter{Poles of the motivic zeta function}\label{ch:poles}

In this chapter, we will prove the following theorem, which has been announced in \cite{Jaspers}.

\begin{theorem*} 
Let $X$ be a $K3$ surface over $K$ with Crauder-Morrison model $\X$ of $X$. Write the special fiber as $\X_k=\sum_{i\in I} N_i E_i$. Let $\omega$ be a volume form on $X$ and let $(N_i,\nu_i)$ be the numerical data of $E_i$. Let $\rho_i = \nu_i/N_i +1$ be the weight of $E_i$ for every $i\in I$.

Define $I^\dagger\subset I$ to be the set of indices $i\in I$, where either
\begin{enumerate}[(i)]
\item $\rho_i$ is minimal, or
\item $E_i$ is the top of a conic-flower.
\end{enumerate}
Define $S^\dagger=\left\{(-\nu_i,N_i)\in \Z\times \Z_{>0}\mid i\in I^\dagger\right\}$. We have
\[Z_{X, \omega}(T)\in \mathcal{M}_k^{\hat{\mu}}\left[T, \frac{1}{1-\L^{a}T^b}\right]_{(a,b)\in S^\dagger}.\]
Moreover, the poles of $Z_{X,\omega}(T)$ are precisely $\left\{-\nu_i/N_i\mid i\in I^\dagger\right\}$.
\end{theorem*}

The main tool in the proof of this theorem is the Denef-Loeser formula (Theorem~\ref{thm:Denef-Loeser-CY}):
\[Z_{X, \omega}(T)=\sum_{\emptyset \neq J \subseteq I} (\L-1)^{|J|-1}[\widetilde{E^\circ_J}]\prod_{j\in J}\frac{\L^{-\nu_j}T^{N_j}}{1-\L^{-\nu_j}T^{N_j}}.\]
So if we define ${S} =\left\{(-\nu_i,N_i)\in \Z\times \Z_{>0}\mid i\in I\right\}$, then it holds that
\[Z_{X, \omega}(T)\in \mathcal{M}_k^{\hat{\mu}}\left[T, \frac{1}{1-\L^{a}T^b}\right]_{(a,b)\in {S}}.\]
In this chapter, we will find all elements that can be omitted from $S$.

The factor $[\widetilde{E_i^\circ}]$ that appears in the Denef-Loeser formula needs to be computed for certain components $E_i$ in the special fiber $\X_k$. We will do so in Section~\ref{sect:E_i^circ}. In Section~\ref{sect:contributions-flowers}, we will define the \emph{contribution of flowers to the motivic zeta function} and explain how we compute these contributions by writing Python code. As a consequence, we get the first part of the theorem:
\begin{equation} \label{eq:poles-zeta-in-ring}
Z_{X, \omega}(T)\in \mathcal{M}_k^{\hat{\mu}}\left[T, \frac{1}{1-\L^{a}T^b}\right]_{(a,b)\in S^\dagger}.
\end{equation}
In the final section, we will first define a \emph{pole} of a rational function over $\mathcal{M}_k^{\hat{\mu}}$, and then we will compute the poles of the motivic zeta function $Z_{X, \omega}(T)$. As a result, we get that $S^\dagger$ is the `smallest' subset of $S$ such that \eqref{eq:poles-zeta-in-ring} holds, i.e., no elements of $S^\dagger$ can be omitted such that \eqref{eq:poles-zeta-in-ring} still holds.

\subsubsection*{Notation}
In this chapter, we fix an algebraically closed field $k$ of characteristic zero. Put $R=k\llbracket t\rrbracket$ and $K=k(\!( t )\!)$. Let $X$ be a $K3$ surface over $K$ and let $\X$ be a Crauder-Morrison model of $X$ with special fiber $\X_k=\sum_{i\in I} N_i E_i$.

 For any local ring $\O$ with maximal ideal $\m$, we denote by $\widehat{\O}$ the completion of $\O$ with respect to $\m$. For any reduced ring $A$, we denote by $Q(A)$ the total ring of fractions of $A$, and by $A'$ the integral closure of $A$ in $Q(A)$.

\section{Computation of \texorpdfstring{$[\widetilde{E_J^\circ}]$ in $\mathcal{M}_{\lowercase{k}}^{\hat{\mu}}$}{[tilde(EJ)]} }\label{sect:E_i^circ}

Recall that for a non-empty subset $J\subseteq I$, we defined 
\[E_J = \bigcap_{j\in J} E_j \qquad \text{and} \qquad E_J^\circ = E_J \setminus \left(\bigcup_{i\in I\setminus J} E_i\right).\]
Notice that if $|J|\geq 3$, then ${E_J}=\emptyset$, because $\X$ is triple-point-free. Furthermore, when $|J|=2$, we have $E_J^\circ = E_J$.

We defined $\widetilde{E^\circ_J}$ as the finite \'etale cover $E_J^\circ \times_\X \mathcal{Y}_J$ of $E_J^\circ$, where $\mathcal{Y}_J$ is the normalization of 
$\X\times_{R} R[\pi]/(\pi^{N_J}- t),$
with $N_J = \lcm_{ j\in J} N_j$. The group $\mu_{N_j}$ acts on $\widetilde{E^\circ_J}$ via its action on $R[\pi]/(\pi^{N_J}-t)$. This $\mu_{N_J}$-action induces a $\hat{\mu}$-action on $\widetilde{E^\circ_J}$.
In \cite[Lemma~3.2.2]{BultotNicaise}, an alternative description of $\widetilde{E^\circ_J}$ is given:  
\begin{lemma}[{\cite[Lemma~3.2.2]{BultotNicaise}}]
\[\widetilde{E^\circ_J} \simeq E_J^\circ \times_\X \mathcal{Y},\]
where $\mathcal{Y}$ is the normalization of 
$\X\times_{R} R[\pi]/(\pi^{n}- t),$
for any multiple $n$ of $N_J$.
\end{lemma}
Moreover, the $\hat{\mu}$-action on $\widetilde{E^\circ_J}$ induced by the $\mu_n$-action on $R[\pi]/(\pi^{n}- t)$ coincides with the $\hat{\mu}$-action induced by the $\mu_{N_J}$-action on $R[\pi]/(\pi^{N_J}- t)$.
In practice, we will often choose $n=\lcm_{i\in I} N_i$.

The aim of this section is to compute the class $[\widetilde{E_J^\circ}]\in \mathcal{M}_k^{\hat{\mu}}$, for certain $E_J$. To be precise, we will prove the following proposition.

\begin{proposition*}[Proposition~\ref{thm:poles-tilde-E_i-top-flower} and Proposition~\ref{thm:poles-tilde-E_i-middle-flower}]
Let $X$ be a $K3$ surface over $K$ and let $\X$ be a Crauder-Morrison model of $X$ with special fiber $\X_k=\sum_{i\in I} N_i E_i$. Let $N_0 F_0+N_1 F_1+\cdots + N_{\ell} F_\ell$ be a flower in $ \X_k$, where $F_i\cap F_j=\emptyset$ if and only if $j\not\in \{i-1,i, i+1\}$. The component $F_\ell$ meets $\Gamma_{min}$ in $F_{\ell+1}$. Denote by $C_j$ the intersection $F_{j-1}\cap F_j$ for $j=1, \ldots, \ell+1$.
 The following relations hold in $\mathcal{M}_k^{\hat{\mu}}$.
\begin{enumerate}[(i)]
\item If $F_0\simeq \P^2$ and $C_1$ is a line, then there exists a $k$-variety $\mathcal{P}$ with a good $\hat{\mu}$-action, such that $[\widetilde{F_0^\circ}]=[\mathcal{P}]\L^2$ and $[\widetilde{C_1}] = [\mathcal{P}](\L+1)$,
\item If $F_0$ is a minimal ruled surface, then  $[\widetilde{F_0^\circ}]=[\widetilde{C_1}]\L$,
\item For any $j=1,\ldots, \ell$, we have $[\widetilde{F_j^\circ}]=[\widetilde{C_j}](\L-1)$ and $[\widetilde{C_j}]=[\widetilde{C_{j+1}}]$
\end{enumerate}
\end{proposition*}

In Subsection~\ref{sect:poles-tilde-local-computation}, we give some local computations and in Subsection~\ref{sect:poles-tilde-structure-Y}, we describe the geometrical structure of $\Y$ and the morphism $\Y\to \X$. Finally, we will prove the proposition in Subsections \ref{sect:poles-tilde-top-flower} and \ref{sect:poles-tilde-middle-flower}.

\subsection{Local computations} \label{sect:poles-tilde-local-computation}

Let $n= \lcm_{i\in I} N_i$ and let $R_n$ be the unique totally ramified extension $R[\pi]/(\pi^{n}- t)$ of $R$ of degree $n$. We then have $R_n\simeq k\llbracket \pi\rrbracket$. \note{Unique since $k$ is algebraically closed.}  Denote by $\mathcal{Y}$ the normalization of $\X\times_R R_n$ and let 
$f\colon \Y \to \X$ be
the induced morphism, which is finite.
\note{
As $\Spec R_n\to \Spec R$ is finite and finiteness is stable under base change, we have that $\X\times_R R_n\to \X$ is a finite morphism. As a normalization morphism is also finite \cite[Exercise II.3.8]{Hartshorne}, we obtain that $f \colon \Y \to \mathcal{X}$ is finite.
}
Let $p'$ be a closed point in the special fiber $\Y_k=\Y\times_{R_n} k$ and let $p=f(p')$.

In this subsection, we will prove the following lemma.
\begin{lemma*}[Lemma~\ref{thm:poles-technical-one-branch} and Lemma~\ref{thm:poles-technical-two-branches}]\hfill 
\begin{enumerate}[(i)]
\item Suppose $p\in E_i^\circ$. Define $n'\in \Z$ such that $n=n'N_i$. Then the morphism $\widehat{\O}_{\mathcal{X},p} \to \widehat{\O}_{\Y,p'}$ is described by
\begin{align*}
R\llbracket x,y,z\rrbracket/(t-x^{N_i}) & \to R_n\llbracket x,y,z\rrbracket/(\pi^{n'}-\xi x)\\
g(x, y ,z, t) & \mapsto g(x,y,z, \pi^n),
\end{align*} 
for some $N_i$-th root of unity $\xi$.
\item Suppose $p \in E_i \cap E_j$ and define $c=\gcd(N_i, N_j)$ and $a,b,n'\in \Z$ such that $N_i=ca$, $N_j=cb$ and $n=cn'$ respectively. Let $d\in \Z$ be the integer such that $n'=abd$.
Then the morphism $\widehat{\O}_{\mathcal{X},p} \to \widehat{\O}_{\Y,p'}$ is described by
\begin{align*}
R\llbracket x,y,z\rrbracket/(t-x^{N_i}y^{N_j})& \to R_n\llbracket x_0,y_0,z_0\rrbracket/(\pi^d-\xi x_0y_0)\\
g(x, y ,z, t) & \mapsto g(x_0^{b},y_0^{a},z_0, \pi^n),
\end{align*} 
for some $c$-th root of unity $\xi$.
\end{enumerate}
\end{lemma*}

To prove this lemma, we first need some results from commutative algebra.

\subsubsection{Some commutative algebra}

\begin{lemma} \label{thm:poles-technical-decomposition-complete-local-ring}
\[\Y\times_\X \Spec(\widehat{\O}_{\X,p})=\coprod_{f(q)=p} \Spec (\widehat{\O}_{\Y,q}).\]
\end{lemma}

\begin{proof}
This is an immediate consequence of \cite[Corollary 7.6]{Eisenbud}.

\note{
Let $\Spec A$ be an affine open neighborhood of $p$. As $f$ is finite, there exists a ring $B$ with $f^{-1}(\Spec A)=\Spec B$. We then have the following carthesian diagram:
\[\xymatrix{
 & \Spec\ \widehat{\O}_{\X,p}\ar[d]\\
\Spec B\ar[r]\ar[d]& \Spec A\ar[d]\\
\Y\ar[r] & \X.
}\]
It is sufficient to show that $\coprod_{f(q)=p} \Spec(\widehat{\O}_{\Y,q})=\Spec B \times_{\Spec A} \Spec(\widehat{\O}_{\X,p}),$
which is equivalent to
\[B \otimes_A \widehat{A}_\p=\prod_{f^{-1}(\q)=\p} \widehat{B}_\q,\]
where $\p$ is the maximal ideal of $A$ corresponding to $p$ and $B$ gets the $A$-algebra structure by $f\colon A \to B$. 

Since $f$ is finite, $B$ is a finite $A$-algebra and hence, by \cite[Corollary 7.6]{Eisenbud}, the ring $B\otimes_A \widehat{A}_\p$ has only finitely many maximal ideals $\m_i$. Moreover,
\[B\otimes_A \widehat{A}_\p=\prod_i (B\otimes_A \widehat{A}_\p)_{\m_i},\]
where $(B\otimes_A \widehat{A}_\p)_{\m_i}$ is the localization of $B\otimes_A \widehat{A}_\p$ with respect to the complement of the maximal ideal $\m_i$.
Since $B\otimes_A \widehat{A}_\p= \widehat{B}\otimes_A A_\p$, \emph{$\widehat{A}_\p\otimes_A B=(\widehat{A}\otimes_A A_\p) \otimes_A B=A_\p \otimes _A (\widehat{A}\otimes_A B)= A_\p \otimes_A \widehat{B}$} the maximal ideals of $B\otimes_A \widehat{A}_\p$ are exactly the maximal ideals of $\widehat{B}_{f(\p)\widehat{B}}$. These are the maximal ideals $\m$ of $\widehat{B}$ with $f(\p)\widehat{B}\subset \q$, which is equivalent to $f^{-1}(\q)=\p$. 
We conclude that
\[B\otimes_A \widehat{A}_\p=\prod_{f^{-1}(\q)=\p} \widehat{B}_\q.\]}
\end{proof}

\begin{lemma} \label{thm:poles-technical-norm-commute-tensor}
\[\Y\times_\X \Spec(\widehat{\O}_{\X,p})=\Spec({R_n\otimes_R\widehat{\O}_{\X,p}})'.\]

\end{lemma}

\begin{proof}
This follows from the fact that the completion of a Noetherian, integrally closed ring is again integrally closed.
\end{proof}

\subsubsection{Point on a single component}

Let $E_i$ be an irreducible component of the special fiber $\X_k$, with multiplicity~$N_i$. Let $p'$ be a closed point of $\Y_k$ such that $p=f(p')\in E_i^\circ$, i.e., $p\in E_i$ is not contained in any other irreducible component $E_j$ with $j\neq i$.

\begin{lemma} \label{thm:poles-technical-one-branch}
The morphism $\widehat{\O}_{\mathcal{X},p} \to \widehat{\O}_{\Y,p'}$ is described by
\begin{align*}
R\llbracket x,y,z\rrbracket/(t-x^{N_i}) & \to R_n\llbracket x,y,z\rrbracket/(\pi^{n'}-\xi x)\colon\\
g(x, y ,z, t) & \mapsto g(x,y,z, \pi^n),
\end{align*} 
for some $N_i$-th root of unity $\xi$  and where $n'\in \Z$ is defined such that $n=n'N_i$.
\end{lemma} 
\begin{proof}
Because $p$ is contained in $E_i$, but not in any other component of the special fiber, we have
\[\widehat{\O}_{\X,p} = R\llbracket x,y,z\rrbracket/(t-x^{N_i}).\]
Recall that $n'\in \Z$ is defined such that $n=n'N_i$.
By Lemma~\ref{thm:poles-technical-decomposition-complete-local-ring} and Lemma~\ref{thm:poles-technical-norm-commute-tensor}, it is clear that we should compute $(R_n\otimes_R \widehat{\O}_{\X,p})'$.
Consider the ring 
\[A=(R_n\otimes_R \widehat{\O}_{\X,p})\simeq R_n\llbracket x,y,z\rrbracket /(\pi^n-x^{N_i}).\]
Define $B$ as
\[B= \prod_{\xi \in \mu_{N_i}} R_n\llbracket x,y,z\rrbracket/(\pi^{n'}-\xi x),\]
where $\mu_{N_i}$ is the group of $N_i$-th roots of unity in $k$.
There is an inclusion
$A  \hookrightarrow B$, because $\pi^n-x^{N_i}=\prod_{\xi\in \mu_{N_i}} (\pi^{n'}-\xi x)$. Since $B$ is integral over $A$, there is also an inclusion $B\hookrightarrow Q(A)$.
\note{
Let $\phi_i\colon B\to B_i$ be a projection for $i\in I$ and let $\phi\colon B\to \prod_{i\in I} B_i$ be an injection. Then $\prod B_i$ is contained in $K(B)$. Indeed we show that $(b_1, 0, \ldots, 0)$ is integral over $B$. Take $b\in B$ such that $\phi_1(b)=b_1$. Then $b_1$ satisfies the equation $X^2-b X$. Indeed $(b_1, 0 , \ldots, 0)^2=(b_1^2, 0 , \ldots, 0)$ and $b\cdot (b_1, 0 , \ldots, 0) = (b_1, b_2 , \ldots, b_n)\cdot (b_1, 0 , \ldots, 0)=(b_1^2, 0 , \ldots, 0)$.
}
Therefore, $A'=B'$.

For every $N_i$-th root of unity, the ring $R_n\llbracket x,y,z\rrbracket/(\pi^{n'}-\xi x)$ is regular, because 
it is isomorphic to $R_n\llbracket y,z \rrbracket$.
\note{
The isomorphism is \begin{align*}
R_n\llbracket x,y,z\rrbracket/(\pi^{n'}-\xi x)&\to R_n\llbracket y,z \rrbracket \colon\\
f(x,y,z)+(\pi^{n'}-\xi x) &\mapsto f(\xi^{N_i-1}\pi^{n'},y,z)
\end{align*}
}
Hence $R_n\llbracket x,y,z\rrbracket/(\pi^{n'}-\xi x)$ is normal and therefore, $B$ is normal too. So
\[(R_n\otimes_R \widehat{\O}_{\X,p})' \simeq \prod_{\xi \in \mu_{N_i}} R_n\llbracket x,y,z\rrbracket/(\pi^{n'}-\xi x).\]
Lemma~\ref{thm:poles-technical-decomposition-complete-local-ring} implies that 
\[\widehat{\O}_{\Y,p'} = R_n\llbracket x,y,z\rrbracket/(\pi^{n'}-\xi x),\]
for some $N_i$-th root of unity $\xi$ and the result follows.
\end{proof}

\subsubsection{Point on the intersection of two components}

Let $E_i$ and $E_j$ be irreducible components of $\X_k$, with respective multiplicities $N_i$ and $N_j$. Let $p'$ be a closed point of $\Y_k$ such that $p=f(p') \in E_i \cap E_j$. Define $c=\gcd(N_i, N_j)$ and $a,b,n'\in \Z$ such that $N_i=ca$, $N_j=cb$ and $n=cn'$.

\begin{lemma}\label{thm:poles-technical-two-branches}
The morphism $\widehat{\O}_{\mathcal{X},p} \to \widehat{\O}_{\Y,p'}$ is described by
\begin{align*}
R\llbracket x,y,z\rrbracket/(t-x^{N_i}y^{N_j})& \to R_n\llbracket x_0,y_0,z_0\rrbracket/(\pi^d-\xi x_0y_0)\colon\\
g(x, y ,z, t) & \mapsto g(x_0^{b},y_0^{a},z_0, \pi^n),
\end{align*} 
for some $c$-th root of unity $\xi$ and where $d\in \Z$ is defined such that $n'=abd$. 
\note{Since $N_i\mid n$ and $N_j\mid n$, we can write $n=eN_i=fN_j$. So $n'=ae=bf$. Since $c=\gcd(N_i,N_j)$, we have $\gcd(a,b)=1$, and hence $a\mid f$ and $b\mid e$. We conclude $n'=dab$ for some $d\in \Z$.} 
\end{lemma}

\begin{proof}
Because $p\in E_i\cap E_j$, we have $\widehat{\O}_{\X,p}=R\llbracket x,y,z\rrbracket/(t-x^{N_i}y^{N_j})$.
By Lemma~\ref{thm:poles-technical-decomposition-complete-local-ring} and Lemma~\ref{thm:poles-technical-norm-commute-tensor}, it is clear that we should compute $(R_n\otimes_R \widehat{\O}_{\X,p})'$. Consider the ring
\[A=R_n\otimes_R \widehat{\O}_{\X,p} \simeq R_n\llbracket x,y,z\rrbracket /(\pi^n-x^{N_i}y^{N_j}).\]
Define $B$ as
\[B = \prod_{\xi \in \mu_c} R_n\llbracket x, y ,z\rrbracket /(\pi^{n'} -\xi x^ay^b),\]
where $\mu_{c}$ is the group of $c$-th roots of unity in $k$.
There is an inclusion
$A  \hookrightarrow B$, because $\pi^n-x^{N_i}y^{N_j}=\prod_{\xi\in \mu_c} (\pi^{n'} -\xi x^ay^b)$. Since $B$ is integral over $A$, there is also an inclusion $B\hookrightarrow Q(A)$.
Therefore, $A'=B'$.

Since $\gcd(a, b) = 1$, there exists an
integer $d\in \Z$ with $n' = abd$.
Consider the $R_n$-algebra homomorphism
\begin{align*}
R_n\llbracket x,y,z\rrbracket/(\pi^{n'}-\xi x^ay^b) & \to R_n\llbracket x_0,y_0,z_0\rrbracket/(\pi^d-\xi x_0y_0)\colon\\
g(x,y,z, \pi) &\mapsto g(x_0^b, y_0^a, z_0,\pi).
\end{align*}

This morphism is injective and integral, and it induces an isomorphism of fraction fields. 
Because the ring $R_n\llbracket x_0,y_0,z_0\rrbracket/(\pi^d-\xi x_0y_0)$ describes a toric singularity, it is normal. Therefore, $R_n\llbracket x_0,y_0,z_0\rrbracket/(\pi^d-\xi x_0y_0)$ is the normalization of $R_n\llbracket x,y,z\rrbracket/(\pi^{n'}-\xi x^ay^b)$. 
So
\[(R_n\otimes_R \widehat{\O}_{\X,p})'  = \prod_{\xi \in \mu_c} R_n\llbracket x_0,y_0,z_0\rrbracket/(\pi^d-\xi x_0y_0).\]

Lemma~\ref{thm:poles-technical-norm-commute-tensor} implies that 
\[\widehat{\O}_{\Y,p'} = R_n\llbracket x_0,y_0,z_0\rrbracket/(\pi^d-\xi x_0y_0),\]
for some $c$-th root of unity $\xi$ and the result follows.
\end{proof}

\subsection{Structure of \texorpdfstring{$\mathcal{Y}_k$}{Yk}} \label{sect:poles-tilde-structure-Y}

\begin{lemma} \label{thm:poles-local-comp-smoothness}
Let $f\colon\Y\to \X$ be the finite morphism constructed before. Let $p'\in \Y_k$ be a closed point, and let $p=f(p')$. Then the following properties hold.
\begin{enumerate}[(i)]
\item If $p$ belongs to a unique irreducible component of $\X_k$, or if $p$ belongs to two irreducible components of $\X_k$ with multiplicities $N_i$ and $N_j$ such that $n=\lcm(N_i, N_j)$, then $\Y$ is regular  at $p'$. The converse also holds.
\item If $p$ belongs to a unique irreducible component of $\X_k$, then $p'$ belongs to a unique irreducible component of $\Y_k$.
\item If $p$ belongs to two distinct irreducible components of $\X_k$, then $p'$ belongs to two distinct irreducible components of $\Y_k$.
\item The irreducible components of $\Y_k$ are smooth. 
\end{enumerate}
\end{lemma}

\begin{proof}
\begin{enumerate}[(i)]
\item By \cite[Lemma~2.26]{Liu}, we have that $\Y$ is regular at $p'$ if and only if $\widehat{\O}_{\Y,p'}$ is regular.

Since $\X$ is triple-point-free, $p$ is either contained in a unique irreducible component $E_i$ of the special fiber $\X_k$, or $p$ is contained in $E_i\cap E_j$ for some $i,j\in I$.

Suppose first that $p$ is contained in a unique irreducible component $E_i$ of the special fiber $\X_k$. Let $N_i$ be the multiplicity of $E_i$ in $\X_k$. By Lemma~\ref{thm:poles-technical-one-branch}, we have that $\widehat{\O}_{\Y,p'}\simeq R_n\llbracket x,y,z\rrbracket/(\pi^{n'}-\xi x)$
for some $N_i$-th root of unity~$\xi$  and where $n'$ is defined as $n=n'N_i$. This is a regular ring, so $\Y$ is regular at $p'$.

Suppose now that $p\in E_i\cap E_j$ for some $i,j\in I$. Let $N_i$ and $N_j$ be the respective multiplicities of $E_i$ and $E_j$ in $\X_k$. Define $c=\gcd(N_i, N_j)$ and $d\in \Z$ such that $nc = N_i N_j d$. By Lemma~\ref{thm:poles-technical-two-branches}, we have $\widehat{\O}_{\Y,p'}\simeq R_n\llbracket x_0,y_0,z_0\rrbracket/(\pi^d-\xi x_0y_0)$ 
for some $c$-th root of unity $\xi$. This ring is regular if and only if $d=1$, which is equivalent to $n=\lcm(N_i, N_j)$.
\note{Let $F=\pi^d-\xi x_0y_0$. Then $\partial F/\partial x = -y$, $\partial F/\partial y = -x$ and $\partial F/\partial z = d z^{d-1}$. So $(0,0,0)$ is a solution of $F=\partial F/\partial x=\partial F/\partial y=\partial F/\partial z=0$ if and only if $d\neq 1$. This means that it is regular if and only if $d=1$.}
\note{$d=1\iff n'=ab \iff n/c = N_i/c \cdot N_j/c \iff cn=N_iN_j$}

\item Suppose $p$ belongs to a unique irreducible component of $\X_k$ of multiplicity~$N_i$. Because of Lemma~\ref{thm:poles-technical-one-branch}, we know that $\widehat{\O}_{\Y,p'}\simeq R_n\llbracket x,y,z\rrbracket/(\pi^{n'}-\xi x)$, for some $N_i$-th root of unity $\xi$ and where $n'\in \Z$ is defined such that $n=n'N_i$. Since $\Y_k \simeq \Y \otimes_{R_n} k \simeq  \Y \otimes_{R_n} R_n/(\pi)$, we have
\[\widehat{\O}_{\Y_k,p'} \simeq \widehat{\O}_{\Y,p'}\otimes_{R_n} R_n/(\pi)\simeq k\llbracket x,y,z\rrbracket /(\xi x) \simeq k\llbracket y,z\rrbracket,\]
which is a domain. 
\note{
The first equality is explained as follows: since $R_n/(\pi)$ is finite over $R_n$, we must have $\widehat{\O}_{\Y,p'}\otimes_{R_n} R_n/(\pi)$ is finite over $\widehat{\O}_{\Y,p'}$. Finite modules over complete rings are complete as well, so the universal property of completion gives $\widehat{\O}_{\Y,p'}\otimes_{R_n} R_n/(\pi) \simeq \left({\O}_{\Y,p'}\otimes_{R_n} R_n/(\pi)\right)^{\widehat{\ }}\simeq \widehat{\O}_{\Y_k,p'}$.
}
Therefore, $\O_{\Y,p'}$ is a domain as well, which implies there is a unique irreducible component of $\Y_k$ passing through $p'$, by \cite[Proposition~2.4.12]{Liu}.

\note{If $\widehat{A}$ is a domain for some local Noetherian ring $A$ with maximal ideal $\m$, then $A$ is a domain as well. Indeed, by Krull intersection theorem~\cite[Corollary~5.4]{Eisenbud}, we have $\bigcap_{i\in \Z_{>0}} \m^i = 0$. This implies that the completion morphism $A\to \widehat{A}=\varprojlim A/\m^i$ is injective. }

\item Suppose $p$ belongs to two irreducible components of $\X_k$. Because of Lemma~\ref{thm:poles-technical-two-branches}, we know that $\widehat{\O}_{\Y,p'}\simeq R_n\llbracket x_0,y_0,z_0\rrbracket/(\pi^d-\xi x_0y_0)$, for some root of unity $\xi$ and where $d\in \Z$ is defined such that $n \gcd(N_i, N_j) = N_i N_j d$. Since $\Y_k \simeq \Y \otimes_{R_n} k \simeq  \Y \otimes_{R_n} R_n/(\pi)$, we have
\[\widehat{\O}_{\Y_k,p'} \simeq \widehat{\O}_{\Y,p'}\otimes_{R_n} R_n/(\pi)\simeq k\llbracket x_0,y_0,z_0\rrbracket /(\xi x_0y_0).\]
This means that there are at most two irreducible components of $\Y_k$ passing through $p'$. 
\note{Intuitively: By looking at $\widehat{\O}_{\Y,p'}$, we can only detect more irreducible components than by looking at $\O_{\Y,p'}$. By \cite[Proposition~2.4.7.(b)]{Liu}, the minimal prime ideals of $\widehat{\O}_{\Y,p'}$, correspond bijectively to the irreducible components of $\Spec \widehat{\O}_{\Y,p'}$. And \cite[Proposition~2.4.12.(a)]{Liu} says that the irreducible components of $\Spec {\O}_{\Y,p'}$ correspond bijectively to the irreducible components of $\Y$ passing through $p'$.}

On the other hand, the going-down theorem~\cite[Lemma 10.4.34]{Liu} says that the irreducible components passing through $p$ can be lifted, via the finite morphism $f$, to irreducible components passing through $p'$. Therefore, we conclude that there are exactly two irreducible components of $\Y_k$ passing through $p'$.

\item Let $E$ be an irreducible component of $\Y_k$, and take $q'\in E$. 
If $E$ is the only irreducible component of $\Y_k$ passing through $q'$, then $\O_{E,q'}=\O_{\Y_k, q'}$. \note{The local equation of $E$ at $q'$ is the same as the local equation of $\Y_k$ at $q'$.} So in that case, we have
\[\widehat{\O}_{E,q'} \simeq k\llbracket y, z\rrbracket,\]
which implies $E$ is smooth at $q'$.
If there is another irreducible component of $\Y_k$ passing through $q'$, then 
\[\widehat{\O}_{\Y,q'} \simeq k \llbracket x_0,y_0,z_0 \rrbracket /(x_0y_0).\]
As there are exactly two components passing through $q'$ by (iii), both components must be smooth at $q'$.
\end{enumerate}
\end{proof}

\begin{lemma} \label{thm:poles-local-comp-etale-branch}
Let $f\colon\Y\to \X$ be the finite morphism constructed before. Write $\X_k = \sum_{i\in I} N_i E_i$. Let 
\[f_{E_i}\colon f^{-1}(E_i)_{red} \to E_i,\]
be the morphism induced by $f$. The following properties hold.
\begin{enumerate}[(i)]
\item The morphism $f_{E_i}$ is \'etale of degree $N_i$ above $E_i^\circ$.
\item Suppose $E_i$ and $E_j$ intersect non-trivially and let $C$ be an irreducible component of $E_i\cap E_j$. Let $C'$ be an irreducible component of $f^{-1}(C)_{red}$. Then $C'$ is smooth, and the ramification index of $f_{E_i}$ at $C'$ is $N_i/\gcd(N_i,N_j)$.
\end{enumerate}
\end{lemma}

\begin{proof}
\begin{enumerate}[(i)]
\item Let $p'\in \Y_k$ be a closed point such that $p=f(p')\in E_i^\circ$. The map $\widehat{\O}_{\X,p} \to \widehat{\O}_{\Y,p'}$ is described by
\begin{align*}
R\llbracket x,y,z\rrbracket/(t-x^{N_i}) & \to R_n\llbracket x,y,z\rrbracket/(\pi^{n'}-\xi x)\colon\\
g(x, y ,z,t) & \mapsto g(x,y,z, \pi^n),
\end{align*} 
where $n'=\frac{\lcm_{j\in I} N_j}{N_i}$.
Therefore the map $\widehat{\O}_{E_i,p}\to \widehat{\O}_{f^{-1}(E_i)_{red}, p'}$ is described by
\begin{align*}
k\llbracket x,y,z\rrbracket/(x) &\to k\llbracket x,y,z \rrbracket / (\xi x) \colon\\
g(x,y,z)&\mapsto g(x,y,z).
\end{align*}

\note{
We tensor the morphisms $\widehat{\O}_{\X,p} \to \widehat{\O}_{\Y,p'}$ and $R/(t)\to R_n/(\pi)$. Moreover we can reduce this morphism since $k\llbracket x,y,z \rrbracket / (\xi x)$ is reduced.
}
We conclude that $f_{E_i}$ is \'etale at $p'$, by \cite[Proposition~4.3.26]{Liu}. Moreover, $f_{E_i}$ has degree $N_i$, since the fiber of $f_{E_i}$ above $p$ consists of $N_i$ points, by Lemma~\ref{thm:poles-technical-decomposition-complete-local-ring} and Lemma~\ref{thm:poles-technical-one-branch}.

\item Let $p'\in \Y_k$ be a closed point such that $p=f(p')\in E_i\cap E_j$. The morphism $\widehat{\O}_{\X,p} \to \widehat{\O}_{\Y,p'}$ is described by
\begin{align*}
R\llbracket x,y,z\rrbracket/(t-x^{N_i}y^{N_j})& \to R_n\llbracket x_0,y_0,z_0\rrbracket/(\pi^d-\xi x_0y_0)\colon\\
g(x, y ,z, t) & \mapsto g(x_0^{N_j/\gcd(N_i,N_j)},y_0^{N_i/\gcd(N_i,N_j)},z_0, \pi^n),
\end{align*}
where $d\in \Z$ is defined such that $(\lcm_{l\in I} N_l) (\gcd(N_i, N_j)) = N_i N_j d$.
Suppose that the ideal $(x)$ describes $E_i$ and $(y)$ describes $E_j$. Then \[\widehat{\O}_{E_i,p} \simeq \widehat{\O}_{\X,p}/(x) \simeq R\llbracket x,y,z\rrbracket /(x, t-x^{N_i}y^{N_j})\simeq k\llbracket y,z\rrbracket.\]
Similarly
\[\widehat{\O}_{f^{-1}(E_i)_{red}, p'}\simeq \left(R_n\llbracket x_0,y_0,z_0\rrbracket /(x_0, \pi^d-\xi x_0 y_0)\right)_{red}\simeq k\llbracket y_0, z_0\rrbracket.\]
Therefore, the morphism  $\widehat{\O}_{E_i,p}\to \widehat{\O}_{f^{-1}(E_i)_{red}, p'}$ is described by
\begin{align*}
k\llbracket y,z\rrbracket &\to k\llbracket y_0,z_0 \rrbracket \colon\\
g(y,z)&\mapsto g(y_0^{N_i/\gcd(N_i,N_j)},z_0).
\end{align*}
We see that the ramification index of $f_{E_i}$ at $p'$ is $N_i/\gcd(N_i,N_j)$.

To prove that $C'$ is smooth at $p'$, we compute 
\[\widehat{\O}_{C',p'}\simeq \left(\widehat{\O}_{f^{-1}(E_i)_{red}, p'}/\left(y_0^{N_i/\gcd(N_i,N_j)}\right)\right)_{red} \simeq k\llbracket z_0\rrbracket,\]
and the result follows.
\end{enumerate}
\end{proof}

The following lemma will be used in the proof of Proposition~\ref{thm:poles-tilde-E_i-top-flower} (ii) and Proposition~\ref{thm:poles-tilde-E_i-middle-flower}.

\begin{lemma} \label{thm:poles-base-change-ruled}
Let $f\colon\Y\to \X$ be the finite morphism constructed before. Write $\X_k = \sum_{i\in I} N_i E_i$. Suppose that $E_j$ is a minimal ruled surface and that any double curve~$C$ on $E_j$ is a section of the ruling. Assume moreover that for any fiber $\fiber$ of the ruling on $E_j$, we have that $f^{-1}(\fiber)_{red}$ is a disjoint union of smooth, rational curves. Then $f^{-1}(E_j)_{red}$ is a disjoint union of minimal ruled surfaces. 

Moreover, let $E'$ be an irreducible component of $f^{-1}(E_j)_{red}$, and let $C$ be a double curve on $E_j$. Exactly one irreducible component of $f^{-1}(C)_{red}$ is contained in $E'$, and it is a section of the ruling on~$E'$.
\end{lemma}

\begin{proof}
It follows from Lemma~\ref{thm:poles-local-comp-smoothness}, that $f^{-1}(E_j)_{red}$ is a disjoint union of smooth surfaces. Let $E'$ be an irreducible component of $f^{-1}(E_j)_{red}$. We will show that $E'$ is a minimal ruled surface. 
Consider the finite morphism $\phi=(f_{E_j})|_{E'}\colon E'\to E_j$. 
Let $\fiber$ be a fiber of the ruling on $E_j$, and let $\fiber'$ be a component of $\phi^{-1}(\fiber)_{red}$. 
\note{
Since $f^{-1}(E_j)_{red}$ is a ramified Galois cover of $E_j$, the Galois group acts on $f^{-1}(E_j)$ such that the quotient of this action is $E_j$. Moreover, on every fiber of $f^{-1}(E_j)\to E_j$ the action is transitive. This means that the curve $C$ has an inverse image in every component of $f^{-1}(E_j)$. Moreover, all components of $f^{-1}(E_j)$ are isomorphic, since the action sends a component to another component and since the action is transitive, every component can be reached by applying the action sufficiently many times.
}
Since $\fiber$ is transversal to the double curves, the generic point $\xi$ of $\fiber$ is contained in $E_j^\circ$. Because $\phi$ is \'etale above $E_j^\circ$ by Lemma~\ref{thm:poles-local-comp-etale-branch}, we can write $\phi^*\fiber = \fiber_1'+\cdots +\fiber_r'$, where $\fiber'_i$ is a smooth, rational curve for $i=1,\ldots, r$ and all $\fiber_i'$ are disjoint. We can assume $\fiber'=\fiber'_1$. On the other hand, we have $\phi_*\fiber' = m\fiber$, where $m$ is the degree of the map $\fiber'\to \fiber$. By the projection formula \cite[Theorem~9.2.12]{Liu}, we get
\begin{align*}
\fiber'\cdot \fiber' &= \fiber' \cdot (\fiber_1'+\cdots + \fiber_r')\\
	&= \fiber' \cdot \phi^*\fiber\\
	&=\phi_*\fiber'\cdot \fiber\\
	&= m\fiber\cdot \fiber\\
	&= 0.
\end{align*}
Since $E'$ contains a smooth, rational curve with self-intersection 0, it is a ruled surface, by \cite[Proposition~V.4.3]{BHPV}.

Since $\fiber'\cdot \fiber'=0$, \cite[Proposition~V.4.3]{BHPV} implies that $\fiber'$ is an irreducible fiber of the ruling on $E'$. As $\fiber$ was chosen arbitrarily, we conclude that any component of the inverse image of a fiber on $E_j$ is an irreducible fiber of the ruling on $E'$.

We will now show that $E'$ is \emph{minimal} ruled. Let $\mathcal{G}$ be a component of a fiber of the ruling on $E'$. By the projection formula
\[\phi_*\mathcal{G} \cdot \fiber = \mathcal{G} \cdot \phi^*\fiber =0.\]
This implies that $\phi(\mathcal{G})$ is a fiber of the ruling on $E_j$. Since any component of the inverse image of a fiber on $E_j$ is an irreducible fiber of the ruling on $E'$, we know that $\mathcal{G}$ is an irreducible fiber, and we conclude that $E'$ is minimal ruled.

To prove the second claim in the statement, let $C$ be a double curve on $E_j$, and let $\mathcal{G}$ be a fiber of the ruling on $E'$. As proven before, $\phi_*\mathcal{G}$ is a fiber of the ruling on $E_j$, and hence the projection formula yields
\[\phi^*C \cdot \mathcal{G} = C \cdot \phi_*\mathcal{G} = 1.\]
This implies that $\phi^{-1}(C)_{red}$ is irreducible, and that it is a section.

\end{proof}

\subsection{Top of a flower} \label{sect:poles-tilde-top-flower}

\begin{proposition}\label{thm:poles-tilde-E_i-top-flower}
Let $X$ be a $K3$ surface over $K$ and let $\X$ be a Crauder-Morrison model of $X$ with special fiber $\X_k=\sum_{i\in I} N_i E_i$. Let $N_0 F_0+N_1 F_1+\cdots + N_{\ell} F_\ell$ be a flower in $ \X_k$, where $F_i\cap F_j=\emptyset$ if and only if $j\not\in \{i-1, i,i+1\}$. Denote by $C_1$ the intersection $F_0\cap F_1$.
 The following relations hold in $\mathcal{M}_k^{\hat{\mu}}$.
\begin{enumerate}[(i)]
\item If $F_0\simeq \P^2$ and $C_1$ is a line, then there exists a $k$-variety $\mathcal{P}$ with a good $\hat{\mu}$-action, such that $[\widetilde{F_0^\circ}]=[\mathcal{P}]\L^2$ and $[\widetilde{C_1}] = [\mathcal{P}](\L+1)$,
\item If $F_0$ is a minimal ruled surface, then  $[\widetilde{F_0^\circ}]=[\widetilde{C_1}]\L$,
\end{enumerate}
\end{proposition}

\begin{proof}
Let $f\colon\Y\to \X$ be the finite morphism constructed before.
\begin{enumerate}[(i)]
\item Assume that $F_0\simeq \P^2$ and that $C_1$ is a line. From Table~\ref{table:flowers-P2-line}, it follows that $\gcd(N_0,N_1)=N_0$.

Since $N_0/\gcd(N_0,N_1)=1$, Lemma~\ref{thm:poles-local-comp-etale-branch} gives that $f$ is \'etale of degree $N_0$ above $F_0$.
Since $F_0\simeq \P^2$ is simply connected, we have
\[f^{-1}(F_0)_{red} = \bigsqcup_{i=1}^{N_0} G_i,\]
where $G_i\simeq \P^2$.

Let $g\colon f^{-1}(F_0)_{red}\to f^{-1}(F_0)_{red}$ be the automorphism of $f^{-1}(F_0)_{red}$ induced by a generator of $\Gal(K(n)/K)$. Since the image of an irreducible component is irreducible, there exists an $\alpha_j\in \{1, \ldots, N_0\}$ for every $j\in \{1, \ldots, N_0\}$, such that $g(G_j) = G_{\alpha_j}$.

Set $\mathcal{P}=\bigsqcup_{i=1}^{N_0} P_i$, where $P_i=\Spec (k)$ is a point.  Define a $\mu_n$-action on $\mathcal{P}$ by $g\cdot P_j = P_{\alpha_j}$. Since every preimage of a point under $f$ consists of exactly $N_0$ points, there are $\mu_n$-equivariant isomorphisms
\[f^{-1}(F_0)_{red} \simeq \mathcal{P} \times_k \P^2,\]
with trivial $\mu_n$-action on $\P^2$, and
\[f^{-1}(C_1)_{red} \simeq \mathcal{P} \times_k \P^1,\]
with trivial $\mu_n$-action on $\P^1$.

We conclude that
\[ [\widetilde{F_0^\circ}]=[\mathcal{P}]\L^2,\]
and 
\[ [\widetilde{C_1}] = [\mathcal{P}](\L+1).\]

\note{definition simply connected: \cite[Example~IV.2.5.3]{Hartshorne}. In \cite[exercise~IV.2.1]{Hartshorne}, it is stated that $\P^n$ is simply connected for any $n$. (answer to the exercise on mathoverflow.) Other reference: Debarre 4.18.}

\item Assume that $F_0$ is a minimal ruled surface. In that case, $C_1$ is a section of the ruling on $F_0$. From Table~\ref{table:flowers-ruled}, it follows that $\gcd(N_0,N_1)=N_0$. Lemma~\ref{thm:poles-local-comp-etale-branch} implies that $f_{F_0}$ is \'etale of degree $N_0$ above $F_0$.

Let $\fiber$ be a fiber of the ruling on $F_0$.
Since $\fiber\simeq \P^1$ is simply connected, $f^{-1}(\fiber)_{red}$ is a disjoint union of smooth, rational curves. Lemma~\ref{thm:poles-base-change-ruled} guarantees that $f^{-1}(F_0)_{red}$ is a disjoint union of minimal ruled surfaces and that every component of $f^{-1}(F_0)_{red}$ contains exactly one component of $f^{-1}(C_1)_{red}$, which is a section. 

Therefore, there is a map
\[\psi\colon f^{-1}(F_0)_{red} \to f^{-1}(C_1)_{red},\]
consistent with the rulings on the components of $f^{-1}(F_0)_{red}$. So, all fibers of this map are isomorphic to $\P^1$. Furthermore, we have a commutative diagram
\[
    \xymatrix{f^{-1}(F_0)_{red} \ar[r] \ar[d]_\psi & F_0 \ar[d] \\
              f^{-1}(C_1)_{red} \ar[r]& C_1 }
\]
where the map $F_0\to C_1$ defines the ruling on $F_0$.
Notice that for any closed point $x\in f^{-1}(C_1)_{red}\subset f^{-1}(F_0)_{red}$, we have $\psi(x)=x$.

We will now show that $\psi$ is $\mu_n$-equivariant. Let $x$ be a closed point of $f^{-1}(F_0)_{red}$ and fix $g\in \mu_n$. Set $y=g\cdot x$. Notice that $f(x)=f(y)$. Let $E_x$ and $E_y$ be the irreducible components of $f^{-1}(F_0)_{red}$ containing $x$ and $y$ respectively. 
Let $\fiber_x$ be the fiber of the ruling on $E_x$, containing $x$ and similarly, let $\fiber_y$ be the fiber of the ruling on $E_y$ containing $y$.
From the proof of Lemma \ref{thm:poles-base-change-ruled}, it follows that $f(\fiber_x)$ and $f(\fiber_y)$ are fibers of the ruling on $F_0$. Since $f(x)=f(y)$, we must have $f(\fiber_x) = f(\fiber_y)$. Therefore, for any point $x'\in \fiber_x$, we have that $g\cdot x'\in \fiber_y$. 

Let $x'$ be the unique point in the intersection of $\fiber_x$ with the component of $f^{-1}(C_1)_{red}$ contained in $E_x$. Since $\psi$ is consistent with the ruling, we have that $\psi(x)=\psi(x')$. Moreover, since $x'\in f^{-1}(C_1)_{red}$, we have $\psi(x')=x'$. So
\[g\cdot \psi(x) = g\cdot x'.\]
On the other hand, we know that $g\cdot x'\in \fiber_y$, and hence $\psi(g\cdot x')=\psi(y)$. Because $g\cdot x' \in f^{-1}(C_1)_{red}$, we know that $\psi(g\cdot x') = g\cdot x'$. Therefore
\[\psi(y) = g\cdot x'.\]
We conclude that $\psi(g\cdot x) = \psi(y)= g\cdot \psi(x)$. This implies that $\psi$ is $\mu_n$-equivariant.

From the fact that $\psi$ is $\mu_n$-equivariant, it follows that $\widetilde{F_0^\circ}=f^{-1}(F_0)\setminus f^{-1}(C_1)$ is an affine bundle over $f^{-1}(C_1)_{red}$ of rank 1. So $[\widetilde{F_0^\circ}]=[\widetilde{C_1}]\L$.
\end{enumerate}
\end{proof}

Ignoring the $\hat{\mu}$-action, we also have a relation for the top of a conic-flower.
\begin{proposition}\label{thm:poles-conic-flower}
Let $X$ be a $K3$ surface over $K$ and let $\X$ be a Crauder-Morrison model of $X$ with special fiber $\X_k=\sum_{i\in I} N_i E_i$. Let $N_0 F_0+N_1 F_1+\cdots + N_{\ell} F_\ell$ be a flower in $ \X_k$, where $F_i\cap F_j=\emptyset$ if and only if $j\not\in \{i-1, i, i+1\}$. Denote by $C_1$ the intersection $F_0\cap F_1$.
We have \[[\widetilde{F_0^\circ}]=\frac{N_0}{2}(\L+1)\L
\qquad \text{and} \qquad  
[\widetilde{C_1}] = \frac{N_0}{2}(\L+1),\]
in $\mathcal{M}_k$.
\end{proposition} 
\begin{proof}

From Table~\ref{table:flowers-P2-conic}, it follows that $N_0$ is even and that $\gcd(N_0,N_1)=N_0/2$.

Let $f\colon\Y\to \X$ be the finite morphism constructed before.
By Lemma~\ref{thm:poles-local-comp-etale-branch}, we know that
$f_{F_0}$
is \'etale of degree $N_0$ above $F_0^\circ=F_0\setminus C_1$, and ramified of index 2 at $f^{-1}(C_1)$. Therefore, $f^{-1}(C_1)_{red}\to C_1$ is \'etale of degree $N_0/2$. Since $C_1\simeq \P^1$ is simply connected, $\widetilde{C_1}$ has $N_0/2$ irreducible components, all isomorphic to $\P^1$. In particular, we have
\[[\widetilde{C_1}] = \frac{N_0}{2}(\L+1).\]
 
Let $E$ be a component of $f^{-1}(F_0)_{red}$, and let $C$ be a component of $f^{-1}(C_1)$ contained in $E$. 
\note{
Since $f^{-1}(E)_{red}$ is a ramified Galois cover of $E$, the Galois group acts on $f^{-1}(E)$ such that the quotient of this action is $E$. Moreover, on every fiber of $f^{-1}(E)\to E$ the action is transitive. This means that the curve $C_1$ has an inverse image in every component of $f^{-1}(E)$. Moreover, all components of $f^{-1}(E)$ are isomorphic, since the action sends a component to another component and since the action is transitive, every component can be reached by applying the action sufficiently many times.
} 
The map $E\to F_0$ is a cyclic covering, branched along the conic~$C_1$. Let $d$ be the degree of this covering.
\note{
definition of cyclic covering of degree $n$ in \cite[Section~16 and 17]{BHPV}
} 
\note{
$d$ is even: take a point $x\in C_1$. since the ramification index is 2, $e_{x/y}=2$ for every $y\in f^{-1}(C_1)$. By \cite[formula~(7.4.8), p.290]{Liu}, $d$ is even.
}
By \cite[Section~1.17]{BHPV}, there exists a unique line bundle $\mathcal{L}$ on $F_0\simeq\P^2$ with $\O_{\P^2}(C_1) = \mathcal{L}^{\otimes d}$. As $C_1$ is a conic, we get $\mathcal{L}^{\otimes d}=\O_{\P^2}(2)$, and hence $d=2$. This means that $E$ contains exactly one component of $f^{-1}(C_1)$ and therefore, $f^{-1}(F_0)_{red}$ has exactly $N_0/2$ irreducible components.

We will now prove that $E\simeq \P^1\times \P^1$. As mentioned before, $E$ is the double cyclic cover of $\P^2$, branched along a conic. Therefore, it is isomorphic to the quadric surface in $\P^3 = \mathrm{Proj}\, k[x:y:z:w]$ defined by the equation $w^2= f(x:y:z)$, where $f$ is a homogeneous polynomial of degree $2$, defining the conic $C_1$ in $\P^2$. Because all quadric surfaces in $\P^3$ are isomorphic, $E$ is isomorphic to the image of $\P^1\times \P^1$ under the Segre embedding $\P^1\times \P^1\to \P^3$. We conclude that $E\simeq \P^1\times \P^1$.

Therefore
\[[\widetilde{F_0^\circ}] = \frac{N_0}{2}([\P^1\times \P^1]-[\P^1]) = \frac{N_0}{2}\L(\L+1).\] 
\end{proof}

\begin{remark} \label{rmk:poles-conic-flower}
We don't know whether $[\widetilde{F_0^\circ}]=[\widetilde{C_1}]\L$ also holds in $\mathcal{M}_k^{\hat{\mu}}$. Let $E\simeq \P^1\times \P^1$ be an irreducible component of $f^{-1}(F_0)_{red}$, and let $C$ be the component of $f^{-1}(C_1)_{red}$ contained in $E$, then it does not seem likely that $E\setminus C$ is an equivariant affine bundle over $C$, because the two fibrations on $E$ are interchanged under the $\hat{\mu}$-action.
\end{remark}

\subsection{Middle component of a flower} \label{sect:poles-tilde-middle-flower}

\begin{proposition}\label{thm:poles-tilde-E_i-middle-flower}
Let $X$ be a $K3$ surface over $K$ and let $\X$ be a Crauder-Morrison model of $X$ with special fiber $\X_k=\sum_{i\in I} N_i E_i$. Let $N_0 F_0+N_1 F_1+\cdots + N_{\ell} F_\ell$ be a flower in $ \X_k$, where $F_i\cap F_j=\emptyset$ if and only if $j\not\in \{i-1,i, i+1\}$. The component $F_\ell$ meets $\Gamma_{min}$ in $F_{\ell+1}$.
Denote by $C_j$ the intersection $F_{j-1}\cap F_j$ for $j=1, \ldots, \ell+1$.
For any $j=1, \ldots, \ell$, we have 
\[ [\widetilde{F_j^\circ}]=[\widetilde{C_j}](\L-1)\quad \text{ and }\quad[\widetilde{C_j}]=[\widetilde{C_{j+1}}],\]
in $\mathcal{M}_k^{\hat{\mu}}$.
\end{proposition}

We first need a lemma.

\begin{lemma}\label{thm:poles-base-change-fiber}
Fix  $j\in \{1, \ldots, \ell\}$.
 Let $\fiber$ be a fiber of the minimal ruled surface $F_j$, and let $f\colon\Y\to \X$ be the finite morphism constructed before. Then we have that $f^{-1}(\fiber)_{red}$ is a disjoint union of smooth, rational curves.
\end{lemma}

\begin{proof}
Let $ p \in \fiber $ be a closed point and let $ E $ be a component of $ f^{-1}(F_j)_{red} $. Let $ p' \in E $ be a point such that $f(p')=p$, and let $\fiber'$ be an irreducible component of $ f^{-1}(\fiber)_{red} $ passing through $p'$.

Suppose first that $ p \in F_j^{\circ} $. From Lemma~\ref{thm:poles-technical-one-branch}, it follows that the map $ {\widehat{\O}_{F_j, p} \to \widehat{\O}_{E, p'}} $ is given by
\begin{align*}
k\llbracket y,z\rrbracket& \to k\llbracket y_0,z_0\rrbracket\colon g(y,z) \mapsto g(y_0,z_0).
\end{align*}
Because $\fiber$ is smooth, \cite[Corollary 4.2.12]{Liu} gives that the ideal defining $\fiber$ at $p$ is generated by a non-zero element $ h \in \mathfrak{m} \setminus \mathfrak{m}^2 $, where $ \mathfrak{m}=(y,z) $ is the maximal ideal of $k\llbracket y,z\rrbracket $. Thus $\fiber'$ is generated by $h'\in k\llbracket y_0, z_0 \rrbracket$ with $h' \in (y_0, z_0)\setminus (y_0,z_0)^2$ and therefore $\fiber'$ is smooth at $p'$. Furthermore, $ f^{-1}(\fiber)_{red} \to \fiber $ is \'etale at $p'$, by \cite[Proposition 4.3.26]{Liu}.

Let $i=j$ or $j+1$. Assume now that $ p \in C_i$.
From Lemma~\ref{thm:poles-technical-one-branch}, it follows that the map $ \widehat{\O}_{F_j, p} \to \widehat{\O}_{E, p'} $ is given by
\begin{align*}
k\llbracket y,z\rrbracket &\to k\llbracket y_0,z_0\rrbracket\colon g(y,z) \mapsto g\left(y_0^{N_j/\gcd(N_i,N_j)},z_0\right).
\end{align*}
Let $(y)$ be the ideal defining $C_i\subset F_j$ at $p$. Let $C'$ be an irreducible component of $f^{-1}(C_i)_{red}$ passing through $p'$.

Because $\fiber$ is smooth, \cite[Corollary 4.2.12]{Liu} gives that the ideal defining $\fiber$ at $p$ is generated by a non-zero element $ h \in \mathfrak{m} \setminus \mathfrak{m}^2 $, where $ \mathfrak{m}=(y,z) $ is the maximal ideal of $k\llbracket y,z\rrbracket $. Since $\fiber$ and $C_i$ meet transversally at $p$, the maximal ideal $(y,z)$ must be generated by $h$ and $y$. Hence $ h \equiv z + ay \mod \mathfrak{m}^2 $ for some $ a \in k$. Thus $\fiber'$ is generated by $h'\in k\llbracket y_0, z_0 \rrbracket$ with $h' \equiv z_0 + a y_0^{N_j/\gcd(N_i,N_j)} \mod (y_0,z_0)^2$ and therefore $\fiber'$ is smooth at $p'$.
As a consequence, $ f^{-1}(\fiber)_{red} $ is a disjoint union of smooth curves. Moreover, $\fiber'$ and $C'$ meet transversally at $p'$, because $C'$ is defined by the ideal $(y_0)$. 

We now show that $\fiber'$ is rational. Consider the morphism $\phi=f|_{\fiber'}\colon\fiber'\to\fiber$. By Lemma~\ref{thm:poles-local-comp-etale-branch}, we know that $\phi$ is only ramified above the points $p_j=\fiber\cap C_j$ and $p_{j+1}=\fiber\cap C_{j+1}$. Hurwitz' formula \cite[Theorem 7.4.16]{Liu} gives that
\[2 g(\fiber') - 2 = n(2g(\fiber) - 2) +(e_{p_{j}}-1) + (e_{p_{j+1}}-1),\]
where $n$ is the degree of $\phi$, and $e_{p_j}$ and $e_{p_{j+1}}$ are the ramification indices of $p_j$ and $p_{j+1}$ respectively. Since $\fiber$ is a fiber of $F_j$, it is rational. Furthermore, \cite[Equation (4.8) p. 290]{Liu} gives that $e_{p_j}=e_{p_{j+1}}=n$. We conclude that $g(\fiber')=0$ and hence $f^{-1}(\fiber)_{red}$ is a disjoint union of smooth, rational curves. 
\end{proof}

\textit{Proof of Proposition~\ref{thm:poles-tilde-E_i-middle-flower}}.
Let $f\colon\Y\to \X$ be the finite morphism constructed before.
By Lemma~\ref{thm:poles-base-change-ruled} and Lemma~\ref{thm:poles-base-change-fiber}, we know that $f^{-1}(F_j)_{red}$ is a disjoint union of minimal ruled surfaces. Moreover, any component of $f^{-1}(F_j)_{red}$  contains exactly one component of $f^{-1}(C_j)_{red}$ and $f^{-1}(C_{j+1})_{red}$ and they are sections of the ruling.

Therefore, there is a map
\[\psi\colon f^{-1}(F_j)_{red} \to f^{-1}(C_j)_{red},\]
consistent with the rulings on the components of $f^{-1}(F_j)_{red}$. So, all fibers of this map are isomorphic to $\P^1$. Furthermore, we have a commutative diagram
\[
    \xymatrix{f^{-1}(F_j)_{red} \ar[r] \ar[d]_\psi & F_j \ar[d] \\
              f^{-1}(C_j)_{red} \ar[r]& C_j }
\]
where the map $F_j\to C_j$ defines the ruling on $F_j$. Notice that for any closed point $x\in f^{-1}(C_j)_{red}\subset f^{-1}(F_j)_{red}$, we have $\psi(x)=x$.

We will now show that $\psi$ is $\mu_n$-equivariant. Let $x$ be a closed point of $f^{-1}(F_j)_{red}$ and denote $y=g\cdot x$, for some $g\in \mu_n$. Notice that $f(x)=f(y)$. Let $E_x$ and $E_y$ be the irreducible components of $f^{-1}(F_j)_{red}$ containing $x$ and $y$ respectively. 
Let $\fiber_x$ be the fiber of the ruling on $E_x$ containing  $x$ and similarly, let $\fiber_y$ be the fiber of the ruling on $E_y$ containing $y$.
From the proof of Lemma~\ref{thm:poles-base-change-ruled}, it follows that $f(\fiber_x)$ and $f(\fiber_y)$ are fibers of the ruling on $F_j$. Since $f(x)=f(y)$, we must have $f(\fiber_x) = f(\fiber_y)$. Therefore, for any point $x'\in \fiber_x$, we have that $g\cdot x'\in \fiber_y$. 

Let $x'$ be the unique point in the intersection of $\fiber_x$ with the component of $f^{-1}(C_j)_{red}$ contained in $E_x$. Since $\psi$ is consistent with the ruling, we have that $\psi(x)=\psi(x')$. Moreover, since $x'\in f^{-1}(C_j)_{red}$, we have $\psi(x')=x'$. So
\[g\cdot \psi(x) = g\cdot x'.\]
On the other hand, we know that $g\cdot x'\in \fiber_y$, and hence $\psi(g\cdot x')=\psi(y)$. Because $g\cdot x' \in f^{-1}(C_j)_{red}$, we know that $\psi(g\cdot x') = g\cdot x'$. Therefore
\[\psi(y) = g\cdot x'.\]
We conclude that $\psi(g\cdot x)=\psi(y)=g\cdot \psi(x)$. 
This implies that $\psi$ is $\mu_n$-equivariant.

In particular, when $x\in f^{-1}(C_{j+1})_{red}$, the previous argument implies that 
\[[\widetilde{C_{j}}]=[\widetilde{C_{j+1}}].\] 
\note{$\widetilde{C_{j+1}}_{red}=f^{-1}(C_{j+1})_{red}$.}

From the fact that $\psi$ is $\mu_n$-equivariant, it follows that $f^{-1}(F_j)_{red}\setminus f^{-1}(C_j)_{red}$ is an affine bundle over $f^{-1}(C_j)_{red}$ of rank 1. So
\[ [f^{-1}(F_j)] - [\widetilde{C_j}] = [\widetilde{C_j}]\L.\]
This means that
\[[\widetilde{F_j^\circ}]=[f^{-1}(F_j)] - [\widetilde{C_j}]-[\widetilde{C_{j+1}}]=[\widetilde{C_j}](\L-1).\]
\qed

\section{Contribution of a flower to the motivic zeta function} \label{sect:contributions-flowers}

\begin{definition}
Let $X$ be a $K3$ surface over $K$ with Crauder-Morrison model~$\X$. Let $\omega$ be a volume form on $X$.
Let $N_0 F_0+N_1 F_1+\cdots + N_{\ell} F_\ell$ be a flower in $ \X_k$, where $F_i\cap F_j=\emptyset$ if and only if $j\not\in \{i-1,i, i+1\}$. The component $F_\ell$ meets $\Gamma_{min}$ in $F_{\ell+1}$. Denote by $C_{j+1}$ the intersection $F_{j}\cap F_{j+1}$, for $j=0, \ldots, \ell$.
Let $(N_j, \nu_j)$ be the numerical data of $F_j$, for every $j=0, \ldots, \ell+1$.

We define the \emph{contribution $Z_F(T)\in \mathcal{M}_k^{\hat{\mu}}\llbracket T \rrbracket$ of the flower $F$} to the motivic zeta function $Z_{X,\omega}(T)$ as
\begin{align*}
Z_F(T) = \sum_{j=0}^{\ell}&\Bigg([\widetilde{F_j^\circ}]\frac{\L^{-\nu_j}T^{N_j}}{1-\L^{-\nu_j}T^{N_j}}\\
&+(\L-1)[\widetilde{C_{j+1}}]\frac{\L^{-\nu_j}T^{N_j}}{1-\L^{-\nu_j}T^{N_j}}\frac{\L^{-\nu_{j+1}}T^{N_{j+1}}}{1-\L^{-\nu_{j+1}}T^{N_{j+1}}}\Bigg).
\end{align*}
\end{definition}

\begin{lemma} \label{thm:poles-sum-contribution-flowers}
Let $X$ be a $K3$ surface over $K$ allowing a Crauder-Morrison model $\X$ with special fiber $\X_k=\sum_{i\in I} N_i E_i$. Let $\omega$ be a volume form on $X$ and let $(N_i, \nu_i)$ be the numerical data of $E_i$, for every $i\in I$. Let $\rho_i = \nu_i/N_i +1$ be the weight of $E_i$, for every $i\in I$.

Define  $I_{min} = \{i\in I\mid \rho_i \text{ is minimal}\}$ and 
\[Z_{min}(T)=\sum_{\emptyset \neq J \subseteq I_{min}} (\L-1)^{|J|-1}[\widetilde{E^\circ_J}]\prod_{j\in J}\frac{\L^{-\nu_j}T^{N_j}}{1-\L^{-\nu_j}T^{N_j}}.\]

We have
\[Z_{X,\omega}(T) = Z_{min}(T) + \sum Z_F(T),\]
where the sum runs over all flowers in $\X_k$.
\end{lemma}

\begin{proof}
This immediately follows from Theorem~\ref{thm:Denef-Loeser-CY} and Theorem~\ref{thm:CM-classification}.
\end{proof}

\begin{theorem} \label{thm:poles-contribution-flowers}
Let $X$ be a $K3$ surface over $K$ allowing a Crauder-Morrison model $\X$. Let $\omega$ be a volume form on $X$.
Let $N_0 F_0+N_1 F_1+\cdots + N_{\ell} F_\ell$ be a flower in $ \X_k$, where $F_i\cap F_j=\emptyset$ if and only if $j\not\in \{i-1,i, i+1\}$. The component $F_\ell$ meets $\Gamma_{min}$ in $F_{\ell+1}$. Denote by $C_{j+1}$ the intersection $F_{j}\cap F_{j+1}$ for $j=0, \ldots, \ell$.
Let $(N_j, \nu_j)$ be the numerical data of $F_j$, for every $j=0, \ldots, \ell+1$.

If $F_0\simeq \P^2$ with $F_0\cap F_1$ a line, or if $F_0$ is a minimal ruled surface, then
\[Z_{F}(T)\in \mathcal{M}_k^{\hat{\mu}}\left[T, \frac{1}{1-\L^{-\nu_{\ell+1}}T^{N_{\ell+1}}}\right].\]

If $F_0\simeq \P^2$ with $F_0\cap F_1$ a conic, then 
\[Z_{F}(T)\in \mathcal{M}_k^{\hat{\mu}}\left[T, \frac{1}{1-\L^{-\nu_{0}}T^{N_{0}}}, \frac{1}{1-\L^{-\nu_{\ell+1}}T^{N_{\ell+1}}}\right].\]
\end{theorem}

\begin{proof}
By Theorem~\ref{thm:CM-classification-flowers} and Lemma~\ref{thm:CM-flowers-numerical-relations}, we can compute $(N_j, \nu_j)$ in function of $(N_0,\nu_0)$, for every $j=1, \ldots, \ell+1$. When $F_0\simeq \P^2$ with $F_0\cap F_1$ a conic, we will express $(N_j, \nu_j)$ in function of $(N_1,\nu_1)$ instead of $(N_0,\nu_0)$ to avoid fractions. By the results in Section~\ref{sect:E_i^circ}, we know that
\begin{enumerate}[(i)]
\item if $F_0\simeq \P^2$ and $F_0\cap F_1$ is a line, then there exists a $k$-variety $\mathcal{P}$ with good $\hat{\mu}$-action, such that $[\widetilde{F_0^\circ}] = [\mathcal{P}] \L^2$ and $[\widetilde{C_1}]=[\mathcal{P}](\L+1)$,
\item if $F_0$ is a minimal ruled surface, then $[\widetilde{F_0^\circ}] = [\widetilde{C_1}](\L+1),$
\item $[\widetilde{F_j^\circ}]=[\widetilde{C_1}](\L-1)$ and $[\widetilde{C_{j+1}}] = [\widetilde{C_1}]$, for every $j=1, \ldots, \ell$.
\end{enumerate}
Therefore, we can explicitely compute $Z_F(T)$.
We will illustrate this for a flower of type $3A$ and for a flower of type $4C$.

Suppose $F$ is a flower of type $3A$, then $\ell=1$ and Theorem~\ref{thm:CM-classification-flowers} and Lemma~\ref{thm:CM-flowers-numerical-relations} imply that the numerical data are given by
\begin{align*}
(N_1, \nu_1) &= (2N_0, 2\nu_0-1),\\
(N_2,\nu_2) &= (3N_0, 3\nu_0-2). 
\end{align*}
Let $\mathcal{P}$ be the $k$-variety with good $\hat{\mu}$-action, such that $[\widetilde{F_0^\circ}] = [\mathcal{P}] \L^2$ and $[\widetilde{C_1}]=[\mathcal{P}](\L+1)$. We then have $[\widetilde{F_1^\circ}]=[\mathcal{P}](\L+1)(\L-1)$ and $[\widetilde{C_2}]=[\mathcal{P}](\L+1)$.

Therefore, we compute
\begin{align*}
Z_F(T) &= \frac{[\widetilde{F_0^\circ}] \L^{-\nu_0}T^{N_0}}{1-\L^{-\nu_0}T^{N_0}} + \frac{[\widetilde{F_1^\circ}] \L^{-\nu_1}T^{N_1}}{1-\L^{-\nu_1}T^{N_1}} \\
&\qquad + \frac{[\widetilde{C_1}](\L-1) \L^{-\nu_0-\nu_1}T^{N_0+N_1}}{(1-\L^{-\nu_0}T^{N_0})(1-\L^{-\nu_1}T^{N_1})} + \frac{[\widetilde{C_2}](\L-1) \L^{-\nu_1-\nu_2}T^{N_1+N_2}}{(1-\L^{-\nu_1}T^{N_1})(1-\L^{-\nu_2}T^{N_2})}\\
&= \frac{[\mathcal{P}]\L^2 \L^{-\nu_0}T^{N_0}}{1-\L^{-\nu_0}T^{N_0}} + \frac{[\mathcal{P}](\L+1)(\L-1) \L^{-2\nu_0+1}T^{2N_0}}{1-\L^{-2\nu_0+1}T^{2N_0}} \\
&\qquad + \frac{[\mathcal{P}](\L+1)(\L-1) \L^{-3\nu_0+1}T^{3N_0}}{(1-\L^{-\nu_0}T^{N_0})(1-\L^{-2\nu_0+1}T^{2N_0})}\\
&\qquad + \frac{[\mathcal{P}](\L+1)(\L-1) \L^{-5\nu_0+3}T^{5N_0}}{(1-\L^{-2\nu_0+1}T^{2N_0})(1-\L^{-3\nu_0+2}T^{3N_0})}\\
&= \frac{[\mathcal{P}] \L^{-\nu_0+2} T^{N_0} (\L^{-2\nu_0+2}T^{2N_0} + \L^{-\nu_0+2}T^{N_0} +1)}{1-\L^{-\nu_{2}}T^{N_{2}}},
\end{align*}
where we computed the last equality with SymPy, a library of the programming language Python for symbolic mathematics.

We conclude that
\[Z_{F}(T)\in \mathcal{M}_k^{\hat{\mu}}\left[T, \frac{1}{1-\L^{-\nu_{\ell+1}}T^{N_{\ell+1}}}\right].\] 

Assume now that $F$ is a flower of type $4C$ of length $\ell$. Theorem~\ref{thm:CM-classification-flowers} and Lemma~\ref{thm:CM-flowers-numerical-relations} imply that the numerical data are given by
\begin{align*}
(N_0, \nu_0) &= (2N_1, 2\nu_1+1),\\
(N_i,\nu_i) &= (N_1, \nu_1 - i + 1) & \text{for } i =1, \ldots, \ell,\\
(N_{\ell+1}, \nu_{\ell+1}) &= (2N_1, 2\nu_1 -2\ell+1).
\end{align*}
We have $[\widetilde{F_j^\circ}]=[\widetilde{C_1}](\L-1)$, and $[\widetilde{C_{j+1}}]=[\widetilde{C_1}]$ for $j=1, \ldots, \ell$.

Therefore, we compute
\begin{align*}
Z_F(T) &= \frac{[\widetilde{F_0^\circ}] \L^{-\nu_0}T^{N_0}}{1-\L^{-\nu_0}T^{N_0}}+ \frac{[\widetilde{C_1}](\L-1) \L^{-\nu_0-\nu_1}T^{N_0+N_1}}{(1-\L^{-\nu_0}T^{N_0})(1-\L^{-\nu_1}T^{N_1})}\\
&\qquad + \sum_{j=1}^{\ell-1}\left(\frac{[\widetilde{F_j^\circ}] \L^{-\nu_j}T^{N_j}}{1-\L^{-\nu_j}T^{N_j}}
+ \frac{[\widetilde{C_{j+1}}](\L-1) \L^{-\nu_j-\nu_{j+1}}T^{N_j+N_{j+1}}}{(1-\L^{-\nu_{j}}T^{N_j})(1-\L^{-\nu_{j+1}}T^{N_{j+1}})}\right)\\
&\qquad +\frac{[\widetilde{F_\ell^\circ}] \L^{-\nu_\ell}T^{N_\ell}}{1-\L^{-\nu_\ell}T^{N_\ell}} + \frac{[\widetilde{C_{\ell+1}}](\L-1) \L^{-\nu_{\ell}-\nu_{\ell+1}}T^{N_{\ell}+N_{\ell+1}}}{(1-\L^{-\nu_\ell}T^{N_\ell})(1-\L^{-\nu_{\ell+1}}T^{N_{\ell+1}})}\\
&= \frac{[\widetilde{F_0^\circ}] \L^{-\nu_0}T^{N_0}}{1-\L^{-\nu_0}T^{N_0}}+ \frac{[\widetilde{C_1}](\L-1) \L^{-3\nu_1-1}T^{3N_1}}{(1-\L^{-2\nu_1-1}T^{2N_1})(1-\L^{-\nu_1}T^{N_1})}\\
&\qquad + \sum_{j=1}^{\ell-1}\Bigg(\frac{[\widetilde{C_1}](\L-1) \L^{-\nu_1+j-1}T^{N_1}}{1-\L^{-\nu_1+j-1}T^{N_1}}\\
&\qquad\qquad + \frac{[\widetilde{C_{1}}](\L-1) \L^{-2\nu_1+2j-1}T^{2N_1}}{(1-\L^{-\nu_{1}+j-1}T^{N_1})(1-\L^{-\nu_{1}+j}T^{N_{1}})}\Bigg)\\
&\qquad +\frac{[\widetilde{C_{1}}](\L-1) \L^{-\nu_1+\ell-1}T^{N_1}}{1-\L^{-\nu_1+\ell-1}T^{N_1}}\\
&\qquad \qquad + \frac{[\widetilde{C_{1}}](\L-1) \L^{-3\nu_{1}+3\ell-2}T^{3N_{1}}}{(1-\L^{-\nu_1+\ell-1}T^{N_1})(1-\L^{-2\nu_{1}+2\ell-1}T^{2N_{1}})}.
\end{align*}

By induction on $m$, one can compute that
\begin{align*}
\sum_{j=1}^{m}&\left(\frac{[\widetilde{C_1}](\L-1) \L^{-\nu_1+j-1}T^{N_1}}{1-\L^{-\nu_1+j-1}T^{N_1}}
 + \frac{[\widetilde{C_{1}}](\L-1) \L^{-2\nu_1-2j-1}T^{2N_1}}{(1-\L^{-\nu_{1}+j-1}T^{N_1})(1-\L^{-\nu_{1}+j}T^{N_{1}})}\right)\\
 &= \frac{[\widetilde{C_1}](\L^m-1)\L^{-\nu_1}T^{N_1}}{(1-\L^{-\nu_1}T^{N_1})(1-\L^{-\nu_1+m}T^{N_1})}
 \end{align*}
 
Therefore, direct computation (using code in SymPy) shows that
 \[Z_F(T) = \frac{[\widetilde{F_0^\circ}] \L^{-\nu_0}T^{N_0}}{1-\L^{-\nu_0}T^{N_0}}+\frac{[\widetilde{C_1}]\L^{-\nu_1}T^{N_1}\Theta(T)}{(1-\L^{-\nu_0}T^{N_0})(1-\L^{-\nu_{\ell+1}}T^{N_{\ell+1}})},\]
 where 
\begin{align*}
\Theta(T) &= (\L^{\ell}-1)+(\L^{2\ell-1}-1)\L^{-\nu_1}T^{N_1}\\
&\quad + (\L^\ell-1)\ell^{-2\nu+l-1}T^{2N_1}
 + (\L-1)\L^{-3\nu_1+2\ell-2}T^{3N_1}.
\end{align*}
We conclude that
\[Z_{F}(T)\in \mathcal{M}_k^{\hat{\mu}}\left[T, \frac{1}{1-\L^{-\nu_{0}}T^{N_{0}}}, \frac{1}{1-\L^{-\nu_{\ell+1}}T^{N_{\ell+1}}}\right].\]

All computations have been implemented in Python, as well as the computations for the other flowers.
The code can be found in Appendix~\ref{app:Sympy-code} and can be downloaded from \mbox{\url{www.github.com/AnneliesJaspers/flowers_contribution}}.
Explicit formulas for the contribution of flowers can be found in Appendix~\ref{app:formula-contribution}.
\end{proof}

\begin{corollary}\label{thm:poles-contributions}
Let $X$ be a $K3$ surface over $K$ allowing a Crauder-Morrison model $\X$ with special fiber $\X_k=\sum_{i\in I} N_i E_i$. Let $\omega$ be a volume form on $X$ and let $(N_i, \nu_i)$ be the numerical data of $E_i$, for every $i\in I$. Let $\rho_i = \nu_i/N_i +1$ be the weight of $E_i$, for every $i\in I$.

Define $I^\dagger\subset I$ to be the set of indices $i\in I$ with either
\begin{enumerate}[(i)]
\item $\rho_i$ is minimal, or
\item $E_i$ is the top of a conic-flower.
\end{enumerate}
Let $S^\dagger=\left\{(-\nu_i,N_i)\in \Z\times \Z_{>0}\mid i\in I^\dagger\right\}$. We have
\[Z_{X,\omega}(T)\in \mathcal{M}_k^{\hat{\mu}}\left[T, \frac{1}{1-\L^{a}T^b}\right]_{(a,b)\in S^\dagger}.\]
\end{corollary}

\begin{proof}
This follows immediately from Lemma~\ref{thm:poles-sum-contribution-flowers} and Theorem~\ref{thm:poles-contribution-flowers}.
\end{proof}

\section{Poles of the motivic zeta function} \label{sect:poles-computation}

\subsection{Definition of a pole of a rational function over \texorpdfstring{$\mathcal{M}_k^{\hat{\mu}}$}{Mk(mu)}} \label{sect:def-poles}

The motivic zeta function $Z_{X,\omega}(T)$ is a rational function over $\mathcal{M}_k^{\hat{\mu}}$, and we are interested in its \emph{poles}. Since $\mathcal{M}_k^{\hat{\mu}}$ is not a domain, it is, a priori, not clear what the definition of a pole of a rational function over $\mathcal{M}_k^{\hat{\mu}}$ should be. We use the definition of pole as defined by Nicaise and Xu in \cite[Remark~3.7]{NicaiseXu}.
\note{Rodrigues and Veys also gave a definition of pole, which is the following:
Write $Z(T)=\frac{F(T)}{G(T)}$, where we allow rational powers of $\L$. Let $n\geq 0$ such that $(T-\L^{-q})^n$ divides $G(T)$ but $(T-\L^{-q})^{n+1}$ does not divide $G(T)$. Then $q$ is a pole of $Z(T)$ iff $(T-\L^{-q})^n$ does not divide $F(T)$. (One needs to work in the appropriate ring.) One has to introduce rational powers of $\L$. We prefer not to do so. Moreover, it is not known whether the two definitions are equivalent. We know that a pole in the definiton of Rodriques and Veys is also a pole in the definition of Nicaise and Xu, but the converse is not clear.} 

%
%

In this thesis, we will only be interested in poles of rational functions over $\mathcal{M}_k^{\hat{\mu}}$ that are elements of $\mathcal{M}_k^{\hat{\mu}}\left[T,\frac{1}{1-\mathbb{L}^{a}T^b}\right]_{(a,b)\in \Z\times\Z_{>0}}$. This ring is the localization of $\mathcal{M}_k^{\hat{\mu}}[T]$ with respect to the multiplicatively closed set $\left\{1-\mathbb{L}^a T^b\mid (a,b)\in \Z\times\Z_{>0}\right\}$. Note that this set does not contain zero-divisors and hence 
\[\frac{F_1}{G_1} = \frac{F_2}{G_2} \in \mathcal{M}_k^{\hat{\mu}}\left[T,\frac{1}{1-\mathbb{L}^{a}T^b}\right]_{(a,b)\in \Z\times\Z_{>0}} \Leftrightarrow F_1\cdot G_2 = F_2 \cdot G_1 \in \mathcal{M}_k^{\hat{\mu}}[T].\]

\begin{definition}\label{def:pole}
Let $Z(T)\in \mathcal{M}_k^{\hat{\mu}}\left[T,\frac{1}{1-\mathbb{L}^{a}T^b}\right]_{(a,b)\in \Z\times\Z_{>0}}$ be a rational function over $\mathcal{M}_k^{\hat{\mu}}$ and let $q\in \mathbb{Q}$ be a rational number. We say that $Z(T)$ has a \emph{pole of order at most $m$ at $q$}, if we find a set $\mathcal{A}$ consisting of multisets in $\Z\times \Z_{>0}$ such that
\begin{enumerate}[(i)]
\item each multiset $A \in \mathcal{A}$ contains at most $m$ elements $(a,b)$ with $\frac{a}{b}=q$, and
\item $Z(T)$ belongs to the sub-$\mathcal{M}_k^{\hat{\mu}}[T]$-module of $\mathcal{M}_k^{\hat{\mu}}\llbracket T\rrbracket$ generated by
\[ \left\{ {\left.\frac{1}{\prod_{(a,b)\in A} (1-\L^aT^b)}\, \right|\,  A\in \mathcal{A}}\right\}.\]
\end{enumerate}
We say that $Z(T)$ has \emph{a pole of order $m$ at $q$}, if it has a pole of order at most $m$, but not of order at most $m-1$ (the latter condition is void for $m = 0$).

We say that $Z(T)$ has \emph{a pole at $q$}, if there exists an integer $m\geq 1$  such that $Z(T)$ has a pole or order $m$ at $q$.
\end{definition}

To explain the intuition behind this definition, we notice that one often formally substitutes $T=\L^{-s}$ and considers $Z(T)$ as a function in the variable $s$. If $\frac{a}{b}$ is a pole of $Z(\L^{-s})$, then we can think of the factor $1-\L^aT^b=1-\L^a\L^{-bs}$ in the denominator as having a zero in $s=\frac{a}{b}$.

\begin{lemma} \label{thm:poles-equivalence-def}
Let $Z(T)\in \mathcal{M}_k^{\hat{\mu}}\left[T,\frac{1}{1-\mathbb{L}^{a}T^b}\right]_{(a,b)\in \Z\times\Z_{>0}}$ be a rational function over $\mathcal{M}_k^{\hat{\mu}}$, and let $q\in \mathbb{Q}$ be a rational number. The following are equivalent:
\begin{enumerate}[(i)]
\item $q$ is a pole of $Z(T)$,
\item for every subset $S \subset \Z\times \Z_{>0}$ with $Z(T)\in \mathcal{M}_k^{\hat{\mu}}\left[T,\frac{1}{1-\mathbb{L}^{a}T^b}\right]_{(a,b)\in S}$, there exists an element~$(a,b)\in S$ with $q=\frac{a}{b}$. 
\end{enumerate}
\end{lemma}

\begin{proof}
Suppose first that there exists a subset~${S\subset \Z\times \Z_{>0}}$ with $Z(T)\in \mathcal{M}_k^{\hat{\mu}}\left[T,\frac{1}{1-\mathbb{L}^{a}T^b}\right]_{(a,b)\in S}$ with the property that $q\neq \frac{a}{b}$ for every $(a,b)\in S$. This means that we can write
\[Z(T) = \frac{F(T)}{\prod_{(a,b)\in S} (1-\L^aT^b)^{m_{(a,b)}}},\]
for some $F(T)\in \mathcal{M}_k^{\hat{\mu}}[T]$ and $m_{(a,b)}\in \Z_{\geq 0}$. Let $A$ be the multiset that contains every $(a,b)\in S$ with multiplicity $m_{(a,b)}$ and define $\mathcal{A}=\{A\}$. It is clear that $Z(T)$ belongs to the sub-$\mathcal{M}_k^{\hat{\mu}}[T]$-module of $\mathcal{M}_k^{\hat{\mu}}\llbracket T\rrbracket$ generated by
\[ \left\{ {\left.\frac{1}{\prod_{(a,b)\in A} (1-\L^aT^b)}\, \right|\,  A\in \mathcal{A}}\right\}.\]
Moreover, $A$ does not contain an element $(a,b)$ with $q=\frac{a}{b}$. Hence $q$ is a pole of order at most 0, which means that it is \emph{not} a pole of $Z(T)$.

Conversely, suppose that $q$ is not a pole of $Z(T)$. So there exists a set $\mathcal{A}$ consisting of multisets in $\Z\times \Z_{>0}$ such that $Z(T)$ belongs to the sub-$\mathcal{M}_k^{\hat{\mu}}[T]$-module of $\mathcal{M}_k^{\hat{\mu}}\llbracket T \rrbracket$ generated by 
\[ \left\{ {\left.\frac{1}{\prod_{(a,b)\in A} (1-\L^aT^b)}\, \right|\,  A\in \mathcal{A}}\right\},\]
and such that no multiset $A\in \mathcal{A}$ contains a couple $(a,b)$ with $\frac{a}{b}=q$.
Define $S = \bigcup_{A\in \mathcal{A}} A$. It is clear that $Z(T)\in \mathcal{M}_k^{\hat{\mu}}\left[T,\frac{1}{1-\mathbb{L}^{a}T^b}\right]_{(a,b)\in S}$, and that there is no couple $(a,b)\in S$ with $q=\frac{a}{b}$.
\end{proof}


Let $P\colon\mathcal{M}_k^{\hat{\mu}}\to \Z[v,v^{-1}]$ be the morphism from Example~\ref{ex:specialization-Grothendieck-ring} (iii) that maps the class $[X]$ of a smooth, proper $k$-variety $X$ to its Poincar\'e polynomial $\sum_{i\geq 0} (-1)^i b_i(X) v^i$. This morphism can be extended to
\[P\colon \mathcal{M}_k^{\hat{\mu}}[T] \to \Z[v,v^{-1}][T].\]
We will use the morphism $P$ to compute poles in practice. More specifically, the following lemma and corollary will turn out to be useful.

\begin{lemma}\label{thm:poles-how-compute}
Let $Z(T)\in \mathcal{M}_k^{\hat{\mu}}[T,\frac{1}{1-\mathbb{L}^{a}T^b}]_{(a,b)\in \Z\times\Z_{>0}}$ be a rational function over $\mathcal{M}_k^{\hat{\mu}}$ and let $q\in \mathbb{Q}$ be a rational number. Write $Z(T)$ as
\[Z(T) = \frac{F(T)}{\prod_{(a,b)\in S} (1-\L^aT^b)^{m_{(a,b)}}},\]
for some finite subset $S\subset \Z\times \Z_{>0}$, and $m_{(a,b)}> 0$ for every $(a,b)\in S$, and where $F(T)\in \mathcal{M}_k^{\hat{\mu}}[T]$. Suppose $q= \frac{a_0}{b_0}$ for some $(a_0,b_0)\in S$. Suppose moreover that $P(F(T))$ does not have a zero in $T=v^{-2a_0/b_0}$, when we consider $P(F(T))$ as an element of $\Z[v^{1/b_0}, v^{-1/b_0}][T]$ by the inclusion $\Z[v,v^{-1}][T] \hookrightarrow \Z[v^{1/b_0}, v^{-1/b_0}][T]$. Then $Z(T)$ has a pole at $q$.
\end{lemma}

\begin{proof}
We will use Lemma~\ref{thm:poles-equivalence-def}. So suppose $Z(T)\in \mathcal{M}_k^{\hat{\mu}}\left[T,\frac{1}{1-\mathbb{L}^{a'}T^{b'}}\right]_{(a',b')\in S'}$ with $S'\subset \Z\times \Z_{>0}$. We have to show that there is a couple $(a'_0,b'_0)\in S'$ with $q=\frac{a'_0}{b'_0}$. Write $Z(T)$ as 
\[Z(T) = \frac{F'(T)}{\prod_{(a',b')\in {S'}} (1-\L^{a'}T^{b'})^{m'_{(a',b')}}}.\]
This means that
\[\frac{F'(T)}{\prod_{(a',b')\in {S'}} (1-\L^{a'}T^{b'})^{m'_{(a',b')}}} = \frac{F(T)}{\prod_{(a,b)\in S} (1-\L^aT^b)^{m_{(a,b)}}}\] 
in $\mathcal{M}_k^{\hat{\mu}}[T,\frac{1}{1-\mathbb{L}^{a}T^b}]_{(a,b)\in \Z\times\Z_{>0}}$. Therefore
\[F'(T)\prod_{(a,b)\in S} (1-\L^aT^b)^{m_{(a,b)}} = F(T)\prod_{(a',b')\in S'} (1-\L^{a'}T^{b'})^{m'_{(a',b')}}.\]
If we apply the morphism $P\colon \mathcal{M}_k^{\hat{\mu}}[T] \to \Z[v,v^{-1}][T]$ to both sides, we get
\[P(F'(T)) \prod_{(a,b)\in S} (1-v^{2a}T^b)^{m_{(a,b)}} = P(F(T)) \prod_{(a',b')\in S'} (1-v^{2a'}T^{b'})^{m'_{(a',b')}}.\]
Let us now work in $\Z[v^{1/b_0}, v^{-1/b_0}][T]$ by the inclusion $\Z[v,v^{-1}][T] \hookrightarrow \Z[v^{1/b_0}, v^{-1/b_0}][T]$. Since $(a_0,b_0)\in S$, the right-hand side of this equation has a zero in $T=v^{-2a_0/b_0}$ and therefore
\[P(F(v^{-2a_0/b_0}))\prod_{(a',b')\in S'} (1-v^{2a'-\frac{2a_0b'}{b_0}})^{m'_{(a',b')}}=0.\]
Since $P(F(v^{-2a_0/b_0}))\neq 0$ by assumption, and $\Z[v^{1/b_0}, v^{-1/b_0}]$ is a domain, there must be a couple $(a_0',b_0')\in S'$ with $2a_0'-\frac{2a_0b_0'}{b_0}=0$. This is equivalent to $\frac{a_0'}{b_0'}=\frac{a_0}{b_0}=q$.
\end{proof}

\begin{remark} \label{rmk:poles-compute-chi}
By taking $v=1$ in Lemma~\ref{thm:poles-how-compute}, we find that, if $\chi(F(T))$ does not have a zero in $T=\exp(-4\pi i a_0/b_0)$, then $Z(T)$ has a pole at $q$. Here, $\chi$ stands for the topological Euler characteristic.
\end{remark}

\begin{corollary} \label{thm:poles-how-compute-sum}
Let $Z(T)\in \mathcal{M}_k^{\hat{\mu}}\left[T,\frac{1}{1-\mathbb{L}^{a}T^b}\right]_{(a,b)\in \Z\times\Z_{>0}}$ be a rational function over $\mathcal{M}_k^{\hat{\mu}}$, and let $q\in \mathbb{Q}$ be a rational number. Suppose we can write
\[Z(T) = Z_1(T)+Z_2(T),\]
with $Z_1(T), Z_2(T)\in \mathcal{M}_k^{\hat{\mu}}\left[T,\frac{1}{1-\mathbb{L}^{a}T^b}\right]_{(a,b)\in \Z\times\Z_{>0}}$, and such that $Z_1(T)$ satisfies the conditions of Lemma~\ref{thm:poles-how-compute} and $Z_2(T)$ does not have a pole at $q$. Then $Z(T)$ has a pole at $q$.
\end{corollary}

\subsection{Computation of the poles of \texorpdfstring{$Z_{X,\omega}(T)$}{Z(T)} }

In this section, we will compute the poles of the motivic zeta function. 

The following result is proven by Halle and Nicaise in \cite[Theorem~3.2.3]{HalleNicaiseKulikov} for general Calabi-Yau varieties.
\begin{theorem}
Let $X$ be a Calabi-Yau variety over $K$.  Let $\X$ be an $snc$-model of $X$ with special fiber $\X_k=\sum_{i\in I} N_i E_i$. Let $\omega$ be a volume form on $X$ and let $(N_i,\nu_i)$ be the numerical data of $E_i$. 

Define
\begin{align*}
lct(X)&=\min\left\{ \nu_i/N_i\mid i\in I\right\},\\
\delta(X)&=\max \left\{ |J| \mid \emptyset \neq J \subset I, E_J \neq \emptyset, \nu_j/N_j = lct(X) \text{ for all } j \in J\right\} 
- 1.
\end{align*}
The motivic zeta function $Z_{X,\omega}(T)$ has a pole at $-lct(X)$ of order $\delta(X)+1$, and this is the largest pole of $Z_{X,\omega}(T)$.
\end{theorem}

As a consequence, we get
\begin{corollary} \label{thm:poles-lct}
Let $X$ be a $K3$ surface over $K$.  Let $\X$ be a Crauder-Morrison model of $X$ with special fiber $\X_k=\sum_{i\in I} N_i E_i$. Let $\omega$ be a volume form on $X$ and let $(N_i, \nu_i)$ be the numerical data of $E_i$ for every $i\in I$.
If $\rho_i$ is minimal, then $-\nu_i/N_i$ is a pole of $Z_{X,\omega}(T)$, and it is the largest pole of $Z_{X,\omega}(T)$. Moreover, it is a pole of order 1 if $\X$ is a flowerpot degeneration, and of order 2 if $\X$ is a chain degeneration.
\end{corollary}

We will now compute the other poles.
The following result has been announced in \cite{Jaspers}.

\begin{theorem} \label{thm:poles}
Let $X$ be a $K3$ surface over $K$.  Let $\X$ be a Crauder-Morrison model of $X$ with special fiber $\X_k=\sum_{i\in I} N_i E_i$. Let $\omega$ be a volume form on $X$ and let $(N_i, \nu_i)$ be the numerical data of $E_i$ for every $i\in I$.

The rational number $q\in \Q$ is a pole of $Z_{X,\omega}(T)$ if and only if there exists an element $i\in I$ with $q=-\nu_i/N_i$ and such that 
\begin{enumerate}[(i)]
\item either $\rho_i$ is minimal,
\item or $E_i$ is the top of a conic-flower.
\end{enumerate}
Moreover, in case (i), $q$ is a pole of order 1 if $\X$ is a flowerpot degeneration, and of order 2 if $\X$ is a chain degeneration. In case (ii), $q$ is a pole of order 1.
\end{theorem}

\begin{proof}

For every $q\in \Q$, we define $I_q=\left\{ i\in I \mid -\nu_i/N_i=q\right\}$. Remark that, because of equation \eqref{def:CM-weight}, we have $\rho_i=\rho_j$ for every $i,j\in I_q$. 

Suppose first that there does \emph{not} exist an $i\in I$ with $q=-\nu_i/N_i$ and such that 
\begin{enumerate}[(i)]
\item either $\rho_i$ is minimal,
\item  or $E_i$ is the top of a conic-flower.
\end{enumerate} 

If $I_q=\emptyset$, then $Z_{X,\omega}(T)$ does not have a pole at $q$ by the Denef-Loeser formula~\ref{thm:Denef-Loeser-CY}.

So from now on we may assume $I_q\neq \emptyset$. 
Let $j\in I_q$. Since $\rho_j$ is not minimal, $E_j$ is a component of a flower. Moreover, if $k\in I_q$ and $k\neq j$, then we know that $E_k$ is component of a different flower than the flower in which $E_j$ is contained. Indeed, the weights $\rho_i=\nu_i/N_i+1$ strictly  decrease within a flower, so no two components in one flower can realize the same $-\nu_i/N_i$. Denote by $F^{(i)}$ the flower containing $E_i$ for every $i\in I_q$. By Lemma~\ref{thm:poles-sum-contribution-flowers}, we can write
\[Z_{X,\omega}(T) = \sum_{i\in I_q}Z_{F^{(i)}}(T) + G(T),\]
where $Z_{F^{(i)}}(T)$ is the contribution of $F^{(i)}$ to the motivic zeta function and  $G(T) \in \mathcal{M}_k^{\hat{\mu}}\left[T,\frac{1}{1-\mathbb{L}^{a}T^b}\right]_{(a,b)\in \Z\times\Z_{>0}}$. By Lemma~\ref{thm:poles-equivalence-def}, we know that $G(T)$ does not have a pole at $q$.
When $i\in I_q$, Theorem~\ref{thm:poles-contribution-flowers} implies that $Z_{F^{(i)}}(T)$ doesn't have a pole at $q$ either.
We conclude that $q$ is not a pole of $Z_{X,\omega}(T)$.

Conversely, suppose that $\rho_i$ is minimal for some $i\in I_q$, and hence for all $i\in I_q$. Corollary~\ref{thm:poles-lct} gives that $-\nu_i/N_i$ is a pole of $Z_{X,\omega}(T)$.

Finally, suppose there is a $j\in I_q$ such that $E_j$ is the top of a conic-flower. For every $i\in I_q$, we have that $E_i$ is a component of a flower $F^{(i)}$, since $\rho_i$ is not minimal. Furthermore, for every $i,k\in I_q$ with $i\neq k$, the flower $F^{(i)}$ and $F^{(k)}$ are distinct. Since the weight strictly decreases within a flower, $q$ is a pole of order at most 1 by the Denef-Loeser formula~\ref{thm:Denef-Loeser-CY}.
 
Define
\[I_q' = \{i\in I_q \mid E_i \text{ is the top of a conic-flower} \}.\]
By assumption $I_q'\neq \emptyset$.
 By Lemma~\ref{thm:poles-sum-contribution-flowers}, we can write
\[Z_{X,\omega}(T) = \sum_{i\in I'_q}Z_{F^{(i)}}(T) + G(T).\]
where $Z_{F^{(i)}}(T)$ is the contribution of $F^{(i)}$ to the motivic zeta function and $G(T) \in \mathcal{M}_k^{\hat{\mu}}\left[T,\frac{1}{1-\mathbb{L}^{a}T^b}\right]_{(a,b)\in \Z\times\Z_{>0}}$. By Lemma~\ref{thm:poles-equivalence-def} and Theorem~\ref{thm:poles-contribution-flowers}, $G(T)$ does not have a pole at $q$ .

Choose $i\in I_q'$.
Theorem~\ref{thm:CM-classification-flowers} and Lemma~\ref{thm:CM-flowers-numerical-relations} explain why $((\nu_{i}-1)/2,N_{i}/2)$ is the numerical data of the component meeting $E_i$ in a conic. 
Define
\[Z_{E_i}(T) = [\widetilde{E_i^\circ}] \frac{\L^{-\nu_i}T^{N_i}}{1-\L^{-\nu_i}T^{N_i}}+(\L-1)[\widetilde{C_i}]\frac{\L^{-\nu_i}T^{N_i}}{1-\L^{-\nu_i}T^{N_i}}\frac{\L^{-(\nu_i-1)/2}T^{N_i/2}}{1-\L^{-(\nu_i-1)/2}T^{N_i/2}},\]
where $C_i$ is the double curve on $E_i$.
 We see
\[Z_{F^{(i)}}(T) = Z_{E_i}(T) + G_i(T),\]
with $G_i(T) \in \mathcal{M}_k^{\hat{\mu}}\left[T,\frac{1}{1-\mathbb{L}^{a}T^b}\right]_{(a,b)\in \Z\times\Z_{>0}}$ and $G_i(T)$ does not have a pole at $q$.
Therefore we can write 
\[Z_{X,\omega}(T) = \sum_{i\in I_q'}Z_{E_{i}}(T) + G'(T),\]
where $G'(T) \in \mathcal{M}_k^{\hat{\mu}}\left[T,\frac{1}{1-\mathbb{L}^{a}T^b}\right]_{(a,b)\in \Z\times\Z_{>0}}$ does not have a pole at $q$.

Define 
\[\Theta_i(T) = [\widetilde{E_i^\circ}] \L^{-\nu_i}T^{N_i}\left(1-\L^{-(\nu_i-1)/2}T^{N_i/2}\right) + (\L-1)[\widetilde{C_i}]\L^{(-3\nu_i+1)/2}T^{3N_i/2},\]
in $\mathcal{M}_k^{\hat{\mu}}[T]$. Then
\[Z_{E_i}(T) = \frac{\Theta_i(T)}{(1-\L^{-\nu_i}T^{N_i})(1-\L^{-(\nu_i-1)/2}T^{N_i/2})}.\] 

Let $N=\lcm\left\{N_{i}\mid i\in I'_q \right\}$ and $\nu=\lcm \left\{\nu_{i}\mid i\in I'_q\right\}$.  Since $-\nu_{i}/N_{i}$ is constant, we have $-\nu/N = -\nu_{i}/N_{i}$ for every $i\in I'_q$.
Hence
\[(1-\L^{-\nu}T^{N})=(1-\L^{-\nu_{i}}T^{N_{i}})\cdot \sum_{k=0}^{N/N_{i}-1}\L^{-k\nu_{i}}T^{kN_{i}}.\]

By bringing all terms in $\sum_{i\in I_q'}Z_{E_i}(T)$ to a common denominator, we get
\[ \sum_{i\in I'_q}Z_{E_i}(T) = \frac{\Theta(T)}{(1-\L^{-\nu}T^{N})\prod_{i\in I_q'}(1-\L^{-(\nu_{i}-1)/2}T^{N_{i}/2})},\]
where
\[\Theta(T) = \sum_{i\in I_q'} \left(\Theta_{i}(T) \cdot  \sum_{k=0}^{N/N_{i}-1} (\L^{-k\nu_{i}}T^{kN_{i}})\cdot \prod_{j\in I_q'\setminus\{i\}} (1-\L^{-(\nu_j-1)/2}T^{N_{j}/2})\right).\]

Let $P\colon\mathcal{M}_k^{\hat{\mu}}\to \Z[v,v^{-1}]$ be the morphism from Example~\ref{ex:specialization-Grothendieck-ring} (iii) that maps the class $[X]$ of a smooth, proper $k$-variety $X$ to its Poincar\'e polynomial $\sum_{i\geq 0} (-1)^i b_i(X) v^i$. This morphism can be extended to
\[P\colon \mathcal{M}_k^{\hat{\mu}}[T] \to \Z[v,v^{-1}][T].\]

Lemma~\ref{thm:poles-how-compute-sum} says that $q$ is a pole of $Z_{X,\omega}(T)$, if $P(\Theta_q(T))$ does not have a zero in $T=v^{-2\nu/N}$, when we consider $P(\Theta_q(T))$ as an element of $\Z[v^{1/N}, v^{-1/N}][T]$ by the inclusion $\Z[v,v^{-1}][T] \hookrightarrow \Z[v^{1/N}, v^{-1/N}][T]$.

Because $P([\widetilde{E_i^\circ}])=\frac{N_i}{2}(v^2+1)v^2$ and $P([\widetilde{C_i}])=\frac{N_0}{2}(v^2+1)$ by Proposition~\ref{thm:poles-conic-flower}, direct computation shows that
\[P(\Theta(v^{-2\nu/N})) = -\frac{\sum_{i\in I_q'} (N_i+N)}{2}v(v^2+1)(1-v)^{|I_q'|},\]
which is clearly non-zero. This concludes the proof.
\note{
Specializing to the Poincar\'e polynomial gives
\[\sum_{i\in I_q'}\left( \frac{N_{i}}{2}(x^2+1)( x^2 - x^{-2\nu_{1i}}T^{N_{1i}})x^{-2\nu_{i}}T^{N_{i}}  \cdot \sum_{k=0}^{N/N_{i}} (x^{-2\nu_{i}}T^{N_{i}})^k \cdot \prod_{j\in I_q'\setminus\{i\}} (1-x^{-\nu_j+1}T^{N_{j}/2})\right).\]
If we substitute $T\to x^{2\nu/N}$, we get an expression in $\C[x, x^{1/N}]$:
\begin{align*}
\sum_{i\in I_q'} & \left( \frac{N_{i}}{2}(x^2+1)( x^2 - x^{-2\nu_{1i}}x^{2\nu N_{1i}/N})x^{-2\nu_{i}}x^{2\nu_{i}}  \cdot \sum_{k=0}^{N/N_{i}} (x^{-2\nu_{i}}x^{2\nu_{i}})^k \cdot \prod_{j\in I_q'\setminus\{i\}} (1-x^{-\nu_j+1}x^{2\nu N_{1j}/N})\right)\\
\sum_{i\in I_q'} & \left( \frac{N_{i}}{2}(x^2+1)( x^2 - x^{-2\nu_{1i}+2\nu N_{1i}/N}) \cdot (N/N_{i}+1) \cdot \prod_{j\in I_q'\setminus\{i\}} (1-x^{-\nu_{j}+1+2\nu N_{1j}/N})\right)\\
\sum_{i\in I_q'} & \left( \frac{N_{i}+N}{2}(x^2+1)( x^2 - x) \cdot \prod_{j\in I_q'\setminus\{i\}} (1-x)\right)\\
\sum_{i\in I_q'} & \left( \frac{N_{i}+N}{2}\right)(x^2+1)( x^2 - x)(1-x)^{|I_q'|-1}.
\end{align*}
This is clearly a non-zero expression.}
\end{proof}

\section{\texorpdfstring{Example of a $K3$ surface where $Z_{X,\omega}(T)$ has two poles}{Example of a K3 surface with a motivic zeta function with two poles}}

Theorem \ref{thm:poles} suggests that the motivic zeta function of a $K3$ surface can have more than one pole, in contrast with abelian varieties \cite[Theorem~8.5]{HalleNicaiseAbelian} and Calabi-Yau varieties with an equivariant Kulikov model \cite[Corollary~5.3.3]{HalleNicaiseKulikov}. We will illustrate this fact with an example. This example first appeared in \cite[Example~5.3.4]{HalleNicaiseKulikov}, and it has also been published in \cite[Example~2]{Jaspers}.

\begin{example} \label{ex:poles-more-than-one}
Let $k=\C$, $R=\C\llbracket t\rrbracket$ and $K=\C(\!(t)\!)$.
Let $X$ be the $K3$ surface defined by the homogeneous equation
\begin{align} \label{eq:example}
 x^2w^2+y^2w^2+z^2w^2+x^4+y^4+z^4+tw^4 = 0
 \end{align}
in $\mathrm{Proj}\, K[x,y,z,w]$.
Let $\mathcal{Y}$ be the closed subscheme of $\mathrm{Proj}\, R[x,y,z,w]$ defined by the homogeneous equation \eqref{eq:example}. It is straightforward to check that $\Y$ is regular. Moreover, $\Y_k$ is a singular irreducible surface with a unique singularity at $P=(0:0:0:1)$, and the singularity is of type $A_1$.
We construct an $snc$-model $\X$ of $X$ by blowing up $\mathcal{Y}$ at $P$. The special fiber of $\X$ is of the form $\X_k = D+2E$, where $D$ is the strict transform of $\mathcal{Y}_k$, and $E\simeq \P^2_k$ is the exceptional divisor. The strict transform $D$ is a smooth $K3$ surface and intersects $E$ transversally along a smooth conic~$C$. \note{
In affine chart, the equation is
\[x^2+y^2+z^2+x^4+y^4=z^4+t=0.\]
When we blow up, in one chart, this is the map $(x,y,z,t)\mapsto(x,xy,xz,xt)$. And hence the equation of $\X$ in that chart is:
\[x^2+x^2y^2+x^2z^2+x^4+x^4y^4+x^4z^4+xt=0.\]
The special fiber can be computed by putting $xt=0$ and hence the equation of $\X_k$ is
\[x^2(1+y^2+z^2+x^2+x^2y^4+x^2z^4)=0.\]
So $E$ is described by $x=0$ and has multiplicity 2 and $D$ is described by $1+y^2+z^2+x^2+x^2y^4+x^2z^4=0$. The intersection $C$ has equation
\[1+y^2+z^2=0,\]
which is a smooth conic.
}
\note{
$D$ is a $K3$ surface. We know that $\Y_k$ is a singular $K3$ surface, so $D$ is for sure birational to a $K3$ surface. But by A'Campo's formula, we find that $\chi(D)=24$ and hence $D$ is a $K3$ surface.
}
For a suitable choice of volume form $\omega$, one has $\nu_D=0$ and $\nu_E=1$. 
\note{
Let $\omega=dx\wedge dy \wedge dz$ be a volume form on $\Y$. It is a generator for the log-canonical bundle away from the singular point $P$. The pull-back of $\omega$ to $\X$ is $x^2 dx\wedge dy \wedge dz$. It has log-discrepancy $2-1=1$ at $x=0$ and $0$ at $D$.
}
 Therefore, $D$ has weight $\rho_D=1$ and $E$ has weight $\rho_E = 3/2$. This means that $\X$ is a flowerpot degeneration, where $D$ is the flowerpot, and $E$ is a conic-flower of type 2B. 
 
 Using the Denef-Loeser formula (Theorem~\ref{thm:Denef-Loeser-CY}), the motivic zeta function can be computed as
\begin{align*}
Z_{X,\omega}(T) &= [\widetilde{D^\circ}]\frac{T}{1-T}+ [\widetilde{E^\circ}]\frac{\L^{-1}T^2}{1-\L^{-1}T^2}+[C] \frac{\L^{-1}T^3}{(1-T)(1-\L^{-1}T^2)}\\
&= \frac{[\widetilde{D^\circ}]T(1-\L^{-1}T^2)+ [\widetilde{E^\circ}]\L^{-1}T^2(1-T)+[C]\L^{-1}T^3}{(1-T)(1-\L^{-1}T^2)}.
\end{align*}

Define $F(T) = [\widetilde{D^\circ}]T(1-\L^{-1}T^2)+ [\widetilde{E^\circ}]\L^{-1}T^2(1-T)+[C]\L^{-1}T^3$, then 
\[\chi(F(T)) = \chi(\widetilde{D^\circ})T(1-T^2)+ \chi(\widetilde{E^\circ})T^2(1-T)+\chi(C)T^3.\]
Since $\chi(F(T))$ does not have a zero in $T=1$, Remark~\ref{rmk:poles-compute-chi} implies that $0$ and $-1/2$ are poles of $Z_{X,\omega}(T)$, which is in agreement with Theorem~\ref{thm:poles}. 
\end{example}

\cleardoublepage


\chapter{\texorpdfstring{$K3$}{K3} surfaces satisfying the monodromy property}\label{ch:GMP-holds}

In this chapter, we will prove that certain types of $K3$ surfaces satisfy the monodromy property. To be precise, we will prove the following theorem, which is slightly more general than the result announced in \cite{Jaspers}.

\begin{theorem*}
Let $X$ be a $K3$ surface over $K$ with Crauder-Morrison model~$\X$ and with special fiber $\X_k=\sum_{i\in I} N_iE_i$. Suppose $\X$ satisfies one of the following:
\begin{enumerate}[(i)]
\item $\X$ is a flowerpot degeneration, or
\item $\X$ is a chain degeneration with chain $V_0 \text{---} V_1 \text{---} \cdots \text{---} V_k \text{---} V_{k+1}$. Set $N=\min\{N(V_i)\mid i=0, \ldots, k+1\}$. Assume moreover that one of the following conditions hold:
\begin{enumerate}[(a)]
\item the components $V_0, V_1, \ldots, V_{k+1}$ all have  multiplicity $N$, or
\item neither $V_0$ nor $V_{k+1}$ is a rational, non-minimal ruled surface, or
\item exactly one end component of the chain is a rational, non-minimal ruled surface and it has multiplicity $N$, or
\item if $V_j$ meets a conic-flower, then $N(V_j)>N$.
\end{enumerate}
\end{enumerate}
Then $X$ satisfies the monodromy property.
\end{theorem*}

In the first section, we will explain the strategy to prove this theorem. We will also give some first results about monodromy eigenvalues. In Section~\ref{sect:holds-flowerpot}, we will prove part (i) of the theorem, i.e., that $K3$ surfaces allowing a flowerpot degeneration satisfy the monodromy property. Finally, in Section~\ref{sect:holds-chain}, we show that the monodromy property holds for $K3$ surfaces with a chain degeneration satisfying one of the extra conditions.

\subsubsection{Notation}
In this chapter, we fix an algebraically closed field $k$ of characteristic zero. Put $R=k\llbracket t\rrbracket$ and $K=k(\!( t )\!)$. Fix an algebraic closure $K^{alg}$ of $K$. Let $\sigma$ be a topological generator of the monodromy group $\Gal(K^{alg}/K)$.

Let $X$ be a $K3$ surface over $K$ and let $\X$ be a Crauder-Morrison model of $X$ with special fiber $\X_k=\sum_{i\in I} N_i E_i$.
Let $\omega$ be a volume form on $X$ and let $(N_i, \nu_i)$ be the numerical data of $E_i$, for every $i\in I$. To  avoid confusion, we will sometimes use the notation $(N(E_i), \nu(E_i))$ instead of $(N_i,\nu_i)$ for the numerical data of a component $E_i$. Let $\rho_i = \nu_i/N_i +1$ be the weight of $E_i$, for every $i\in I$.

 We define
\[\xi(E_i) = \exp\left(-2\pi i \, \frac{\nu(E_i)}{N(E_i)}\right),\] 
the `candidate' monodromy eigenvalue associated with $E_i$.

For any flower $F$ in $\X_k$, we denote by $F_0$ the top component of the flower, i.e., the component in the flower with maximal weight $\rho(F_0)$. By $F_{\ell+1}$, we mean the component in $\Gamma_{min}$ that intersects with the flower $F$.

\section{Strategy and first results} \label{sect:strategy}
\subsection{Strategy}

In Corollary~\ref{thm:poles-contributions}, we found that
\[Z_{X,\omega}(T)\in \mathcal{M}_k^{\hat{\mu}}\left[T, \frac{1}{1-\L^{a}T^b}\right]_{(a,b)\in S^\dagger},\]
with $S^\dagger=\left\{(-\nu_i,N_i)\in \Z\times \Z_{>0}\mid i\in I^\dagger\right\}$ and $I^\dagger\subset I$ the set of indices $i\in I$ with either
\begin{enumerate}[(i)]
\item $\rho_i$ is minimal, or
\item $E_i$ is the top of a conic-flower.
\end{enumerate}

In order to prove that $X$ satisfies the monodromy property, we need to verify that, for every $i\in I^\dagger$, the complex number~$\exp(-2\pi i\, \nu_i/N_i)$ is an eigenvalue of $\sigma$ on $H^m(X\times_K K^{alg}, \Q_\ell)$, for some $m\geq 0$ and for every embedding of $\Q_\ell$ into~$\C$.
\note{This embedding is needed, because $\exp(-2\pi i\, \nu_i/N_i)\in \C$ and $H^m(X\times_K K^{alg}, \Q_\ell)\subset \Q_\ell$.}

Remember that we defined the monodromy zeta function in Definition~\ref{def:GMP-monodromy-zeta} as
\[\zeta_X(T) = \prod_{m\geq 0} \Bigl(\det \left( T\cdot Id - \sigma\mid H^m(X\times_K K^{alg}, \Q_\ell)\right)\Bigr)^{(-1)^{m+1}} \in \Q_\ell(T).\]
This rational function has the property that all zeroes and poles are monodromy eigenvalues, but in general there may be more monodromy eigenvalues than can be seen from the monodromy zeta function. The reason is that cancellations may occur in the alternating product. 
However, for $K3$ surfaces we have the following result.

\begin{proposition} \label{thm:holds-form-zeta-function-K3}
Let $X$ be a $K3$ surface. The monodromy zeta function can be written as
\[
\zeta_X(T) = \frac{1}{Q(T)},\]
for some $Q\in \Q_\ell[T]$ of degree 24.
Moreover, the poles of $\zeta_X(T)$ are \emph{exactly} the monodromy eigenvalues of $X$.
\end{proposition}
\begin{proof}
Because $X$ is a $K3$ surface, the cohomology spaces of $X$ are trivial in odd degree and $\chi(X)=24$, by Proposition~\ref{thm:pre-K3-cohomology}. Therefore, no cancellations can occur, and the poles of the monodromy zeta function are {exactly} the monodromy eigenvalues. 
\end{proof}

In this chapter, we will compute the poles of $\zeta_X(T)$. A helpful tool will be the A'Campo formula, discussed in Proposition~\ref{thm:GMP-ACampo-dvr}.
\[\zeta_X(T) = \prod_{i\in I} \left(T^{N_i}-1\right)^{-\chi(E_i^\circ)},\]
where $\chi(E_i^\circ)$ is the topological Euler characteristic of $E_i^\circ= E_i\setminus \left( \cup_{j\in I\setminus \{ i\}} E_j\right)$.
In Lemma~\ref{thm:CM-euler-characteristics}, we computed Euler characteristics $\chi(E_i^\circ)$ for some relevant components $E_i$ of the special fiber $\X_k=\sum_{i\in I} N_i E_i$.

\subsection{Minimal weight}
Let $E_i$ be a component of $\X_k$ with numerical data $(N_i, \nu_i)$ and such that the weight $\rho_i$ is minimal. Halle and Nicaise proved that $\xi(E_i)=\exp(-2\pi i \nu_i/N_i)$ is a monodromy eigenvalue of $X$.

\begin{theorem} \label{thm:holds-lct}
Let $X$ be a $K3$ surface over $K$ with Crauder-Morrison model~$\X$. Let $\omega$ be a volume form on $X$.
Let $E_i$ be a component of $\X_k$ with numerical data $(N_i, \nu_i)$ and such that the weight $\rho_i$ is minimal. Then $\xi(E_i)=\exp(-2\pi i \nu_i/N_i)$ is a monodromy eigenvalue of $X$.
\end{theorem}
\begin{proof}
This is a special case of {\cite[Theorem~3.3.3]{HalleNicaiseKulikov}}. 
\end{proof}

\subsection{Conic-flowers}

Let $F$ be a conic-flower. The conic-flowers are classified in types $2B$, $2C$, $4C$, $6C$, $6D$ and $6E$. Notice that we proved in Corollary~\ref{thm:CM-K3-no-4D} that flowers of type $4D$ cannot occur.
The following lemma gives some information about the `candidate monodromy eigenvalue' induced by the top of a flower of type $2$ or $6$.

\begin{lemma} \label{thm:holds-conic-flower-xi}
Let $X$ be a $K3$ surface over $K$ with Crauder-Morrison model~$\X$.
If $F$ is a conic-flower in $\X_k$ of type $2$ or $6$, i.e., of type $2B, 2C, 6C, 6D$ or $6E$, then $\xi(F_0)$ is not an $N(F_{\ell+1})$-th root of unity.
\end{lemma}

\begin{proof}
Let $\xi$ be a primitive $N(F_0)$-th root of unity, such that
\[\xi(F_0) = \xi ^{\nu(F_0)}.\]
\begin{enumerate}[(i)]

\item If $F$ is a conic-flower of type $2$, then $N(F_0)=2N(F_{\ell+1})$, by Table~\ref{table:flowers-P2-conic}.
Suppose first that $F$ is a flower of type $2B$. In Corollary~\ref{thm:CM-flowers-relation-mu}, we computed that $\nu(F_0) = 2 \nu(F_{\ell+1}) +1 $. Therefore
\[\xi(F_0) ^{N(F_{\ell+1})} = \xi^{\nu(F_0) N(F_{\ell+1})} = \xi^{\nu(F_{\ell+1}) N(F_0)} \xi^{N(F_{\ell+1})} = \xi^{N(F_0)/2}  \neq 1,\]
since $\xi$ is a primitive $N(F_0)$-th root of unity.

Analogously, when $F$ is a flower of type $2C$, we have $\nu(F_0) = 2\nu(F_{\ell+1})+2\ell-1$, and it follows that $\xi(F_0)^{N(F_{\ell+1})}  \neq 1$.

\item If $F$ is a conic-flower of type $6$, then $3N(F_0)=2N(F_{\ell+1})$, by Table~\ref{table:flowers-P2-conic}.
Suppose first that $F$ is a flower of type $6C$. In Corollary~\ref{thm:CM-flowers-relation-mu}, we computed that $3\nu(F_0) = 2 \left(\nu(F_{\ell+1}) +3\ell-2\right)-1 $. Therefore,
\begin{align*}
\xi(F_0)^{N(F_{\ell+1})} &= \xi^{\nu(F_0) N(F_{\ell+1})} = \xi^{(\nu(F_{\ell+1})+3\ell-2) N(F_0)} \xi^{-N(F_{\ell+1})/3}\\ &= \xi^{-N(F_0)/2}  \neq 1,
\end{align*}
since $\xi$ is a primitive $N(F_0)$-th root of unity.

Analogously, when $F$ is a flower of type $6D$, we have $3\nu(F_0) = 2\nu(F_{\ell+1})+6\ell - 1$, and when $F$ is a flower of type $6E$, we have $3\nu(F_0) = 2\nu(F_{\ell+1})+ 1$. In both cases, it follows immediately that $\xi(F_0)^{N(F_{\ell+1})}\neq 1$.
%
\end{enumerate}
\vspace{-18pt}
\end{proof}

\begin{lemma}
Let $X$ be a $K3$ surface over $K$ with Crauder-Morrison model~$\X$.
If $F$ is a conic-flower in $\X_k$ of type $4C$, then $N(F_{\ell+1})$ is even and $\xi(F_0)$ is not an $N(F_{\ell+1})/2$-th root of unity. 
\end{lemma}

\begin{proof}
In Table~\ref{table:flowers-P2-conic}, we can see that $N(F_{\ell+1})=N(F_0)$, and that $N(F_{\ell+1})$ is even.
Let $\xi$ be a primitive $N(F_0)$-th root of unity, such that
\[\xi(F_0) = \xi ^{\nu(F_0)}.\]
In Lemma~\ref{thm:CM-flowers-numerical-relations}, we computed that $\nu(F_0) = 2 \nu(F_{1}) +1 $. Therefore
\[\xi(F_0) ^{N(F_{\ell+1})/2} = \xi^{\nu(F_0) N(F_{0})/2} = \xi^{\nu(F_{1}) N(F_0)} \xi^{N(F_{0})/2} = \xi^{N(F_0)/2}  \neq 1,\]
since $\xi$ is a primitive $N(F_0)$-th root of unity.

\end{proof}

\section{Flowerpot degenerations}\label{sect:holds-flowerpot}

The aim of this section is to prove the following theorem.
\begin{theorem} \label{thm:holds-pot}
Let $X$ be a $K3$ surface over $K$ with Crauder-Morrison model~$\X$.
If $\X$ is a flowerpot degeneration, then $X$ satisfies the monodromy property. 
\end{theorem}

This is a consequence of Theorems~\ref{thm:holds-easy-pot} and \ref{thm:holds-difficult-pot}, combined with Proposition~\ref{thm:CM-geometry-pot}.

Throughout this section, we will assume that $\X$ is a flowerpot degeneration, i.e., there is a unique component $P$ of $\X_k$ with minimal weight $\rho(P)$.

In Subsection~\ref{sect:holds-4C}, we discuss why the existence of flowers of type $4C$ causes an extra difficulty. In Subsection~\ref{sect:holds-easy-pot}, we show that $X$ satisfies the monodromy property if $P\simeq\P^2$, if $P$ is minimal ruled, or if $K_P\equiv 0$.
Finally in Subsection~\ref{sect:holds-ruled-pot}, we show that $X$ also satisfies the monodromy property in the only remaining case: when $P$ is a rational, non-minimal ruled surface.

\subsection{Flowers of type \texorpdfstring{$4C$}{4C}} \label{sect:holds-4C}

\begin{proposition} \label{thm:holds-no-4C}
Let $X$ be a $K3$ surface over $K$ with Crauder-Morrison model~$\X$. Suppose that $\X$ is a flowerpot degeneration. If $F$ is a conic-flower of type 2 or 6, then $\xi(F_0)$ is a monodromy eigenvalue for $X$.
\end{proposition}

\begin{proof}
We have to prove that $\xi(F_0)$ is a pole of $\zeta_X(T) = \prod_{i\in I} (T^{N_i}-1)^{-\chi(E_i^\circ)}$. From Lemma~\ref{thm:CM-euler-characteristics}, we see that the only factors in this formula with a \emph{positive} exponent are either 
\[(T^{N(P)}-1)^{-\chi(P^\circ)},\]
when $\chi(P^\circ)<0$, or
\[(T^{N(F_0')}-1)^{2g(F')-2},\]
where $F'$ is a flower of genus $g(F')>1$. Notice that for  any non-rational flower $F'$, the multiplicity~$N(F'_0)$ is a divisor of $N(P)$, by Table~\ref{table:flowers-ruled}. Therefore, $T^{N(F_0')}-1$ is a divisor of $T^{N(P)}-1,$ for any non-rational flower $F'$.

Since $F$ is a flower of type $2$ or $6$, Corollary~\ref{thm:holds-conic-flower-xi} implies that $\xi(F_0)$ is not a zero of $T^{N(P)}-1$ and $T^{N(F_0')}-1$. For this reason, $\xi(F_0)$ is a pole of $\zeta_X(T)$.
\end{proof}

If a conic-flower in a flowerpot degeneration of a $K3$ surface is not of type 2 or 6, then it is a flower of type 4C. The proof of the previous proposition cannot be adapted to prove that $\xi(F_0)$ is a monodromy eigenvalue for flowers $F$ of type $4C$. The reason is that from Table~\ref{table:flowers-P2-conic}, we see that $N(F_0)=N(P)$ and therefore we need to argue why there is no cancellation with the term $(T^{N(P)}-1)^{-\chi(P^\circ)}$.

\subsection{Easy cases} \label{sect:holds-easy-pot}

In Proposition~\ref{thm:CM-geometry-pot}, we proved that the pot $P$ satisfies one of the following:
\begin{enumerate}[(i)]
\item $P\simeq \P^2$, or
\item $K_P\equiv 0$, or
\item $P$ is a minimal ruled surface, or
\item $P$ is a rational, non-minimal ruled surface.
\end{enumerate}
In the first three cases, it is not difficult to show that $X$ satisfies the monodromy property.

\begin{theorem} \label{thm:holds-easy-pot}
Let $X$ be a $K3$ surface over $K$ with Crauder-Morrison model~$\X$. Suppose $\X$ is a flowerpot degeneration with pot $P$. If $P$ satisfies one of the following:
\begin{enumerate}[(i)]
\item $P\simeq \P^2$, or
\item $K_P\equiv 0$, or
\item $P$ is a minimal ruled surface,
\end{enumerate}
then $X$ satisfies the monodromy property.
\end{theorem}

\begin{proof}
Define $I^\dagger\subset I$ to be the set of indices $i\in I$ with either
$E_i=P$, or
$E_i$ is the top of a conic-flower.
Let $S^\dagger=\left\{(-\nu_i,N_i)\in \Z\times \Z_{>0}\mid i\in I^\dagger\right\}$. Corollary~\ref{thm:poles-contributions} gives that
\[Z_{X,\omega}(T)\in \mathcal{M}_k^{\hat{\mu}}\left[T, \frac{1}{1-\L^{a}T^b}\right]_{(a,b)\in S^\dagger}.\]
We still have to prove that $\xi(F_0)$ is a monodromy eigenvalue for any conic-flower $F$, since it has already been proven in Theorem~\ref{thm:holds-lct} that $\xi(P)$ is a monodromy eigenvalue.

Suppose first $P\simeq \P^2$. By Proposition~\ref{thm:CM-pot-P2}, there are no conic-flowers. Therefore $X$ certainly satisfies the monodromy property.

Suppose now that $K_P\equiv 0$. By Proposition~\ref{thm:CM-pot-NK=0}, any conic-flower is of type $2$. Proposition~\ref{thm:holds-no-4C} shows that $X$ satisfies the monodromy property.

Suppose now that $P$ is a minimal ruled surface. Let $F$ be a conic-flower. By Proposition~\ref{thm:CM-pot-minimal-ruled}, there is exactly one other flower $F'$, and it has genus $g\geq 2$. Let $C$ and $ C'$ be the flowercurves on $P$ of $F$ and $F'$ respectively. Then 
\[\chi(P^\circ) = \chi(P) - \chi(C) - \chi(C') = 4 - 2 - (2-2g) = 2g.\]
Lemma~\ref{thm:CM-euler-characteristics} gives that $\chi\left((F'_0)^\circ\right)=2-2g$ and $\chi(F_0^\circ)=1$. Moreover, $\chi(E_i^\circ)=0$ for $E_i\neq P,F_0$ or $F_0'$. Therefore
\[\zeta_X(T) = \frac{(T^{N(F_0')}-1)^{2g-2}}{(T^{N(P)}-1)^{2g}(T^{N(F_0)}-1)}.\]
Because $N(F_0')$ is a divisor of $N(P)$ by Table~\ref{table:flowers-ruled}, it is clear that $\xi(F_0)$ is a pole of $\zeta_X(T)$, and hence it is a monodromy eigenvalue. This concludes that $X$ satisfies the monodromy property.
\end{proof}

\subsection{Rational, non-minimal ruled pot} \label{sect:holds-ruled-pot}

From Theorem~\ref{thm:holds-easy-pot} and Proposition~\ref{thm:CM-geometry-pot}, it is clear that the only case we still need to prove in this section, is that of $K3$ surfaces allowing a flowerpot degeneration where the pot $P$ is a rational, non-minimal ruled surface. This is not trivial and will require some work. So assume $P$ is a rational, non-minimal ruled surface. One of the main tools will be the fact that $\xi(P)$ is certainly a monodromy eigenvalue by Theorem~\ref{thm:holds-lct}.

To start, we will derive some information on possible flowers on a rational, non-minimal ruled pot. Then, we will prove some lemmas about $\xi(P)$ and finally, we will prove that the monodromy property holds for a $K3$ surface allowing a flowerpot degeneration with a rational, non-minimal ruled pot.

\subsubsection{Flowers on $P$}

We will show the following:
\begin{itemize}
\item If there is a flower of type $4C$, then there is no flower of type $2A$, $4A$, $4B$, $6A$ or $6B$. (Lemma~\ref{thm:holds-4C-no-other-flowers})
\item There are no flowers of type $8\alpha$, $8\beta$, $12\alpha$ or $12\beta$. (Corollary~\ref{thm:CM-pot-no-flower-8-12})
\item If there is a flower of type $4C$, then any flower of genus $g>1$ is of type~$4\alpha$. (Proposition~\ref{thm:holds-pot-only-4alpha})
\end{itemize}

\begin{lemma} \label{thm:holds-4C-no-other-flowers}
Let $X$ be a $K3$ surface over $K$ with a Crauder-Morrison model~$\X$. Suppose that $\X$ is a flowerpot degeneration with pot $P$. If there exists a flower of type $4C$, then there are no flowers of type $2A$, $4A$, $4B$, $6A$ and $6B$.
\end{lemma}
\begin{proof}
Let $F$ be a flower of type $4C$. Lemma~\ref{thm:CM-flowers-numerical-relations} implies that $\nu(F_0) = 2\nu(F_1)+1$,
and hence $\nu(F_0)$ is odd. This yields that $\nu(P)$ is odd as well, by Corollary~\ref{thm:CM-flowers-relation-mu}. On the other hand, if there exists a flower of type $2A$, $4A$, $4B$, $6A$ or $6B$, then Corollary~\ref{thm:CM-flowers-relation-mu} implies that $\nu(P)$ is even, which is a contradiction.
\end{proof}

\begin{lemma} \label{thm:holds-pot-curve-self-intersection}
Let $X$ be a $K3$ surface over $K$ with a Crauder-Morrison model~$\X$. Suppose that $\X$ is a flowerpot degeneration with pot $P$. For any smooth, rational curve $C$ on $P$, the following properties hold.
\begin{enumerate}[(i)]
\item If $C^2\leq -3$, then $C$ is a flowercurve.
\item If $C^2=-2$, then $C$ does not meet any flowercurves of type $\geq 3$.

\end{enumerate}
\end{lemma}
\begin{proof}

We know that a flowercurve on $P$ does not meet any other flowercurves, because there are no triple points. So let $C$ be a smooth, rational curve on $P$ with $C^2\leq -2$, and suppose it is \emph{not} a flowercurve. We will derive a contradiction for $C^2\leq -3$. Moreover, for $C^2=-2$, we will find that $C$ does not meet any flowercurves of type $\geq 3$.
Let $C_1, \ldots, C_n$ be the flowercurves on $P$ of type $\geq 3$.
By \eqref{eq:canonical-bundle-component-general} and Remark~\ref{rmk:CM-canonical-bundle-min-weight}, we can write
\[K_{P/k}\equiv \sum_{i=1}^n a_i C_i,\]
with $-1<a_i<0$. 
The adjunction formula gives
\[-2 = (K_P+C)\cdot C = C^2 + \sum_{i=1}^n a_i(C_i\cdot C).\]
By assumption $C^2\leq -2$ and therefore
\[\sum_{i=1}^n a_i(C_i\cdot C) \geq 0, \]
where the inequality is strict if and only if $C^2\leq -3$.

On the other hand, since $C$ is not a flowercurve, we have $C_i\cdot C\geq 0$ for $i=1, \ldots, n$. Because $a_i<0$, we find
\[\sum_{i=1}^n a_i(C_i\cdot C) \leq 0. \]
This is a contradiction for $C^2\leq -3$. For $C^2=-2$, this is only possible if $C_i\cdot C=0$ for all $i=1, \ldots, n$.
\end{proof}

\begin{lemma}\label{thm:CM-pot-flowercurve-meets-exceptional}
Let $X$ be a $K3$ surface over $K$ with a Crauder-Morrison model~$\X$. Suppose that $\X$ is a flowerpot degeneration with a rational, non-minimal ruled pot $P$. Let $C$ be a flowercurve of type $M\geq 3$. Then $C$ meets a smooth, rational curve $E \subset P$ with $E^2=-1$.
\end{lemma}

\begin{proof}
Let $C$ be a flowercurve of type $ M \geq 3$. The flowercurve~$C$ meets some component $E$ of a reducible fiber of the ruling on $P$. Indeed, if $g(C)\geq 1$, then $C$ is a horizontal curve, and if $g(C)=0$, then $C^2=-M\leq -3$ by Lemma~\ref{thm:CM-self-intersection-flower-curve}, and hence $C$ is either a horizontal curve or a component of a reducible fiber of the ruling on $P$. Since $E$ is a component of a reducible fiber, it is a smooth, rational curve with $E^2 < 0$. Lemma~\ref{thm:holds-pot-curve-self-intersection} implies that $E^2=-1$.
\end{proof}


\begin{lemma} \label{thm:CM-pot-exceptional-meets-flowercurve}
Let $X$ be a $K3$ surface over $K$ with a Crauder-Morrison model~$\X$. Suppose that $\X$ is a flowerpot degeneration with pot $P$. Let $E$ be a smooth, rational curve on $P$ with $E^2=-1$. There are three possibilities:
\begin{enumerate}[(i)]
\item $E$ meets only flowercurves of type 2 and 3, 
\item $E$ meets only flowercurves of type 2 and 4,
\item $E$ meets exactly one flowercurve of type 3 and exactly one of type 6. All other flowercurves meeting $E$ are of type 2.
\end{enumerate}
\end{lemma}

\begin{proof}

Let $E$ be a smooth, rational curve on $P$ with $E^2=-1$. Let $C_1, \ldots, C_n$ be the flowercurves of type $M_i\geq 3$ on $P$. By adjunction, we must have
\[ \sum_{i=1}^n a_i (C_i \cdot E) = 1,\]
where $a_i=  1 - \frac{2}{M_i}$.
Since $M_i\geq 3$, we have $ 1/3\leq a_i <1$, and therefore ${ 1 <\sum_{i=1}^n (C_i \cdot E)\leq 3}$. Moreover, Lemma~\ref{thm:CM-self-intersection-flower-curve} implies that $E\neq C_i$ and therefore $C_i\cdot E\geq 0$.

If $ \sum_{i=1}^n (C_i \cdot E)= 3$, then we must have $M_i=3$ for all $1\leq i\leq n$. This is case~(i). 

If $ \sum_{i=1}^n (C_i \cdot E)=2$, then we have
\begin{align*}
\frac{2}{M_1}\cdot 2&=1 &\text{when } n=1,\\ 
\frac{2}{M_1}+\frac{2}{M_2}&=1& \text{when } n=2.
\end{align*}
For $n=1$, it is clear that $M_1=4$, and hence $E$ meets only flowercurves of type 2 and 4.
For $n=2$, it is straightforward to see that the only possibilities are $M_1=M_2=4$, or $M_1=3$ and $M_2=6$. 
\end{proof}

\begin{corollary}\label{thm:CM-pot-no-flower-8-12}
Let $P$ be a flowerpot that is ruled, but not minimal ruled. There are no flowers of type $8\alpha, 8\beta, 12\alpha$ and $ 12\beta$.
\end{corollary}

\begin{proof}
This follows from Lemma~\ref{thm:CM-pot-flowercurve-meets-exceptional} and Lemma~\ref{thm:CM-pot-exceptional-meets-flowercurve}.
\end{proof}

\begin{lemma} \label{thm:CM-pot-ruled-horizontal-flowers}
Let $X$ be a $K3$ surface over $K$ with a Crauder-Morrison model~$\X$. Suppose that $\X$ is a flowerpot degeneration with a rational, non-minimal ruled pot $P$. Let $\fiber$ be a general fiber of the ruling on $P$. If there is a flowercurve $C$ of type 6 with $C\cdot \fiber \geq 2$, then there is no horizontal flowercurve of type 4.
\end{lemma}

\begin{proof}
Let $C_1, \ldots, C_n$ be the flowercurves of flowers of type $M_i\geq 3$. After renumbering, we can assume $C_1=C$.
Since $\fiber$ is a general fiber on $P$, we have $\fiber^2= 0$ and $g(\fiber)=0$ and hence, the adjunction formula gives
\[
 -2 = K_{P/k} \cdot \fiber + \fiber^2 =  -\sum_{i=1}^n a_i C_i \cdot \fiber,
\]
where  $a_i=  1 - \frac{2}{M_i}$.
By Lemma~\ref{thm:CM-self-intersection-flower-curve}, either $g(C_i)\geq 1$ or $C_i^2<0$ for every $1\leq i\leq n$. Therefore, $ C_i \neq \fiber $ and hence $ C_i \cdot \fiber \geq 0 $. In particular, $C_i \cdot \fiber > 0$ if and only if $C_i$ is horizontal.
 
 Since $C_1$ is a flowercurve of type 6, we have $a_1=2/3$. Suppose $C_2$ is a horizontal flowercurve of type 4, then $a_2 = 1/2$. Thus we get
 \[\sum_{i=1}^n a_i C_i \cdot \fiber \geq \frac{2}{3}\cdot 2 + \frac{1}{2} >2,\]
where the first inequality holds because $C_1\cdot\fiber \geq 2 $ and $ a_i C_i\cdot \fiber\geq 0$ for all $i=3,\ldots,n$. This is a contradiction.
\end{proof}

\begin{lemma} \label{thm:holds-pot-technical-fiber}
Let $V$ be a rational, non-minimal ruled surface. Let $\fiber$ be a reducible fiber of the ruling on $V$ with at least $r+1$ components, for some integer $r\geq 2$. Let $\fiber_1, \ldots, \fiber_r$ be some components of $\fiber$ such that 
$\fiber_1^2 = -1$ and $\fiber_i^2 =-2$ for  $i=2, \ldots, r-1$. Assume moreover that $\fiber_{i} \cdot \fiber_{i+1}\geq 1$ for $i=1, \ldots, r-1$, and otherwise $\fiber_{i} \cdot \fiber_{j}=0$.
Then $\fiber_r^2\leq  -2$.
\end{lemma}

\begin{proof}
We prove this by induction on $r$. Let $\fiber_1, \fiber_2$ be two components of a reducible fiber with at least $3$ components and suppose 
\[\fiber_1^2=-1 \qquad \text{ and }\qquad \fiber_1 \cdot \fiber_2\geq 1.\]
Let $V\to \overline{V}$ be the contraction of $\fiber_1$, and let $\overline{\fiber}$ and $\overline{\fiber}_2$ be the images of $\fiber$ and $\fiber_2$ respectively. Then $\overline{V}$ is a rational, ruled surface, and $\overline{\fiber}$ is a reducible fiber of the ruling on $\overline{V}$. Since $\overline{\fiber}_2$ is a component of $\overline{\fiber}$, we must have
\[\overline{\fiber}_2^2\leq -1.\]
On the other hand, \cite[Proposition~3.2]{Hartshorne} and \cite[Proposition~3.6]{Hartshorne} give that
\[\fiber_2^2<\overline{\fiber}_2^2.\]
This implies $\fiber_2^2\leq  -2$.

Suppose now that there exists a reducible fiber $\fiber$ with at least $r+1$ components, and let $\fiber_1, \ldots, \fiber_r$ be some components of $\fiber$. Assume
$\fiber_1^2 = -1$ and $\fiber_i^2 =-2$ for  $i=2, \ldots, r-1$. Assume moreover that $\fiber_{i} \cdot \fiber_{i+1}\geq 1$ for $i=1, \ldots, r-1$, and otherwise $\fiber_{i} \cdot \fiber_{j}=0$.
Suppose that
\[\fiber_r^2>-2.\]
Since $\fiber_r$ is a component of a reducible fiber, its self-intersection must be strictly negative, and therefore $\fiber_r^2=-1$.
Let $V\to \overline{V}$ be the contraction of $\fiber_r$ and let $\overline{\fiber}$ and $\overline{\fiber_i}$ be the images of $\fiber$ and $\fiber_i$ respectively, for $i=1, \ldots, r-1$. Then $\overline{V}$ is a rational, ruled surface, and $\overline{\fiber}$ is a reducible fiber of the ruling on $\overline{V}$ with at least $r$ components. Moreover, $\overline{\fiber}_1, \ldots, \overline{\fiber}_{r-1}$ are components of $\overline{\fiber}$.
By \cite[Proposition~3.2]{Hartshorne} and \cite[Proposition~3.6]{Hartshorne}, we know that $\overline{\fiber}_j^2=\fiber_j^2$ for $j=1, \ldots, r-2$ and
\[\overline{\fiber}_{r-1}^2> \fiber_{r-1}^2=-2.\]

On the other hand, the surface $\overline{V}$ with the curves $\overline{\fiber}_1, \ldots, \overline{\fiber}_{r-1}$ satisfy the induction hypothesis, and therefore $\overline{\fiber}_{r-1}^2\leq -2$.
We derived a contradiction, so 
\[\fiber_r^2<-1.\]
\end{proof}

\begin{proposition} \label{thm:holds-pot-only-4alpha}
Let $X$ be a $K3$ surface over $K$ with a Crauder-Morrison model~$\X$. Suppose that $\X$ is a flowerpot degeneration with a rational, non-minimal ruled pot $P$. If there exists a flower of type~$4C$, then any flower of genus $g>1$ is of type $4\alpha$.
\end{proposition}
\begin{proof} 
Let $F$ be a flower of type $4C$ with flowercurve $C$. 
Let $F'$ be a flower with flowercurve $C'$ such that the genus $g(C')>1$. Suppose $F'$ is not of type $4\alpha$. We will find a contradiction.

\textbf{Step 1:} \emph{$C$ is a component of a reducible fiber $\fiber^\dagger$ and $C'$ is horizontal.}\\
Since $P$ is a rational, ruled surface and $g(C')>1$, we must have $C'\cdot \fiber \geq 2$ for any fiber $\fiber$ of the ruling on $P$, so $C'$ is horizontal. Since $F'$ is not of type $4\alpha$, Table~\ref{table:flowers-ruled} and Corollary~\ref{thm:CM-pot-no-flower-8-12} give that $F'$ is of type $6$. 
By Lemma~\ref{thm:CM-pot-ruled-horizontal-flowers}, the curve $C$ is vertical. Since $C^2=-4$ by Lemma~\ref{thm:CM-self-intersection-flower-curve}, the curve $C$ is a component of a reducible fiber~$\fiber^\dagger$.

\textbf{Step 2:} \emph{For any component $\fiber_1$ of $\fiber^\dagger$, we have 
\[\fiber_1^2\in \{-1, -2, -4\}.\]}
Let $\fiber_1$ be a component of $\fiber^\dagger$. Since any component of a reducible fiber has strict negative self-intersection, we have $\fiber_1^2\leq -1$. Suppose $\fiber_1^2\leq -3$, and we will prove that $\fiber_1^2=-4$. By Lemma~\ref{thm:holds-pot-curve-self-intersection}, the curve $\fiber_1$ is a flowercurve.
Because $\fiber^\dagger$ is connected, it is path-connected, and hence there are irreducible components $E_0, E_1, \ldots, E_m$ of $\fiber^\dagger$ such that
\[E_i\cdot E_{i+1} >0,\]  
for $i=0, \ldots, m-1$ and where $E_0=C$ and $E_{m}=\fiber_1$. By removing any cycles, we can assume that $E_i\cdot E_j=0$ for all $i,j$ with $j\not\in \{i-1, i, i+1\}$. 

We show, by induction, that $E_i^2\in \{-1, -4\}$ for $i=0, \cdots, m$. It is given that $E_0^2=-4$. Furthermore, by Lemma~\ref{thm:holds-pot-curve-self-intersection}, we also have $E_1^2=-1$. So suppose $E_i^2\in \{-1, -4\}$ for $i\leq k$, we will show that $E_{k+1}^2\in \{-1, -4\}$. 

If $E_k^2=-4$, then $E_k$ is a flowercurve and hence, $E_{k+1}^2=-1$, by Lemma~\ref{thm:holds-pot-curve-self-intersection}.
If $E_k^2=-1$, then  $E_{k-1}^2\neq -1$ by Lemma~\ref{thm:holds-pot-technical-fiber}, so $E_{k-1}^2=-4$. Since Lemma~\ref{thm:holds-pot-technical-fiber} also implies that $E_{k+1}^2\neq -1$, Lemma~\ref{thm:CM-pot-exceptional-meets-flowercurve} gives that $E_{k+1}^2=-2$ or $-4$. We will exclude $E_{k+1}^2=-2$ by contradiction. So suppose $E_{k+1}^2=-2$. Lemma~\ref{thm:holds-pot-technical-fiber} implies that $E_{k+2}^2\leq -2$ and therefore by Lemma~\ref{thm:holds-pot-curve-self-intersection}, we find $E_{k+2}^2= -2$ as well. Applying again Lemma~\ref{thm:holds-pot-curve-self-intersection} together with Lemma~\ref{thm:holds-pot-technical-fiber}, we find that for all $j\geq k+1$, we have $E_j^2=-2$. This is in contradiction with the fact that $E_m^2=\fiber_1^2\leq -3$. We conclude that $E_{k+1}^2=-4$ which ends the induction argument. In particular, we find that $\fiber_1^2=-4$.

\textbf{Step 3:} \emph{For any component $\fiber_1$ of $\fiber^\dagger$ with $\fiber_1^2=-1$, there exists a component~$\fiber_2$ of $\fiber^\dagger$ with $\fiber_2^2=-4$ and $\fiber_1\cdot \fiber_2>0$.}

Let $\fiber_1$ be a component of $\fiber^\dagger$ with $\fiber_1^2=-1$.
Because $\fiber^\dagger$ is connected, it is path-connected, and hence there are irreducible components $E_0, E_1, \ldots, E_m$ of $\fiber^\dagger$ such that
\[E_i\cdot E_{i+1} >0,\]  
for $i=0, \ldots, m-1$ and where $E_0=C$ and $E_{m}=\fiber_1$. By removing any cycles, we can assume that $E_i\cdot E_j=0$ for all $i,j$ with $j\not\in \{i-1,i, i+1\}$. The same argument as in step 2, shows that $E_i^2\in \{-1, -4\}$. Since no two components with self-intersection $-1$ can meet each other by Lemma~\ref{thm:holds-pot-technical-fiber}, we conclude $E_{m-1}^2=-4$. Therefore $E_{m-1}$ satisfies the conditions of the statement.

\textbf{Step 4:} \emph{Contradiction.}\\
Let $\fiber_1$ be a component of $\fiber^\dagger$ such that $C'\cdot \fiber_1>0$. Such a component certainly exists, since $C'$ is horizontal by Step 1. In Step 2, we showed that  $\fiber_1^2\in \{-1, -2, -4\}$. By Lemma~\ref{thm:holds-pot-curve-self-intersection}, we have that $\fiber_1^2=-1$, since $\fiber_1$ meets the flowercurve $C'$.

By Step 3, there exists a component $\fiber_2$ of $\fiber^\dagger$ with $\fiber_2^2=-4$ and $\fiber_1\cdot \fiber_2>0$. Lemma~\ref{thm:holds-pot-curve-self-intersection} gives that $\fiber_2$ is a flowercurve. This is a contradiction with Lemma~\ref{thm:CM-pot-exceptional-meets-flowercurve}. 
\end{proof}

\begin{remark}
If we additionally assume that $X$ has a rational point, then Proposition~\ref{thm:holds-pot-only-4alpha} is more easily proven. Indeed, let $F$ be a flower of type $4C$ and let $F'$ be a flower of genus $g>1$, which is not of type $4\alpha$. Lemma~\ref{thm:CM-pot-no-flower-8-12} implies that $F'$ is of type $6\alpha$ or $6\beta$. Let $N_P$ be the multiplicity of the pot $P$. Because $F$ is a flower of type $4\alpha$, Table~\ref{table:flowers-P2-conic} gives that $N_P$ is divisible by 2. Moreover, Table~\ref{table:flowers-ruled} shows that $N_P$ is divisible by 3, because $F'$ is a flower of type $6\alpha$ or $6\beta$. Therefore, $N_P$ is divisible by 6. Since there are no flowers of type $6A$, $6B$, $12\alpha$ and $12\beta$ by Lemma~\ref{thm:holds-4C-no-other-flowers} and Corollary~\ref{thm:CM-pot-no-flower-8-12}, Tables~\ref{table:flowers-P2-line}, \ref{table:flowers-P2-conic} and \ref{table:flowers-ruled} give that there is no component of $\X_k$ with multiplicity 1. Therefore $X$ does not have a rational point.
\end{remark}

\subsubsection{The monodromy eigenvalue $\xi(P)$}

\begin{lemma} \label{thm:holds-pot-xi-P-root-unity}
Let $X$ be a $K3$ surface over $K$ with a Crauder-Morrison model~$\X$. Suppose that $\X$ is a flowerpot degeneration with pot $P$. The following statements hold.
\begin{enumerate}[(i)]
\item If there exists a flower $F$ of type $4C$, then $N(P)$ is even and $\xi(P)$ is not an $N(P)/2$-th root of unity.
\item If there exists a flower $F$ of type $3A$, $3B$, $6\alpha$ or $6\beta$, then $N(P)$ is divisible by 3 and $\xi(P)$ is not an $N(P)/3$-th root of unity.
\item If there exists a flower $F$ of type $6C$, $6D$ or $6E$, then $N(P)$ is divisible by 3 and $\xi(P)$ is not an $2N(P)/3$-th root of unity.
\end{enumerate}
\end{lemma}

\begin{proof}
Let $\xi$ be a primitive $N(P)$-th root of unity such that
\[\xi(P) = \xi^{\nu(P)}.\]

\begin{enumerate}[(i)]
\item Let $F$ be a flower of type  $4C$. From Table~\ref{table:flowers-P2-conic}, we can see that $N(P)=N(F_0)$, and moreover, $N(P)$ is even. From Lemma~\ref{thm:CM-flowers-numerical-relations}, it follows that $\nu(F_0)$ is odd, and hence by Corollary~\ref{thm:CM-flowers-relation-mu}, $\nu(P)$ is odd as well. Therefore
\[\xi(P) ^{N(P)/2} = \xi^{\nu(P) N(P)/2} = \xi^{N(P)/2}  \neq 1,\]
since $\xi$ is a primitive $N(P)$-th root of unity.

\item Let $F$ be a flower of type  $3A$, $3B$, $6\alpha$ or $6\beta$. From Tables~\ref{table:flowers-P2-line} and \ref{table:flowers-ruled}, we see that $N(P)$ is divisible by 3.
If $F$ is a flower of type $3A$, then $\nu(P) = 3\nu(F_0)-4$, by Corollary~\ref{thm:CM-flowers-relation-mu}. Therefore
\[\xi(P)^{N(P)/3} = \xi^{\nu(P)N(P)/3} = \xi^{(\nu(F_0)-1)N(P)}\xi^{-N(P)/3}=\xi^{-N(P)/3}\neq 1,\]
since $\xi$ is a primitive $N(P)$-the root of unity. Completely analogously, we compute $\xi(P)^{N(P)/3}\neq 1$, when $F$ is a flower of type  $3B$, $6\alpha$ or $6\beta$.

\item Let $F$ be a flower of type  $6C$, $6D$ or $6E$. From Table~\ref{table:flowers-P2-conic}, we see that $N(P)$ is divisible by 3. If $F$ is a flower of type $6C$, then $\nu(P) = (3\nu(F_0)-6\ell+5)/2$, by Corollary~\ref{thm:CM-flowers-relation-mu}. Therefore
\begin{align*}
\xi(P)^{2N(P)/3} &= \xi^{2\nu(P)N(P)/3} = \xi^{(\nu(F_0)-2\ell+2)N(P)}\xi^{-N(P)/3}\\
&=\xi^{-N(P)/3}\neq 1,
\end{align*}
since $\xi$ is a primitive $N(P)$-the root of unity.  Completely analogously, we compute $\xi(P)^{2N(P)/3}\neq 1$, when $F$ is a flower of type $6D$ or $6E$.

%
\end{enumerate}
\end{proof}

\subsubsection{The monodromy property holds}

\begin{theorem} \label{thm:holds-difficult-pot}
Let $X$ be a $K3$ surface over $K$ with Crauder-Morrison model~$\X$. Suppose $\X$ is a flowerpot degeneration with pot $P$. If $P$ is a rational, non-minimal ruled surface, then $X$ satisfies the monodromy property.
\end{theorem}

\begin{proof}
Define $I^\dagger\subset I$ to be the set of indices $i\in I$ with either
$E_i=P$, or
$E_i$ is the top of a conic-flower.
Let $S^\dagger=\left\{(-\nu_i,N_i)\in \Z\times \Z_{>0}\mid i\in I^\dagger\right\}$. Corollary~\ref{thm:poles-contributions} says that
\[Z_{X,\omega}(T)\in \mathcal{M}_k^{\hat{\mu}}\left[T, \frac{1}{1-\L^{a}T^b}\right]_{(a,b)\in S^\dagger}.\]
Let $F$ be a conic-flower. We need to prove that $\xi(F_0)$ is a monodromy eigenvalue, since it has already been proven in Theorem~\ref{thm:holds-lct} that $\xi(P)$ is a monodromy eigenvalue.

To simplify notation, let $N=N(P)$ be the multiplicity of the pot~$P$.
If $F$ is not a flower of type $4C$, then Proposition~\ref{thm:holds-no-4C} concludes the proof. So suppose that $F$ is a flower of type $4C$. Since $N(F_0)=N$, we have that $\xi(F_0)$ is an $N$-th root of unity.

Denote by $B_*$ the number of flowers meeting $P$ in a \emph{rational} curve of type $*$. For instance, $B_{6\alpha\beta}$ denotes the number of flowers of type $6\alpha$ and $6\beta$ meeting $P$ in a \emph{rational} curve.

By Lemma~\ref{thm:CM-pot-hodge-index}, one of the following two possibilities hold.
\begin{enumerate}[(i)]
\item There is a unique flower meeting $P$ in a curve of genus $g>1$. All other flowercurves on $P$ are rational.
\item All flowercurves on $P$ have genus 0 or 1.
\end{enumerate}

Suppose first that all flowercurves on $P$ have genus 0 or 1. Define
\[ Q_1(T) = (T^{N/3}-1)^{B_{3AB}+2B_{6\alpha\beta}}(T^{2N/3}-1)^{B_{6CDE}}.\]
The monodromy zeta function can be written as 
\[\zeta_X(T) = \frac{(T^N-1)^{-\chi(P^\circ)}}{(T^N-1)^{B_{4C}}(T^{2N}-1)^{B_{2BC}}(T^{N/2}-1)^{2B_{4\alpha}}Q_1(T)},\]
which can be rewritten to
\[\zeta_X(T) = \frac{1}{(T^N-1)^{B_{4C}+B_{2BC}+\chi(P^\circ)}(T^{N}+1)^{B_{2BC}}(T^{N/2}-1)^{2B_{4\alpha}}Q_1(T)}.\]
Since $\xi(P)$ is a monodromy eigenvalue, it is a pole of $\zeta_X(T)$ by Theorem~\ref{thm:holds-form-zeta-function-K3}.
Theorem~\ref{thm:holds-pot-xi-P-root-unity} (i) gives that $\xi(P)$ is not an $N/2$-th root of unity. Moreover, $\xi(P)$ is not a zero of $Q_1(T)$, by Theorem~\ref{thm:holds-pot-xi-P-root-unity} (ii) and (iii). Therefore 
\[B_{4C}+B_{2BC}+\chi(P^\circ)>0,\]
which implies that $\xi(F_0)$ is a pole of the monodromy zeta function, and hence it is a monodromy eigenvalue.

Suppose now that there is a (necessarily unique) flower with genus $g>1$. By Proposition~\ref{thm:holds-pot-only-4alpha}, this flower is of type $4\alpha$.
The monodromy zeta function can be written as 
\[\zeta_X(T) = \frac{(T^N-1)^{-\chi(P^\circ)}(T^{N/2}-1)^{2g-2}}{(T^N-1)^{B_{4C}}(T^{2N}-1)^{B_{2BC}}(T^{N/2}-1)^{2B_{4\alpha}}Q_1(T)},\]
where $Q_1(T)$ is defined as before.
This can be rewritten to
\[\zeta_X(T) = \frac{1}{(T^N-1)^{B_{4C}+B_{2BC}+\chi(P^\circ)}(T^{N}+1)^{B_{2BC}}(T^{N/2}-1)^{2B_{4\alpha}-2g+2}Q_1(T)}.\]
Since $\xi(P)$ is a monodromy eigenvalue, it is a pole of $\zeta_X(T)$ by Theorem~\ref{thm:holds-form-zeta-function-K3}.
Theorem~\ref{thm:holds-pot-xi-P-root-unity} (i) gives that $\xi(P)$ is not an $N/2$-th root of unity. Moreover, $\xi(P)$ is not a zero of $Q_1(T)$, by Theorem~\ref{thm:holds-pot-xi-P-root-unity} (ii) and (iii). Therefore 
\[B_{4C}+B_{2BC}+\chi(P^\circ)>0,\]
which implies that $\xi(F_0)$ is a pole of the monodromy zeta function, and hence it is a monodromy eigenvalue.
\end{proof}

\section{Chain degenerations with an extra condition} \label{sect:holds-chain}

In this section, we will prove the following theorem.  
\begin{theorem} \label{thm:holds-chain}
Let $X$ be a $K3$ surface over $K$ with Crauder-Morrison model~$\X$ with special fiber $\X_k=\sum_{i\in I} N_iE_i$.
Suppose $\X$ is a chain degeneration with chain $V_0 \text{---} V_1 \text{---} \cdots \text{---} V_k \text{---} V_{k+1}$. Set $N=\min\left\{N(V_i)\mid i=0, \ldots, k+1\right\}$. Assume moreover that one of the following conditions hold:
\begin{enumerate}[(i)]
\item the components $V_0, V_1, \ldots, V_{k+1}$ all have  multiplicity $N$, or
\item neither $V_0$ nor $V_{k+1}$ is a rational, non-minimal ruled surface, or
\item exactly one end component of the chain is a rational, non-minimal ruled surface, and it has multiplicity $N$, or
\item if $V_j$ meets a conic-flower, then $N(V_j)>N$.
\end{enumerate}
Then $X$ satisfies the monodromy property.
\end{theorem}

\begin{proof}
Define $I^\dagger\subset I$ to be the set of indices $i\in I$ where either
$\rho_i$ is minimal, or
$E_i$ is the top of a conic-flower.
Let $S^\dagger=\left\{(-\nu_i,N_i)\in \Z\times \Z_{>0}\mid i\in I^\dagger\right\}$. Corollary~\ref{thm:poles-contributions} says that
\[Z_{X,\omega}(T)\in \mathcal{M}_k^{\hat{\mu}}\left[T, \frac{1}{1-\L^{a}T^b}\right]_{(a,b)\in S^\dagger}.\]
Let $F$ be a conic-flower. Denote by $F_0$ the top of the flower $F$ and by $F_{\ell+1}$ the component of the chain meeting the flower $F$. By Proposition~\ref{thm:CM-chain-flowers}, $F$ is a flower of type $2B$ or $2C$. Moreover, Theorem~\ref{thm:CM-classification-flowers} and Proposition~\ref{thm:CM-chain-good} imply that $N(F_0)=2N$ or $4N$.
 We need to prove that $\xi(F_0)$ is a monodromy eigenvalue, since it has already been proven in Theorem~\ref{thm:holds-lct} that $\xi(E_i)$ is a monodromy eigenvalue if $\rho_i$ is minimal.

\begin{enumerate}[(i)]
\item
Suppose first that $V_0, V_1, \ldots, V_{k+1}$ all have multiplicity $N$. Table~\ref{table:flowers-P2-conic} gives that $N(F_0)=2N$. By Proposition~\ref{thm:CM-chain-flowers}, all flowers are of type 2A, 2B, 2C or 4$\alpha$, and moreover flowers of type 4$\alpha$ can only meet a component of the chain in $V_0$ or $V_{k+1}$.

To show that $\xi(F_0)$ is a monodromy eigenvalue, we will prove that it is a pole of the monodromy zeta function $\zeta_X(T)$. Proposition~\ref{thm:GMP-ACampo-dvr} gives that
\[\zeta_X(T)=\prod_{i\in I}\left(T^{N_i}-1\right)^{-\chi(E_i^\circ)},\]
where $\chi(E_i^\circ)$ is the topological Euler characteristic of $E_i^\circ$. 
From Lemma~\ref{thm:CM-euler-characteristics}, we see that the only factors in this formula with a \emph{positive} exponent are either 
\[(T^{N}-1)^{-\sum_{i=0}^{k+1}\chi(V_i^\circ)},\]
when $\sum_{i=0}^{k+1}\chi(V_i^\circ)<0$, or
\[(T^{N(F_0')}-1)^{2g(F')-2},\]
if there exists a flower $F'$ of type $4\alpha$ and of genus $g(F')>1$. Notice that $N(F'_0)=N/2$, by Table~\ref{table:flowers-ruled}. Therefore, $T^{N(F_0')}-1$ is a divisor of $T^{N}-1$.
Since $F$ is a flower of type $2$, Corollary~\ref{thm:holds-conic-flower-xi} implies that $\xi(F_0)$ is not a zero of $T^{N}-1$ and $T^{N(F_0')}-1$. For this reason, $\xi(F_0)$ is a pole of $\zeta_X(T)$. 

\item
Assume now that neither $V_0$ nor $V_{k+1}$ is a rational, non-minimal ruled surface. By Proposition~\ref{thm:CM-chain-flowers}, all flowers are of type 2A, 2B, 2C or 4$\alpha$, and moreover flowers of type 4$\alpha$ can only meet a component of the chain in $V_0$ or $V_{k+1}$.

To show that $\xi(F_0)$ is a monodromy eigenvalue, we will prove that it is a pole of the monodromy zeta function $\zeta_X(T)$. Proposition~\ref{thm:GMP-ACampo-dvr} gives that
\[\zeta_X(T)=\prod_{i\in I}\left(T^{N_i}-1\right)^{-\chi(E_i^\circ)},\]
where $\chi(E_i^\circ)$ is the topological Euler characteristic of $E_i^\circ$. 
From Lemma~\ref{thm:CM-euler-characteristics}, we see that the only factor in this formula with a \emph{positive} exponent is
\[(T^{N(F_0')}-1)^{2g(F')-2},\]
if there exists a flower $F'$ of type $4\alpha$ and of genus $g(F')>1$. Hence, if there is no such flower, then we are done. So suppose there is a flower $F'$ of type $4\alpha$ and of genus $g(F')>1$. Since only $V_0$ or $V_{k+1}$ can meet flowers of type $4\alpha$, we can assume, without loss of generality, that $V_0$ meets $F'$. From Proposition~\ref{thm:CM-chain-flowers}, we see that $V_0$ is an elliptic, ruled surface and hence $N(V_0)=N$ or $2N$, by Proposition~\ref{thm:CM-chain-good}(iv). Therefore, Theorem~\ref{thm:CM-classification-flowers} gives that $N(F_0')=N/2$ or $N$.
Since $F$ is a flower of type $2$ and $N(F_0)=2N$ or $4N$, Corollary~\ref{thm:holds-conic-flower-xi} implies that $\xi(F_0)$ is not a zero of $T^{N(F_0')}-1$. For this reason, $\xi(F_0)$ is a pole of $\zeta_X(T)$.

\item Assume that there is exactly one end component of the chain that is a rational, non-minimal ruled surface. Without loss of generality, we can assume that this component is $V_0$. Assume that $N(V_0)=N$.

To show that $\xi(F_0)$ is a monodromy eigenvalue, we will prove that it is a pole of the monodromy zeta function $\zeta_X(T)$. Proposition~\ref{thm:GMP-ACampo-dvr} gives that
\[\zeta_X(T)=\prod_{i\in I}\left(T^{N_i}-1\right)^{-\chi(E_i^\circ)},\]
where $\chi(E_i^\circ)$ is the topological Euler characteristic of $E_i^\circ$. 
From Lemma~\ref{thm:CM-euler-characteristics}, we see that the only factor in this formula with a \emph{positive} exponent are either
\[(T^{N}-1)^{-\sum_{i=0}^{\beta}\chi(V_i^\circ)},\]
when $\sum_{i=0}^{\beta}\chi(V_i^\circ)<0$ and where $\beta$ is defined as in Proposition~\ref{thm:CM-chain-good}, or
\[(T^{N(F_0')}-1)^{2g(F')-2},\]
if there exists a flower $F'$ of type $4\alpha$ and of genus $g(F')>1$. 

Suppose there is a flower $F'$ of type $4\alpha$ and of genus $g(F')>1$. From Proposition~\ref{thm:CM-chain-flowers}, it follows that only end components that are elliptic ruled can meet flowers of type $4\alpha$, and therefore we know that $V_{k+1}$ meets $F'$, and that $V_{k+1}$ is an elliptic, ruled surface. Hence $N(V_{k+1})=N$ or $2N$, by Proposition~\ref{thm:CM-chain-good}(iv). Therefore, Theorem~\ref{thm:CM-classification-flowers} gives that $N(F_0')=N/2$ or $N$.

Since $F$ is a flower of type $2$ and $N(F_0)=2N$ or $4N$, Corollary~\ref{thm:holds-conic-flower-xi} implies that $\xi(F_0)$ is not a zero of $T^N-1$ and $T^{N(F_0')}-1$. For this reason, $\xi(F_0)$ is a pole of $\zeta_X(T)$.

\item Finally, assume that every component $V_j$ meeting a conic-flower has multiplicity $N(V_j)>N$. Therefore, $N(F_{\ell+1})>N$, since $F$ is a conic-flower. Statements (ii) and (iii) of Proposition~\ref{thm:CM-chain-good} imply that we can assume, without loss of generality, that $F_{\ell+1}=V_0$. Moreover, from Proposition~\ref{thm:CM-chain-good} (iv), it follows that $V_0$ is a ruled surface and that $N(V_0)=2N$. Therefore, $N(F_0)=4N$, by Theorem~\ref{thm:CM-classification-flowers}.

To show that $\xi(F_0)$ is a monodromy eigenvalue, we will prove that it is a pole of the monodromy zeta function $\zeta_X(T)$. Proposition~\ref{thm:GMP-ACampo-dvr} gives that
\[\zeta_X(T)=\prod_{i\in I}\left(T^{N_i}-1\right)^{-\chi(E_i^\circ)},\]
where $\chi(E_i^\circ)$ is the topological Euler characteristic of $E_i^\circ$.

Suppose first that $N(V_{k+1})=2N$. 
From Lemma~\ref{thm:CM-euler-characteristics}, we see that the only factors in this formula with a \emph{positive} exponent are either 
\[(T^{2N}-1)^{-\chi(V_0^\circ)-\chi(V_{k+1}^\circ)},\]
when $-\chi(V_0^\circ)-\chi(V_{k+1}^\circ)<0$, or
\[(T^{N(F_0')}-1)^{2g(F')-2},\]
if there exists a flower $F'$ of type $4\alpha$ and of genus $g(F')>1$. 
If there is a flower $F'$ of type $4\alpha$ and of genus $g(F')>1$, then $N(F_0')=N$, by Theorem~\ref{thm:CM-classification-flowers} and Proposition~\ref{thm:CM-chain-flowers}.
Since $F$ is a flower of type $2$ and $N(F_0)=4N$, Corollary~\ref{thm:holds-conic-flower-xi} implies that $\xi(F_0)$ is not a zero of $T^N-1$ and $T^{N(F_0')}-1$. For this reason, $\xi(F_0)$ is a pole of $\zeta_X(T)$.

Suppose now that $N(V_{k+1})\neq 2N$. 
From Lemma~\ref{thm:CM-euler-characteristics}, we see that the only factors in the formula for $\zeta_X(T)$ with a \emph{positive} exponent are either 
\[(T^{2N}-1)^{-\chi(V_0^\circ)},\]
when $-\chi(V_0^\circ)<0$, or
\[(T^{N(F_0')}-1)^{2g(F')-2},\]
if there exists a flower $F'$ of type $4\alpha$ and of genus $g(F')>1$. 
If there is a flower $F'$ of type $4\alpha$ and of genus $g(F')>1$, then $N(F_0')=N/2$ or $N$, by Theorem~\ref{thm:CM-classification-flowers} and Proposition~\ref{thm:CM-chain-flowers}.
Since $F$ is a flower of type $2$ and $N(F_0)=4N$, Corollary~\ref{thm:holds-conic-flower-xi} implies that $\xi(F_0)$ is not a zero of $T^N-1$ and $T^{N(F_0')}-1$. For this reason, $\xi(F_0)$ is a pole of $\zeta_X(T)$.
\end{enumerate}
\end{proof}

\cleardoublepage


\chapter{Future research: a proof or a counterexample?}\label{ch:future}

Let $X$ be a $K3$ surface over $K$ admitting a triple-point-free model, and let $\X$ be a Crauder-Morrison model of $X$. Instead of proving directly that $X$ satisfies the monodromy property, we can ask ourselves the following question:

\emph{If $X$ does \emph{not} satisfy the monodromy property, then what should $\X_k$ look like?} 

Descriptions of such special fibers are called \emph{combinatorial countercandidates}.
A combinatorial countercandidate is a combinatorial characterization of a special fiber, such that, \emph{if} it occurs as the special fiber of a $K3$ surface $X$, then $X$ does not satisfy the monodromy property. In other words, if there exists a $K3$ surface $X$ admitting an $snc$-model with a combinatorial countercandidate as special fiber, then $X$ is a counterexample of the monodromy property.

The interest of combinatorial countercandidates is twofold: on the one hand, with a concrete special fiber in mind that would lead to a counterexample of the monodromy property, one can try to prove that there indeed exists a $K3$ surface with a model with this particular special fiber. So if we can smoothen a combinatorial countercandidate to a $K3$ surface, then we will have constructed a $K3$ surface not satisfying the monodromy property. On the other hand, it is clearly beneficial to have deep knowledge about the properties of possible $K3$ surfaces not satisfying the monodromy property. If these properties are contradictory, no such $K3$ surface exists. Therefore, a strategy to prove that $K3$ surfaces admitting a triple-point-free model satisfy the monodromy property, could be to exclude all combinatorial countercandidates, i.e., find reasons why combinatorial countercandidates do \emph{not} occur as special fibers of $K3$ surfaces. 

In the first section, we will give some more useful results that will be used throughout this chapter.
In Section~\ref{sect:future-comb-count}, we will explain a strategy to produce \emph{all} combinatorial countercandidates.
One application of combinatorial countercandidates would be to smoothen one to a $K3$ surface. A strategy to do so, is described in Section~\ref{sect:future-exist}.
Another application would be to rule out all combinatorial countercandidates, and we illustrate this strategy in Section~\ref{sect:future-exclude}, where we will indeed exclude some combinatorial countercandidates. 
The final section of this chapter is devoted to discussing some questions for further research.

\subsection*{Notation}

In this chapter, we fix an algebraically closed field $k$ of characteristic zero. Put $R=k\llbracket t\rrbracket$ and $K=k(\!( t )\!)$. Fix an algebraic closure $K^{alg}$ of $K$.

For any ruled surface $V$ over $k$, we define $L$ to be the number of blow-ups in the contraction $V\to \overline{V}$ to the minimal ruled surface $\overline{V}$. In particular, when $V$ is minimal ruled, we have $L=0$.

If $X$ is a $K3$ surface over $K$ with $\X$ a Crauder-Morrison model of $X$ over $R$ and if $\omega$ is a volume form on $X$, then we write the special fiber as $\X_k=\sum_{i\in I} N_i E_i$. We denote by $(N_i, \nu_i)$ the numerical data of $E_i$, for every $i\in I$. To avoid confusion, we will sometimes use the notation $(N(E_i), \nu(E_i))$ instead of $(N_i,\nu_i)$. Set $\rho_i = \nu_i/N_i +1$ to be the weight of $E_i$, for every $i\in I$.
For every component~$E_i$ of $\X_k$, we define
\[\xi(E_i) = \exp\left(-2\pi i \, \frac{\nu(E_i)}{N(E_i)}\right),\] 
the `candidate' monodromy eigenvalue associated with $E_i$.

\section{Some more results} \label{sect:future-results}

\subsection{Index of a variety}

\begin{definition} \label{def:index}
Let $X$ be a variety over $K$.
The index $\iota(X)$ of $X$ is defined as the greatest common divisor of the degrees $[L : K]$ of
all finite extensions~$L/K$ such that $X(L) \neq \emptyset$.
\end{definition}
In particular, if $X$ has a $K$-rational point, then $\iota(X)=1$. The converse is not true.

The following characterization of the index is well known.
\begin{proposition} \label{thm:future-index-gcd}
Let $X$ be variety over $K$. Let $\X$ be an $snc$-model of $X$ with $\X_k=\sum_{i\in I} N_i E_i$. Then
\[\iota(X) = \gcd_{i\in I} N_i.\]
\end{proposition}

The following statement is a special case of a theorem proved by Esnault, Levine and Wittenberg.
\begin{theorem} \label{thm:future-index-K3}
Let $X$ be a $K3$ surface over $K$. The index $\iota(X)$ is either 1 or 2.
\end{theorem}
\begin{proof}
This follows from \cite[Theorem~2.1]{ELW}, where it is shown that $\iota(X)$ is a divisor of the holomorphic Euler characteristic of $X$, and the fact that the holomorphic Euler characteristic of a $K3$ surface is~2.
\end{proof}

\begin{lemma}
Let $X$ be a $K3$ surface over $K$ of index $\iota(X)=2$. 
Denote by $K(2)$ the unique totally ramified extension of degree $2$ of $K$ in $K^{alg}$, and define $X(2)=X\times_K K(2)$. Then $X(2)$ is a $K3$ surface with index $\iota(X(2))=1$.
\end{lemma}
\begin{proof}
The surface $X(2)$ is proper, geometrically connected and smooth, as these properties are stable under base change. Moreover, \cite[Proposition~II.8.10]{Hartshorne} guarantees that $X(2)$ has trivial canonical bundle, and \cite[Tag~02KH]{stacks-project} that $H^1\left(X(2), \O_{X(2)}\right) = 0$. Therefore, $X(2)$ is a $K3$ surface.


Every finite extension $L$ of $K$ in $K^{alg}$ of even degree contains $K(2)$, and therefore we have that $X$ has an $L$-rational point if and only if $X(2)$ has an $L$-rational point. Since $[L:K]=2\cdot [L:K(2)]$, we conclude that $\iota(X(2))=1$.
\end{proof}

\begin{lemma} \label{thm:future-index}
Let $X$ be a $K3$ surface over $K$ of index $\iota(X)=2$. Denote by $K(2)$ the unique totally ramified extension of degree $2$ of $K$ in $K^{alg}$, and define $X(2) = X\times_K K(2)$. If $X(2)$ satisfies the monodromy property, then $X$ satisfies the monodromy property as well.
\end{lemma}

\begin{proof}
Let $\omega$ be a volume form on $X$, and let $\omega(2)$ be the pullback of $\omega$ to $X(2)$. Let $\sigma$ be a topological generator of $\Gal(K^{alg}/K)$. Then $\sigma^2$ is a topological generator of $\Gal(K^{alg}/K(2))$. Note that obviously, we have $X\times_K K^{alg} = X(2) \times_{K(2)} K^{alg}$.

For every integer $d\geq 1$, denote by $K(d)$ the unique totally ramified extension of degree $d$ of $K$ in $K^{alg}$, and define $X(d) = X\times_K K(d)$.  Let $\omega(d)$ be the pullback of $\omega$ to $X(d)$.
Because the index $\iota(X)=2$, we have $X(K(d))=\emptyset$, whenever $d$ is odd. Therefore, when $d$ is odd, $X(d)$ itself is a weak N\'eron model of $X(d)$, and hence $\int_{X(d)}|\omega(d)|=0$.
\note{
Let $\Y$ be a weak N\'eron model of $X(d)$ over $R(d)$. Then $\Y(R(d))=X(d)(K(d))=\emptyset$. Since all components of $\Y$ are smooth, there is a surjection $\Y(R(d))\to \Y_k$ (something about that through every point on a component of $\Y_k$, there is a section.) This means that the special fiber is empty.
}
 From Definition~\ref{def:motivic-zeta-function-CY}, we see that
\begin{equation}
Z_{X,\omega}(T) = Z_{X(2),\omega(2)}(T^2). \label{eq:future-zeta-extension}
\end{equation}
Assume that $X(2)$ satisfies the monodromy property. Define $S^{(2)}\subset \Z\times \Z_{\geq 0}$ such that
\[Z_{X(2),\omega(2)}(T)\in \mathcal{M}_k^{\hat{\mu}}\left[T, \frac{1}{1-\L^{a}T^b}\right]_{(a,b)\in S^{(2)}},\] 
and such that for each $(a,b)\in S^{(2)}$, we have that $\exp(2\pi i a/b)$ is an eigenvalue of $\sigma^2$ on $H^m(X\times_K K^{alg}, \Q_\ell)$, for some $m\geq 0$ and for every embedding of $\Q_\ell$ into $\C$. 
From \eqref{eq:future-zeta-extension}, it follows that
\[Z_{X,\omega}(T)\in \mathcal{M}_k^{\hat{\mu}}\left[T, \frac{1}{1-\L^{a}T^b}\right]_{(a,b)\in S},\] 
where $S=\left\{(a,2b)\in \Z\times \Z_{\geq 0}\mid (a,b)\in S^{(2)}\right\}$.

Take $(a,2b)\in S$, we have to prove that $\alpha=\exp(2\pi i a/(2b))$ is a monodromy eigenvalue. By definition of $S$, we know that $(a,b)\in S^{(2)}$, and therefore $\alpha^2=\exp(2\pi i a/b)$ is an eigenvalue of $\sigma^2$ on $H^m(X\times_K K^{alg}, \Q_\ell)$, for some $m\geq 0$. This implies that $\alpha$ or $-\alpha$ is an eigenvalue of $\sigma$ on $H^m(X\times_K K^{alg}, \Q_\ell)$. 

Suppose $-\alpha$ is an monodromy eigenvalue. 
By Lemma~\ref{thm:holds-form-zeta-function-K3}, we know that $-\alpha$ is a pole of the monodromy zeta function~$\zeta_X(T)$.
Let $\X$ be an $snc$-model of $X$, and write $\X_k=\sum_{i\in I} N_i E_i$. Because the index $\iota(X)=2$, we have that $N_i$ is even, for every $i\in I$. As a consequence, A'Campo's formula~\ref{thm:GMP-ACampo-dvr} gives that $\zeta_X(T)$ is a rational function in $T^2$. Therefore, $\alpha$ is a pole of $\zeta_X(T)$ as well, and hence it is a monodromy eigenvalue. We conclude that $X$ satisfies the monodromy property.
\end{proof}

\subsection{Number of \texorpdfstring{$(-2)$-curves}{(-2)-curves} on a rational, ruled surface}

In Chapter~\ref{ch:GMP-holds}, we proved that $K3$ surfaces allowing a triple-point-free model satisfy the monodromy property, except in the following case: if the Crauder-Morrison model is a chain degeneration where the chain has at least one rational, ruled end component $V_0$.
The main difficulty to prove the monodromy property in this case, is that the Euler characteristic $\chi(V_0^\circ)$ may be negative. Therefore, the factor $T^{N(V_0)}-1$ may cause cancellations in the A'Campo formula of Theorem~\ref{thm:GMP-ACampo-dvr}. Negative Euler characteristics do not occur for other surfaces in the chain, as can be seen from Lemma~\ref{thm:CM-euler-characteristics}.

Let us explore this problem of negative Euler characteristics in more detail. Let $X$ be a $K3$ surface over $K$ with a Crauder-Morrison model $\X$ such that $\X_k$ is a chain degeneration. Let $V_0$ be an end component of the chain that is a rational, ruled surface. Let $D$ be the elliptic double curve on $V_0$, where $V_0$ meets the next component in the chain. 

Let $C_1, \ldots, C_\lambda$ be the flowercurves on $V_0$. By Lemma~\ref{thm:CM-chain-flowers} and Lemma~\ref{thm:CM-self-intersection-flower-curve}, we know that $C_i$ is a rational curve with $C_i^2=-2$, for every $i=1, \ldots, \lambda$. Moreover, they are all disjoint, and don't meet $D$.

Computing the Euler characteristic $\chi(V_0^\circ)$ gives
\[\chi(V_0^\circ) = \chi(V_0) - \chi(D) - \sum_{i=1}^\lambda \chi(C_i) = \chi(V_0) - 2\lambda,\]
since $D$ is elliptic, and all $C_i$ are rational. 

Denote by $L$ the number of blow-ups in the contraction $\phi\colon V_0\to \overline{V}_0$ to the minimal ruled surface $\overline{V}_0$. 
Then we have that $\chi(V_0)=4+L$, and hence
\[\chi(V_0^\circ) = 4 + L - 2\lambda.\]
Moreover, from the proof of Theorem~\ref{thm:CM-chain-structure-fibers}, we know that $\phi$ consists of contractions of rational curves with self-intersection $-1$ meeting $D$ (or the image of $D$).

Therefore, the following question arises: Let $\overline{V}_0$ be a rational, minimal ruled surface with elliptic curve $\overline{D}\equiv - K_{\overline{V}_0}$, and let $V_0$ be the surface obtained by blowing up $\overline{V}_0$ in $L$ points on $\overline{D}$ or the intermediate strict transforms of $\overline{D}$. Let $D\subset V_0$ be the strict transform of $\overline{D}$. 
\emph{What is the maximal number $\lambda$ of pairwise disjoint, rational curves on $V_0$ with self-intersection $-2$ that are disjoint from~$D$?}

If $L\geq 2\lambda-4$, then $\chi(V_0^\circ)\geq 0$, but we cannot hope for this inequality in general, as Example~\ref{ex:future-Wim} shows. This example has $L=\lambda=6$, and hence $\chi(V_0^\circ)=-2$.
However, the following lemma gives a bound on $\lambda$ anyway.

\begin{lemma} \label{thm:future-relation-L-flower}
Let $V$ be a rational, ruled surface such that there is an elliptic curve $D$ with $D\equiv -K_V$. Let $L$ be the number of blow-ups in the contraction $\phi\colon V\to \overline{V}$ to the minimal ruled surface $\overline{V}$. Assume that $\phi$ is a composition of blow-ups with center on $D$, or strict transforms of $D$. Suppose there exist  pairwise disjoint, rational curves $C_1, \ldots, C_\lambda$ with $C_i^2=-2$ and $C_i\cdot D=0$, for every $i=1, \ldots, \lambda$. Then 
\[\lambda\leq L+1.\]

\end{lemma}

\begin{proof}
We know by \cite[Proposition V.2.3 and Exercise II.8.5]{Hartshorne} that
\[\Pic(V)\simeq \Z^{L+2}.\]

We will show that $[C_1], \ldots, [C_\lambda]$ and $[D]$ are linearly independent classes in $\Pic(V)\otimes_\Z \Q$, where $[\cdot]$ denotes the class in $\Pic(V)\otimes_\Z \Q$. From this, it follows that $\lambda\leq L+1$.

We start by proving that $[C_1],\ldots,[C_\lambda]$ are linearly independent in the vector space $\Pic(V)\otimes_\Z \Q$.  Let $a_i\in \Q$ such that
$\sum_{i=1}^\lambda a_i[C_i] = 0$.
Then for every $j=1,\ldots, \lambda$, we have
\[-2a_j=\sum_{i=1}^\lambda a_i C_i\cdot C_j =  0,\]
because $C_j^2=-2$ and $C_i \cdot C_j=0$ for every $i\neq j$. Hence $a_j=0$ for every $j$ and therefore, we have that $[C_1],\ldots,[C_\lambda]$ are linearly independent classes in $\Pic(V)\otimes_\Z \Q$.

We now prove that $[D]$ is linearly independent from $[C_1],\ldots, [C_\lambda]$ in $\Pic(V)\otimes_\Z \Q$. Suppose the contrary, then we can write
\[[D]=\sum_{i=1}^\lambda a_i [C_i],\]
for some $a_i\in \Q$. As linear equivalence implies numerical equivalence, we must have
\[D\cdot C_j = \sum_{i=1}^\lambda a_i C_i\cdot C_j,\]
for every $j=1,\ldots, \lambda$. Because $D\cdot C_j=0$ and $\sum_{i=1}^\lambda a_i C_i\cdot C_j=-2a_j$, we have $a_j=0$ for all $j=1, \ldots, \lambda$. Therefore $[D]=0$, which is a contradiction, since $D\equiv -K_V$.
\end{proof}

We don't know whether the inequality in Lemma~\ref{thm:future-relation-L-flower} is sharp.
In any case, it cannot be improved a lot because of the following example, which was shown to us by W. Veys.

\begin{example}[due to W. Veys] \label{ex:future-Wim}
We will construct a rational, ruled surface~$V$ with the following properties: there is an elliptic curve $D$ with $D\equiv -K_V$, and there are six pairwise disjoint, smooth, rational curves $C_1, \ldots, C_6$ with $C_i^2=-2$ and $C_i\cdot D=0$, for every $i=1, \ldots, 6$. The contraction $\phi\colon V\to \overline{V}$  to the minimal ruled surface $\overline{V}$ is a composition of six contractions of rational $(-1)$-curves meeting $D$. 

We will start from $\overline{V}\simeq \Sigma_1$, the first Hirzebruch surface. Then we will define six blow-ups
\[V=V^{(6)}\to V^{(5)}\to V^{(4)}\to V^{(3)}\to V^{(2)}\to V^{(1)}\to V^{(0)}=\overline{V}.\]
The strict transform of any curve $C^{(i)} \subset V^{(i)}$ will be denoted by $C^{(i+1)}\subset V^{(i+1)}$. For any point $P^{(i)}\in V^{(i)}$ that is not the center of the blow-up $V^{(i+1)}\to V^{(i)}$, we denote its inverse image by $P^{(i+1)}\in V^{(i+1)}$.

Since $\overline{V}=\Sigma_1$ is the blow-up of $\P^2$ in one point, we will define some relevant curves on $\overline{V}$ by defining them on $\P^2$.
So let $V^{(-1)}\simeq \P^2$ and consider three lines $\ell, C_2^{(-1)}, C_3^{(-1)}$ in $V^{(-1)}$, such that $\ell\cap C_2^{(-1)}\cap C_3^{(-1)}= \{P_1^{(-1)}\}$ for some point $P_1^{(-1)}\in V^{(-1)}$. Let $D^{(-1)}$ be a cubic through $P_1^{(-1)}$, tangent to the lines $\ell, C_2^{(-1)}, C_3^{(-1)}$, in the points $P_\ell, P_2^{(-1)}, P_3^{(-1)}$ respectively. None of the points $P_\ell, P_2^{(-1)}, P_3^{(-1)}$ equals $P_1^{(-1)}$. It is not difficult to check that such a cubic indeed exists. A sketch of this configuration can be found in Figure~\ref{fig:future-example-wim-V1}, where the dotted line denotes $D^{(-1)}$.

\begin{figure}[H]
\begin{center}
\begin{tikzpicture}[scale=0.6]
\draw (0,-1) -- (0,4.5);
\draw (-4,1) -- (4,-1);
\draw (-4,-1) -- (4,1);

\node [above] at (0, 4.5) {$\ell$};
\node [right] at (4,1) {$C_2^{(-1)}$};
\node [right] at (4,-1) {$C_3^{(-1)}$};

\node [draw,shape=circle, fill, inner sep=0pt,minimum size=5pt] at (0, 3) {};
\node [below left] at (0,3) {$P_\ell$};

\node [draw,shape=circle, fill, inner sep=0pt,minimum size=5pt] at (-3, 0.75) {};
\node [below] at (-3,0.75) {$P_3^{(-1)}$};

\node [draw,shape=circle, fill, inner sep=0pt,minimum size=5pt] at (3, 0.75) {};
\node [below] at (3,0.75) {$P_2^{(-1)}$};

\node [draw,shape=circle, fill, inner sep=0pt,minimum size=5pt] at (0, 0) {};
\node [below left] at (0,0) {$P_1^{(-1)}$};

\draw[dotted, thick] plot [smooth] coordinates{(0.5, 4) (0.25, 3.75) (0,3) (0.25, 2.25) (0.5, 2)};

\draw[dotted, thick] (4,2) parabola bend (2.875,47/64) (2,1.5);

\draw[dotted, thick] (-4,2) parabola bend (-2.875,47/64) (-2,1.5);

\draw[dotted, thick] plot [smooth] coordinates{(-1, 1) (-0.5, 0.75)  (0,0) (0.5, -0.75) (1,-1)};
\end{tikzpicture}
\end{center}
\caption{Sketch of $V^{(-1)}$} \label{fig:future-example-wim-V1}
\end{figure}
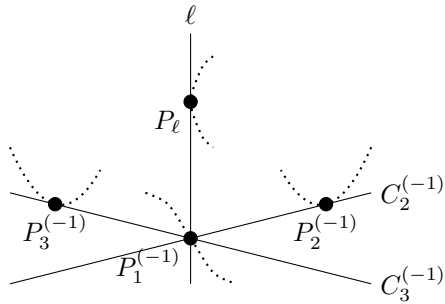

Let $\overline{V} = V^{(0)}$ be the blow-up of $V^{(-1)}$ in $P_\ell$. So we have $\overline{V}\simeq\Sigma_1$, the first Hirzebruch surface. Let $C_1^{(0)}$ be the exceptional divisor. Define the curves $D^{(0)}, \fiber_1^{(0)}, C_2^{(0)}$ and $C_3^{(0)}$ as the strict transforms of $D^{(-1)}, \ell, C_2^{(-1)}$ and $C_3^{(-1)}$ respectively.
Define the points $P_1^{(0)}, P_2^{(0)}$ and $P_3^{(0)}$ as the inverse images of $P_1^{(-1)}, P_2^{(-1)}$ and $P_3^{(-1)}$ respectively under the  morphism $V^{(0)}\to V^{(-1)}$. Define $P^{(0)}_4 = C_1^{(0)}\cap \fiber_1^{(0)}$. The curve $D^{(-1)}$ meets $\fiber_1^{(0)}$ transversally in $P^{(0)}_1$ and $P^{(0)}_4$, and it is tangent to  $C_2^{(0)}$ in $P^{(0)}_2$, and to $C_3^{(0)}$ in $P^{(0)}_3$.

Some important lines and points on $V^{(0)}$ are sketched in Figure~\ref{fig:future-example-wim-Vbar}, where the dashed lines are fibers of the ruling on ${V}_0$, and the solid lines are horizontal curves.

\begin{figure}[H]
\begin{center}
\begin{tikzpicture}[scale=0.6]
\draw (-4,4) -- (4,4);
\draw (-4,1) -- (4,-1);
\draw (-4, -1) -- (4,1);

\draw[dashed] (0, -1) -- (0,4.5);
\draw[dashed] (-3, -1) -- (-3,4.5);
\draw[dashed] (3, -1) -- (3,4.5);

\node [above] at (0, 4.5) {$\fiber_1^{(0)}$};
\node [right] at (4,4) {$C_1^{(0)}$};
\node [right] at (4,1) {$C_2^{(0)}$};
\node [right] at (4,-1) {$C_3^{(0)}$};

\node [draw,shape=circle, fill, inner sep=0pt,minimum size=5pt] at (-3, 0.75) {};
\node [below left] at (-3,0.75) {$P_3^{(0)}$};

\node [draw,shape=circle, fill, inner sep=0pt,minimum size=5pt] at (3, 0.75) {};
\node [below left] at (3,0.75) {$P_2^{(0)}$};

\node [draw,shape=circle, fill, inner sep=0pt,minimum size=5pt] at (0, 4) {};
\node [below left] at (0,4) {$P_4^{(0)}$};

\node [draw,shape=circle, fill, inner sep=0pt,minimum size=5pt] at (0, 0) {};
\node [below left] at (0,0) {$P_1^{(0)}$};
\end{tikzpicture}
\end{center}
\caption{Sketch of $V^{(0)}$}
\label{fig:future-example-wim-Vbar}
\end{figure}
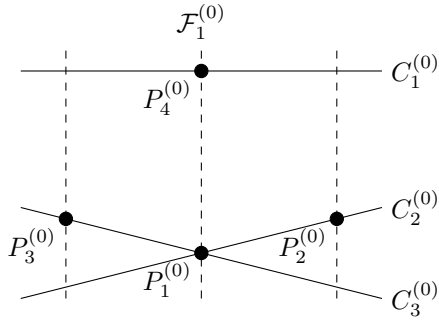

Because $C_1^{(0)}$ is the exceptional curve, and because $\left(C_2^{(-1)}\right)^2=\left(C_3^{(-1)}\right)^2=1$, we have that 
\[\left(C_1^{(0)}\right)^2 = -1 \text{ and } \left(C_2^{(0)}\right)^2 = \left(C_3^{(0)}\right)^2 = 1.\]

Let $V^{(2)}\to V^{(0)}$ be the birational morphism obtained by blowing up $V^{(0)}$ in the points $P_1^{(0)}$ and $P_4^{(0)}$. Some important lines and points on $V^{(2)}$ are sketched in Figure~\ref{fig:future-example-wim-V2}.

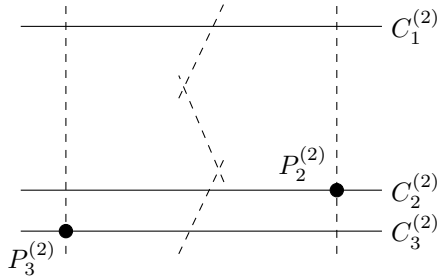
\begin{figure}[H]
\begin{center}
\begin{tikzpicture}[scale=0.6]
\draw (-4,4) -- (4,4);
\draw (-4,0.4) -- (4,0.4);
\draw (-4, -0.5) -- (4,-0.5);

\draw[dashed] 	(-0.5,-1) -- (0.5, 1.08)
				(0.5, 0.58) -- (-0.5, 2.92)
				(-0.5, 2.42) -- (0.5,4.5)							
				;

\draw[dashed] (-3, -1) -- (-3,4.5);
\draw[dashed] (3, -1) -- (3,4.5);

\node [right] at (4,4) {$C_1^{(2)}$};
\node [right] at (4,0.4) {$C_2^{(2)}$};
\node [right] at (4,-0.5) {$C_3^{(2)}$};

\node [draw,shape=circle, fill, inner sep=0pt,minimum size=5pt] at (3, 0.4) {};
\node [above left] at (3,0.4) {$P_2^{(2)}$};

\node [draw,shape=circle, fill, inner sep=0pt,minimum size=5pt] at (-3, -0.5) {};
\node [below left] at (-3,-0.5) {$P_3^{(2)}$};

\end{tikzpicture}
\end{center}
\caption{Sketch of $V^{(2)}$}
\label{fig:future-example-wim-V2}
\end{figure}

The curves $C_2^{(2)}$ and $C_3^{(2)}$ are the strict transforms of $C_2^{(0)}$ and $C_3^{(0)}$ respectively. We have that 
\[\left(C_1^{(2)}\right)^2 = -2 \text{ and } \left(C_2^{(2)}\right)^2 =\left(C_3^{(2)}\right)^2 = 0.\]
Furthermore, $C_1^{(2)}, C_2^{(2)}$ and $C_3^{(2)}$ are disjoint.
Moreover, the curve $D^{(2)}$ is tangent to $C_2^{(2)}$ in $P_2^{(2)}$, and to $C_3^{(2)}$ in $P_3^{(2)}$.

There is one reducible fiber, with three irreducible components as seen in Figure~\ref{fig:future-example-wim-V2}. The middle component has self-intersection $-2$ and the other two components have self-intersection $-1$. Furthermore, the component with self-intersection $-2$ is disjoint from $C_1^{(2)}, C_2^{(2)}$ and $C_3^{(2)}$.

Because we blew up smooth points on $D^{(0)}$, we have that $ D^{(2)}\equiv -K_{V^{(2)}} $ by \cite[Proposition~V.3.3 and Proposition~V.3.6]{Hartshorne}.

Let $V^{(3)}\to V^{(2)}$ be the birational morphism obtained by blowing up $V^{(2)}$ in the point $P_3^{(2)}$. Some important lines and points on $V^{(3)}$ are sketched in Figure~\ref{fig:future-example-wim-V3}.

\begin{figure}[H]
\begin{center}
\begin{tikzpicture}[scale=0.6]
\draw (-4,4) -- (4,4);
\draw plot [smooth] coordinates{(-4,3) (-2, 2) (0.24, 0.4) (3,0.4) (4,0.4)};
\draw (-4, -0.5) -- (4,-0.5);

\draw[dashed] 	(-0.5,-1) -- (0.5, 1.08)
				(0.5, 0.58) -- (-0.5, 2.92)
				(-0.5, 2.42) -- (0.5,4.5)
				;

\draw[dashed]	(-2.5, -1) -- (-3.5,2)
				(-3.5, 1.5) -- (-2.5, 4.5)
;
\draw[dashed] (3, -1) -- (3,4.5);

\node [right] at (4,4) {$C_1^{(3)}$};
\node [right] at (4,0.4) {$C_2^{(3)}$};
\node [right] at (4,-0.5) {$C_3^{(3)}$};

\node [draw,shape=circle, fill, inner sep=0pt,minimum size=5pt] at (3, 0.4) {};
\node [above left] at (3,0.4) {$P_2^{(3)}$};

\node [draw,shape=circle, fill, inner sep=0pt,minimum size=5pt] at (-2.66, -0.5) {};
\node [below left] at (-2.66,-0.5) {$P_3^{(3)}$};

\end{tikzpicture}
\end{center}
\caption{Sketch of $V^{(3)}$}
\label{fig:future-example-wim-V3}
\end{figure}
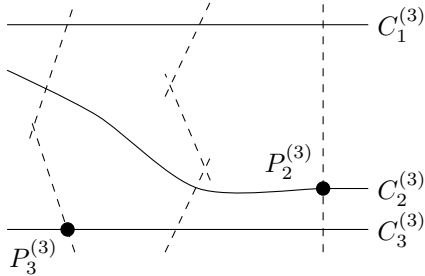

We have that $D^{(3)}$ intersects $C_3^{(3)}$ transversally, and we denote this point by~$P_3^{(3)}$.

Let $V^{(4)}\to V^{(3)}$ be the birational morphism obtained by blowing up $V^{(3)}$ in the point $P_3^{(3)}$. Some important lines and points on $V^{(4)}$ are sketched in Figure~\ref{fig:future-example-wim-V4}.

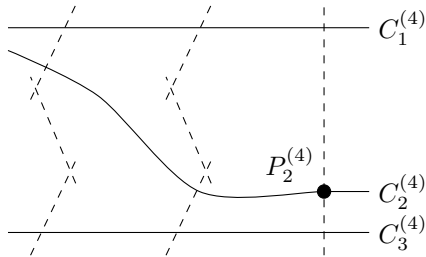
\begin{figure}[H]
\begin{center}
\begin{tikzpicture}[scale=0.6]
\draw (-4,4) -- (4,4);
\draw plot [smooth] coordinates{(-4,3.5) (-2, 2.5) (0.24, 0.4) (3,0.4) (4,0.4)};
\draw (-4, -0.5) -- (4,-0.5);

\draw[dashed] 	(-0.5,-1) -- (0.5, 1.08)
				(0.5, 0.58) -- (-0.5, 2.92)
				(-0.5, 2.42) -- (0.5,4.5)
				;

\draw[dashed] 	(-3.5,-1) -- (-2.5, 1.08)
				(-2.5, 0.58) -- (-3.5, 2.92)
				(-3.5, 2.42) -- (-2.5,4.5)
				;
\draw[dashed] (3, -1) -- (3,4.5);

\node [right] at (4,4) {$C_1^{(4)}$};
\node [right] at (4,0.4) {$C_2^{(4)}$};
\node [right] at (4,-0.5) {$C_3^{(4)}$};

\node [draw,shape=circle, fill, inner sep=0pt,minimum size=5pt] at (3, 0.4) {};
\node [above left] at (3,0.4) {$P_2^{(4)}$};

\end{tikzpicture}
\end{center}
\caption{Sketch of $V^{(4)}$}
\label{fig:future-example-wim-V4}
\end{figure}

The curves $C_2^{(4)}$ and $C_3^{(4)}$ are the strict transforms of $C_2^{(2)}$ and $C_3^{(2)}$ respectively. We have that 
\[\left(C_1^{(4)}\right)^2 =\left(C_3^{(4)}\right)^2= -2 \text{ and } \left(C_2^{(4)}\right)^2  = 0.\]
Furthermore, $C_1^{(4)}, C_2^{(4)}$ and $C_3^{(4)}$ are disjoint.
Moreover, the curve $D^{(4)}$ is tangent to $C_2^{(4)}$ in $P_2^{(4)}$.

There are two reducible fibers, each with three irreducible components as seen in Figure~\ref{fig:future-example-wim-V4}. The middle component of each of these fibers has self-intersection~$-2$, and the other two components have self-intersection $-1$. Furthermore, the two components with self-intersection $-2$ are disjoint from $C_1^{(4)}, C_2^{(4)}$ and $C_3^{(4)}$.

Because we blew up smooth points on $D^{(2)}$ and $D^{(3)}$, we have that $ D^{(4)}\equiv -K_{V^{(4)}} $ by \cite[Proposition~V.3.3 and Proposition~V.3.6]{Hartshorne}.

Let $V^{(5)}\to V^{(4)}$ be the birational morphism obtained by blowing up $V^{(4)}$ in the point $P_2^{(4)}$. We have that $D^{(5)}$ intersects $C_2^{(5)}$ transversally, and we denote this point by $P_2^{(5)}$.
Let $V=V^{(6)}\to V^{(5)}$ be the birational morphism obtained by blowing up $V^{(5)}$ in the point $P_2^{(5)}$.

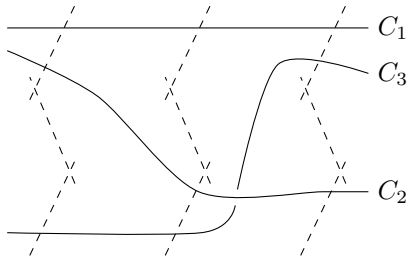
\begin{figure}[H]
\begin{center}
\begin{tikzpicture}[scale=0.6]
\draw (-4,4) -- (4,4);

\draw plot [smooth] coordinates{(-4,3.5) (-2, 2.5) (0.24, 0.4) (3,0.4) (4,0.4)};

\draw plot [smooth] coordinates {(-4,-0.5) (0.3, -0.5) (1.05, 0.1)}; 
\draw plot [smooth] coordinates {(1.1, 0.45) (2,3.2) (4, 3)};

\draw[dashed] 	(-0.5,-1) -- (0.5, 1.08)
				(0.5, 0.58) -- (-0.5, 2.92)
				(-0.5, 2.42) -- (0.5,4.5)
				;

\draw[dashed] 	(-3.5,-1) -- (-2.5, 1.08)
				(-2.5, 0.58) -- (-3.5, 2.92)
				(-3.5, 2.42) -- (-2.5,4.5)
				;
				
\draw[dashed] 	(2.5,-1) -- (3.5, 1.08)
				(3.5, 0.58) -- (2.5, 2.92)
				(2.5, 2.42) -- (3.5,4.5)
				;

\node [right] at (4,4) {$C_1$};
\node [right] at (4,0.4) {$C_2$};
\node [right] at (4,3) {$C_3$};

\end{tikzpicture}
\end{center}
\caption{Sketch of $V$}
\label{fig:future-example-wim-V6}
\end{figure}

The curves $C_2$ and $C_3$ are the strict transforms of $C_2^{(4)}$ and $C_3^{(4)}$ respectively. We have that 
\[C_1^2 =C_2^2=C_3^2= -2  .\]
Furthermore, $C_1, C_2$ and $C_3$ are disjoint.

There are three reducible fibers, each with three irreducible components. In each reducible fiber, the middle component has self-intersection $-2$, and the other two components have self-intersection $-1$. Furthermore, these vertical curves with self-intersection $-2$ are pairwise disjoint, and disjoint from $C_1, C_2$ and $C_3$.

Because we blew up smooth points on $D^{(4)}$ and $D^{(5)}$, we have that $ D\equiv -K_{V} $ by \cite[Proposition~V.3.3 and Proposition~V.3.6]{Hartshorne}.

To conclude, we have constructed a rational, ruled surface $V$ such that there is an elliptic curve $D$ with $D\equiv -K_V$. There are six pairwise disjoint, smooth, rational curves $C_1, \ldots, C_6$ with $C_i^2=-2$ and $C_i\cdot D=0$, for every $i=1, \ldots, 6$. The contraction $\phi\colon V\to \overline{V}$  to the minimal ruled surface $\overline{V}$ is a composition of six blow-ups with center on $D$.
\end{example}

\subsection{Number of blow-ups in the chain}

Let $X$ be a $K3$ surface over $K$ with Crauder-Morrison model $\X$. Suppose $\X$ is a chain degeneration, and let $V_0, V_1, \ldots, V_k, V_{k+1}$ be the components in the chain, such that $V_i\cap V_j=\emptyset$ if and only if $j\not \in \{i-1,i,i+1\}$. For every $i=0, \ldots, k+1$, if $V_i$ is a ruled surface, then we define $L_i$ to be the number of blow-ups in the contraction $V_i\to \overline{V}_i$ to the minimal ruled surface $\overline{V}_i$. If $V_i\simeq \P^2$, then we set $L_i=0$ and $\overline{V_i}=V_i$.
The following formula restricts the possibilities of the values~$L_i$.

\begin{lemma} \label{thm:future-relation-L}
Let $X$ be a $K3$ surface over $K$ with Crauder-Morrison model~$\X$ and with special fiber $\X_k=\sum_{i\in I} N_iE_i$.
Suppose $\X$ is a chain degeneration, and let $V_0, V_1, \ldots, V_k, V_{k+1}$ be the components in the chain, such that $V_i\cap V_j=\emptyset$ if and only if $j\not \in \{i-1,i,i+1\}$. 
We have
\[\sum_{i=0}^{k+1}  \frac{L_i}{N(V_i)} = \frac{K_{\overline{V}_0}^2}{N(V_0)}+ \frac{K_{\overline{V}_{k+1}}^2}{N(V_{k+1})}.\]
\end{lemma}

\begin{proof}
For $i=0,\ldots, k$, set $C_i=V_i\cap V_{i+1}$ . We have that 
\begin{align*}
K_{V_i} & \equiv -C_{i-1}-C_i,
\end{align*}
for every $i=1, \ldots, k$, because of \eqref{eq:canonical-bundle-component-general} and the fact that $V_i$ only meets flowers of type 2 by Proposition~\ref{thm:CM-chain-flowers}.
On the other hand, since ${V_i}$ is elliptic ruled for $i=1,\ldots, k$, we know that
\[K_{V_i}^2=-L_i,\]
by \cite[Corollary V.2.11 and Proposition V.3.3]{Hartshorne}.
As there are no triple points, this gives
\begin{equation}
(C_{i-1})^2_{V_i}+(C_i)^2_{V_i}=-L_i. \label{eq:C^2OnSameSurf}
\end{equation}

Using techniques similar as in the proof of the triple point formula \cite[Corollary~2.4.2]{Persson}, one can prove that
\begin{equation}
(C_i)^2_{V_i}=-\frac{N(V_i)}{N(V_{i+1})}(C_i)^2_{V_{i+1}}, \label{eq:C^2OnDifferentSurf}
\end{equation}
for $i=0,\ldots, k$.
\note{
$(C_i)^2_{V_i}=V_j\cdot V_j \cdot V_i = -\sum_{k\neq j } N_k/N_j V_k\cdot V_j\cdot V_i = -N_i/N_j V_i\cdot V_i\cdot V_j=-\frac{N_i}{N_{i+1}}(C_i)^2_{V_{i+1}}$
}

Combining formulas \eqref{eq:C^2OnSameSurf} and \eqref{eq:C^2OnDifferentSurf}, we obtain
\begin{align*}
(C_0)^2_{V_0} 	&= -\frac{N(V_0)}{N(V_1)}(C_0)^2_{V_1} \\
				&= \frac{N(V_0)}{N(V_1)}\left((C_1)^2_{V_1}+L_1\right)\\
				&= -\frac{N(V_0)}{N(V_2)}(C_1)^2_{V_2}+ \frac{N(V_0)}{N(V_1)} L_1 \\
				&= \frac{N(V_0)}{N(V_2)}(C_2)^2_{V_2}+\sum_{i=1}^2\frac{N(V_0)}{N(V_i)} L_i \\
				&= \cdots \\
				&= -\frac{N(V_0)}{N(V_{k+1})}(C_{k})^2_{V_{k+1}}+ \sum_{i=1}^{k} \frac{N(V_0)}{N(V_i)}L_i.
\end{align*}
On the other hand
\[(C_0)^2_{V_0}		= K_{V_0}^2 = K_{\overline{V}_0}^2-L_0,\]
and
\[(C_k)^2_{V_{k+1}} 	=K_{V_{k+1}}^2 =K_{\overline{V}_{k+1}}^2-L_{k+1},\]
by \cite[Proposition~V.3.3]{Hartshorne}.

Therefore, 
\[K_{\overline{V}_0}^2-L_0 = -\frac{N(V_0)}{N(V_{k+1})}\left(K_{\overline{V}_{k+1}}^2-L_{k+1}\right)+ \sum_{i=1}^{k} \frac{N(V_0)}{N(V_i)}L_i,\]
which is equivalent to
\[\sum_{i=0}^{k+1}  \frac{L_i}{N(V_i)} = \frac{K_{\overline{V}_0}^2}{N(V_0)}+ \frac{K_{\overline{V}_{k+1}}^2}{N(V_{k+1})}.\]
\end{proof}

\begin{remark}
The equation in the previous lemma can be simplified if one assumes that $\X$ satisfies the properties of Proposition~\ref{thm:CM-chain-good}. If we define $\alpha$ and $\beta$ as in that proposition, and if we denote $N=N(V_\alpha)$, then we have the equality
\[\frac{L_0}{N(V_0)} +  \frac{\sum_{i=\alpha}^{\beta} L_i}{N} + \frac{L_{k+1}}{N(V_{k+1})} = \frac{K_{\overline{V}_0}^2}{N(V_0)}+ \frac{K_{\overline{V}_{k+1}}^2}{N(V_{k+1})}.\]
\end{remark}

\section{Combinatorial countercandidates}
\label{sect:future-comb-count}

In this section, we will discuss combinatorial countercandidates and a strategy to produce them all. 
We will start by explaining the concept of combinatorial countercandidates. In Subsection~\ref{sect:future-count-ex}, we produce all combinatorial countercandidates for a specific kind of chain degenerations. Then we will explain how to generalize this strategy to all chain degenerations in Subsection~\ref{sect:future-count-strategy}. In the final subsection, we will use this strategy to prove that certain $K3$ surfaces with a chain degeneration satisfy the monodromy property.

\subsection{Combinatorial countercandidates}
Although the concept of a combinatorial countercandidate is not mathematically precise, we will explain this notion in more detail. To do so, we first need the definition of a strict normal crossings surface.

\begin{definition}
A $k$-surface $D$ is called a \emph{strict normal crossings surface}, if for every closed point $p\in D$, we have $\widehat{\O}_{D,p} \simeq k\llbracket x,y,z \rrbracket /(x^ay^bz^c)$ for some non-negative integers $a,b,c\in \Z_{\geq 0}$ and $a\neq 0$, and if moreover for every irreducible component $E$ of $D$, the reduction $E_{red}$ is smooth.

If for every closed point $p\in D$, we have $\widehat{\O}_{D,p} \simeq k\llbracket x,y,z \rrbracket /(x^ay^b)$ for some integers $a,b\in \Z_{\geq 0}$ with $a\neq 0$, then we say $D$ is \emph{triple-point-free}.

We write $D=\sum_{i\in I} N_i E_i$, where $E_i$ are the reduced irreducible components of $D$ and where $N_i$ is the multiplicity of $E_i$.
\end{definition}

 We define the \emph{index} of a strict normal crossings surface $D$ as $\iota(D) = \gcd_{i\in I} N_i$. Note that, when $X$ is a $K$-variety with an $snc$-model $\X$, then $\iota(X)=\iota(\X_k)$, by Proposition~\ref{thm:future-index-gcd}. We also define the \emph{dual graph} of a triple-point-free, strict normal crossings surface in the obvious way.

A \emph{combinatorial countercandidate} is a combinatorial description of a triple-point-free, strict normal crossings surface (for example with a given dual graph, multiplicities of the components, some geometrical properties of the components,~\ldots) such that 
\begin{itemize}
\item \emph{if} there exists a strict normal crossings surface $D$ with the described properties, and \emph{if} there exists a $K3$ surface $X$ admitting an $snc$-model $\X$ such that $\X_k\simeq D$ (i.e., if $D$ can be smoothened to a $K3$ surface $X$), then $X$ does not satisfy the monodromy property,
\item we don't know whether there exists a $K3$ surface $X$ over $K$ with an $snc$-model $\X$ such that the special fiber $\X_k$ satisfies the combinatorial characterization, i.e., there is no reason (yet) to rule it out.
\end{itemize}
We think of a combinatorial countercandidate as a description of a candidate to produce a $K3$ surface not satisfying the monodromy property. 
As soon as we find a reason why this candidate will not produce a counterexample of the monodromy property, it is obviously not a candidate anymore. This explains the non-precise nature of the second requirement.

If there exists a $K3$ surface over $K$ allowing a triple-point-free model that does not satisfy the monodromy property, then it has a Crauder-Morrison model. Therefore, we are mainly interested in combinatorial countercandidates that satisfy the properties of the Crauder-Morrison classification. For this reason, we can use the terminology from this classification in the setting of combinatorial countercandidates as well. We will for example talk about chains, flowerpots and flowers in this context. Moreover, from Theorem~\ref{thm:holds-pot}, it is clear that combinatorial countercandidates are chain degenerations.


Note that by Lemma~\ref{thm:future-index}, we are mainly interested in combinatorial countercandidates $D$ of index 1. 
Indeed, if we can exclude all combinatorial countercandidates of index 1, then we have proven that the monodromy property holds for all $K3$ surface with a triple-point-free model, also for those of index 2.

From the results in this section and in Appendix~\ref{app:combinatorial-countercandidate}, we will see that we have a list of 63 combinatorial countercandidates. If there exists a $K3$ surface $X$ allowing a triple-point-free model, not satisfying the monodromy property, then the special fiber of its Crauder-Morrison model has the properties of one of these 63 combinatorial countercandidates.

\subsection{Constructing combinatorial countercandidates: example}\label{sect:future-count-ex}

Suppose there exists a $K3$ surface $X$ over $K$ of index $\iota(X)=1$ allowing a triple-point-free model that does not satisfy the monodromy property. Let $\X$ be a Crauder-Morrison model of $X$. By Theorem~\ref{thm:holds-pot}, we have that $\X$ is a chain degeneration. Let $V_0, V_1, \ldots, V_{k}, V_{k+1}$ be the components in the chain, such that $V_i\cap V_j=\emptyset$ if and only if $j\not \in \{i-1,i,i+1\}$. Because of Theorem~\ref{thm:holds-chain}~(ii), we can assume that $V_0$ is a rational, ruled surface.

In this subsection, we will construct combinatorial countercandidates where the chain has the following properties: $V_{k+1}$ is a rational, ruled surface and there exists an integer $N\geq 1$ such that
\[N(V_i) = \begin{cases}
2N & \text{if } i=0,\\
N&\text{otherwise.}
\end{cases}\] 

Let $\phi$ be the total number of flowers of type $2B$ and $2C$ meeting $V_1, \ldots, V_{k+1}$, and let $\gamma$ be the total number of flowers of type $2A$ meeting $V_1, \ldots, V_{k+1}$. We have that $\phi\geq 1$,
by Theorem~\ref{thm:holds-chain} (iv).
 Moreover, let $\phi'$ be the total number of flowers of type $2B$ and $2C$ meeting $V_0$, and let $\gamma'$ be the total number of flowers of type $2A$ meeting $V_0$.

Define $I^\dagger\subset I$ to be the set of indices $i\in I$, where either
$\rho_i$ is minimal, or
$E_i$ is the top of a conic-flower.
Let $S^\dagger=\left\{(-\nu_i,N_i)\in \Z\times \Z_{>0}\mid i\in I^\dagger\right\}$. Corollary~\ref{thm:poles-contributions} gives that
\[Z_{X,\omega}(T)\in \mathcal{M}_k^{\hat{\mu}}\left[T, \frac{1}{1-\L^{a}T^b}\right]_{(a,b)\in S^\dagger}.\]
By Theorem~\ref{thm:holds-lct}, we know that $\xi(E_i)$ is a monodromy eigenvalue, if $\rho_i$ is minimal. Therefore, because $X$ does not satisfy the monodromy property, there must be a conic-flower $F$ with top $F_0$, such that $\xi(F_0)$ is \emph{not} a monodromy eigenvalue. By Proposition~\ref{thm:CM-chain-flowers}, the flower $F$ is of type $2B$ or $2C$.

Because of Lemma~\ref{thm:CM-euler-characteristics}, the monodromy zeta function is given by
\begin{equation}
\zeta_X(T) = \frac{(T^{2N} - 1)^{- \chi(V_0^{\circ})}}{(T^{4N} - 1)^{\phi'} (T^{2N} - 1)^{\phi}(T^{N} - 1)^{\gamma'+\sum_{i=1}^{k+1}\chi(V_i^{\circ})}(T^{N/2} - 1)^{\gamma} }. \label{eq:zetafunction'}
\end{equation}

Suppose $F$ meets $V_0$. Then $N(F_0)=4N$, and by Lemma~\ref{thm:holds-conic-flower-xi}, we have that $\xi(F_0)$ is not a $2N$-th root of unity. We see from \eqref{eq:zetafunction'} that $\xi(F_0)$ is a pole of $\zeta_X(T)$, and hence a monodromy eigenvalue, contrary to our assumption.

This means that $F$ meets one of the components $V_1, \ldots, V_{k+1}$ and $N(F_0)=2N$. Since $\xi(F_0)$ is not a monodromy eigenvalue, it is neither a zero nor a pole of $\zeta_X(T)$.
By Lemma~\ref{thm:holds-conic-flower-xi}, we know that $\xi(F_0)$ is a $2N$-th root of unity, but not an $N$-th root of unity, and hence the factor $(T^{2N}-1)$ must cancel completely in $\zeta_X(T)$. Therefore we must have
$-\chi(V_0^\circ) = \phi+\phi'$.
Since $\chi(V_0^{\circ}) = 4 + L_0 -2\phi' - 2\gamma'$, we get
\begin{equation}
\phi'+2\gamma'-\phi-L_0 = 4. \label{eq:ineq-phi2}
\end{equation}
Therefore, equation~\eqref{eq:zetafunction'} can be rewritten as
\begin{equation*}
\zeta_X(T) = \frac{1}{(T^{2N} + 1)^{\phi'}(T^{N} - 1)^{\gamma'+\sum_{i=1}^{k+1}\chi(V_i^{\circ})}  (T^{N/2} - 1)^{\gamma} }. 
\end{equation*}
Because $X$ has index $\iota(X)=1$, Proposition~\ref{thm:future-index-gcd} gives that $\gcd_{i\in I} N_i = 1$. This means that
\begin{equation*}
N = \begin{cases}
1 & \text{if } \gamma=0,\\
2 & \text{if } \gamma>0.\\
\end{cases} \label{eq:index'}
\end{equation*}

Proposition~\ref{thm:holds-form-zeta-function-K3} gives that
\[\sum_{i=1}^{k+1}\chi(V_i^{\circ}) + 2\phi'+\gamma'+ \frac{\gamma}{2}=\frac{24}{N}. \]
In particular, $\gamma$ is even.

Since $\sum_{i=1}^{k+1}\chi(V_i^\circ)=4+\sum_{i=1}^{k+1}L_i-2\phi-2\gamma$, we get
\begin{equation}
\sum_{i=1}^{k+1}L_i+2\phi'+\gamma' -2\phi - \frac{3\gamma}{2} = \frac{24}{N} - 4. \label{eq:degree-zeta-function}
\end{equation}

Lemma~\ref{thm:future-relation-L} implies that 
\begin{equation}
L_0+2 \sum_{i=1}^{k+1}L_i = 24. \label{eq:relation-L}
\end{equation}

Note that, since $V_0$ is a non-minimal ruled surface by Theorem~\ref{thm:holds-chain}, we have that $L_0\geq 1$.

The variables $\phi,\gamma, \phi', \gamma', N, L_0$ and $\sum_{i=1}^{k+1} L_i$ are non-negative integers, and therefore, we have obvious lower bounds for all of these variables. We will now also derive upper bounds. 
Equation \eqref{eq:relation-L} immediately gives upper bounds for $L_0$ and $\sum_{i=1}^{k+1}L_i$.

From Lemma~\ref{thm:future-relation-L-flower}, it follows that 
\begin{equation*}
\phi'+\gamma'\leq L_0+1. \label{eq:bound-flowers'}
\end{equation*}
Moreover, Lemma~\ref{thm:CM-relation-L-flower-elliptic} combined with Lemma~\ref{thm:future-relation-L-flower} gives that
\[\phi + \gamma \leq \Big(\sum_{i=1}^k L_i\Big)/2 + (L_{k+1}+1),\]
and therefore
\begin{equation*}
\phi+\gamma\leq \sum_{i=1}^{k+1} L_i+1. \label{eq:bound-flowers}
\end{equation*}

As a consequence, finding combinatorial countercandidates with this particular chain,
comes down to finding non-negative integer solutions of
\[\begin{cases}
\phi'+2\gamma'-\phi-L_0 = 4,\\
\sum_{i=1}^{k+1} L_i + 2\phi'+\gamma'-2\phi-\frac{3\gamma}{2}=\frac{24}{N}-4,\\
L_0 + 2\sum_{i=1}^{k+1} L_i =24,\\
\phi'+\gamma'\leq L_0+1,\\
\phi+\gamma \leq \sum_{i=1}^{k+1} L_i +1,\\
\phi\geq 1,\\
L_0\geq 1,\\
N = \begin{cases}
1 & \text{if } \gamma=0,\\
2 & \text{if } \gamma>0.
\end{cases}
\end{cases}
\]
Since all variables $N, \phi, \gamma, \phi', \gamma', L_0, \sum_{i=1}^{k+1}L_i$ are bounded, it is easy to let a computer check every possibility whether it satisfies these relations. We have implemented this in Python and the code and output can be found in section~\ref{sect:app-case1} of the appendix. The code can also be downloaded from \mbox{\url{www.github.com/AnneliesJaspers/combinatorial-counterexamples}}.

In this way, we find 68 combinatorial countercandidates (modulo the length $k+1$ of the chain), of which we will now give two as an illustration. In Section~\ref{sect:future-exclude}, we will exclude 26 of these combinatorial countercandidates, including Combinatorial Countercandidate~\ref{count:exclude}.

\begin{countercandidate} \label{count:exclude}
Let $D$ be a strict normal crossings surface with the following dual graph, where the labels denote the multiplicity of the components:

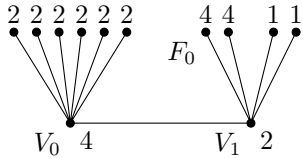
\begin{center}
\begin{tikzpicture}[scale=0.6]

\node [draw,shape=circle, fill, inner sep=0pt,minimum size=3pt] at (0,0) {};
\node [below left] at (0,0) {$V_0$};
\node [below right] at (0,0) {$4$};
\node [draw,shape=circle, fill, inner sep=0pt,minimum size=3pt] at (4,0) {};
\node [below left] at (4,0) {$V_1$};
\node [below right] at (4,0) {$2$};

\node [draw,shape=circle, fill, inner sep=0pt,minimum size=3pt] at (-1.25,2) {};
\node [above] at (-1.25,2) {$2$};
\node [draw,shape=circle, fill, inner sep=0pt,minimum size=3pt] at (-0.75,2) {};
\node [above] at (-0.75,2) {$2$};
\node [draw,shape=circle, fill, inner sep=0pt,minimum size=3pt] at (-0.25,2) {};
\node [above] at (-0.25,2) {$2$};
\node [draw,shape=circle, fill, inner sep=0pt,minimum size=3pt] at (0.25,2) {};
\node [above] at (0.25,2) {$2$};
\node [draw,shape=circle, fill, inner sep=0pt,minimum size=3pt] at (0.75,2) {};
\node [above] at (0.75,2) {$2$};
\node [draw,shape=circle, fill, inner sep=0pt,minimum size=3pt] at (1.25,2) {};
\node [above] at (1.25,2) {$2$};

\node [draw,shape=circle, fill, inner sep=0pt,minimum size=3pt] at (3,2) {};
\node [above] at (3,2) {$4$};
\node [below left] at (3,2) {$F_0$};
\node [draw,shape=circle, fill, inner sep=0pt,minimum size=3pt] at (3.5,2) {};
\node [above] at (3.5,2) {$4$};

\node [draw,shape=circle, fill, inner sep=0pt,minimum size=3pt] at (4.5,2) {};
\node [above] at (4.5,2) {$1$};
\node [draw,shape=circle, fill, inner sep=0pt,minimum size=3pt] at (5,2) {};
\node [above] at (5,2) {$1$};

\draw (0,0) -- (4,0)
(0,0) -- (-1.25,2)
(0,0) -- (-0.75,2)
(0,0) -- (-0.25, 2)
(0,0) -- (0.25,2)
(0,0) -- (0.75,2)
(0,0) -- (1.25,2)
(4,0) -- (3,2)
(4,0) -- (3.5, 2)
(4,0) -- (4.5, 2)
(4,0) -- (5,2);
\end{tikzpicture}\end{center}

The components $V_0$ and $V_1$ are rational, ruled surfaces. Moreover, if we define $L_i$ as the number of blow-ups in the birational morphism $V_i\to \overline{V}_i$ to the minimal ruled surface $\overline{V}_i$ for $i=0,1$, then $L_0=6$ and $L_1=9$. All other ten components are isomorphic to $\P^2$. 

The six flowers meeting $V_0$ are all of type $2A$. There are two flowers of type $2A$ meeting $V_1$ and two flowers of type $2B$.
Note that the surface $V_0$ exists as Example~\ref{ex:future-Wim} shows.

Suppose there exists a $K3$ surface $X$ over $K$ with an $snc$-model $\X$ such that ${\X_k\simeq D}$. Let $\omega$ be a volume form on $V_0$. Then the monodromy property does not hold for $X$. Indeed, Theorem~\ref{thm:poles} guarantees that $-\nu_0/4$ is a pole of the motivic zeta function of $X$, where $(4, \nu_0)$ is the numerical data of $F_0$, a conic-flower meeting $V_1$ with $N(F_0)=4$. 

On the other hand, the monodromy zeta function is
\begin{align*}
\zeta_X(T) &= \frac{(T^4-1)^{2}}{(T^2-1)^5(T^2-1)^6(T^4-1)^2(T-1)^2}\\
&= \frac{1}{(T^2-1)^{11}(T-1)^{2}},
\end{align*}
because
\begin{align*}
\chi(V_0^\circ) &= 4+L_0-2\lambda_0 = -2,\\
\chi(V_1^\circ) &= 4+L_1-2\lambda_1 = 5,
\end{align*}
where $\lambda_i$ is the number of flowers on $V_i$ for $i=0,1$.
Since $\xi(F_0)=\exp(\pi i \nu_0/2)$ is not a $2$-th root of unity by Lemma~\ref{thm:holds-conic-flower-xi}, we conclude that $\xi(F_0)$ is not a monodromy eigenvalue. Therefore, $X$ does not satisfy the monodromy property.

It is important to note that, a priori, we don't know whether such a strict normal crossings surface $D$ exists. And if such a strict normal crossings surface exist, we also don't know if $D$ can be smoothened to a $K3$ surface $X$. However, we will show in Section~\ref{sect:future-exclude}, that there does \emph{not} exist such a $K3$ surface $X$.
\end{countercandidate}
 
\begin{countercandidate} \label{count:exist}
Let $D$ be a strict normal crossings surface with the following dual graph, where the labels denote the multiplicity of the components:

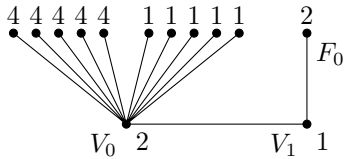
\begin{center}
\begin{tikzpicture}[scale=0.6]

\node [draw,shape=circle, fill, inner sep=0pt,minimum size=3pt] at (0,0) {};
\node [below left] at (0,0) {$V_0$};
\node [below right] at (0,0) {$2$};
\node [draw,shape=circle, fill, inner sep=0pt,minimum size=3pt] at (4,0) {};
\node [below left] at (4,0) {$V_1$};
\node [below right] at (4,0) {$1$};

\node [draw,shape=circle, fill, inner sep=0pt,minimum size=3pt] at (-2.5,2) {};
\node [above] at (-2.5,2) {$4$};
\node [draw,shape=circle, fill, inner sep=0pt,minimum size=3pt] at (-2,2) {};
\node [above] at (-2,2) {$4$};
\node [draw,shape=circle, fill, inner sep=0pt,minimum size=3pt] at (-1.5,2) {};
\node [above] at (-1.5,2) {$4$};
\node [draw,shape=circle, fill, inner sep=0pt,minimum size=3pt] at (-1,2) {};
\node [above] at (-1,2) {$4$};
\node [draw,shape=circle, fill, inner sep=0pt,minimum size=3pt] at (-0.5,2) {};
\node [above] at (-0.5,2) {$4$};
\node [draw,shape=circle, fill, inner sep=0pt,minimum size=3pt] at (0.5,2) {};
\node [above] at (0.5,2) {$1$};
\node [draw,shape=circle, fill, inner sep=0pt,minimum size=3pt] at (1,2) {};
\node [above] at (1,2) {$1$};
\node [draw,shape=circle, fill, inner sep=0pt,minimum size=3pt] at (1.5,2) {};
\node [above] at (1.5,2) {$1$};
\node [draw,shape=circle, fill, inner sep=0pt,minimum size=3pt] at (2,2) {};
\node [above] at (2,2) {$1$};
\node [draw,shape=circle, fill, inner sep=0pt,minimum size=3pt] at (2.5,2) {};
\node [above] at (2.5,2) {$1$};

\node [draw,shape=circle, fill, inner sep=0pt,minimum size=3pt] at (4,2) {};
\node [above] at (4,2) {$2$};
\node [below right] at (4,2) {$F_0$};

\draw (0,0) -- (4,0)
(0,0) -- (-2.5,2)
(0,0) -- (-2,2)
(0,0) -- (-1.5,2)
(0,0) -- (-1,2)
(0,0) -- (-0.5, 2)
(0,0) -- (0.5,2)
(0,0) -- (1,2)
(0,0) -- (1.5,2)
(0,0) -- (2,2)
(0,0) -- (2.5,2)
(4,0) -- (4,2);

\end{tikzpicture}\end{center}

The components $V_0$ and $V_1$ are rational, ruled surfaces. Moreover, if we define $L_i$ as the number of blow-ups in the birational morphism $V_i\to \overline{V}_i$ to the minimal ruled surface $\overline{V}_i$ for $i=0,1$, then $L_0=10$ and $L_1=7$. All other eleven components are isomorphic to $\P^2$.

There are five flowers of type $2A$ meeting $V_0$, and five of type $2B$. The flower meeting $V_1$ is of type $2B$, and we call the unique component of this flower $F_0$.

Suppose there exists a $K3$ surface $X$ over $K$ with an $snc$-model $\X$ such that ${\X_k\simeq D}$. Let $\omega$ be a volume form on $X$. Then the monodromy property does not hold for $X$. Indeed, Theorem~\ref{thm:poles} guarantees that $-\nu_0/2$ is a pole of the motivic zeta function of $X$, where $(2, \nu_0)$ is the numerical data of $F_0$. 

On the other hand, the monodromy zeta function is
\begin{align*}
\zeta_X(T) &= \frac{(T^2-1)^{6}}{(T-1)^9(T^4-1)^5(T-1)^5(T^2-1)}\\
&= \frac{1}{(T^2+1)^{5}(T-1)^{14}},
\end{align*}
because
\begin{align*}
\chi(V_0^\circ) &= 4+L_0-2\lambda_0 = -6,\\
\chi(V_1^\circ) &= 4+L_1-2\lambda_1 = 9,
\end{align*}
where $\lambda_i$ is the number of flowers on $V_i$ for $i=0,1$.
Since $\xi(F_0)$ is a primitive $2$-th root of unity by Lemma~\ref{thm:holds-conic-flower-xi}, we conclude that $\xi(F_0)$ is not a monodromy eigenvalue. Therefore, $X$ does not satisfy the monodromy property.

It is important to note that, a priori, we don't know whether  such a strict normal crossings surface $D$ exists. Actually, A. H\"oring suggested a strategy which seems to rule out the existence of $V_0$, and therefore excludes this particular combinatorial countercandidate.
Even if such a strict normal crossings surface $D$ existed, we would still need to prove that it can be smoothened to a $K3$ surface $X$. In Section~\ref{sect:future-exist}, we will give a strategy of how one might try to smoothen such a strict normal crossings surface to a $K3$ surface.
\end{countercandidate}

\subsection{Constructing combinatorial countercandidates: strategy} \label{sect:future-count-strategy}

Suppose $X$ is a $K3$ surface over $K$ of index $\iota(X)=1$ with a Crauder-Morrison model $\X$, not satisfying the monodromy property.
Theorem~\ref{thm:holds-pot} implies that $\X$ is a chain degeneration. 
Let $V_0, \ldots, V_{\alpha-1}, V_\alpha, \ldots, V_\beta, V_{\beta+1}, \ldots, V_{k+1}$ be the components in the chain, 
where $\alpha$ and $\beta$ are defined as in Proposition~\ref{thm:CM-chain-good} and $V_i\cap V_j=\emptyset$, except when $j\in \{i-1, i, i+1\}$. Set $N=N(V_\alpha)$.

As $X$ does not satisfy the monodromy property, Theorem~\ref{thm:holds-chain} and Proposition~\ref{thm:CM-chain-good} (iv) imply that (possibly after renumbering) the component $V_0$ is a rational, non-minimal ruled surface, that $N(V_0)=2N$, and that ${\alpha=1}$. Moreover, Theorem~\ref{thm:holds-chain} (iv) also implies that there exists a conic-flower $F$ meeting one of the components $V_1, \ldots, V_\beta$. Notice that this flower must be of type $2B$ or $2C$, because of Theorem~\ref{thm:CM-chain-flowers}.

From Table~\ref{table:CM-chain-classification} in Proposition~\ref{thm:CM-chain-good}, we can deduce that one of the following ten cases hold:
\begin{enumerate}
\label{future-cases}
\item $\beta=k+1$, and $V_{k+1}$ is a rational, ruled surface,
\item $\beta=k+1$, and $V_{k+1}$ is an elliptic, ruled surface,
\item $\beta=k+1$, and $V_{k+1}\simeq \P^2$,
\item $V_{k+1}\simeq \P^2$, and $N(V_{k+1})=3N$,
\item $V_{k+1}\simeq \P^2$, and $N(V_{k+1})=\frac{3}{2} N$,
\item $V_{k+1}\simeq \Sigma_2$, and $N(V_{k+1})=4N$,
\item $V_{k+1}\simeq \Sigma_2$, $\beta=k-1$, and $N(V_{k+1})=2N$,
\item $V_{k+1}\simeq \Sigma_2$, and $N(V_{k+1})=\frac{4}{3}N$,
\item $V_{k+1}$ is a rational, ruled surface, $\beta =k$, and $N(V_{k+1})=2N$,
\item $V_{k+1}$ is an elliptic, ruled surface, and $N(V_{k+1})=2N$.
\end{enumerate}

Define $\phi$ to be the number of conic-flowers meeting $V_1, \ldots, V_{\beta}$, and define $\gamma$ to be the number of flowers of type $2A$ meeting $V_1, \ldots, V_{\beta}$. We have that $\phi\geq 1$, because of the conic-flower $F$ meeting one of the components $V_1, \ldots,V_\beta$.

Define $\phi'$ and $\gamma'$ as follows: 
\begin{itemize}
\item in cases 9 and 10, $\phi'$ is defined as the number of conic-flowers meeting $V_0$ and $V_{k+1}$, and $\gamma'$ is defined as the number of flowers of type $2A$ meeting $V_0$ and $V_{k+1}$.
\item otherwise, we define $\phi'$ as the number of conic-flowers meeting $V_0$, and $\gamma'$ as the number of flowers of type $2A$ meeting $V_0$.
\end{itemize}

If $V_{k+1}$ is an elliptic, ruled surface (cases 2 and 10), then either there is exactly one non-rational flower of type $4\alpha$, or there are two non-rational flowers of type $4\alpha$, as seen in Proposition~\ref{thm:CM-chain-flowers}. In the latter case, both flowers are elliptic. If there is exactly one non-rational flower of type $4\alpha$, we denote by $g$ its genus. If there are two elliptic flowers of type $4\alpha$, we set $g=1$.

In cases 1 to 8, we will consider the variables $L_0$ and $\sum_{i=1}^\beta L_i$. Note that in cases 3 to 8, we have $L_{k+1}=0$. In cases 9 and 10, we will consider the variables $L_0+L_{k+1}$ and $\sum_{i=1}^\beta L_i$.

This means that we will work with the variables $\phi,\gamma, \phi', \gamma', N, \sum_{i=1}^{\beta} L_i$, and $L_0$ or $L_0+L_{k+1}$. In cases 2 and 10, we also consider the variable $g$. All these variables are non-negative integers.
We will describe how to find three equalities, similar to equations \eqref{eq:ineq-phi2}, \eqref{eq:degree-zeta-function} and \eqref{eq:relation-L}, that need to be satisfied in order to find a combinatorial countercandidate. We will also describe how to deduce upper bounds for all of these variables.

Define $I^\dagger\subset I$ to be the set of indices $i\in I$, where either
$\rho_i$ is minimal, or
$E_i$ is the top of a conic-flower.
Set $S^\dagger=\left\{(-\nu_i,N_i)\in \Z\times \Z_{>0}\mid i\in I^\dagger\right\}$. Corollary~\ref{thm:poles-contributions} says that
\[Z_{X,\omega}(T)\in \mathcal{M}_k^{\hat{\mu}}\left[T, \frac{1}{1-\L^{a}T^b}\right]_{(a,b)\in S^\dagger}.\]
By Theorem~\ref{thm:holds-lct}, we know that $\xi(E_i)$ is a monodromy eigenvalue, if $\rho_i$ is minimal. Therefore, since $X$ does not satisfy the monodromy property, there must be a conic-flower $F$ such that $\xi(F_0)$ is \emph{not} a monodromy eigenvalue. By Proposition~\ref{thm:CM-chain-flowers}, the flower $F$ is of type $2B$ or $2C$.
An argument similar as in the proof of Theorem~\ref{thm:holds-chain} (iii), gives that $F$ meets the chain in one of the components $V_1, \ldots, V_\beta$, and therefore $N(F_0)=2N$.

Recall that we have the following formula for the monodromy zeta-function:
\[\zeta_X(T) = \prod_{i\in I} (T^{N_i}-1)^{-\chi(E_i^\circ)}.\]
Using Lemma~\ref{thm:CM-euler-characteristics}, we can compute all relevant $\chi(E_i^\circ)$ in terms of the variables $\phi,\gamma, \phi', \gamma', N, \sum_{i=1}^{\beta} L_i, g$, and $L_0$ or $L_0+L_{k+1}$. For example, in all cases, except for case 2, we have
\[\sum_{i=1}^\beta \chi(V_i^\circ) = \sum_{i=1}^\beta L_i - 2\phi -2\gamma.\]


Since $\xi(F_0)$ is not a monodromy eigenvalue, it is not a zero nor a pole of $\zeta_X(T)$. Lemma~\ref{thm:holds-conic-flower-xi} gives that $\xi(F_0)$ is not an $N$-th root of unity, and it is not difficult to see that the factor $(T^{2N}-1)$ must cancel completely, by writing down the monodromy zeta function explicitely in all ten cases. This induces the \textbf{first equality}. 
Moreover, the monodromy zeta function must be of the form $\frac{1}{Q(T)}$ with $Q(T)$ a polynomial of degree 24. Computing the degree of the monodromy zeta function gives the \textbf{second equality}. 
The \textbf{third equality} is given by Lemma~\ref{thm:future-relation-L}.

Since all variables $\phi,\gamma, \phi', \gamma', N, g, \sum_{i=1}^{\beta} L_i$, and $L_0$ or $L_0+L_{k+1}$ are non-negative integers, we have obvious lower bounds. We can deduce upper bounds for all variables as well. 
Lemma~\ref{thm:future-relation-L} gives an upper bound for $\sum_{i=1}^\beta L_i$, and $L_0$ or $L_0+ L_{k+1}$. Lemma~\ref{thm:CM-relation-L-flower-elliptic} and Lemma~\ref{thm:future-relation-L-flower} give upper bounds for $\phi, \gamma, \phi', \gamma'$. Since we require the index $\iota(X)$ to be 1, the $\gcd$ of the multiplicities of the components in the special fiber must be $1$, which gives information on the variable $N$. Finally, also the variable $g$ in cases 2 and 10, is bounded, because for $g$ too big, it is not possible for the monodromy zeta function to have degree~24.

Because of the upper and lower bounds on integer variables, one can verify for every possibility whether the equalities described above are satisfied. This results in a finite list of combinatorial countercandidates. All 63 combinatorial countercandidates are listed in Appendix~\ref{app:combinatorial-countercandidate}, where we omitted those that will be excluded in Section~\ref{sect:future-exclude}.

\subsection{Partial results}\label{sect:future-count-result}

The strategy described in Subsection~\ref{sect:future-count-strategy} can also be used to prove the monodromy property for certain chain degenerations. We will see that there are no combinatorial countercandidates in cases 6 to 10, which implies that the monodromy property holds for $K3$ surfaces of index 1 with such chain degenerations.

\begin{proposition}\label{thm:future-chain-extra}
Let $X$ be a $K3$ surface over $K$ of index ${\iota(X)=1}$, with Crauder-Morrison model $\X$ and with special fiber $\X_k=\sum_{i\in I} N_iE_i$.
Suppose $\X$ is a chain degeneration satisfying Proposition~\ref{thm:CM-chain-good}. 
Let $V_0, V_{1}, \ldots, V_{k}, V_{k+1}$ be the components in the chain, 
where $V_i\cap V_j=\emptyset$, except when $j\in \{i-1, i, i+1\}$. Suppose there exists an integer $N\geq 1$, such that
\[N(V_i) = \begin{cases} 2N & \text{if } i=0 \text{ or } k+1,\\
N &\text{otherwise.} \end{cases}\]
Then $X$ satisfies the monodromy property.
\end{proposition}

\begin{proof}

Let $X$ be a $K3$ surface over $K$ of index $\iota(X)=1$. Suppose $X$ admits a Crauder-Morrison model $\X$ with special fiber $\X_k=\sum_{i\in I} N_iE_i$ as described in the statement of the proposition. Suppose $X$ does \emph{not} satisfy the monodromy property.

From Theorem~\ref{thm:holds-chain} (ii), it follows that $V_0$ and/or $V_{k+1}$ is a rational, non-minimal ruled surface. Without loss of generality, we can assume $V_0$ is rational, non-minimal ruled. Because $N(V_{k+1})=2N$ and $\beta=k$, Proposition~\ref{thm:CM-chain-good} gives that $V_{k+1}$ is a rational or elliptic, ruled surface.

Let $\phi$ be the total number of conic-flowers meeting $V_1, \ldots, V_k$, and let $\gamma$ be the total number of flowers of type $2A$ meeting $V_1, \ldots, V_k$. We have that 
$\phi\geq 1$,
by Theorem~\ref{thm:holds-chain} (iv). Moreover, let $\phi'$ be the total number of conic-flowers meeting $V_0$ and $V_{k+1}$, and let $\gamma'$ be the total number of flowers of type $2A$ meeting $V_0$ and $V_{k+1}$.

Define $I^\dagger\subset I$ to be the set of indices $i\in I$, where either
$\rho_i$ is minimal, or
$E_i$ is the top of a conic-flower.
Set $S^\dagger=\left\{(-\nu_i,N_i)\in \Z\times \Z_{>0}\mid i\in I^\dagger\right\}$. Corollary~\ref{thm:poles-contributions} says that
\[Z_{X,\omega}(T)\in \mathcal{M}_k^{\hat{\mu}}\left[T, \frac{1}{1-\L^{a}T^b}\right]_{(a,b)\in S^\dagger}.\]
By Theorem~\ref{thm:holds-lct}, we know that $\xi(E_i)$ is a monodromy eigenvalue, when $\rho_i$ is minimal. Therefore, there must be a conic-flower $F$ with top $F_0$, such that $\xi(F_0)$ is \emph{not} a monodromy eigenvalue, since $X$ does not satisfy the monodromy property. By Proposition~\ref{thm:CM-chain-flowers}, the flower $F$ is of type $2B$ or $2C$.

Suppose first that $V_{k+1}$ is an elliptic, ruled surface.
Proposition~\ref{thm:CM-chain-flowers} implies that either there is exactly one non-rational flower, which is of type $4\alpha$, or there are exactly two non-rational flowers, both elliptic and of type $4\alpha$. In the first case, set $g$ to be the genus of this unique non-rational flower, and in the latter case, set $g=1$.
We then have
\[\chi(V_{k+1}^\circ) = L_{k+1} -(2-2g) - 2\lambda,\]
where $\lambda$ is the number of rational flowers meeting $V_{k+1}$.

Because of Lemma~\ref{thm:CM-euler-characteristics}, the monodromy zeta function is given by
\begin{equation}
\zeta_X(T) = \frac{(T^{2N} - 1)^{- \chi(V_0^{\circ}) - \chi(V_{k+1}^{\circ})}(T^{N}-1)^{2g-2}}{(T^{4N} - 1)^{\phi'}(T^{2N} - 1)^{\phi}(T^{N} - 1)^{\sum_{i=1}^k\chi(V_i^{\circ})+\gamma'} (T^{N/2} - 1)^{\gamma} }. \label{eq:zetafunction}
\end{equation}
Suppose $F$ meets $V_0$ or $V_{k+1}$. Then $N(F_0)=4N$, and by Lemma~\ref{thm:holds-conic-flower-xi}, we have that $\xi(F_0)$ is not a $2N$-th root of unity. Therefore, $\xi(F_0)$ is a pole of $\zeta_X(T)$, and hence a monodromy eigenvalue, contrary to our assumption.

This means that $F$ meets one of the components $V_1, \ldots, V_k$, and $N(F_0)=2N$. Since $\xi(F_0)$ is not a monodromy eigenvalue, it is not a pole nor a zero of $\zeta_X(T)$.
By Lemma~\ref{thm:holds-conic-flower-xi}, we have that $\xi(F_0)$ is not an $N$-th root of unity, and therefore the factor $(T^{2N}-1)$ must cancel completely in $\zeta_X(T)$. This means that
$-\chi(V_0^\circ)-\chi(V_{k+1}^\circ) = \phi+\phi'$.
Since $\chi(V_0^{\circ}) + \chi(V_{k+1}^{\circ}) = 2+L_0+L_{k+1}-2\phi'-2\gamma'+2g$, we get
\begin{equation}\phi'+2\gamma'-\phi-2g-L_0-L_{k+1} = 2.
\label{eq:GMP-not3}
\end{equation}

Moreover, equation~\eqref{eq:zetafunction} can be written as
\begin{equation*}
\zeta_X(T) = \frac{1}{(T^{2N} + 1)^{\phi'}(T^{N} - 1)^{\sum_{i=1}^{k}\chi(V_i^{\circ})+\gamma'-2g+2}  (T^{N/2} - 1)^{\gamma} }. 
\end{equation*}
Proposition~\ref{thm:holds-form-zeta-function-K3} gives that
\[2\phi'+\sum_{i=1}^{k}\chi(V_i^{\circ}) +\gamma'-2g+2+ \frac{\gamma}{2}=\frac{24}{N}. \]
Since $\sum_{i=1}^{k}\chi(V_i^\circ)=\sum_{i=1}^{k}L_i-2\phi-2\gamma$, we get
\begin{equation}
\sum_{i=1}^{k}L_i+2\phi'+\gamma' -2g -2\phi - \frac{3\gamma}{2} = \frac{24}{N}-2. \label{eq:geen-inspiratie3}
\end{equation}

Furthermore, Lemma~\ref{thm:future-relation-L} implies that 
\begin{equation}
L_0+2\sum_{i=1}^k L_i + L_{k+1} = 8. \label{eq:relation-L-2N-N-2N-rational3}
\end{equation} 

We will now derive upper bounds for the variables $\phi, \gamma, \phi', \gamma', g, N, L_0+L_{k+1}$ and $\sum_{i=1}^k L_i$.
First of all, \eqref{eq:relation-L-2N-N-2N-rational3} gives obvious upper bounds for $L_0+L_{k+1}$ and $\sum_{i=1}^k L_i$.

Combining Lemma~\ref{thm:future-relation-L-flower} and Corollary~\ref{thm:CM-relation-L-flower-elliptic}, we get that 
\begin{equation}
\phi'+\gamma'\leq L_0+L_{k+1}+1. \label{eq:number-flowers-L3}
\end{equation}

Finally, Corollary~\ref{thm:CM-relation-L-flower-elliptic} implies that
\begin{equation}
2\phi+2\gamma \leq \sum_{i=1}^k L_i. \label{eq:relation-L3}
\end{equation}

We also have a bound on $g$, because by Lemma~\ref{thm:holds-form-zeta-function-K3}, we know that
\[\sum_{i=1}^k \chi(V_i^\circ) + \gamma' - 2g+2 \geq 0,\]
and therefore
\begin{equation}
g\leq \Big(\sum_{i=1}^k L_i - 2\phi -2\gamma + \gamma' +2\Big)/2.
\end{equation}

Because $X$ has index $\iota(X)=1$, Proposition~\ref{thm:future-index-gcd} gives that $\gcd_{i\in I} N_i = 1$. It is not difficult to see that
\[ \gcd_{i\in I} N_i =
\begin{cases}
N & \text{if } \gamma=0,\\
N/2 & \text{if } \gamma>0,
\end{cases}
\]
and hence
\begin{equation}
N = \begin{cases}
1 & \text{if } \gamma=0,\\
2 & \text{if } \gamma>0.
\end{cases} \label{eq:index3}
\end{equation}

One can verify that there are no integer solutions of equations \eqref{eq:GMP-not3}, \eqref{eq:geen-inspiratie3} and \eqref{eq:relation-L-2N-N-2N-rational3} bounded by \eqref{eq:number-flowers-L3} to \eqref{eq:index3}. This is a contradiction, and we conclude that $X$ satisfies the monodromy property. 
We have implemented this in Python, and the code can be found in Section~\ref{sect:app-case10} of the appendix.
The code can also be downloaded from \mbox{\url{www.github.com/AnneliesJaspers/combinatorial-counterexamples}}.

Suppose now that $V_{k+1}$ is a rational, ruled surface.
A similar strategy gives the following equations:
\[\begin{cases}
\phi'+2\gamma'-\phi-(L_0+L_{k+1}) = 8,\\
\sum_{i=1}^k L_i + 2\phi' + \gamma'-2\phi-\frac{3\gamma}{2} =  \frac{24}{N},\\
L_0 + 2\sum_{i=1}^k L_i + L_{k+1} = 16,\\
\phi'+\gamma'\leq L_0+L_{k+1}+2,\\
2\phi + 2\gamma \leq \sum_{i=1}^k L_i,\\
N = \begin{cases}
1 & \text{if } \gamma = 0,\\
2 & \text{if } \gamma > 0.
\end{cases}
\end{cases}\]

One can verify that there are no non-negative integer solutions satifying these relations. This is a contradiction, and we conclude that $X$ satisfies the monodromy property. 
We have implemented this in Python, and the code can be found in Section~\ref{sect:app-case9} of the appendix.
The code can also be downloaded from \mbox{\url{www.github.com/AnneliesJaspers/combinatorial-counterexamples}}.
\end{proof}

\begin{proposition}\label{thm:future-unique-countercandidate}
Let $X$ be a $K3$ surface over $K$ of index ${\iota(X)=1}$, with Crauder-Morrison model $\X$ and with special fiber $\X_k=\sum_{i\in I} N_iE_i$.
Suppose $\X$ is a chain degeneration. 
Let $V_0, V_{1}, \ldots, V_{k}, V_{k+1}$ be the components in the chain, 
where $V_i\cap V_j=\emptyset$, except when $j\in \{i-1,i,  i+1\}$. Suppose there exists an integer $N\geq 1$, such that
\[N(V_i) = \begin{cases} 2N & \text{if } i=0,\\
N &\text{otherwise.} \end{cases}\]
Assume moreover that $V_0$ is a rational, ruled surface, and that $V_{k+1}$ is an elliptic, ruled surface.
If $\X_k\not\simeq D$, where $D$ is the strict normal crossings surface described in Combinatorial Countercandidate~\ref{count:elliptic}, then $X$ satisfies the monodromy property. 
\end{proposition}

\begin{countercandidate} \label{count:elliptic}
Let $D$ be a strict normal crossings surface as follows: it is a chain degeneration with components $V_0, V_{1}, \ldots, V_{k}, V_{k+1}$, where $V_i\cap V_j=\emptyset$, except when $j\in \{i-1,i,  i+1\}$. We have $N(V_0)=4$, and $N(V_i)=2$, for $i=1, \ldots, k+1$.

There is a unique flower of type $2B$ or $2C$ meeting $V_0$, and four flowers of type $2A$ are meeting $V_0$.
Moreover, there is a unique flower of type $2B$ or $2C$ meeting one of the components $V_1, \ldots, V_{k+1}$, and denote the top of this flower by $F_0$. Finally, there is a unique flower of type $4\alpha$ meeting $V_{k+1}$, and it has genus 7. 

Furthermore, we have $L_0=4$ and $\sum_{i=1}^{k+1} L_i=2$.

For $k=0$ and when there are no flowers of type $2C$, the dual graph of $D$ is given in Figure~\ref{fig:countercandidate-elliptic}.
\begin{figure}[H]
\begin{center}
\begin{tikzpicture}[scale=0.6]

\node [draw,shape=circle, fill, inner sep=0pt,minimum size=3pt] at (0,0) {};
\node [below left] at (0,0) {$V_0$};
\node [below right] at (0,0) {$4$};
\node [draw,shape=circle, fill, inner sep=0pt,minimum size=3pt] at (4,0) {};
\node [below left] at (4,0) {$V_1$};
\node [below right] at (4,0) {$2$};

\node [draw,shape=circle, fill, inner sep=0pt,minimum size=3pt] at (-1.25,2) {};
\node [above] at (-1.25,2) {$8$};

\node [draw,shape=circle, fill, inner sep=0pt,minimum size=3pt] at (-0.25,2) {};
\node [above] at (-0.25,2) {$2$};
\node [draw,shape=circle, fill, inner sep=0pt,minimum size=3pt] at (0.25,2) {};
\node [above] at (0.25,2) {$2$};
\node [draw,shape=circle, fill, inner sep=0pt,minimum size=3pt] at (0.75,2) {};
\node [above] at (0.75,2) {$2$};
\node [draw,shape=circle, fill, inner sep=0pt,minimum size=3pt] at (1.25,2) {};
\node [above] at (1.25,2) {$2$};

\node [draw,shape=circle, fill, inner sep=0pt,minimum size=3pt] at (3.5,2) {};
\node [above] at (3.5,2) {$1$};

\node [draw,shape=circle, fill, inner sep=0pt,minimum size=3pt] at (4.5,2) {};
\node [above] at (4.5,2) {$4$};
\node [below right] at (4.5, 2) {$F_0$};

\draw (0,0) -- (4,0)
(0,0) -- (-1.25,2)
(0,0) -- (-0.25, 2)
(0,0) -- (0.25,2)
(0,0) -- (0.75,2)
(0,0) -- (1.25,2)
(4,0) -- (3.5,2)
(4,0) -- (4.5, 2);
\end{tikzpicture}\end{center}
\caption{} \label{fig:countercandidate-elliptic}
\end{figure}

Suppose there exists a $K3$ surface $X$ over $K$ with an $snc$-model $\X$ such that ${\X_k\simeq D}$. Let $\omega$ be a volume form on $X$. Then $X$ does not satisfy the monodromy property. Indeed, Theorem~\ref{thm:poles} guarantees that $-\nu_0/4$ is a pole of the motivic zeta function of $X$, where $(4,\nu_0)$ is the numerical data of $F_0$.

Since $\chi(V_0^\circ) = 4 + 4 - 2 \cdot (1+4) = -2$, and $\sum_{i=1}^{k+1} \chi(V_i^{\circ}) = 2 - 2 + (2\cdot 7 - 2) = 12$, Lemma~\ref{thm:CM-euler-characteristics} gives that
\begin{align*}
\zeta_X(T) &= \frac{(T^4-1)^2 (T-1)^{12}}{(T^2-1)^{12}(T^8-1)(T^2-1)^4(T^4-1)}\\
&= \frac{1}{(T^4+1)(T^2-1)^4(T+1)^{12}}.
\end{align*}

Since $\xi(F_0)$ is a primitive $4$-th root of unity by Lemma~\ref{thm:holds-conic-flower-xi}, we conclude that $\xi(F_0)$ is not a monodromy eigenvalue. Therefore, $X$ does not satisfy the monodromy property.

It is important to note that, a priori, we don't know whether  such a strict normal crossings surface $D$ exists. 
And even if such a strict normal crossings surface $D$ exists, we don't know whether it can be smoothened to a $K3$ surface~$X$. 
\end{countercandidate}

\emph{Proof of Proposition~\ref{thm:future-unique-countercandidate}.}
Let $X$ be a $K3$ surface over $K$ of index $\iota(X)=1$. Suppose $X$ admits a Crauder-Morrison model $\X$ with special fiber $\X_k=\sum_{i\in I} N_iE_i$ as described in the statement of the proposition. Suppose $X$ does \emph{not} satisfy the monodromy property.

Let $\phi$ be the total number of conic-flowers meeting $V_1, \ldots, V_{k+1}$, and let $\gamma$ be the total number of flowers of type $2A$ meeting $V_1, \ldots, V_{k+1}$. Moreover, let $\phi'$ be the total number of flowers of type $2B$ and $2C$ meeting $V_0$, and let $\gamma'$ be the total number of flowers of type $2A$ meeting $V_0$.

Following the strategy described in Subsection~\ref{sect:future-count-strategy}, we find the following equalities
\[\begin{cases}
\phi'+2\gamma'-\phi-L_0 = 4,\\
4\phi'+2\gamma'-4\phi-3\gamma + 2g+2\sum_{i=1}^{k+1} L_i = 26,\\
L_0+2\sum_{i=1}^{k+1} L_i = 8,\\
2\phi+2\gamma \leq \sum_{i=1}^{k+1} L_i,\\
\phi'+\gamma' \leq L_0+1.
\end{cases}\]
One can verify that there is a unique integer solution, which is $L_0 = 4, \sum_{i=1}^{k+1}L_i = 2, \phi = 1, \gamma = 0, \phi'=1, \gamma'=4, g=7$. We have implemented this in Python, and the code can be found in Section~\ref{sect:app-case2} of the appendix. The code can also be downloaded from\\ \mbox{\url{www.github.com/AnneliesJaspers/combinatorial-counterexamples}}.
\qed

The following proposition shows that a $K3$ surface of index 1 with a chain described by the cases 6 to 8, satisfies the monodromy property.
\begin{proposition}
Let $X$ be a $K3$ surface over $K$ of index $\iota(X)=1$, with Crauder-Morrison model $\X$ and with special fiber $\X_k=\sum_{i\in I} N_i E_i$. Suppose $\X$ is a chain degeneration.
Let $V_0, V_1, \ldots, V_k, V_{k+1}$ be the components in the chain, where $V_i\cap V_j = \emptyset$, except when $j\in \{i-1, i,i+1\}$. Set $N=\min \{N(V_i)\mid i=0, \ldots, k+1\}$.
Suppose that $V_{k+1}$ is a rational, minimal ruled surface, and that $N(V_{k+1})> N$. Then $X$ satisfies the monodromy property.
\end{proposition}

\begin{proof}
From Theorem~\ref{thm:holds-chain} (ii) and (iii), we know that $V_0$ is a rational, non-minimal ruled surface, and $N(V_0)>N$. From Proposition~\ref{thm:CM-chain-good}, it follows that one of the following cases hold:
\begin{itemize}
\item $V_{k+1}\simeq \Sigma_2$, and $N(V_{k+1})=4N$,
\item $V_{k+1}\simeq \Sigma_2$, $\beta=k-1$, and $N(V_{k+1})=2N$,
\item $V_{k+1}\simeq \Sigma_2$, and $N(V_{k+1})=\frac{4}{3}N$,
\item $V_{k+1}$ is a rational, ruled surface, $\beta =k$, and $N(V_{k+1})=2N$.
\end{itemize}
The result now follows from the Python code in Sections~\ref{sect:app-case6}, \ref{sect:app-case7}, \ref{sect:app-case8}, and \ref{sect:app-case9}, where we have proven that there are no combinatorial countercandidates.
\end{proof}

\section{Existence of a counterexample: a strategy} \label{sect:future-exist}

In this section, we will describe a possible strategy to produce a $K3$ surface~$X$ that does not satisfy the monodromy property. The idea is to smoothen a combinatorial countercandidate to a $K3$ surface, which then will not satisfy the monodromy property. In the literature, there are several smoothening results for \emph{semi-stable} normal crossings varieties to $K3$ surfaces, for example \cite[Theorem~5.10]{Friedman} and \cite[Corollary~5.4.15]{Hanke}. However, we would like to smoothen a strict normal crossings variety that is \emph{not} semi-stable. In this section, we will describe how it may be possible to use the results for semi-stable normal crossings varieties anyway.

\textbf{Step 1: Choose a combinatorial countercandidate.}

The aim is to prove that there exists a $K3$ surface admitting a model with the chosen combinatorial countercandidate as special fiber. So we keep one of the 63 combinatorial countercandidates in mind. 

\textbf{Step 2: Check that there exists a strict normal crossings surface $D_\X$ as described by the combinatorial countercandidate.}

A priori, it is not clear whether there exists a strict normal crossings surface $D_\X$ as described by the combinatorial countercandidate. We need to verify the existence of the irreducible components, and whether they can be glued in the described way. Part of the difficulty is to check that on every irreducible component, the desired double curves exist.

For example, to check the existence of the surface $V_0$ in Combinatorial Countercandidate~\ref{count:exclude}, we need to verify that there exists a rational, ruled surface, which is constructed by blowing up six points on a minimal ruled surface, and that allows seven disjoint, irreducible curves, where one curve is elliptic and numerically equivalent to the anticanonical bundle, and the other six curves are rational of self-intersection $-2$.
This surface exists by Example~\ref{ex:future-Wim}.

Proving the existence of the components is not at all trivial. For example, to prove that the component $V_0$ in Combinatorial Countercandidate~\ref{count:exist} exists,
one needs to verify the existence of a rational, ruled surface, which is constructed by blowing up ten points on a minimal ruled surface, with ten disjoint, rational curves of self-intersection $-2$. A. H\"oring suggested a strategy to prove that such a surface does \emph{not} exist, which would exclude this combinatorial countercandidate as a candidate to produce a counterexample of the monodromy property.

If the desired strict normal crossings surface $D_\X$ exists, we write $D_\X=\sum_{i\in I} N_i E_i$.

\textbf{Step 3: Compute the semi-stable strict normal crossings surface $D_\Y$ obtained by normalizing the semi-stable reduction.}

Suppose there exists a $K3$ surface $X$ over $K$ with an $snc$-model $\X$ such that the special fiber $\X_k\simeq D_\X$. Set $n=\lcm_{i\in I} N_i$ and $R_n=R[\pi]/(\pi^{n}- t)$. Let $\Y$ be the normalization of $\X\times_{R} R_n$. It is possible to compute some properties of the special fiber~$\Y_k$, using the results in Section~\ref{sect:poles-tilde-structure-Y} and general results of cyclic coverings, for example found in \cite[Sections~I.16 and I.17]{BHPV}.

Note that we don't know whether $\Y$ actually exists, but it should be possible to compute some properties of the special fiber $\Y_k$, \emph{if} $\Y$ exists. We can then ignore $\Y$ and define $D_\Y$ as a strict normal crossings surface that satisfies these properties.  Similar as in Step 2, we should verify that such a strict normal crossings surface indeed exists.

\textbf{Step 4: Check the existence of the $\hat{\mu}$-action on $D_\Y$.}

If $\Y$ exists, there is a $\hat{\mu}$-action on $\Y$, and therefore there should be a $\hat{\mu}$-action on $\Y_k$ satisfying the properties found in the proofs of the propositions in Subsections~\ref{sect:poles-tilde-top-flower} and \ref{sect:poles-tilde-middle-flower}.

Therefore, we need to verify that $D_\Y$ allows a $\hat{\mu}$-action with the desired properties.

\textbf{Step 5: Compute the semi-stable strict normal crossings variety~$D_\mathcal{Z}$ obtained by $\hat{\mu}$-equivariantly resolving singularities.}

If $\Y$ exists, Theorem~\ref{thm:poles-local-comp-smoothness} (i) gives the singular locus of $\Y$. Let $\mathcal{Z}$ be the resolution of singularities of $\Y$, where the blow-ups are done in a $\hat{\mu}$-equivariant way.
Note that we don't know whether $\mathcal{Z}$ actually exists, but it should be possible to compute some properties of the special fiber $\mathcal{Z}_k$, \emph{if} $\mathcal{Z}$ exists. We can then ignore $\mathcal{Z}$ and define $D_\mathcal{Z}$ as a strict normal crossings surface equipped with a $\hat{\mu}$-action, that satisfies these properties.

Similar as in Step 2, we should verify that such a strict normal crossings surface $D_\mathcal{Z}$ indeed exists.

\textbf{Step 6: Use results on smoothening semi-stable normal crossings varieties to get an algebraic space $\widetilde{\X}$ over $R_n$. This needs to be done $\hat{\mu}$-equivariantly.}

We would like to prove that there exists a smooth algebraic space $\widetilde{\X}$ over $R_n$ with a $\hat{\mu}$-action, such that there exists a $\hat{\mu}$-equivariant isomorphism $\widetilde{X}_k\simeq D_\mathcal{Z}$. Possibly useful results are \cite[Theorem~5.10]{Friedman} and \cite[Corollary~5.4.15]{Hanke}.

\textbf{Step 7: Define $\X$ to be the quotient space $\widetilde{\X}/\hat{\mu}$. Prove that $X\simeq \X\times_R K$ is a $K3$ surface that does not satisfy the monodromy property.}

Define $\X$ to be the quotient space $\widetilde{\X}/\hat{\mu}$. This is an algebraic space over $R$, and $X\simeq \X\times_R K$ is a $K3$ surface over $K$. In general, we won't have $\X_k\simeq D_\X$, but we expect that it should be possible to contract components in $\X_k$ such that we obtain $D_\X$ as a special fiber. This means that $X$ admits a model with $D_\X$ as a special fiber, and hence is does not satisfy the monodromy property.

\section{Excluding some combinatorial countercandidates} \label{sect:future-exclude}

By executing the steps described in Section~\ref{sect:future-exist}, we may come across an ad hoc reason to exclude a combinatorial countercandidate. For example, for Combinatorial Countercandidate~\ref{count:exist}, we are stuck at Step 2. We are not able to prove the existence of the surface $V_0$, and we ask ourselves whether such a surface exists. If there is indeed no such surface $V_0$, then there does not exist a strict normal crossings surface as described in Combinatorial Countercandidate~\ref{count:exist}, which would rule out this combinatorial countercandidate. A strategy shown to us by A. H\"oring, indeed seems to prove that the surface $V_0$ does not exist, and therefore it excludes Combinatorial Countercandidate~\ref{count:exist}.

We also tried to execute the strategy for Combinatorial Countercandidate~\ref{count:exclude}, but failed to complete Step 3. This results in proving that there does not exist an algebraic space $\Y$ over $R_n$, as defined in Step 3. Therefore, Combinatorial Countercandidate~\ref{count:exclude} cannot be smoothened to a $K3$ surface. This section is devoted to explaining this reasoning.
We will prove that there does not exist a $K3$ surface $X$ over $K$ with an $snc$-model $\X$ such $\X_k\simeq D$, where $D$ is a strict normal crossings surface with the properties described in Combinatorial Countercandidate~\ref{count:exclude}. Therefore, we can rule out \ref{count:exclude} as a combinatorial countercandidate. 

\begin{proposition}
There does \emph{not} exist a $K3$ surface $X$ over $K$ with Crauder-Morrison model $\X$ such that the special fiber $\X_k$ has the following properties: the dual graph is described in Figure~\ref{fig:countercandidate-exclude}, where the labels denote the multiplicities of the components.

\begin{figure}[H]
\begin{center}
\begin{tikzpicture}[scale=0.6]

\node [draw,shape=circle, fill, inner sep=0pt,minimum size=3pt] at (0,0) {};
\node [below left] at (0,0) {$V_0$};
\node [below right] at (0,0) {$4$};
\node [draw,shape=circle, fill, inner sep=0pt,minimum size=3pt] at (4,0) {};
\node [below left] at (4,0) {$V_1$};
\node [below right] at (4,0) {$2$};

\node [draw,shape=circle, fill, inner sep=0pt,minimum size=3pt] at (-1.25,2) {};
\node [above] at (-1.25,2) {$2$};
\node [draw,shape=circle, fill, inner sep=0pt,minimum size=3pt] at (-0.75,2) {};
\node [above] at (-0.75,2) {$2$};
\node [draw,shape=circle, fill, inner sep=0pt,minimum size=3pt] at (-0.25,2) {};
\node [above] at (-0.25,2) {$2$};
\node [draw,shape=circle, fill, inner sep=0pt,minimum size=3pt] at (0.25,2) {};
\node [above] at (0.25,2) {$2$};
\node [draw,shape=circle, fill, inner sep=0pt,minimum size=3pt] at (0.75,2) {};
\node [above] at (0.75,2) {$2$};
\node [draw,shape=circle, fill, inner sep=0pt,minimum size=3pt] at (1.25,2) {};
\node [above] at (1.25,2) {$2$};

\node [draw,shape=circle, fill, inner sep=0pt,minimum size=3pt] at (3,2) {};
\node [above] at (3,2) {$4$};
\node [draw,shape=circle, fill, inner sep=0pt,minimum size=3pt] at (3.5,2) {};
\node [above] at (3.5,2) {$4$};

\node [draw,shape=circle, fill, inner sep=0pt,minimum size=3pt] at (4.5,2) {};
\node [above] at (4.5,2) {$1$};
\node [draw,shape=circle, fill, inner sep=0pt,minimum size=3pt] at (5,2) {};
\node [above] at (5,2) {$1$};

\draw (0,0) -- (4,0)
(0,0) -- (-1.25,2)
(0,0) -- (-0.75,2)
(0,0) -- (-0.25, 2)
(0,0) -- (0.25,2)
(0,0) -- (0.75,2)
(0,0) -- (1.25,2)
(4,0) -- (3,2)
(4,0) -- (3.5, 2)
(4,0) -- (4.5, 2)
(4,0) -- (5,2);
\end{tikzpicture}\end{center}
\caption{} \label{fig:countercandidate-exclude}
\end{figure}

Denote by $V_0$ and $V_1$ the two components in the chain, where $V_0$ has multiplicity~4 and $V_1$ has multiplicity 2.
The surfaces $V_0$ and $V_1$ are rational, ruled surfaces, and if we define $L_i$ as the number of blow-ups in the birational morphism $V_i\to \overline{V}_i$ to the minimal ruled surface $\overline{V}_i$ for $i=0,1$, then $L_0=6$ and $L_1=9$.
The six flowers meeting $V_0$ are all of type $2A$. There are two flowers of type $2A$ meeting $V_1$, and two flowers of type $2B$ meeting $V_1$. 
\end{proposition}

\begin{proof}
Suppose there exists a $K3$ surface $X$ with the properties described in the statement.
Set $C_1=V_0\cap V_1$ and, denote by $C_2, \ldots, C_7$ the flowercurves on $V_0$. We know that $C_1$ is an elliptic curve and $C_2, \ldots, C_7$ are rational curves. Moreover, they are pairwise disjoint as there are no triple points.

Let $\Y$ be the normalization of $\X\times_R R[\pi]/(\pi^4-t)$, and let $f\colon\Y\to\X$ be the induced morphism. Denote by $W$ an irreducible component of $f^{-1}(V_0)_{red}$.
Proposition \ref{thm:poles-local-comp-etale-branch} gives that $f_{V_0}\colon f^{-1}(V_0)_{red}\to V_0$ is \'etale above $V_0^\circ$ of degree 4, and that it is ramified of index 2 above $C_1, \ldots, C_7$. 
 Therefore, $W\to V_0$ is an $m$-cyclic cover, with $m=2$ or $4$, ramified above $C_1+C_2+\cdots + C_7$.

In \cite[Section I.17]{BHPV}, it is explained that the existence of this $m$-cyclic cover implies the existence of a line bundle $\mathcal{L}$ on $V_0$ with
\[\mathcal{L}^{\otimes m} = \O_{V_0}(C_1+C_2+\cdots + C_7).\]

Let $D$ be the divisor on $V_0$ associated with the line bundle $\mathcal{L}$. We have $mD \equiv C_1+C_2+\cdots+ C_7$, and therefore, we find
\begin{align*}
m^2 D^2 &= 	(C_1+C_2+\cdots + C_7)^2\\
	&=		 (C_1^2+C_2^2+\cdots +C_7^2),
\end{align*}
where the second equality is explained by the fact that the curves $C_1, \cdots , C_7$ are pairwise disjoint.
From Lemma~\ref{thm:CM-self-intersection-flower-curve}, we know that $C_2^2=\cdots = C_7^2=-2$. Moreover, from equation~\eqref{eq:canonical-bundle-component-general} and Remark~\ref{rmk:CM-canonical-bundle-min-weight}, it follows that
\[K_{V_0}\equiv -C_1.\]
Because $L_0=6$, we have $C_1^2= 8-6 =2$, by \cite[Corollary~V.2.11 and Proposition~V.3.3]{Hartshorne}

It follows that $m^2D^2=10$, which implies that $D^2$ is not an integer, which is a contradiction. As a consequence, an $m$-cyclic cover of $V_0$ ramified along $C_1+\cdots + C_7$ does not exist. 
We conclude that there does not exist a \mbox{$K3$ surface $X$} with the desired properties.
\end{proof}

We can repeat this argument in more generality. Let $X$ be a $K3$ surface over $K$ of index $\iota(X)=1$ with a Crauder-Morrison model $\X$ that does not satisfy the monodromy property. Write $\X_k=\sum_{i\in I} N_i E_i$ and set $n=\lcm_{i\in I} N_i$. Since $X$ does not satisfy the monodromy property, $\X$ is a chain degeneration. Let $V_0, \ldots, V_{k+1}$ be the surfaces in the chain, where $V_i\cap V_j=\emptyset$, except when $j\in \{i-1,i,i+1\}$. Assume that $V_0$ is a rational, ruled surface, and suppose moreover that $N(V_0)=2N(V_1)$, and that $N(V_1)=\min \{N(V_i)\mid i=0, \ldots, k+1\}$. 
Since $X$ has index $\iota(X)=1$, it is possible to conclude that $N(V_1)$ is either 1 or~2.

Set $C_1=V_0\cap V_1$, and let $C_2, \ldots, C_{\gamma'+1}$ be the flowercurves of flowers of type $2A$ on $V_0$.
Let $\Y$ be the normalization of $\X\times_R R[\pi]/(\pi^n-t)$, and let $f\colon\Y\to\X$ be the induced morphism. Denote by $W$ an irreducible component of $f^{-1}(V_0)_{red}$.
The morphism $W\to V_0$ is an $m$-cyclic cover, with $m=2$ or 4, ramified along $C_1+C_2+\cdots+C_{\gamma'+1}$. Therefore, there exists a divisor $D$ on $V_0$ with the property that
\begin{align*}
m^2D^2 &= (C_1+C_2+\cdots+C_{\gamma'+1})^2\\
&= 8-L_0 -2\gamma'.
\end{align*}
Since $m$ is 2 or 4, we conclude that $L_0+2\gamma'$ is divisible by 4. Therefore, we can exclude all combinatorial countercandidates with $L_0+2\gamma' \not\equiv 0 \mod 4$.

By adding this requirement to our Python-code, we are able to exclude 26 of the 68 combinatorial countercandidates where the chain has components $V_0, \ldots, V_{k+1}$ and
\[N(V_i) = 
\begin{cases}
2N(V_1) & \text{if } i=0,\\
N(V_1) & \text{otherwise.}
\end{cases}\]

The method described above does not rule out Combinatorial Countercandidate~\ref{count:exist}, but A. H\"oring suggested a strategy to also exclude this combinatorial countercandidate. We believe this method could maybe be generalized to rule out even more combinatorial countercandidates.

\section{Some questions for further research}
\label{sect:future-questions}

As is probably clear by now, we don't know yet whether all $K3$ surfaces allowing a triple-point-free model satisfy the monodromy property. We produced a finite list of combinatorial countercandidates, but we don't know if one or more of them can be smoothened to a $K3$ surface. This is the first obvious question we suggest for further research:

\begin{itemize}
\item Is it possible to smoothen a combinatorial countercandidate to a $K3$ surface? If so, then there exists a $K3$ surface not satisfying the monodromy property.
\item Is it possible to exclude all combinatorial countercandidates? If so, then all $K3$ surfaces allowing a triple-point-free model satisfy the monodromy property. A possible strategy would be to generalize a method suggested by A. H\"oring to exclude Combinatorial Countercandidate~\ref{count:exist}.
\end{itemize}

If it is possible to construct a $K3$ surface not satisfying the monodromy property, then it would be very interesting to investigate the poles of the motivic zeta function that do not give rise to monodromy eigenvalues. Some possible questions are:
\begin{itemize}
\item If there is a $K3$ surface $X$ not satisfying the monodromy property, can we classify the extra poles of $Z_{X,\omega}(T)$ that do not induce monodromy eigenvalues? Do these poles have special characteristics?
\item Is there a property $P$ such that all $K3$ surfaces admitting a triple-point-free model satisfy the following: all poles of the motivic zeta function with property $P$, induce monodromy eigenvalues?
\end{itemize}

Although there are several differences between the monodromy conjecture for hypersurface singularities and the context of Calabi-Yau varieties, there are also many mysterious similarities. Therefore, it would be quite interesting if we could find close links between the two settings. Some questions, we ask ourselves are:
\begin{itemize}
\item If there is a $K3$ surface $X$ not satisfying the monodromy property, what is the influence of the existence of $X$ on the monodromy conjecture for hypersurface singularities? Could $X$ serve as an inspiration to construct a counterexample of this famous conjecture?
\item Are there techniques used in the setting of Calabi-Yau varieties that could turn out to be fruitful in the context of hypersurface singularities as well?
\end{itemize}

If all combinatorial countercandidates can be excluded, then we know that every $K3$ surface allowing a triple-point-free model satisfies the monodromy property. So far, we distinguished many cases for which we proved the monodromy property separately, and there were a lot of ad hoc arguments. It would be nice to see a more general principle that summarizes all these individual cases. Furthermore, we would like to know whether \emph{all} $K3$ surfaces satisfy the monodromy property, also those without a triple-point-free model. And what about general Calabi-Yau varieties? To summarize our questions:
\begin{itemize}
\item Is there a general argument that shows that $K3$ surfaces allowing a triple-point-free model satisfy the monodromy property?
\item Do all $K3$ surfaces satisfy the monodromy property?
\item Do all Calabi-Yau varieties satisfy the monodromy property?
\end{itemize}

\cleardoublepage



\appendix

\chapter{Formulas for the contribution of a flower and Python code}\label{ch:appendix1}

In Section~\ref{sect:contributions-flowers}, we defined the contribution $Z_F(T)$ of a flower $F$ to the motivic zeta function. In the proof of Theorem~\ref{thm:poles-contribution-flowers}, we explained how to compute the contribution for every type of flower. In Section~\ref{app:formula-contribution}, we will list the explicit formulas for the contribution of flowers. In Section~\ref{app:Sympy-code}, we give the Python code that was used to compute the contributions. This implementation is written in Python 3.6 and SymPy 1.0. The code can be downloaded from \mbox{\url{www.github.com/AnneliesJaspers/flowers_contribution}}. 

\section{Formulas for the contribution of a flower}
\label{app:formula-contribution}

In this section, we list the explicit formulas for the contribution of a flower to the motivic zeta fucntion. For flowers where the top is isomorphic to $\P^2$ and where the double curve is a line, the contribution can be found in Table~\ref{table:contr_P^2_line}. In Table~\ref{table:contr_P^2_conic}, we give the contribution of a conic-flower. Finally, the contribution of a flower with a ruled top is given in Table~\ref{table:contr_ruled}.

\def\arraystretch{2}
\begin{sidewaystable}
\centering
\begin{tabular}{|c|l|}
\hline
 & Contribution $Z_F(T)$ \\ 
\hline  \hline
$2A$ & $\frac{[\mathcal{P}] \L^{-\nu_0+2} T^{N_0} (1+\L^{-\nu_0} T^{N_0})}{1-\L^{-\nu_{\ell+1}}T^{N_{\ell+1}}}$ \\
\hline 
$3A$ & $\frac{[\mathcal{P}] \L^{-\nu_0+2} T^{N_0} (\L^{-2\nu_0+2}T^{2N_0} + \L^{-\nu_0+2}T^{N_0} +1)}{1-\L^{-\nu_{\ell+1}}T^{N_{\ell+1}}} $\\
\hline
$3B$ & $\frac{[\mathcal{P}] \L^{-\nu_0+2} T^{N_0} (\L^{-2\nu_0}T^{2N_0} + \L^{-\nu_0}T^{N_0} +1)}{1-\L^{-\nu_{\ell+1}}T^{N_{\ell+1}}} $\\
\hline
$4A$ & $\frac{[\mathcal{P}] \L^{-\nu_0+2} T^{N_0} (\L^{-3\nu_0+4}T^{3N_0} + \L^{-2\nu_0+4}T^{2N_0} + \L^{-\nu_0+2}T^{N_0} +1)}{1-\L^{-\nu_{\ell+1}}T^{N_{\ell+1}}} $\\
\hline
$4B$ & $\frac{[\mathcal{P}] \L^{-\nu_0+2} T^{N_0} (\L^{-\nu_0}T^{N_0}+1)(\L^{-2\nu_0}T^{2N_0}+1)}{1-\L^{-\nu_{\ell+1}}T^{N_{\ell+1}}} $\\
\hline
$6A$ & $\frac{[\mathcal{P}] \L^{-\nu_0+2} T^{N_0} (\L^{-5\nu_0+8}T^{5N_0}+\L^{-4\nu_0+8}T^{4N_0}+\L^{-3\nu_0+6}T^{3N_0} + \L^{-2\nu_0+4}T^{2N_0} + \L^{-\nu_0+2}T^{N_0} +1)}{1-\L^{-\nu_{\ell+1}}T^{N_{\ell+1}}} $\\
\hline
$6B$ & $\frac{[\mathcal{P}] \L^{-\nu_0+2} T^{N_0} (\L^{-\nu_0}T^{N_0}+1)(\L^{-2\nu_0}T^{2N_0}-\L^{-\nu_0}T^{N_0}+1)( \L^{-2\nu_0}T^{2N_0} + \L^{-\nu_0}T^{N_0} +1)}{1-\L^{-\nu_{\ell+1}}T^{N_{\ell+1}}} $\\
\hline
\end{tabular}
\caption{Contribution to the motivic zeta function of a flower with $F_0\simeq\P^2$ and $C_1$ a line.}
\label{table:contr_P^2_line}
\end{sidewaystable}

\begin{sidewaystable}
\centering
\begin{tabular}{|c|l|}
\hline 
 & Contribution $Z_F(T)$\\ 
\hline 
\hline
$2B$ & $\frac{[\widetilde{F_0^\circ}] \L^{-\nu_0}T^{N_0}}{1-\L^{-\nu_0}T^{N_0}}+(\L-1)[\widetilde{C_{1}}]\frac{\L^{-\nu_0}T^{N_0}}{1-\L^{-\nu_0}T^{N_0}}\frac{\L^{-\nu_{\ell+1}}T^{N_{\ell+1}}}{1-\L^{-\nu_{\ell+1}}T^{N_{\ell+1}}}$\\
\hline
$2C$ & $\frac{[\widetilde{F_0^\circ}] \L^{-\nu_0}T^{N_0}}{1-\L^{-\nu_0}T^{N_0}} 
+
\frac{
[\widetilde{C_1}]\L^{-\nu_1-1}T^{N_1}\left((\L-1)\L^\ell T^{2N_1}+(\L^\ell-1)\L^{-\nu_1+1}T^{N_1} + (\L^\ell-1)\L\right)
}
{(1-\L^{-\nu_0}T^{N_0})(1-\L^{-\nu_{\ell+1}}T^{N_{\ell+1}})}$\\
\hline
$4C$ & $\frac{[\widetilde{F_0^\circ}] \L^{-\nu_0}T^{N_0}}{1-\L^{-\nu_0}T^{N_0}} 
+
\frac{
[\widetilde{C_1}]\L^{-\nu_1-1}T^{N_1}\left((\L-1)\L^{-3\nu_1+2\ell-1} T^{3N_1}+(\L^\ell-1)\L^{-2\nu_1+\ell}T^{2N_1} + (\L^{2\ell-1}-1)\L^{-\nu_1+1}T^{N_1} + (\L^\ell-1)\L\right)
}
{(1-\L^{-\nu_0}T^{N_0})(1-\L^{-\nu_{\ell+1}}T^{N_{\ell+1}})}$\\
\hline
$6C$ & $\frac{[\widetilde{F_0^\circ}] \L^{-\nu_0}T^{N_0}}{1-\L^{-\nu_0}T^{N_0}}+
\Big(\frac{
[\widetilde{C_1}]\L^{-\nu_1-1}T^{N_1}
}
{(1-\L^{-\nu_0}T^{N_0})(1-\L^{-\nu_{\ell+1}}T^{N_{\ell+1}})}$\\
&${\scriptstyle \times \left((\L-1)\L^{-4\nu_1+3\ell-4} T^{4N_1}+(\L^{\ell-1}-1)\L^{-3\nu_1+2\ell-2}T^{3N_1} + (\L^{2\ell-2}-1)\L^{-2\nu_1+\ell-1}T^{2N_1} + (\L^{2\ell-2}-1)\L^{-\nu_1+1}T^{N_1} + (\L^{\ell-1}-1)\L\right)\Big)}$\\
\hline
$6D$ & $\frac{[\widetilde{F_0^\circ}] \L^{-\nu_0}T^{N_0}}{1-\L^{-\nu_0}T^{N_0}}+
\Big(\frac{
[\widetilde{C_1}]\L^{-\nu_1-1}T^{N_1}
}
{(1-\L^{-\nu_0}T^{N_0})(1-\L^{-\nu_{\ell+1}}T^{N_{\ell+1}})}$\\
&${\scriptstyle \times \left((\L-1)\L^{-4\nu_1+3\ell-2} T^{4N_1}+(\L^{\ell}-1)\L^{-3\nu_1+2\ell-1}T^{3N_1} + (\L^{2\ell-1}-1)\L^{-2\nu_1+\ell}T^{2N_1} + (\L^{2\ell-1}-1)\L^{-\nu_1+1}T^{N_1} + (\L^{\ell}-1)\L\right)\Big)}$\\
\hline
$6E$ & $\frac{[\widetilde{F_0^\circ}] \L^{-\nu_0}T^{N_0}}{1-\L^{-\nu_0}T^{N_0}}+(\L-1)[\widetilde{C_{1}}]\frac{\L^{-\nu_0}T^{N_0}}{1-\L^{-\nu_0}T^{N_0}}\frac{\L^{-\nu_{\ell+1}}T^{N_{\ell+1}}}{1-\L^{-\nu_{\ell+1}}T^{N_{\ell+1}}}$\\
\hline
\end{tabular}
\caption{Contribution to the motivic zeta function of a flower with $F_0\simeq\P^2$ meeting $F_1$ in a conic.}
\label{table:contr_P^2_conic}
\end{sidewaystable}

\begin{sidewaystable}
\centering
\begin{tabular}{|c|l|}
\hline 
 & Contribution $Z_F(T)$\\ 
\hline 
\hline
$4\alpha$ & $\frac{[\widetilde{C_1}] \L^{-\nu_0+1} T^{N_0} (1+\L^{-\nu_0} T^{N_0})}{1-\L^{-\nu_{\ell+1}}T^{N_{\ell+1}}}$ \\
\hline 
$6\alpha$ & $\frac{[\widetilde{C_1}] \L^{-\nu_0+1} T^{N_0} (\L^{-2\nu_0+1}T^{2N_0} + \L^{-\nu_0+1}T^{N_0} +1)}{1-\L^{-\nu_{\ell+1}}T^{N_{\ell+1}}} $\\
\hline
$6\beta$ & $\frac{[\widetilde{C_1}] \L^{-\nu_0+1} T^{N_0} (\L^{-2\nu_0}T^{2N_0} + \L^{-\nu_0}T^{N_0} +1)}{1-\L^{-\nu_{\ell+1}}T^{N_{\ell+1}}} $\\
\hline
$8\alpha$ & $\frac{[\widetilde{C_1}] \L^{-\nu_0+1} T^{N_0} (\L^{-3\nu_0+2}T^{3N_0} + \L^{-2\nu_0+2}T^{2N_0} + \L^{-\nu_0+1}T^{N_0} +1)}{1-\L^{-\nu_{\ell+1}}T^{N_{\ell+1}}} $\\
\hline
$8\beta$ & $\frac{[\widetilde{C_1}] \L^{-\nu_0+1} T^{N_0} (\L^{-\nu_0}T^{N_0}+1)(\L^{-2\nu_0}T^{2N_0}+1)}{1-\L^{-\nu_{\ell+1}}T^{N_{\ell+1}}} $\\
\hline
$12\alpha$ & $\frac{[\widetilde{C_1}] \L^{-\nu_0+1} T^{N_0} (\L^{-5\nu_0+4}T^{5N_0}+\L^{-4\nu_0+4}T^{4N_0}+\L^{-3\nu_0+3}T^{3N_0} + \L^{-2\nu_0+2}T^{2N_0} + \L^{-\nu_0+1}T^{N_0} +1)}{1-\L^{-\nu_{\ell+1}}T^{N_{\ell+1}}} $\\
\hline
$12\beta$ & $\frac{[\widetilde{C_1}] \L^{-\nu_0+1} T^{N_0} (\L^{-\nu_0}T^{N_0}+1)(\L^{-2\nu_0}T^{2N_0}-\L^{-\nu_0}T^{N_0}+1)( \L^{-2\nu_0}T^{2N_0} + \L^{-\nu_0}T^{N_0} +1)}{1-\L^{-\nu_{\ell+1}}T^{N_{\ell+1}}} $\\
\hline
\end{tabular}
\caption{Contribution to the motivic zeta function of a flower with $F_0$ a minimal ruled surface.} \label{table:contr_ruled}
\end{sidewaystable}

\pagebreak
\section{Python code} \label{app:Sympy-code}

In this section, we explain and give the Python code written to compute the contribution of a flower to the motivic zeta function. This implementation is written in Python 3.6 and  SymPy 1.0. Note that the code can be downloaded from \mbox{\url{www.github.com/AnneliesJaspers/flowers_contribution}}.

\subsection{Implementation of relevant functions}

The following code defines some preliminary functions used to compute the contributions of flowers to the motivic zeta function.

\begin{lstlisting}[caption = {}]
from sympy import * 
from enum import Enum

L, T, P, c, F = symbols('L T P c F', positive=True)
N, nu = symbols('N nu', real=True)

# L is class of affine line
# P is class of variety such that F_0^\circ = P L^2

class Group(Enum):
    LINE = 1
    CONIC = 2
    RULED = 3

def compute_nu_vector(N_vector, group):
    """Compute the multiplicities nu_i of div (omega) along E_i.
    
    
    Tool is Lemma 3.3.10 in thesis 'Monodromy property for K3 surfaces
        allowing a triple-point-free model'.
    Define nu such that 
        nu_0=nu (if F_0 ruled or F_0 isomorphic to P^2 with C_1 a line), or
        nu_0=2nu+1 (if F_0 isomorphic to P^2 with C_1e a conic).
    
    Arguments:
        N_vector is a vector with the multiplicities of a flower, 
            from N_0 to N_{ell+1}
        group is Group.LINE, Group.CONIC, or Group.RULED,
            depending on the type of the flower
    """
    
    if group == Group.LINE:
        # Define nu_0
        nu_vector = [nu] 
        
        # Define nu_1
        nu_vector.append(expand(nu_vector[0]*N_vector[1]/N_vector[0]-2))
        
    elif group == Group.CONIC:
        # Define nu_0 (always odd, when C_1 a conic)
        nu_vector = [2*nu+1] 
        
        # Define nu_1
        nu_vector.append(expand(nu_vector[0]*N_vector[1]/N_vector[0]-1/2))
    
    elif group == Group.RULED:
        # Define nu_0
        nu_vector = [nu]
        
        # Define nu_1
        nu_vector.append(expand(nu_vector[0]*N_vector[1]/N_vector[0]-1))

        
    for i in range(2,len(N_vector)):
        
        # Define nu_i
        nu_i = nu_vector[-1]*(N_vector[i-2]+N_vector[i])/N_vector[i-1]\
                 - nu_vector[-2]
        nu_vector.append(expand(nu_i))
        
    assert len(N_vector) == len(nu_vector)
        
    return nu_vector


    
def compute_fraction(N_1, nu_1):
    return L**(-nu_1) * T**(N_1)/ (1 - L**(-nu_1) * T**(N_1))

    
def compute_contribution(N_vector, group):
    """Computes contribution of the flower to the motivic zeta function.
    
    For Definition contribution see thesis, Definition 4.2.1.
    
    
    Arguments:
        N_vector is a vector with the multiplicities of a flower, 
            from N_0 to N_{ell+1}
        group is Group.LINE, Group.CONIC, or Group.RULED,
            depending on the type of the flower
    """
    
    expr = 0
    
    if group == Group.LINE:
        C = P * (L+1) # when F_0=P^2, then [C] = P * (L+1)
        F0 = P * L**2 # when F_0=P_2 with C a line, then [F_0^circ]= P * L^2
    elif group == Group.CONIC:
        C = c  
        F0 = F
    elif group == Group.RULED:
        C = c # when F_0 is ruled, we don't know what [C] is, so variable c
        F0 = c * L # when F_0 is ruled, C is section, so [F_0^circ] = c * L
        
    nu_vector = compute_nu_vector(N_vector, group)
    
    #First we add all the terms of the sum with |J| = 1
    for j in range(0,len(N_vector)-1):
        term = compute_fraction(N_vector[j],nu_vector[j])
        if j==0:
            term = F0 * term
        else:
            term = (L-1) * C * term # for all i, we have [C_i]=[C_{i+1}]

        expr += term
    
    # Now, we add the terms of the sum with |J| = 2.
    
    for j in range(0, len(N_vector)-1):
        # for all i, we have [C_i]=[C_{i+1}]
        term = (L-1) * C  * compute_fraction(N_vector[j], nu_vector[j])\
            * compute_fraction(N_vector[j+1], nu_vector[j+1]) 
        expr += term
    
    return expr

def compute_contribution_middle_conic(N_vector, m):
    """ Compute terms of the contribution of conic-flowers
    to the motivic zeta function, corresponding to components F_1, ..., F_m.
    
    
    Arguments:
        N_vector is a vector with the multiplicities of a flower, 
            from N_0 to N_{ell+1}
        m for how many terms we add: 
            those corresponding to components F_1, ..., F_m.
    """
    
    expr = 0
    
    C = c # when F_0=P^2, then [C] = permutation * (L+1) 
        
    nu_vector = compute_nu_vector(N_vector, Group.CONIC)
    
    #First we add all the terms of the sum with |J| = 1
    for j in range(1,m+1):
        term = compute_fraction(N_vector[j],nu_vector[j])
        term = (L-1) * C * term # for all i, we have [C_i]=[C_{i+1}]

        expr += term
    
    # Now, we add the terms of the sum with |J| = 2.   
    for j in range(1, m+1):
        # for all i, we have [C_i]=[C_{i+1}]
        term = (L-1) * C  * compute_fraction(N_vector[j], nu_vector[j])\
            * compute_fraction(N_vector[j+1], nu_vector[j+1]) 
        expr += term
    
    return expr
\end{lstlisting}

\subsection{Computation of the contribution of a flower, that is not a conic-flower}

Let $F$ be a flower that is not a conic-flower. From Table~\ref{table:flowers-P2-line} and \ref{table:flowers-ruled}, we see that, if we know the type of the flower, then we know its length $\ell$. For example, when $F$ is of type $4A$, we have $\ell=2$. Therefore, it is straightforward to compute the contribution of $F$ to the motivic zeta function. The following Python code computes these contributions.

\begin{lstlisting}[caption = {}]
from sympy import * 
from functions import *

##### Flowers of with top isomorphic to P^2 with C_1 a line #####
N_2A = [N, 2*N]
zeta_2A = factor(cancel(compute_contribution(N_2A, Group.LINE)))

N_3A = [N, 2*N,3*N]
zeta_3A = factor(cancel(compute_contribution(N_3A, Group.LINE)))

N_3B = [N, 3*N]
zeta_3B = factor(cancel(compute_contribution(N_3B, Group.LINE)))

N_4A = [N, 2*N, 3*N, 4*N]
zeta_4A = factor(cancel(compute_contribution(N_4A, Group.LINE)))

N_4B = [N, 4*N]
zeta_4B = factor(cancel(compute_contribution(N_4B, Group.LINE)))

N_6A = [N, 2*N, 3*N, 4*N, 5*N, 6*N]
zeta_6A = factor(cancel(compute_contribution(N_6A, Group.LINE)))

N_6B = [N, 6*N]
zeta_6B = factor(cancel(compute_contribution(N_6B, Group.LINE)))

print("The following are the contributions to the motivic zeta function\
          of flowers with top isomorphic to P^2 and C_1 a line: \n" )

print("Contribution of flower 2A:", zeta_2A, "\n")
print("Contribution of flower 3A:", zeta_3A, "\n")
print("Contribution of flower 3B:", zeta_3B, "\n")
print("Contribution of flower 4A:", zeta_4A, "\n")
print("Contribution of flower 4B:", zeta_4B, "\n")
print("Contribution of flower 6A:", zeta_6A, "\n")
print("Contribution of flower 6B:", zeta_6B, "\n")

##### Flowers with ruled top #####
N_4alpha = [N, 2*N]
zeta_4alpha = factor(cancel(compute_contribution(N_4alpha, Group.RULED)))

N_6alpha = [N, 2*N,3*N]
zeta_6alpha = factor(cancel(compute_contribution(N_6alpha, Group.RULED)))

N_6beta = [N, 3*N]
zeta_6beta = factor(cancel(compute_contribution(N_6beta, Group.RULED)))

N_8alpha = [N, 2*N, 3*N, 4*N]
zeta_8alpha = factor(cancel(compute_contribution(N_8alpha, Group.RULED)))

N_8beta = [N, 4*N]
zeta_8beta = factor(cancel(compute_contribution(N_8beta , Group.RULED)))

N_12alpha = [N, 2*N, 3*N, 4*N, 5*N, 6*N]
zeta_12alpha = factor(cancel(compute_contribution(N_12alpha, Group.RULED)))

N_12beta = [N, 6*N]
zeta_12beta = factor(cancel(compute_contribution(N_12beta, Group.RULED)))

print("The following are the contributions to the motivic zeta function\
          of flowers with top a ruled surface: \n" )

print("Contribution of flower 4alpha:", zeta_4alpha, "\n")
print("Contribution of flower 6alpha:", zeta_6alpha, "\n")
print("Contribution of flower 6beta:", zeta_6beta, "\n")
print("Contribution of flower 8alpha:", zeta_8alpha, "\n")
print("Contribution of flower 8beta:", zeta_8beta, "\n")
print("Contribution of flower 12alpha:", zeta_12alpha, "\n")
print("Contribution of flower 12beta:", zeta_12beta, "\n")
\end{lstlisting}

\subsection{Computation of the contribution of a conic-flower}

Let $F$ be a conic-flower. If $F$ is a flower of type $2B$ or $6E$, then $\ell=0$, otherwise, $\ell$ is chosen arbitrarily.
In the first case, it is immediately clear that
\[Z_{F}(T)\in \mathcal{M}_k^{\hat{\mu}}\left[T, \frac{1}{1-\L^{-\nu_{0}}T^{N_{0}}}, \frac{1}{1-\L^{-\nu_{\ell+1}}T^{N_{\ell+1}}}\right].\]

Therefore, we may assume $F$ has type $2C, 4C, 6C$ or $6D$. We will need an induction argument on $\ell$ to obtain a closed formula for $Z_F(T)$.

Note that for these flowers, we can find an integer $m$ with $1\leq m \leq \ell$ such that $N(F_i)=N(F_1)$ for all $1\leq i \leq m+1$. Note that $N(F_0)=2N(F_1)$. 
Define 
\begin{align*}
Z_{F,m}(T) = \sum_{j={1}}^{m}&\Bigg([\widetilde{F_j^\circ}]\frac{\L^{-\nu_j}T^{N_j}}{1-\L^{-\nu_j}T^{N_j}}\\
&+(\L-1)[\widetilde{C_{j+1}}]\frac{\L^{-\nu_j}T^{N_j}}{1-\L^{-\nu_j}T^{N_j}}\frac{\L^{-\nu_{j+1}}T^{N_{j+1}}}{1-\L^{-\nu_{j+1}}T^{N_{j+1}}}\Bigg).
\end{align*}

By induction on $m$, one can compute that
\begin{align*}
Z_{F,m}(T) = \frac{[\widetilde{C_1}](\L^m-1)\L^{-\nu_1}T^{N_1}}{(1-\L^{-\nu_1}T^{N_1})(1-\L^{-\nu_1+m}T^{N_1})}.
 \end{align*}

Let $m_0$ be the maximal $m$ with $1\leq m \leq \ell$ such that $N(F_i)=N(F_1)$ for all $1\leq i \leq m+1$.
Then
\begin{align*}
Z_F(T) = &\frac{[\widetilde{F_0^\circ}] \L^{-\nu_0}T^{N_0}}{1-\L^{-\nu_0}T^{N_0}}+(\L-1)[\widetilde{C_{1}}]\frac{\L^{-\nu_0}T^{N_0}}{1-\L^{-\nu_0}T^{N_0}}\frac{\L^{-\nu_{1}}T^{N_{1}}}{1-\L^{-\nu_{1}}T^{N_{1}}}\\
&+Z_{F,m_0}(T)\\
& + \sum_{j={m_0+1}}^{\ell}\Bigg([\widetilde{F_j^\circ}]\frac{\L^{-\nu_j}T^{N_j}}{1-\L^{-\nu_j}T^{N_j}}\\
&\qquad \qquad +(\L-1)[\widetilde{C_{j+1}}]\frac{\L^{-\nu_j}T^{N_j}}{1-\L^{-\nu_j}T^{N_j}}\frac{\L^{-\nu_{j+1}}T^{N_{j+1}}}{1-\L^{-\nu_{j+1}}T^{N_{j+1}}}\Bigg).
\end{align*} 

The following code first computes $Z_F(T)$ for flowers of type $2B$ and for flowers of type $6C$. Then it computes $Z_{F,m}(T)$ for $m=2, 3$ and $4$. Finally, we compute 
\begin{align*}
(\L-1)[\widetilde{C_{1}}]&\frac{\L^{-\nu_0}T^{N_0}}{1-\L^{-\nu_0}T^{N_0}}\frac{\L^{-\nu_{1}}T^{N_{1}}}{1-\L^{-\nu_{1}}T^{N_{1}}} +Z_{F,m_0}(T)\\
& + \sum_{j={m_0+1}}^{\ell}\Bigg([\widetilde{F_j^\circ}]\frac{\L^{-\nu_j}T^{N_j}}{1-\L^{-\nu_j}T^{N_j}}\\
&\qquad \qquad +(\L-1)[\widetilde{C_{j+1}}]\frac{\L^{-\nu_j}T^{N_j}}{1-\L^{-\nu_j}T^{N_j}}\frac{\L^{-\nu_{j+1}}T^{N_{j+1}}}{1-\L^{-\nu_{j+1}}T^{N_{j+1}}}\Bigg),
\end{align*} 
for flowers of type $2C, 4C, 6C$ and $6D$. To get the contribution of the flower, we need to add $\frac{[\widetilde{F_0^\circ}] \L^{-\nu_0}T^{N_0}}{1-\L^{-\nu_0}T^{N_0}}$ to the result.

\begin{lstlisting}[caption = {}]
from sympy import * 
from functions import *

l = symbols('l', real=True)

##### Computation of Z_{F,m} #####

print("\nWe first compute Z_{F,m}(T) for some small m.\n")                      
N_v = [2*N, N, N, N, N, N, N, N, N, N]

#m=2
Z_2 = compute_contribution_middle_conic(N_v, 2)
print('For m=2:', factor(cancel(Z_2)), "\n")
    
#m=3
Z_3 = compute_contribution_middle_conic(N_v, 3)
print('For m=3:', factor(cancel(Z_3)), "\n")

#m=4
Z_4 = compute_contribution_middle_conic(N_v, 4)
print('For m=4:', factor(cancel(Z_4)), "\n")
    
##### Flower of type 2C #####
 
# Declare Z_{F, m_0}. 
# m_0 = ell -1
Z_F_m = c * (L**(l-1)-1) * L**(-nu) * T**N \
                  / ((1-L**(-nu) * T**N)*(1-L**(-nu+l-1) * T**N))
                  
Z_2C = Z_F_m + c*(L-1)*compute_fraction(2*N,2*nu+1)*compute_fraction(N,nu)\
                 + c*(L-1)*compute_fraction(N,nu-l+1)\
                   + c*(L-1)*compute_fraction(N,nu-l+1)\
                      *compute_fraction(N,nu-l) 
                   
print('Contribution of flower 2C:', factor(cancel(Z_2C)), '\n')  

##### Flower of type 4C #####

# Declare Z_{F, m_0}. 
# m_0 = ell - 1
Z_F_m = c * (L**(l-1)-1) * L**(-nu) * T**N \
                  / ((1-L**(-nu) * T**N)*(1-L**(-nu+l-1) * T**N))
             
Z_4C = Z_F_m + c*(L-1)*compute_fraction(2*N,2*nu+1)*compute_fraction(N,nu) \
                   + c*(L-1)*compute_fraction(N,nu-l+1)\
                   + c*(L-1)*compute_fraction(N,nu-l+1)\
                      *compute_fraction(2*N,2*nu-2*l+1)       
    
print('Contribution of flower 4C:', factor(cancel(Z_4C)), '\n') 
    
##### Flower of type 6C #####

# Declare Z_{F, m_0}. 
# m_0 = ell - 2
Z_F_m = c * (L**(l-2)-1) * L**(-nu) * T**N \
          / ((1-L**(-nu) * T**N)*(1-L**(-nu+l-2) * T**N))
             
Z_6C = Z_F_m + c*(L-1)*compute_fraction(2*N,2*nu+1)*compute_fraction(N,nu) \
                   +c*(L-1)*compute_fraction(N,nu-l+2)\
                   +c*(L-1)*compute_fraction(N,nu-l+2)\
                      *compute_fraction(2*N,2*nu-2*l+3)\
                   +c*(L-1)*compute_fraction(2*N,2*nu-2*l+3)\
                   +c*(L-1)*compute_fraction(2*N,2*nu-2*l+3)\
                      *compute_fraction(3*N, 3*nu-3*l+4)
    
print('Contribution of flower 6C:', factor(cancel(Z_6C)), '\n')

##### Flower of type 6D #####

# Declare Z_{F, m_0}. 
# m_0 = ell - 1
Z_F_m = c * (L**(l-1)-1) * L**(-nu) * T**N \
                  / ((1-L**(-nu) * T**N)*(1-L**(-nu+l-1) * T**N))
             
Z_6D = Z_F_m + c*(L-1)*compute_fraction(2*N,2*nu+1)*compute_fraction(N,nu) \
              +c*(L-1)*compute_fraction(N,nu-l+1)\
              +c*(L-1)*compute_fraction(N,nu-l+1)\
                      *compute_fraction(3*N,3*nu-3*l+2)
    
print('Contribution of flower 6D:', factor(cancel(Z_6D)), '\n')      
\end{lstlisting}

\chapter{List of combinatorial countercandidates and Python code} \label{app:combinatorial-countercandidate}

In this appendix, we will give the Python code to compute all combinatorial countercandidates. We will see that we have a total of 63 combinatorial countercandidates, and we will list them all.

Suppose $X$ is a $K3$ surface of index $\iota(X)=1$ with a Crauder-Morrison model $\X$, not satisfying the monodromy property.
Theorem~\ref{thm:holds-pot} implies that $\X$ is a chain degeneration. 
Let $V_0, \ldots, V_{\alpha-1}, V_\alpha, \ldots, V_\beta, V_{\beta+1}, \ldots, V_{k+1}$ be the components in the chain, 
where $\alpha$ and $\beta$ are defined as in Proposition~\ref{thm:CM-chain-good}, and $V_i\cap V_j=\emptyset$, except when $j\in \{i-1,i, i+1\}$. Set $N=N(V_\alpha)$.

As explained in Section~\ref{sect:future-count-strategy}, one of the following ten cases hold:

\begin{enumerate}
\item $\beta=k+1$, and $V_{k+1}$ is a rational, ruled surface,
\item $\beta=k+1$, and $V_{k+1}$ is an elliptic, ruled surface,
\item $\beta=k+1$, and $V_{k+1}\simeq \P^2$,
\item $V_{k+1}\simeq \P^2$, and $N(V_{k+1})=3N$,
\item $V_{k+1}\simeq \P^2$, and $N(V_{k+1})=\frac{3}{2} N$,
\item $V_{k+1}\simeq \Sigma_2$, and $N(V_{k+1})=4N$,
\item $V_{k+1}\simeq \Sigma_2$, $\beta=k-1$, and $N(V_{k+1})=2N$,
\item $V_{k+1}\simeq \Sigma_2$, and $N(V_{k+1})=\frac{4}{3}N$,
\item $V_{k+1}$ is a rational, ruled surface, $\beta =k$, and $N(V_{k+1})=2N$,
\item $V_{k+1}$ is an elliptic, ruled surface, and $N(V_{k+1})=2N$.
\end{enumerate}

For every case, we will compute the combinatorial countercandidates and give the Python code.

This code has been written in Python 3.6. It can be downloaded from \mbox{\url{www.github.com/AnneliesJaspers/combinatorial-counterexamples}}.

\section{\texorpdfstring{$\beta=\lowercase{k}+1$, and $V_{\lowercase{k}+1}$ is a rational, ruled surface}{beta=k+1, and V(k+1) is a rational, ruled surface}}\label{sect:app-case1}
Suppose $\beta=k+1$, and $V_{k+1}$ is a rational, ruled surface.
The strategy described in Section~\ref{sect:future-count-strategy} explains that we have to find integer solutions satisfying the following conditions:

\[\begin{cases}
\phi'+2\gamma'-\phi-L_0 = 4,\\
\sum_{i=1}^{k+1} L_i + 2\phi'+\gamma'-2\phi-\frac{3\gamma}{2}=\frac{24}{N}-4,\\
L_0 + 2\sum_{i=1}^{k+1} L_i =24,\\
\phi'+\gamma'\leq L_0+1,\\
\phi+\gamma \leq \sum_{i=1}^{k+1} L_i +1,\\
N = \begin{cases}
1 & \text{if } \gamma=0,\\
2 & \text{if } \gamma>0,
\end{cases}\\
\sum_{i=1}^{k+1} L_i, \gamma, \phi' \text{ and } \gamma' \geq 0,\\
L_0 \text{ and } \phi \geq 1.
\end{cases}
\]

Moreover, in Section~\ref{sect:future-exclude}, we found the extra equation $L_0+2\gamma' \not\equiv 0 \mod 4$.

One can verify that there are \textbf{42 integer solutions} satisfying these relations, for example by using the Python code below.
\begin{lstlisting}[caption = {}]
count=0

# L1 for \sum_{i=1}^{k+1} L_i and L0 for L_0
for L1 in range(0, (24//2)+1): 	# L0 +2L1 = 24
	L0 = 24 - 2*L1				
	if L0 == 0:
		#V_0 is non-minimal ruled, so L_0 \neq 0.
		continue

	for phi in range(1,(L1+1)+1):	# phi+gamma \leq L1+1 and phi\geq 1
		for gamma in range(0, (L1-phi+1)+1, 2): # gamma is even because of degree zeta function
			
			# index = 1
			if gamma == 0:
				N=1				
			else:
				N=2

			for phi_ in range(0,(L0+1)+1):	# phi'+gamma' \leq L0+1
				for gamma_ in range(0,(L0-phi_+1)+1):

					# checking nothing went wrong
					assert phi_>=0 and gamma_>=0 and phi>=1 and gamma>=0 and L0>=1 and L1>=0 and (N==1 or N==2)
					assert gamma % 2 ==0
					assert (L0+2*L1) == 24
					assert (phi_+gamma_) <= (L0+1)
					assert (phi+gamma) <= (L1+1)
					
					# checking the two conditions 
					if (phi_+2*gamma_-phi-L0) != 4:
						# making sure the mondromy property does not hold
						continue
					if (L1+2*phi_+gamma_-2*phi-3*gamma//2) != ((24//N)-4):
						# maling sure the degree zeta function is 24
						continue
					
					if (L0+2*gamma_)%4 != 0:	
						# asserting the existence of line bundle such that cyclic cover exists
						continue
					
					# if we arrived here, we have found a combinatorial counterexample
					count += 1
					print('N =', N, ' L0 =',L0, ' sum L_i =',L1, ' phi =',phi, ' gamma =', gamma, ' phi\' =', phi_, ' gamma\' =', gamma_)
							
print('\nWe have produced', count, 'combinatorial counterexamples.')
\end{lstlisting}

The output of this code is the following.
{\footnotesize
\begin{verbatim}
N = 1  L0 = 24  sum L_i = 0  phi = 1  gamma = 0  phi' = 5  gamma' = 12
N = 1  L0 = 22  sum L_i = 1  phi = 1  gamma = 0  phi' = 5  gamma' = 11
N = 1  L0 = 22  sum L_i = 1  phi = 2  gamma = 0  phi' = 6  gamma' = 11
N = 1  L0 = 20  sum L_i = 2  phi = 1  gamma = 0  phi' = 5  gamma' = 10
N = 1  L0 = 20  sum L_i = 2  phi = 2  gamma = 0  phi' = 6  gamma' = 10
N = 1  L0 = 20  sum L_i = 2  phi = 3  gamma = 0  phi' = 7  gamma' = 10
N = 1  L0 = 18  sum L_i = 3  phi = 1  gamma = 0  phi' = 5  gamma' = 9
N = 1  L0 = 18  sum L_i = 3  phi = 2  gamma = 0  phi' = 6  gamma' = 9
N = 1  L0 = 18  sum L_i = 3  phi = 3  gamma = 0  phi' = 7  gamma' = 9
N = 1  L0 = 18  sum L_i = 3  phi = 4  gamma = 0  phi' = 8  gamma' = 9
N = 1  L0 = 16  sum L_i = 4  phi = 1  gamma = 0  phi' = 5  gamma' = 8
N = 2  L0 = 16  sum L_i = 4  phi = 1  gamma = 4  phi' = 1  gamma' = 10
N = 1  L0 = 16  sum L_i = 4  phi = 2  gamma = 0  phi' = 6  gamma' = 8
N = 1  L0 = 16  sum L_i = 4  phi = 3  gamma = 0  phi' = 7  gamma' = 8
N = 1  L0 = 16  sum L_i = 4  phi = 4  gamma = 0  phi' = 8  gamma' = 8
N = 1  L0 = 16  sum L_i = 4  phi = 5  gamma = 0  phi' = 9  gamma' = 8
N = 1  L0 = 14  sum L_i = 5  phi = 1  gamma = 0  phi' = 5  gamma' = 7
N = 2  L0 = 14  sum L_i = 5  phi = 1  gamma = 4  phi' = 1  gamma' = 9
N = 1  L0 = 14  sum L_i = 5  phi = 2  gamma = 0  phi' = 6  gamma' = 7
N = 2  L0 = 14  sum L_i = 5  phi = 2  gamma = 4  phi' = 2  gamma' = 9
N = 1  L0 = 14  sum L_i = 5  phi = 3  gamma = 0  phi' = 7  gamma' = 7
N = 1  L0 = 14  sum L_i = 5  phi = 4  gamma = 0  phi' = 8  gamma' = 7
N = 1  L0 = 12  sum L_i = 6  phi = 1  gamma = 0  phi' = 5  gamma' = 6
N = 2  L0 = 12  sum L_i = 6  phi = 1  gamma = 4  phi' = 1  gamma' = 8
N = 1  L0 = 12  sum L_i = 6  phi = 2  gamma = 0  phi' = 6  gamma' = 6
N = 2  L0 = 12  sum L_i = 6  phi = 2  gamma = 4  phi' = 2  gamma' = 8
N = 1  L0 = 12  sum L_i = 6  phi = 3  gamma = 0  phi' = 7  gamma' = 6
N = 2  L0 = 12  sum L_i = 6  phi = 3  gamma = 4  phi' = 3  gamma' = 8
N = 1  L0 = 10  sum L_i = 7  phi = 1  gamma = 0  phi' = 5  gamma' = 5
N = 2  L0 = 10  sum L_i = 7  phi = 1  gamma = 4  phi' = 1  gamma' = 7
N = 1  L0 = 10  sum L_i = 7  phi = 2  gamma = 0  phi' = 6  gamma' = 5
N = 2  L0 = 10  sum L_i = 7  phi = 2  gamma = 4  phi' = 2  gamma' = 7
N = 2  L0 = 10  sum L_i = 7  phi = 3  gamma = 4  phi' = 3  gamma' = 7
N = 2  L0 = 10  sum L_i = 7  phi = 4  gamma = 4  phi' = 4  gamma' = 7
N = 1  L0 = 8  sum L_i = 8  phi = 1  gamma = 0  phi' = 5  gamma' = 4
N = 2  L0 = 8  sum L_i = 8  phi = 1  gamma = 4  phi' = 1  gamma' = 6
N = 2  L0 = 8  sum L_i = 8  phi = 1  gamma = 8  phi' = 5  gamma' = 4
N = 2  L0 = 8  sum L_i = 8  phi = 2  gamma = 4  phi' = 2  gamma' = 6
N = 2  L0 = 8  sum L_i = 8  phi = 3  gamma = 4  phi' = 3  gamma' = 6
N = 2  L0 = 6  sum L_i = 9  phi = 1  gamma = 4  phi' = 1  gamma' = 5
N = 2  L0 = 6  sum L_i = 9  phi = 2  gamma = 4  phi' = 2  gamma' = 5
N = 2  L0 = 4  sum L_i = 10  phi = 1  gamma = 4  phi' = 1  gamma' = 4

We have produced 42 combinatorial counterexamples.
\end{verbatim}
}

\section{\texorpdfstring{$\beta=\lowercase{k}+1$, and $V_{\lowercase{k}+1}$ is an elliptic, ruled surface}{beta=k+1, and V(k+1) is an elliptic, ruled surface}}
\label{sect:app-case2}

Suppose $\beta=k+1$, and $V_{k+1}$ is an elliptic, ruled surface.
The strategy described in Section~\ref{sect:future-count-strategy} explains that we have to find integer solutions satisfying the following conditions:

\[\begin{cases}
\phi'+2\gamma'-\phi-L_0 = 4,\\
4\phi'+2\gamma'-4\phi-3\gamma + 2g+2\sum_{i=1}^{k+1} L_i = 26,\\
L_0+2\sum_{i=1}^{k+1} L_i = 8,\\
2\phi+2\gamma \leq \sum_{i=1}^{k+1} L_i,\\
\phi'+\gamma' \leq L_0+1,\\
\sum_{i=1}^{k+1} L_i, \gamma, \phi' \text{ and } \gamma' \geq 0,\\
L_0 \text{ and } \phi \geq 1.
\end{cases}\]

One can verify that there is exactly \textbf{1 integer solution} satisfying these relations, for example by using the Python code below.
\begin{lstlisting}[caption = {}]
count=0

# L1 for \sum_{i=1}^{k+1} L_i and L0 for L_0
for L1 in range(0, (8//2)+1): 	# L0 +2L1 = 8
	L0 = 8-2*L1
	if L0 == 0:
		#V_0 is non-minimal ruled, so L_0 \neq 0.
		continue

	for phi in range(1,(L1//2)+1):	# 2*phi+2*gamma \leq sum L_i and phi\geq 1
		for gamma in range(0, ((L1-2*phi)//2)+1):

			for phi_ in range(0,(L0+1)+1):	# phi'+gamma' \leq L0+1
				for gamma_ in range(0, (L0-phi_+1)+1):

					# checking nothing went wrong
					assert phi_>=0 and gamma_>=0 and phi>=1 and gamma>=0 and L0>=1 and L1>=0
					assert (L0+2*L1) == 8
					assert (phi_+gamma_) <= (L0+1)
					assert (2*phi+2*gamma) <= L1



					# checking the two conditions 
					if (phi_+2*gamma_-phi -L0) != 4:
						# making sure the mondromy property does not hold
						continue

					if (26 -4*phi_-2*gamma_+4*phi-3*gamma-2*L1)%2 != 0:
						#making sure g is an integer
						continue

					g = (26 -4*phi_-2*gamma_+4*phi-3*gamma-2*L1)//2
						# degree zeta function is 24

					if g<1:
						# g \geq 1
						continue
										
					# if we arrived here, we have found a combinatorial counterexample
					count += 1
					print('N = 1', ' L0 =',L0, ' sum L_i =',L1, ' phi =',phi,\
						  ' gamma =', gamma,' phi\' =', phi_, ' gamma\' =', gamma_,\
						  ' g =', g)
							
print('\nWe have produced', count, 'combinatorial counterexamples.')
\end{lstlisting}

The output of this code is the following.
{\footnotesize
\begin{verbatim}
N = 1  L0 = 4  sum L_i = 2  phi = 1  gamma = 0  phi' = 1  gamma' = 4  g = 7

We have produced 1 combinatorial counterexamples.
\end{verbatim}
}

\section{\texorpdfstring{$\beta=\lowercase{k}+1$, and $V_{\lowercase{k}+1}\simeq \P^2$}{beta=k+1, and V(k+1) isomorphic to P2}}

Suppose $\beta=k+1$, and $V_{k+1}\simeq \P^2$.
The strategy described in Section~\ref{sect:future-count-strategy} explains that we have to find integer solutions satisfying the following conditions:

\[\begin{cases}
\phi'+2\gamma'-\phi-L_0 = 4,\\
\sum_{i=1}^{k+1} L_i+ \gamma'+2\phi'-2\phi-\frac{3\gamma}{2}= \frac{24}{N} -3,\\
L_0+2\sum_{i=1}^{k+1} L_i = 26,\\
2\phi+2\gamma \leq \sum_{i=1}^{k+1} L_i,\\
\phi'+\gamma' \leq L_0+1,\\
N = \begin{cases}
1 & \text{if } \gamma=0,\\
2 & \text{if } \gamma>0,
\end{cases}\\
\sum_{i=1}^{k+1} L_i, \gamma, \phi' \text{ and } \gamma' \geq 0,\\
L_0 \text{ and } \phi \geq 1.
\end{cases}\]

Moreover, in Section~\ref{sect:future-exclude}, we found the extra equation $L_0+2\gamma' \not\equiv 0 \mod 4$.

One can verify that there are \textbf{17 integer solutions} satisfying these relations, for example by using the Python code below.
\begin{lstlisting}[caption = {}]
count=0

# L1 for \sum_{i=1}^{k+1} L_i and L0 for L_0
for L1 in range(0, (26//2)+1): 	# L0+2*L1=26
	L0 = 26-2*L1	
	if L0 == 0:
		#V_0 is non-minimal ruled, so L_0 \neq 0
		continue

	for phi in range(1,(L1//2)+1):	# 2*phi+2*gamma \leq sum L_i and phi\geq 1
		for gamma in range(0, ((L1-2*phi)//2)+1, 2): #gamma even because of degree zeta function

			if gamma == 0:
				N=1
			else:
				N=2

			for phi_ in range(0,(L0+1)+1):	# phi'+gamma' \leq L0+1
				for gamma_ in range(0,(L0-phi_+1)+1):

					# checking nothing went wrong
					assert phi_>=0 and gamma_>=0 and phi>=1 and gamma>=0 and L0>=1 and L1>=0 and (N==1 or N==2)
					assert gamma%2 == 0
					assert (L0+2*L1) == 26
					assert (phi_+gamma_) <= L0+1
					assert (2*phi+2*gamma) <= L1
					
					# checking the two conditions 
					if (phi_+2*gamma_-phi-L0) != 4:
						# making sure the mondromy property does not hold
						continue
					if (L1+gamma_+2*phi_-2*phi-3*gamma//2) != (24//N-3):
						# making sure degree of zeta function is 24
						continue

					if (L0+2*gamma_)%4 != 0:	
						# asserting the existence of line bundle such that cyclic cover exists
						continue
					
					# if we arrived here, we have found a combinatorial counterexample
					count += 1
					print('N =', N, ' L0 =',L0, ' sum L_i =',L1, ' phi =',phi, ' gamma =', gamma, ' phi\' =', phi_, ' gamma\' =', gamma_)
							
print('\nWe have produced', count, 'combinatorial counterexamples.')
\end{lstlisting}

The output of this code is the following.
{\footnotesize
\begin{verbatim}
N = 1  L0 = 22  sum L_i = 2  phi = 1  gamma = 0  phi' = 5  gamma' = 11
N = 1  L0 = 20  sum L_i = 3  phi = 1  gamma = 0  phi' = 5  gamma' = 10
N = 1  L0 = 18  sum L_i = 4  phi = 1  gamma = 0  phi' = 5  gamma' = 9
N = 1  L0 = 18  sum L_i = 4  phi = 2  gamma = 0  phi' = 6  gamma' = 9
N = 1  L0 = 16  sum L_i = 5  phi = 1  gamma = 0  phi' = 5  gamma' = 8
N = 1  L0 = 16  sum L_i = 5  phi = 2  gamma = 0  phi' = 6  gamma' = 8
N = 1  L0 = 14  sum L_i = 6  phi = 1  gamma = 0  phi' = 5  gamma' = 7
N = 1  L0 = 14  sum L_i = 6  phi = 2  gamma = 0  phi' = 6  gamma' = 7
N = 1  L0 = 14  sum L_i = 6  phi = 3  gamma = 0  phi' = 7  gamma' = 7
N = 1  L0 = 12  sum L_i = 7  phi = 1  gamma = 0  phi' = 5  gamma' = 6
N = 1  L0 = 12  sum L_i = 7  phi = 2  gamma = 0  phi' = 6  gamma' = 6
N = 1  L0 = 12  sum L_i = 7  phi = 3  gamma = 0  phi' = 7  gamma' = 6
N = 1  L0 = 10  sum L_i = 8  phi = 1  gamma = 0  phi' = 5  gamma' = 5
N = 1  L0 = 10  sum L_i = 8  phi = 2  gamma = 0  phi' = 6  gamma' = 5
N = 1  L0 = 8  sum L_i = 9  phi = 1  gamma = 0  phi' = 5  gamma' = 4
N = 2  L0 = 6  sum L_i = 10  phi = 1  gamma = 4  phi' = 1  gamma' = 5
N = 2  L0 = 4  sum L_i = 11  phi = 1  gamma = 4  phi' = 1  gamma' = 4

We have produced 17 combinatorial counterexamples.
\end{verbatim}
}

\section{\texorpdfstring{$V_{\lowercase{k}+1}\simeq \P^2$, and $N(V_{\lowercase{k}+1})=3N$}{V(k+1) isomorphic to P2, and N(V(k+1))=3N}}

Suppose $V_{k+1}\simeq \P^2$, and $N(V_{k+1})=3N$.
The strategy described in Section~\ref{sect:future-count-strategy} explains that we have to find integer solutions satisfying the following conditions:

\[\begin{cases}
\phi'+2\gamma'-\phi-L_0 = 4,\\
\sum_{i=1}^{k+1} L_i+ \gamma'+2\phi'-2\phi-\frac{3\gamma}{2}= \frac{24}{N} -9,\\
L_0+2\sum_{i=1}^{k+1} L_i = 14,\\
2\phi+2\gamma \leq \sum_{i=1}^{k+1} L_i,\\
\phi'+\gamma' \leq L_0+1,\\
N = \begin{cases}
1 & \text{if } \gamma=0,\\
2 & \text{if } \gamma>0,
\end{cases}\\
\sum_{i=1}^{k+1} L_i, \gamma, \phi' \text{ and } \gamma' \geq 0,\\
L_0 \text{ and } \phi \geq 1.
\end{cases}\]

One can verify that there are \textbf{2 integer solutions} satisfying these relations, for example by using the Python code below.
\begin{lstlisting}[caption = {}]
count=0

# L1 for \sum_{i=1}^beta L_i and L0 for L_0
for L1 in range(0, (14//2)+1):		# L0+2* sum L_i=14
	L0 = 14 - 2*L1
	if L0 == 0:
		#V_0 is non-minimal ruled, so L_0 \neq 0
		continue		

	for phi in range(1,(L1//2)+1):	# 2*phi+2*gamma \leq sum L_i and phi\geq 1
		for gamma in range(0, ((L1-2*phi)//2)+1, 2): #gamma even because of degree zeta function
			
			# index = 1
			if gamma == 0:
				N=1				
			else:
				N=2

			for phi_ in range(0,(L0+1)+1):	# phi'+gamma' \leq L0+1
				for gamma_ in range(0,(L0-phi_+1)+1):

					# checking nothing went wrong
					assert phi_>=0 and gamma_>=0 and phi>=1 and gamma>=0 and L0>=1 and L1>=0 and (N==1 or N==2)
					assert gamma%2 == 0
					assert (L0+2*L1) == 14
					assert (phi_+gamma_) <= L0+1
					assert (2*phi+2*gamma) <= L1
					
					# checking the two conditions 
					if (phi_+2*gamma_-phi-L0) != 4:
						# making sure the mondromy property does not hold
						continue
					if (L1+gamma_+2*phi_-2*phi-3*gamma//2) != ((24//N)-9):
						# making sure degree of zeta function is 24
						continue
					
					# if we arrived here, we have found a combinatorial counterexample
					count += 1
					print('N =', N, ' L0 =',L0, ' sum L_i =',L1, ' phi =',phi, ' gamma =', gamma, ' phi\' =', phi_, ' gamma\' =', gamma_)
							
print('\nWe have produced', count, 'combinatorial counterexamples.')
\end{lstlisting}

The output of this code is the following.
{\footnotesize
\begin{verbatim}
N = 1  L0 = 10  sum L_i = 2  phi = 1  gamma = 0  phi' = 5  gamma' = 5
N = 1  L0 = 8  sum L_i = 3  phi = 1  gamma = 0  phi' = 5  gamma' = 4

We have produced 2 combinatorial counterexamples.
\end{verbatim}
}

\section{\texorpdfstring{$V_{\lowercase{k}+1}\simeq \P^2$, and $N(V_{\lowercase{k}+1})=\frac{3}{2} N$}{V(k+1) isomorphic to P2, and N(V(k+1))=3N/2}}

Suppose $V_{k+1}\simeq \P^2$, and $N(V_{k+1})=\frac{3}{2} N$.
The strategy described in Section~\ref{sect:future-count-strategy} explains that we have to find integer solutions satisfying the following conditions:

\[\begin{cases}
\phi'+2\gamma'-\phi-L_0 = 4,\\
2\sum_{i=1}^{k+1} L_i+ 2\gamma'+4\phi'-4\phi-3\gamma=15,\\
L_0+2\sum_{i=1}^{k+1} L_i = 20,\\
2\phi+2\gamma \leq \sum_{i=1}^{k+1} L_i,\\
\phi'+\gamma' \leq L_0+1,\\
\sum_{i=1}^{k+1} L_i, \gamma, \phi' \text{ and } \gamma' \geq 0,\\
L_0 \text{ and } \phi \geq 1.
\end{cases}\]

Moreover, in Section~\ref{sect:future-exclude}, we found the extra equation $L_0+2\gamma' \not\equiv 0 \mod 4$.

One can verify that there is exactly \textbf{1 integer solution} satisfying these relations, for example by using the Python code below.
\begin{lstlisting}[caption = {}]
count=0

#L1 for \sum_{i=1}^beta L_i and L0 for L_0
for L1 in range(0, (20//2)+1): 		# L0+2* sum L_i=20
	L0 = 20 - 2*L1			
	if L0 == 0:
		#V_0 is non-minimal ruled, so L_0 \neq 0
		continue	

	for phi in range(1,(L1//2)+1):	# 2*phi+2*gamma \leq sum L_i and phi\geq 1
		for gamma in range(0, ((L1-2*phi)//2)+1):

			for phi_ in range(0,(L0+1)+1):	# phi'+gamma' \leq L0+1
				for gamma_ in range(0,(L0-phi_+1)+1):

					# checking nothing went wrong
					assert phi_>=0 and gamma_>=0 and phi>=1 and gamma>=0 and L0>=1 and L1>=0
					assert (L0+2*L1) == 20
					assert (phi_+gamma_) <= (L0+1)
					assert (2*phi+2*gamma) <= L1

					# checking the two conditions 
					if (phi_+2*gamma_-phi-L0) != 4:
						# making sure the mondromy property does not hold
						continue
					if (2*L1+2*gamma_+4*phi_-4*phi-3*gamma) != 15:
						# making sure degree of zeta function is 24
						continue
					
					if (L0+2*gamma_)%4 != 0:	
						# asserting the existence of line bundle such that cyclic cover exists
						continue
					
					# if we arrived here, we have found a combinatorial counterexample
					count += 1
					print('N = 1  L0 =',L0, ' sum L_i =',L1, ' phi =',phi, ' gamma =', gamma, ' phi\' =', phi_, ' gamma\' =', gamma_)
							
print('\nWe have produced', count, 'combinatorial counterexamples.')
\end{lstlisting}

The output of this code is the following.
{\footnotesize
\begin{verbatim}
N = 1  L0 = 4  sum L_i = 8  phi = 1  gamma = 3  phi' = 1  gamma' = 4

We have produced 1 combinatorial counterexamples.
\end{verbatim}
}

\section{\texorpdfstring{$V_{\lowercase{k}+1}\simeq \Sigma_2$, and $N(V_{\lowercase{k}+1})=4N$}{V(k+1) isomorphic to Sigma2, and N(V(k+1))=4N}} \label{sect:app-case6}

Suppose $V_{k+1}\simeq \Sigma_2$, and $N(V_{k+1})=4N$.
The strategy described in Section~\ref{sect:future-count-strategy} explains that we have to find integer solutions satisfying the following conditions:

\[\begin{cases}
\phi'+2\gamma'-\phi-L_0 = 7,\\
\sum_{i=1}^{k+1} L_i+ \gamma'+2\phi'-2\phi-\frac{3\gamma}{2}=\frac{24}{N} - 4,\\
L_0+2\sum_{i=1}^{k+1} L_i = 12,\\
2\phi+2\gamma \leq \sum_{i=1}^{k+1} L_i,\\
\phi'+\gamma' \leq L_0+1,\\
N = \begin{cases}
1 & \text{if } \gamma=0,\\
2 & \text{if } \gamma>0,
\end{cases}\\
\sum_{i=1}^{k+1} L_i, \gamma, \phi' \text{ and } \gamma' \geq 0,\\
L_0 \text{ and } \phi \geq 1.
\end{cases}\]

One can verify that there are \textbf{no integer solutions} satisfying these relations, for example by using the Python code below.
\begin{lstlisting}[caption = {}]
count=0

# L1 for \sum_{i+1}^beta L_i and L0 for L_0
for L1 in range(0, (12//2)+1): 	# L0+2* sum L_i=12
	L0 = 12 - 2*L1				
	if L0 == 0:
		#V_0 is non-minimal ruled, so L_0 \neq 0
		continue

	for phi in range(1,(L1//2)+1):	# 2*phi+2*gamma \leq sum L_i and phi\geq 1
		for gamma in range(0, ((L1-2*phi)//2)+1, 2): #gamma even because of degree zeta function
			
			# index = 1
			if gamma == 0:
				N=1				
			else:
				N=2

			for phi_ in range(0,(L0+1)+1):	# phi'+gamma' \leq L0+1
				for gamma_ in range(0, (L0-phi_+1)+1):

					# checking nothing went wrong
					assert phi_>=0 and gamma_>=0 and phi>=1 and gamma>=0 and L0>=1 and L1>=0 and (N==1 or N==2)
					assert gamma%2 == 0
					assert (L0+2*L1) == 12
					assert (phi_+gamma_) <= L0+1
					assert (2*phi+2*gamma) <= L1
					
					# checking the two conditions 
					if (phi_+2*gamma_-phi-L0) != 7:
						# making sure the mondromy property does not hold
						continue
					if (L1+2*phi_+gamma_-2*phi-3*gamma//2) != ((24//N)-4):
						# making sure degree of zeta function is 24
						continue
					
					# if we arrived here, we have found a combinatorial counterexample
					count += 1
					print('N =', N, ' L0 =',L0, ' sum L_i =',L1, ' phi =',phi, ' gamma =', gamma, ' phi\' =', phi_, ' gamma\' =', gamma_)
							
print('\nWe have produced', count, 'combinatorial counterexamples.')
\end{lstlisting}

\section{\texorpdfstring{$V_{\lowercase{k}+1}\simeq \Sigma_2$, $\beta=\lowercase{k}-1$, and $N(V_{\lowercase{k}+1})=2N$}{V(k+1) isomorphic to Sigma2, beta=k-1, and N(V(k+1))=2N}}\label{sect:app-case7}

Suppose $V_{k+1}\simeq \Sigma_2$, $\beta=k-1$, and $N(V_{k+1})=2N$.
The strategy described in Section~\ref{sect:future-count-strategy} explains that we have to find integer solutions satisfying the following conditions:

\[\begin{cases}
\phi'+2\gamma'-\phi-L_0 = 6,\\
2\sum_{i=1}^{k+1} L_i+ 2\gamma'+4\phi'-4\phi-3\gamma=22,\\
L_0+2\sum_{i=1}^{k+1} L_i = 16,\\
2\phi+2\gamma \leq \sum_{i=1}^{k+1} L_i,\\
\phi'+\gamma' \leq L_0+1,\\
\sum_{i=1}^{k+1} L_i, \gamma, \phi' \text{ and } \gamma' \geq 0,\\
L_0 \text{ and } \phi \geq 1.
\end{cases}\]

Moreover, in Section~\ref{sect:future-exclude}, we found the extra equation $L_0+2\gamma' \not\equiv 0 \mod 4$.

One can verify that there are \textbf{no integer solutions} satisfying these relations, for example by using the Python code below.
\begin{lstlisting}[caption = {}]
count=0

# L1 for $\sum_{i=1}^beta L_i and L0 for L_0
for L1 in range(0, (16//2)+1): 	# L0+2* sum L_i=16
	L0 = 16-2*L1			
	if L0 == 0:
		# V_O is non-minimal ruled, so L_0 \neq 0
		continue

	for phi in range(1, (L1//2)+1):	# 2*phi+2*gamma \leq sum L_i and phi\geq 1
		for gamma in range(0, ((L1-2*phi)//2)+1):

			for phi_ in range(0, (L0+1)+1):	# phi'+gamma' \leq L0+1
				for gamma_ in range(0, (L0-phi_+1)+1):

					# checking nothing went wrong
					assert phi_>=0 and gamma_>=0 and phi>=1 and gamma>=0 and L0>=1 and L1>=0
					assert (L0+2*L1) == 16
					assert (phi_+gamma_) <= (L0+1)
					assert (2*phi+2*gamma) <= L1
					
					# checking the conditions 
					if (phi_+2*gamma_-phi-L0) != 6:
						# making sure the mondromy property does not hold
						continue
					if (2*L1+4*phi_+2*gamma_-4*phi-3*gamma) != 22:
						# making sure degree of zeta function is 24
						continue

					if (L0+2*gamma_)%4 != 0:
						#necessary condition for the existence of line bundle such that cyclic cover exists
						continue

					# if we arrived here, we have found a combinatorial counterexample
					count += 1
					print('N = 1  L0 =',L0, ' sum L_i =',L1, ' phi =',phi, ' gamma =', gamma, ' phi\' =', phi_, ' gamma\' =', gamma_)
							
print('\nWe have produced', count, 'combinatorial counterexamples.')
\end{lstlisting}

\section{\texorpdfstring{$V_{\lowercase{k}+1}\simeq \Sigma_2$, and $N(V_{\lowercase{k}+1})=\frac{4}{3}N$}{V(k+1) isomorphic to Sigma2, and N(V(k+1))=4N/3}} \label{sect:app-case8}

Suppose $V_{k+1}\simeq \Sigma_2$, and $N(V_{k+1})=\frac{4}{3}N$.
The strategy described in Section~\ref{sect:future-count-strategy} explains that we have to find integer solutions satisfying the following conditions:

\[\begin{cases}
\phi'+2\gamma'-\phi-L_0 = 4,\\
3\sum_{i=1}^{k+1} L_i+ 3\gamma'+6\phi'-6\phi-\frac{9\gamma}{2}=\frac{24}{N} - 10,\\
L_0+2\sum_{i=1}^{k+1} L_i = 20,\\
2\phi+2\gamma \leq \sum_{i=1}^{k+1} L_i,\\
\phi'+\gamma' \leq L_0+1,\\
N = \begin{cases}
1 & \text{if } \gamma=0,\\
2 & \text{if } \gamma>0,
\end{cases}\\
\sum_{i=1}^{k+1} L_i, \gamma, \phi' \text{ and } \gamma' \geq 0,\\
L_0 \text{ and } \phi \geq 1.
\end{cases}\]

One can verify that there are \textbf{no integer solutions} satisfying these relations, for example by using the Python code below.
\begin{lstlisting}[caption = {}]
count=0

for L1 in range(0, (20//2)+1): 	# L0+2* sum L_i=20
	L0 = 20 - 2*L1		
	if L0 == 0:
		#V_0 is non-minimal ruled, so L_0 \neq 0
		continue		

	for phi in range(1,(L1//2)+1):	# 2*phi+2*gamma \leq sum L_i and phi\geq 1
		for gamma in range(0, ((L1-2*phi)//2)+1, 2): #gamma even because of degree zeta function
			
			# index = 1
			if gamma == 0:
				N=1				
			else:
				N=2

			for phi_ in range(0,(L0+1)+1):	# phi'+gamma' \leq L0+1
				for gamma_ in range(0, (L0-phi_+1)+1):

					# checking nothing went wrong
					assert phi_>=0 and gamma_>=0 and phi>=1 and gamma>= 0 and L0>=1 and L1>=0 and (N==1 or N==2)
					assert (L0+2*L1) == 20
					assert (phi_+gamma_) <= (L0 + 1)
					assert (2*phi+2*gamma) <= L1
					
					# checking the two conditions 
					if (phi_+2*gamma_-phi-L0) != 4:
						# making sure the mondromy property does not hold
						continue
					if (3*L1+6*phi_+3*gamma_-6*phi-9*gamma//2) != ((24//N)-10):
						# making sure degree of zeta function is 24
						continue
					
					# if we arrived here, we have found a combinatorial counterexample
					count += 1
					print('N =', N, ' L0 =',L0, ' sum L_i =',L1, ' phi =',phi, ' gamma =', gamma, ' phi\' =', phi_, ' gamma\' =', gamma_)
							
print('\nWe have produced', count, 'combinatorial counterexamples.')
\end{lstlisting}

\section{\texorpdfstring{$V_{\lowercase{k}+1}$ is a rational, ruled surface, $\beta =\lowercase{k}$, and $N(V_{\lowercase{k}+1})=2N$}{V(k+1) is a rational, ruled surface, beta =k, and N(V(k+1))=2N}}\label{sect:app-case9}

Suppose $V_{k+1}$ is a rational, ruled surface, $\beta =k$, and $N(V_{k+1})=2N$.
The strategy described in Section~\ref{sect:future-count-strategy} explains that we have to find integer solutions satisfying the following conditions:

\[\begin{cases}
\phi'+2\gamma'-\phi-(L_0+L_{k+1}) = 8,\\
\sum_{i=1}^k L_i + 2\phi' + \gamma'-2\phi-\frac{3\gamma}{2} =  \frac{24}{N},\\
L_0 + 2\sum_{i=1}^k L_i + L_{k+1} = 16,\\
\phi'+\gamma'\leq L_0+L_{k+1}+2,\\
2\phi + 2\gamma \leq \sum_{i=1}^k L_i,\\
N = \begin{cases}
1 & \text{if } \gamma = 0,\\
2 & \text{if } \gamma > 0,
\end{cases}\\
\sum_{i=1}^k L_i, \gamma, \phi' \text{ and } \gamma' \geq 0,\\
L_0+L_{k+1} \text{ and } \phi \geq 1.
\end{cases}\]

One can verify that there are \textbf{no integer solutions} satisfying these relations, for example by using the Python code below.
\begin{lstlisting}[caption = {}]
count=0

# L1 for \sum_{i=1}^k L_i and L0 for L_0+L_{k+1}.
for L1 in range(0, (16//2)+1): 	# L0+2L1 = 16.
	L0 = 16 - 2*L1
	if L0 == 0:
		# V_0 is non-minimal ruled, so L_0 \neq 0.
		continue

	for phi in range(1,(L1//2)+1):	# 2*phi+2*gamma \leq sum L_i and phi\geq 1
		for gamma in range(0, ((L1-2*phi)//2)+1):
			
			# index = 1
			if gamma == 0:
				N=1				
			else:
				N=2

			for phi_ in range(0,(L0+2)+1):	# phi'+gamma' \leq L_0+L_{k+1}+2
				for gamma_ in range(0, (L0-phi_+2)+1):

					# checking nothing went wrong
					assert phi_>=0 and gamma_>=0 and phi>=1 and gamma>=0 and L0>=1 and L1>=0 and (N==1 or N==2)
					assert (L0+2*L1) == 16
					assert (phi_+gamma_) <= (L0+2)
					assert (2*phi+2*gamma) <= L1
					
					# checking the two conditions 
					if (phi_+2*gamma_-phi-L0) != 8:
						# making sure the mondromy property does not hold
						continue
					if (L1+2*phi_+gamma_-2*phi-3*gamma//2) != (24//N):
						# making sure the degree zeta function is 24
						continue
					
					# if we arrived here, we have found a combinatorial counterexample
					count += 1
					print('N =', N, ' L0 =',L0, ' sum L_i =',L1, ' phi =',phi, ' gamma =', gamma, ' phi\' =', phi_, ' gamma\' =', gamma_)
							
print('\nWe have produced', count, 'combinatorial counterexamples.')
\end{lstlisting}

\section{\texorpdfstring{$V_{\lowercase{k}+1}$ is an elliptic, ruled surface, and $N(V_{\lowercase{k}+1})=2N$}{V(k+1) is an elliptic, ruled surface, and N(V(k+1))=2N}}
\label{sect:app-case10}

Suppose $V_{k+1}$ is an elliptic, ruled surface, and $N(V_{k+1})=2N$.
The strategy described in Section~\ref{sect:future-count-strategy} explains that we have to find integer solutions satisfying the following conditions:

\[\begin{cases}
\phi'+2\gamma'-\phi-2g-(L_0+L_{k+1}) = 2,\\
\sum_{i=1}^k L_i + 2\phi' + \gamma'-2g-2\phi-\frac{3\gamma}{2} =  \frac{24}{N}-2,\\
L_0 + 2\sum_{i=1}^k L_i + L_{k+1} = 8,\\
\phi'+\gamma'\leq L_0+L_{k+1}+1,\\
2\phi + 2\gamma \leq \sum_{i=1}^k L_i,\\
g\leq (\sum_{i=1}^k L_i -2\phi -2\gamma + \gamma' +2) /2\\
N = \begin{cases}
1 & \text{if } \gamma = 0,\\
2 & \text{if } \gamma > 0,
\end{cases}\\
\sum_{i=1}^k L_i, \gamma, \phi' \text{ and } \gamma' \geq 0,\\
L_0+L_{k+1}, g \text{ and } \phi \geq 1.
\end{cases}\]

One can verify that there are \textbf{no integer solutions} satisfying these relations, for example by using the Python code below.
\begin{lstlisting}[caption = {}]
count=0

# L1 for $\sum_{i=1}^k L_i and L0 for L_0+L_{k+1}
for L1 in range(0, (8//2)+1):
	L0 = 8-2*L1
	if L0 == 0:
		#V_0 is non-minimal ruled, so L_0 \neq 0.
		continue

	for phi in range(1,(L1//2)+1):	# 2*phi+2*gamma \leq sum L_i and phi\geq 1
		for gamma in range(0, ((L1-2*phi)//2)+1, 2): #gamma is even because of degree zeta function
			
			# index = 1
			if gamma == 0:
				N=1				
			else:
				N=2

			for phi_ in range(0,(L0+1)+1):	# phi'+gamma' \leq L0+1
				for gamma_ in range(0, (L0-phi_+1)+1):

					for g in range(1, ((L1-2*phi-2*gamma+gamma_+2)//2)+1):

						# checking nothing went wrong
						assert phi_>=0 and gamma_>=0 and phi>=1 and gamma>=0 and (N==1 or N==2)
						assert g>=1 and L0>=1 and L1>=0
						assert gamma % 2 ==0
						assert (L0+2*L1) == 8
						assert (phi_+gamma_) <= (L0+1)
						assert (2*phi+2*gamma) <= L1
						assert g <= (L1-2*phi-2*gamma+gamma_+2)//2
					
						# checking the two conditions 
						if (phi_+2*gamma_) != (2+L0+2*g+phi):
							# making sure the mondromy property does not hold
							continue
						if (L1+2*phi_+gamma_-2*g-2*phi-3*gamma//2) != (24//N-2):
							# making sure the degree zeta function is 24
							continue
					
						# if we arrived here, we have found a combinatorial counterexample
						count += 1
						print('N =', N, ' L0 =',L0, ' sum L_i =',L1, ' L_k+1 =',L2, ' phi =',phi, ' gamma =', gamma,\
							' phi\' =', phi_, ' gamma\' =', gamma_, ' g =', g)
							
print('\nWe have produced', count, 'combinatorial counterexamples.')
\end{lstlisting}

\cleardoublepage


\backmatter
\includebibliography
\bibliographystyle{mystyle}
\bibliography{algbib}

\begin{thebibliography}{rTSPA17}
\providecommand{\url}[1]{\texttt{#1}}
\providecommand{\urlprefix}{URL }
\expandafter\ifx\csname urlstyle\endcsname\relax
  \providecommand{\doi}[1]{doi:\discretionary{}{}{}#1}\else
  \providecommand{\doi}{doi:\discretionary{}{}{}\begingroup
  \urlstyle{rm}\Url}\fi
\providecommand{\eprint}[2][]{\url{#2}}

\bibitem[SGA7]{SGA7}
\emph{Groupes de monodromie en g\'eometrie alg\'ebrique. I}, volume 288 of
  \emph{Lecture {N}otes in {M}athematics}.
\newblock Springer-Verlag, Berlin.
\newblock S{\'e}minaire de G{\'e}om{\'e}trie Alg{\'e}brique du Bois-Marie
  1967--1969 (SGA 7 I), Dirig{\'e} par Alexandre Grothendieck. Avec la
  collaboration de {M}. Raynaud et D.S. Rim (1972).

\bibitem[{A'C}75]{Acampo}
N.~{A'Campo}.
\newblock La fonction zeta d'une monodromie.
\newblock \emph{Commentarii mathematici Helvetici}, volume~50: pp. 233--248
  (1975).

\bibitem[ACLM02]{Artal-Cassou-et-al}
E.~{Artal Bartolo}, P.~{Cassou-Nogu\`es}, I.~Luengo, and A.~{Melle
  Hern\'andez}.
\newblock Monodromy conjecture for some surface singularities.
\newblock \emph{{A}nnales {S}cientifiques de l’{\'E}cole {N}ormale
  {S}up\'erieure}, volume~35(4): pp. 605--640 (2002).

\bibitem[Art70]{Artin}
M.~Artin.
\newblock Algebraization of formal moduli {II}: Existence of modifications.
\newblock \emph{Annals of Mathematics}, volume~91(1): pp. 88--135 (1970).

\bibitem[BHPV95]{BHPV}
W.~P. Barth, K.~Hulek, C.~A. Peters, and A.~{Van de Ven}.
\newblock \emph{Compact complex surfaces}, volume~4 of \emph{Ergebnisse der
  Mathematik und ihrer Grenzgebiete. 3. Folge. A Series of Modern Surveys in
  Mathematics}.
\newblock Springer-Verlag, Berlin, second edition (1995).

\bibitem[BS66]{BorevichShafarevich}
Z.~I. Borevich and I.~R. Shafarevich.
\newblock \emph{Number theory}, volume~20 of \emph{Pure and Applied
  Mathematics}.
\newblock Academic Press, New York-London (1966).

\bibitem[Bor13]{Bart-thesis}
B.~Bories.
\newblock \emph{Zeta functions, {B}ernstein–{S}ato polynomials, and the
  monodromy conjecture}.
\newblock Ph.D. thesis, KU Leuven in Leuven, Belgium (2013).

\bibitem[Bor17]{Borisov15}
L.~Borisov.
\newblock Class of the affine line is a zero divisor in the {G}rothendieck
  ring.
\newblock To appear in \emph{Journal of Algebraic Geometry} (2017).

\bibitem[BN16]{BultotNicaise}
E.~Bultot and J.~Nicaise.
\newblock Computing motivic zeta functions on log smooth models.
\newblock {a}rXiv: 1610.00742v1 (2016).

\bibitem[Cau16]{Cauwbergs}
T.~Cauwbergs.
\newblock Splicing motivic zeta functions.
\newblock \emph{Revista Matem\'atica Complutense}, volume~29(2): pp. 455--483
  (2016).

\bibitem[CNS17]{MotivicIntegration}
A.~{Chambert-Loir}, J.~Nicaise, and J.~Sebag.
\newblock \emph{Motivic integration}, volume 571 of \emph{Progress in
  Mathematics}.
\newblock Birkha\"user (2017).

\bibitem[Cra83]{Crauder}
B.~Crauder.
\newblock Threefold birational morphisms and degenerations without triple
  points.
\newblock In: R.~Friedman and D.~R. Morrison (eds.), \emph{The Birational
  Geometry of Degenerations}, volume~29 of \emph{Progress in Mathematics}, pp.
  299--352. Birkh\"auser Verlag, Basel (1983).

\bibitem[CM83]{CrauMor}
B.~Crauder and D.~R. Morrison.
\newblock Triple-point-free degenerations of surfaces with {K}odaira number
  zero.
\newblock In: R.~Friedman and D.~R. Morrison (eds.), \emph{The Birational
  Geometry of Degenerations}, volume~29 of \emph{Progress in Mathematics}, pp.
  353--386. Birkh\"auser Boston, Boston (1983).

\bibitem[CM94]{CrauMor2}
B.~Crauder and D.~R. Morrison.
\newblock Minimal models and degenerations of surfaces with {K}odaira number
  zero.
\newblock \emph{Transactions of the American Mathematical Society}, volume
  343(2): pp. 525--558 (1994).

\bibitem[Den84]{Denef84}
J.~Denef.
\newblock The rationality of the {P}oincaré series associated to the $p$-adic
  points on a variety.
\newblock \emph{Inventiones Mathematicae}, volume~77: pp. 1--24 (1984).

\bibitem[Den93]{Denef93}
J.~Denef.
\newblock Degree of local zeta functions and monodromy.
\newblock \emph{Compositio Mathematica}, volume~89(2): pp. 207--216 (1993).

\bibitem[DL98]{DenefLoeser98}
J.~Denef and F.~Loeser.
\newblock Motivic {I}gusa zeta functions.
\newblock \emph{Journal of {A}lgebraic {G}eometry}, volume~7: pp. 505--537
  (1998).

\bibitem[DL01]{DenLoe01}
J.~Denef and F.~Loeser.
\newblock Geometry on arc spaces of algebraic varieties.
\newblock In: \emph{Proceedings of the third ECM, Barcelona}, volume 201 of
  \emph{Progress in Mathematics}, pp. 327--348 (2001).

\bibitem[Eis95]{Eisenbud}
D.~Eisenbud.
\newblock \emph{Commutative algebra with a view toward algebraic geometry},
  volume 150 of \emph{Graduate Texts in Mathematics}.
\newblock Springer-Verlag, New York (1995).

\bibitem[Eke09]{Ekedahl}
T.~Ekedahl.
\newblock The {G}rothendieck group of algebraic stacks.
\newblock {a}rXiv:0903.3143v2 (2009).

\bibitem[ELW15]{ELW}
H.~Esnault, M.~Levine, and O.~Wittenberg.
\newblock Index of varieties over {H}enselian fields and {E}uler characteristic
  of coherent sheaves.
\newblock \emph{Journal of {A}lgebraic {G}eometry}, volume~24: pp. 693--718
  (2015).

\bibitem[Fri83]{Friedman}
R.~Friedman.
\newblock Global smoothings of varieties with normal crossings.
\newblock \emph{Annals of Mathematics}, volume 118(1): pp. 75--114 (1983).

\bibitem[Ful84]{Fulton-intersection}
W.~Fulton.
\newblock \emph{Intersection theory}.
\newblock Ergebnisse der Mathematik und ihrer Grenzgebiete. 3. Folge; 2.
  Springer-Verlag (1984).

\bibitem[HN11]{HalleNicaiseAbelian}
L.~H. Halle and J.~Nicaise.
\newblock Motivic zeta functions of abelian varieties, and the monodromy
  conjecture.
\newblock \emph{Advances in Mathematics}, volume 227: pp. 610--653 (2011).

\bibitem[HN12]{HalleNicaiseMotivicCY}
L.~H. Halle and J.~Nicaise.
\newblock Motivic zeta functions for degenerations of abelian varieties and
  {C}alabi-{Y}au varieties.
\newblock In: A.~C. et~al. (ed.), \emph{Recent trends on Zeta functions in
  algebra and geometry}, volume 566 of \emph{Contemporary Mathematics}, pp.
  233--259. American Mathematical Society (2012).

\bibitem[HN17]{HalleNicaiseKulikov}
L.~H. Halle and J.~Nicaise.
\newblock Motivic zeta functions of degenerating {C}alabi-{Y}au varieties.
\newblock {a}rXiv:1701.09155v1 (2017).

\bibitem[Han15]{Hanke}
D.~A. Hanke.
\newblock \emph{Logarithmically smooth deformations of strict normal crossing
  logarithmically symplectic varieties}.
\newblock Ph.D. thesis, Johannes Gutenberg-Universit\"at in Mainz, Germany
  (2015).

\bibitem[Har15]{Hartmann_motivic}
A.~Hartmann.
\newblock Equivariant motivic integration on formal schemes and the motivic
  zeta function.
\newblock {a}rXiv: 1511.08656v1 (2015).

\bibitem[Har16]{Hartmann_quotient}
A.~Hartmann.
\newblock The quotient map on the equivariant {G}rothendieck ring of varieties.
\newblock \emph{Manuscripta Mathematica}, volume 151(3): pp. 419--451 (2016).

\bibitem[Har77]{Hartshorne}
R.~Hartshorne.
\newblock \emph{Algebraic geometry}, volume~52 of \emph{Graduate Texts in
  Mathematics}.
\newblock Springer-Verlag, New York (1977).

\bibitem[Hir64]{Hironaka}
H.~Hironaka.
\newblock Resolution of singularities of an algebraic variety over a field of
  characteristic zero {I}.
\newblock \emph{Annals of Mathematics}, volume~79(1): pp. 109--203 (1964).

\bibitem[Huy16]{HuybrechtsK3}
D.~Huybrechts.
\newblock \emph{Lectures on $K3$ surfaces}, volume 158 of \emph{Cambridge
  Studies in Advanced Mathematics}.
\newblock Cambridge University Press (2016).

\bibitem[Igu74]{Igusa74}
J.-I. Igusa.
\newblock Complex powers and asymptotic expansions {I}.
\newblock \emph{Journal f\"ur die {R}eine und {A}ngewandte {M}athematik},
  volume 268/269: pp. 110--130 (1974).

\bibitem[Igu75]{Igusa75}
J.-I. Igusa.
\newblock Complex powers and asymptotic expansions {II}.
\newblock \emph{Journal f\"ur die {R}eine und {A}ngewandte {M}athematik},
  volume 278/279: pp. 307--321 (1975).

\bibitem[Igu77]{Igusa77}
J.-I. Igusa.
\newblock Some observations on higher degree characters.
\newblock \emph{American Journal of Mathematics}, volume~99(2): pp. 393--417
  (1977).

\bibitem[Jas17]{Jaspers}
A.~Jaspers.
\newblock The global monodromy property for {$K3$} surfaces allowing a
  triple-point-free model.
\newblock \emph{Comptes Rendus - Math\'ematique}, volume 355(2): pp. 200--204
  (2017).

\bibitem[Kat94]{Katz-ell-adic-cohomology}
N.~M. Katz.
\newblock Review of $\ell$-adic cohomology.
\newblock In: U.~Jannsen, S.~Kleiman, and J.-P. Serre (eds.), \emph{Motives},
  volume 55, part 1 of \emph{Proceedings of Symposia in Pure Mathematics}, pp.
  21--30. American Mathematical Society (1994).

\bibitem[KKMS73]{SemistableReductionTheorem}
G.~Kempf, F.~Knudsen, D.~Mumford, and B.~{Saint-Donat}.
\newblock \emph{Toroidal Embeddings 1}, volume 339 of \emph{Lecture notes in
  Mathematics}.
\newblock Springer-Verlag (1973).

\bibitem[Knu71]{Knutson}
D.~Knutson.
\newblock \emph{Algebraic Spaces}, volume 203 of \emph{Lecture Notes in
  Mathematics}.
\newblock Springer (1971).

\bibitem[KM98]{KollarMori}
J.~Koll\'ar and S.~Mori.
\newblock \emph{Birational geometry of algebraic varieties}, volume 134 of
  \emph{Cambridge tracts in mathematics}.
\newblock Cambridge University Press (1998).

\bibitem[KX16]{DualComplex}
J.~Koll\'ar and C.~Xu.
\newblock The dual complex of {C}alabi-{Y}au pairs.
\newblock \emph{Inventiones Mathematicae}, volume 205(3): pp. 527--557 (2016).

\bibitem[Kul77]{Kulikov1}
V.~Kulikov.
\newblock Degenerations of {K}3 surfaces and {E}nriques surfaces.
\newblock \emph{Math. USSR Izvestija}, volume~11: pp. 957--989 (1977).

\bibitem[Lau81]{Laumon}
G.~Laumon.
\newblock Comparaison de caract\'eristiques d'{E}uler-{P}oincar\'e en
  cohomologie $\ell$-adique.
\newblock \emph{Comptes Rendus - Math\'ematique}, volume 292(3): pp. 209--212
  (1981).

\bibitem[Liu02]{Liu}
Q.~Liu.
\newblock \emph{Algebraic geometry and arithmetic curves}, volume~6 of
  \emph{Oxford Graduate Texts in Mathematics}.
\newblock Oxford University Press, Oxford (2002).

\bibitem[Loe88]{Loeser88}
F.~Loeser.
\newblock Fonctions d'{I}gusa $p$-adiques et polynomes de {B}ernstein.
\newblock \emph{American Journal of Mathematics}, volume 110(1): pp. 1--21
  (1988).

\bibitem[LS03]{LoeserSebag}
F.~Loeser and J.~Sebag.
\newblock Motivic integration on smooth rigid varieties and invariants of
  degenerations.
\newblock \emph{Duke Mathematical Journal}, volume 119(2): pp. 315--344 (2003).

\bibitem[Nag62]{Nagata1}
M.~Nagata.
\newblock Imbedding of an abstract variety in a complete variety.
\newblock \emph{Journal of Mathematics of Kyoto University}, volume~2(1): pp.
  1--10 (1962).

\bibitem[Nag63]{Nagata2}
M.~Nagata.
\newblock A generalization of the imbedding problem of an abstract variety in a
  complete variety.
\newblock \emph{Journal of Mathematics of Kyoto University}, volume~3(1): pp.
  89--102 (1963).

\bibitem[Nic09]{Nicaise_trace-formula-rigid-varieties}
J.~Nicaise.
\newblock A trace formula for rigid varieties, and motivic weil generating
  series for formal schemes.
\newblock \emph{Mathematische Annalen}, volume 343(2): pp. 285--349 (2009).

\bibitem[Nic10]{NicaisePadicMotivicMonConj}
J.~Nicaise.
\newblock An introduction to $p$-adic and motivic zeta functions and the
  monodromy conjecture.
\newblock In: K.~M. G.~Bhowmik and H.~Tsumura (eds.), \emph{Algebraic and
  analytic aspects of zeta functions and L-functions}, volume~21 of \emph{MSJ
  Memoirs}, pp. 115--140. Mathematical Society of Japan (2010).

\bibitem[Nic11]{NicaiseTrace}
J.~Nicaise.
\newblock A trace formula for varieties over a discretely valued field.
\newblock \emph{Journal f\"ur die Reine und Angewandte Mathematik}, volume 650:
  pp. 193--238 (2011).

\bibitem[Nic13]{NicaiseTameRam}
J.~Nicaise.
\newblock Geometric criteria for tame ramification.
\newblock \emph{Mathematische Zeitschrift}, volume 273(3): pp. 839--868 (2013).

\bibitem[NS07]{NicaiseSebag_motivic_serre}
J.~Nicaise and J.~Sebag.
\newblock Motivic serre invariants, ramification, and the analytic {M}ilnor
  fiber.
\newblock \emph{Inventiones Mathematicae}, volume 168(1): pp. 133--173 (2007).

\bibitem[NS11]{NicaiseSebag_Grothendieck_ring}
J.~Nicaise and J.~Sebag.
\newblock The {G}rothendieck ring of varieties.
\newblock In: R.~Cluckers, J.~Nicaise, and J.~Sebag (eds.), \emph{Motivic
  integration and its interactions with model theory and non-archimedean
  geometry}, volume 383 of \emph{London Mathematical Society Lecture Notes
  Series}, pp. 145--188. Cambridge University Press (2011).

\bibitem[NX16]{NicaiseXu}
J.~Nicaise and C.~Xu.
\newblock Poles of maximal order of motivic zeta functions.
\newblock \emph{Duke Mathematical Journal}, volume 165(2): pp. 217--243 (2016).

\bibitem[Per77]{Persson}
U.~Persson.
\newblock \emph{On degenerations of algebraic surfaces}, volume~11 of
  \emph{Memoirs of the American Mathematical Society} (1977).

\bibitem[PP81]{Persson-Pinkham}
U.~Persson and H.~Pinkham.
\newblock Degeneration of surfaces with trivial canonical bundle.
\newblock \emph{Annals of Mathematics}, volume 113(1): pp. 45--66 (1981).

\bibitem[Poo02]{Poonen02}
B.~Poonen.
\newblock The {G}rothendieck ring of varieties is not a domain.
\newblock \emph{Mathematical Research Letters}, volume~9(4): pp. 493--498
  (2002).

\bibitem[Stacks]{stacks-project}
\relax The Stacks Project~Authors.
\newblock Stacks project.
\newblock \url{http://stacks.math.columbia.edu} (2017).

\bibitem[RV01]{Rodrigues-Veys01}
B.~Rodrigues and W.~Veys.
\newblock Holomorphy of {I}gusa's and topological zeta functions for
  homogeneous polynomials.
\newblock \emph{Pacific Journal of Mathematics}, volume 201(2): pp. 429--440
  (2001).

\bibitem[Seg11]{Segers}
D.~Segers.
\newblock The asymptotic behaviour of the number of solutions of polynomial
  congruences.
\newblock \emph{Analele \c{S}tiin\c{t}ifice ale Universitatii “Ovidius”
  Constan\c{t}a, Seria Matematica}, volume~19(1): pp. 255--262 (2011).

\bibitem[SV11]{Stewart-Vologodsky}
A.~J. Stewart and V.~Vologodsky.
\newblock Motivic integral of {$K3$} surfaces over a nonarchimedean field.
\newblock \emph{Advances in Mathematics}, volume 228(5): pp. 2688--2730 (2011).

\end{thebibliography}

\makebackcoverXII

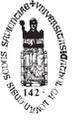
\end{document}